\documentclass[a4paper,11pt]{article}
\usepackage[all]{xy}
\usepackage[T1]{fontenc}
\usepackage[dvips]{graphicx}
\usepackage{amsmath,amssymb,amsfonts,amsthm,latexsym,mathrsfs,textcomp,verbatim, enumerate}
\xyoption{web}
\begin{document}

\title{Extensions of flat functors\\ and theories of presheaf type}

\author{Olivia Caramello\thanks{Supported by a CARMIN IH\'ES-IHP post-doctoral position (as from 1/12/2013) and by a visiting position of the Max Planck Institute for Mathematics (in the period 1/10/2013 - 30/11/2013).}\vspace{3 mm}\\ {\small Institut des Hautes \'Etudes Scientifiques}\\{\small 35 route de Chartres 91440, Bures-sur-Yvette, France}\\{\small olivia@ihes.fr}  }

\bgroup           
\let\footnoterule\relax  
\date{June 20, 2014}
\maketitle

\begin{abstract}
We develop a general theory of extensions of flat functors along geometric morphisms of toposes, and apply it to the study of the class of theories whose classifying topos is equivalent to a presheaf topos. As a result, we obtain a characterization theorem providing necessary and sufficient semantic conditions for a theory to be of presheaf type. This theorem subsumes all the previous partial results obtained on the subject and has several corollaries which can be used in practice for testing whether a given theory is of presheaf type as well as for generating new examples of theories belonging to this class. Along the way, we establish a number of other results of independent interest, including developments about colimits in the context of indexed categories, expansions of geometric theories and methods for constructing theories classified by a given presheaf topos.    
\end{abstract}

\egroup
\tableofcontents


\def\Monthnameof#1{\ifcase#1\or
   January\or February\or March\or April\or May\or June\or
   July\or August\or September\or October\or November\or December\fi}
\def\today{\number\day~\Monthnameof\month~\number\year}

%
%
%
\def\pushright#1{{
   \parfillskip=0pt            
   \widowpenalty=10000         
   \displaywidowpenalty=10000  
   \finalhyphendemerits=0      
  %
   \leavevmode                 
   \unskip                     
   \nobreak                    
   \hfil                       
   \penalty50                  
   \hskip.2em                  
   \null                       
   \hfill                      
   {#1}                        
  %
   \par}}                      

\def\qed{\pushright{$\square$}\penalty-700 \smallskip}

\newtheorem{theorem}{Theorem}[section]

\newtheorem{proposition}[theorem]{Proposition}

\newtheorem{scholium}[theorem]{Scholium}

\newtheorem{lemma}[theorem]{Lemma}

\newtheorem{corollary}[theorem]{Corollary}

\newtheorem{conjecture}[theorem]{Conjecture}

\newenvironment{proofs}%
 {\begin{trivlist}\item[]{\bf Proof }}%
 {\qed\end{trivlist}}

  \newtheorem{rmk}[theorem]{Remark}
\newenvironment{remark}{\begin{rmk}\em}{\end{rmk}}

  \newtheorem{rmks}[theorem]{Remarks}
\newenvironment{remarks}{\begin{rmks}\em}{\end{rmks}}

  \newtheorem{defn}[theorem]{Definition}
\newenvironment{definition}{\begin{defn}\em}{\end{defn}}

  \newtheorem{eg}[theorem]{Example}
\newenvironment{example}{\begin{eg}\em}{\end{eg}}

  \newtheorem{egs}[theorem]{Examples}
\newenvironment{examples}{\begin{egs}\em}{\end{egs}}


\mathcode`\<="4268  
\mathcode`\>="5269  
\mathcode`\.="313A  
\mathchardef\semicolon="603B 
\mathchardef\gt="313E
\mathchardef\lt="313C

\newcommand{\cod}
 {{\rm cod}}

\newcommand{\comp}
 {\circ}

\newcommand{\Cont}
 {{\bf Cont}}

\newcommand{\dom}
 {{\rm dom}}

\newcommand{\empstg}
 {[\,]}

\newcommand{\epi}
 {\twoheadrightarrow}

\newcommand{\hy}
 {\mbox{-}}

\newcommand{\im}
 {{\rm im}}

\newcommand{\imp}
 {\!\Rightarrow\!}

\newcommand{\Ind}[1]
 {{\rm Ind}\hy #1}

\newcommand{\mono}
 {\rightarrowtail}
 
\newcommand{\name}[1]
 {\mbox{$\ulcorner #1 \urcorner$}}

\newcommand{\ob}
 {{\rm ob}}

\newcommand{\op}
 {^{\rm op}}

\newcommand{\Set}
 {{\bf Set }}

\newcommand{\Sh}
 {{\bf Sh}}

\newcommand{\sh}
 {{\bf sh}}

\newcommand{\Sub}
 {{\rm Sub}}

\newcommand{\fu}[2]
{[#1,#2]}

\newcommand{\st}[2]
 {\mbox{$#1 \to #2$}}

\section{Introduction}

Following \cite{beke}, we say that a geometric theory is of \emph{of presheaf type} if it is classified by a presheaf topos.  

A geometric theory $\mathbb T$ is of presheaf type if and only if it is classified by the topos $[\textrm{f.p.} {\mathbb T}\textrm{-mod}(\Set), \Set]$, where $\textrm{f.p.} {\mathbb T}\textrm{-mod}(\Set)$ is (a skeleton of) the full subcategory of ${\mathbb T}\textrm{-mod}(\Set)$ on the finitely presentable $\mathbb T$-models (cf. \cite{OC}). 

The subject of theories of presheaf type has a long history, starting with the book \cite{Hakim} by Hakim, which first introduced the point of view of classifying toposes in the context of the theory of commutative rings with unit and its quotients. The subsequent pionereeing work \cite{Lawvere} by Lawvere led to the discovery that any finitary algebraic theory is of presheaf type, classified by the topos of presheaves on the opposite of its category of finitely presentable models (cf. \cite{johwra}).  This result was later generalized to cartesian (or essentially algebraic) theories as well as to universal Horn theories (cf. \cite{blasce}). At the same time, new examples of non-cartesian theories of presheaf type were discovered (cf. for instance \cite{beke} for a long, but by no means exhaustive, list of examples), and partial results in connection to the problem of characterizing the class of theories of presheaf type emerged; for instance, \cite{diers}, \cite{beke} and \cite{vickers} contain different sets of sufficient conditions for a theory to be of presheaf type. 

Theories of presheaf type occupy a central role in Topos Theory for a number of reasons:

\begin{enumerate}[(i)]

\item Every small category $\cal C$ can be seen, up to Cauchy-completion, as the category of finitely presentable models of a theory of presheaf type (namely, the theory of flat functors on ${\cal C}^{\textrm{op}}$);  

\item As every Grothendieck topos is a subtopos of some presheaf topos, so every geometric theory is a quotient of some theory of presheaf type (cf. the duality theorem of \cite{OC6} between subtoposes of the classifying topos of a geometric theory and quotients of the theory);

\item Every finitary algebraic theory (and more generally, any cartesian theory) is of presheaf type;

\item The class of theories of presheaf type contains, besides all cartesian theories, many other interesting mathematical theories pertaining to different fields of mathematics (for instance, the coherent theory of linear orders or the geometric theory of algebraic extensions of a given field);  

\item The `bridge technique' of \cite{OC10} can be fruitfully applied in the context of theories of presheaf type due to the fact that the classifying topos of any such theory admits (at least) two quite different representations, one of semantic nature (namely, set-valued functors on the category of finitely presentable models of the theory) and one of syntactic nature (namely, sheaves on the syntactic site of the theory).
\end{enumerate}

It is therefore important to dispose of effective criteria for testing whether a theory is of presheaf type, as well as of methods for generating new theories of presheaf type.     

In this paper, we carry out a systematic investigation of this class of theories, obtaining in particular a characterization theorem providing necessary and sufficient conditions for a theory to be of presheaf type, expressed in terms of the models of the theory in arbitrary Grothendieck toposes. This theorem, whose general statement is quite abstract, admits several ramifications and simpler corollaries which can be effectively applied in practice to test whether a given theory is of presheaf type as well as for generating new examples of theories of presheaf type, also through appropriate `modifications' of given geometric theories. All the partial results and recognition criteria previously obtained on the subject are naturally subsumed by this general result; moreover, the constructive nature of the characterization theorem allows to replace the requirements that the theory should have enough set-based models in the sufficient criteria of \cite{beke} and \cite{diers} with explicit semantic conditions which can be directly verified without having to invoke any form of the axiom of choice.  

In order to establish our characterization theorem, we embark, in the first two sections of the paper, in a general analysis of indexed colimits in toposes and extensions of flat functors along geometric morphisms. In fact, the necessary and sufficient conditions for a theory $\mathbb T$ to be of presheaf type provided by the characterization theorem arise precisely from the requirement that for any Grothendieck topos $\cal E$, the operation of extension of flat functors with values in $\cal E$ from the opposite of the category of finitely presentable models of $\mathbb T$ to the geometric syntactic category of $\mathbb T$ should define an equivalence of categories onto the category of $\mathbb T$-models in $\cal E$, naturally in $\cal E$. We then investigate the preservation, by `faithful interpretations' of theories, of each of the conditions in the characterization theorem, obtaining results of the form `under appropriate conditions, a geometric theory in which a theory of presheaf type faithfully interprets is again of presheaf type'. Finally, we discuss known and new examples of theories of presheaf type in light of the theory developed in the paper.

More specifically, the contents of the paper can be summarized as follows.   

In section \ref{indexed}, we investigate $\cal E$-indexed colimits of internal diagrams in Grothendieck toposes $\cal E$, analyzing in particular their behaviour with respect to final $\cal E$-indexed functors (cf. section \ref{finsub}) and establishing explicit characterizations for a $\cal E$-indexed cocone to be colimiting (cf. section \ref{charindcol}). In section \ref{seccol}, we exploit the abstract interpretation of colimits as kinds of tensor products to derive commutation results which play an important role in the subsequent parts of the paper as they allow us to interpret certain set-indexed colimits arising in the context of our main characterization theorem as special kinds of filtered indexed colimits.    

In section 3, we investigate the properties of the operation on flat functors induced by a geometric morphism of toposes via Diaconescu's equivalence. We focus in particular on geometric morphisms between presheaf toposes induced by embeddings between small categories, and on geometric morphisms to the classifying topos of a geometric theory induced by a small category of set-based models of the theory. We also establish, in section \ref{generaladjunction}, a general `hom-tensor' adjunction between categories of ${\cal E}$-valued functors (for $\cal E$ a Grothedieck topos) $[{\cal C}, {\cal E}]$ and $[{\cal D}, {\cal E}]$ induced by a functor $P:{\cal C}\to [{\cal D}^{\textrm{op}}, \Set]$, which generalizes the well-known adjunction induced by Kan extensions along a given functor.   

In section 4, in order to set up the field for the statement and proof of the characterization theorem, we identify some notable properties of theories of presheaf type, notably including the fact that every finitely presentable model of such a theory is finitely presented - in a strong sense which we make precise in section \ref{strong} - and admits an entirely syntactic description in terms of the signature of the theory and the notion of provability of geometric sequents over it in the theory (cf. section \ref{syntacticdescription}). We also show, in section \ref{objhom}, that for  any geometric theory $\mathbb T$ and any Grothendieck topos $\cal E$, there exists for each pair of $\mathbb T$-models $M$ and $N$ in $\cal E$, an `object of $\mathbb T$-model homomorphisms' in $\cal E$ from $M$ to $N$ which classifies the $\mathbb T$-model homomorphisms in slices of $\cal E$ between the localizations of $M$ and $N$ in it.  

In section 5, we establish our main characterization theorem providing necessary and sufficient conditions for a geometric theory to be classified by a presheaf topos. We first state the result abstractly and then proceed to obtain explicit reformulations of each of the conditions. We also derive some corollaries which allow to verify the satisfaction of the conditions of the theorem in specific situations which naturally arise in practice. Lastly, we show that, once recast in the language of indexed colimits of internal diagrams in toposes, the conditions of the characterization theorem for a given geometric theory $\mathbb T$ amount precisely to requirement that every model of $\mathbb T$ in any Grothendieck topos $\cal E$ should be a canonical $\cal E$-indexed colimit of a certain $\cal E$-filtered diagram of `constant' finitely presentable models of $\mathbb T$ which are $\cal E$-finitely presentable.     

In section 6, we introduce the notion of faithful interpretation of geometric theories and investigate to what extent the satisfaction of the conditions of the characterization theorem is preserved by this kind of interpretations. As applications of the general results that we obtain on this topic, we consider in particular the case of quotients of a given geometric theory $\mathbb T$ and that of injectivizations (i.e., theories obtained by adding, for each sort over the signature of the theory, a binary predicate which is provably complemented to the equality relation relative to that sort), providing various sufficient conditions for these theories to be of presheaf type if $\mathbb T$ is. To this end, we carry out in section \ref{fpfgmodels} a general analysis of the relationship between finitely presentable and finitely generated models of a given geometric theory. In section \ref{findingpresheaf}, we treat the problem of finding a geometric theory classified by a given presheaf topos $[{\cal K}, \Set]$, and prove a general theorem ensuring that if the category $\cal K$ can be realized as a full subcategory of the category of finitely presentable models of a theory of presheaf type $\mathbb T$, there exists a quotient of $\mathbb T$ classified by the topos $[{\cal K}, \Set]$, which can be described in most explicit ways in terms of $\mathbb T$ and $\cal K$ provided that some natural conditions are satisfied. We also discuss, in section \ref{rigidtopologies}, the relationship between rigid topologies on the opposite of the category of finitely presentable models of a theory of presheaf type $\mathbb T$ and the presheaf-type quotients of $\mathbb T$.        

In section 7, we investigate expansions of geometric theories from the point of view of the geometric morphisms that they induce between the respective classifying toposes. In particular, we introduce the notion of a localic (resp. hyperconnected) expansion, and show that it naturally corresponds to the notion of localic (resp. hyperconnected) geometric morphism at the level of classifying toposes; as a result, we obtain a logical characterization of the hyperconnected-localic factorization of a geometric morphism. Next, we address the problem of expanding a given geometric theory $\mathbb T$ to a theory classified by a presheaf topos of the form $[{\cal K}, \Set]$, where $\cal K$ is a small category of set-based models of $\mathbb T$, and describe a general method for defining such expansions. 

In section 8, we discuss classical, as well as new, examples of theories of presheaf type from the perspective of the theory developed in the paper. We revisit in particular well-known examples of theories of presheaf type whose finitely presentable models are all finite, notably including the geometric theory of finite sets, and give fully constructive proofs of the fact that Moerdijk's theory of abstract circles and Johnstone's theory of \emph{Diers fields} are of presheaf type. Next, we introduce new examples of theories of presheaf type, including the theory of algebraic extensions of a given field, the theory of locally finite groups, the theory of vector spaces with linear independence predicates and the theory of abelian $\ell$-groups with strong unit. We also show that the injectivization of the algebraic theory of groups is not of presheaf type and explicitly describe a presheaf completion for it.

\section{Indexed colimits in toposes}\label{indexed}

\subsection{Background on indexed categories}\label{backgr}

Before proceeding further, we need to recall some standard notions and facts from the theory of indexed categories; we refer the reader to \cite{El} (especially sections B1.2, B2.3 and B3) and to \cite{PS} for the background.   

Given an internal category $\mathbb C$ in a topos $\cal E$, we denote by ${\mathbb C}_{1}$ its object of arrows, by ${\mathbb C}_{0}$ its object of objects and by $d^{\mathbb C}_{0}, d^{\mathbb C}_{1}:{\mathbb C}_{1}\to {\mathbb C}_{0}$ the domain and codomain arrows.

Given a small category $\cal{C}$ and a topos $\cal{E}$ defined over $\Set$, we can always internalize $\cal{C}$ into $\cal{E}$ by means of the inverse image $\gamma_{\cal{E}}^{\ast}$ of the unique geometric morphism $\gamma_{\cal E}:{\cal E}\to \Set$ from $\cal E$ to $\Set$; the resulting internal category in $\cal{E}$ will be denoted by the symbol $\mathbb{C}$.

Indexed categories will be denoted by underlined letters, with possibly a subscript indicating the indexing category; the fibre at an object $E$ of a $\cal E$-indexed category $\underline{{\cal A}}$ will be denoted by the symbol ${\cal A}_{E}$, and the functor ${\cal A}_{E}\to {\cal A}_{E'}$ corresponding to an arrow $e:E\to E'$ in $\cal E$ will be denoted by the symbol ${\cal A}_{e}$.

A $\cal E$-indexed subcategory ${\underline{\cal B}}_{\cal E}$ of a $\cal E$-indexed category ${\underline{\cal A}}_{\cal E}$ consists, for each object every $E\in {\cal E}$, of a subcategory ${\cal B}_{E}$ of the category ${\cal A}_{E}$ such that for any arrow $\alpha:E'\to E$ in $\cal E$ the functor $ {\cal A}_{\alpha}:{\cal A}_{E}\to {\cal A}_{E'}$ restricts to the subcategories ${\cal B}_{E}$ and ${\cal B}_{E'}$. This clearly defines a $\cal E$-indexed category ${\underline{\cal B}}_{\cal E}$ with a $\cal E$-indexed inclusion ${\underline{\cal B}}_{\cal E} \hookrightarrow {\underline{\cal A}}_{\cal E}$.  

A $\cal E$-indexed functor $F:{\underline{\cal B}}_{\cal E} \to {\underline{\cal A}}_{\cal E}$ is said to be full if for every object $E$ of $\cal E$ the functor $F_{E}:{\cal A}_{E}\to {\cal B}_{E}$ is full. A $\cal E$-indexed subcategory ${\underline{\cal B}}_{\cal E}$ of a $\cal E$-indexed category ${\underline{\cal A}}_{\cal E}$ is said to be a full $\cal E$-indexed category of ${\underline{\cal A}}_{\cal E}$ if the associated $\cal E$-indexed inclusion functor is full.

Every Grothendieck topos $\cal{E}$ gives rise to a $\cal{E}$-indexed category $\underline{\cal{E}}_{\cal E}$ obtained by indexing $\cal{E}$ over itself. 

Recall that if $\mathbb C=(d^{\mathbb C}_{0}, d^{\mathbb C}_{1}:{\mathbb C}_{1}\to {\mathbb C}_{0})$ is an internal category in a topos $\cal E$, a diagram of shape $\mathbb C$ in $\cal E$ is a pair $(f:F\to {\mathbb C}_{0}, \phi:{\mathbb C}_{1}\times_{{\mathbb C}_{0}} F \to F)$ of arrows in $\cal E$ satisfying appropriate conditions, where the pullback ${\mathbb C}_{1}\times_{{\mathbb C}_{0}} F$ is taken relatively to the arrow $d^{\mathbb C}_{0}:{\mathbb C}_{1}\to {\mathbb C}_{0}$ and the arrow $f:F\to {\mathbb C}_{0}$:
\[  
\xymatrix {
{\mathbb C}_{1}\times_{{\mathbb C}_{0}} F  \ar[r]^{\pi_{F}} \ar[d]^{\pi_{1}} & F \ar[d]^{f} \\
{\mathbb C}_{1}  \ar[r]_{d^{\mathbb C}_{0}} & {\mathbb C}_{0}. } 
\]

For a Grothendieck topos $\cal E$ and an internal category $\mathbb C$ in $\cal E$, we have a $\cal{E}$-indexed category $\underline{\fu{\mathbb{C}}{\cal{E}}}$, whose underlying category is the category $\fu{\mathbb{C}}{\cal{E}}$ of diagrams of shape $\mathbb{C}$ in $\cal E$ and morphisms between them. 

Any internal category $\mathbb C$ in $\cal E$ gives naturally rise to a $\cal E$-indexed category, which we call the \emph{$\cal E$-externalization} of $\mathbb C$ and denote by the symbol $\underline{{\mathbb C}}_{\cal E}$. The category $\fu{\mathbb{C}}{\cal{E}}$ is equivalent (naturally in $\cal{E}$) to the category $[\underline{{\mathbb C}}_{\cal E}, \underline{\cal{E}}_{\cal E}]_{\cal E}$ of $\cal{E}$-indexed functors $\st{\underline{{\mathbb C}}_{\cal E}}{\underline{\cal{E}}_{\cal E}}$ and indexed natural transformations between them (by Lemma B2.3.13 in \cite{El}) and also to the category $\fu{\cal{C}}{\cal{E}}$ (by Corollary B2.3.14 in \cite{El}). The equivalence between $\fu{\mathbb{C}}{\cal{E}}$ and $\fu{\cal{C}}{\cal{E}}$ restricts to an equivalence between the full subcategories $\bf{Tors}(\mathbb{C},\cal{E})$ of $\mathbb{C}$-torsors in $\cal{E}$ (as in section B3.2 of \cite{El}) and $\bf Flat(\cal{C}, \cal{E})$ of flat functors $\st{\cal{C}}{\cal{E}}$ (as in chapter VII of \cite{MM}). For any internal functor between internal categories in $\cal E$ or internal diagram $F$ in $\cal E$, we denote the corresponding $\cal E$-indexed functor by the symbol $\underline{F}_{{\cal E}}$. For any functor $F:{\cal C}\to {\cal E}$, we denote the internal diagram in $[{\mathbb C}, {\cal E}]$ corresponding to it under the equivalence $[{\cal C}, {\cal E}]\simeq [{\mathbb C}, {\cal E}]$ by the symbol $\overline{F}$.

The discrete opfibration $p:{\mathbb F} \to {\mathbb C}$ over $\mathbb C$ corresponding to a diagram $(f:F\to {\mathbb C}_{0}, \phi:{\mathbb C}_{1}\times_{{\mathbb C}_{0}} F \to F)$ of shape $\mathbb C$ in $\cal E$ is defined as follows: ${\mathbb F}_{0}=F$, ${\mathbb F}_{1}={\mathbb C}_{1}\times_{{\mathbb C}_{0}} F$, $d_{0}^{\mathbb F}=\pi_{F}:F_{1}={\mathbb C}_{1}\times_{{\mathbb C}_{0}} F \to F=F_{0}$, $d_{1}^{\mathbb F}=\phi :F_{1}={\mathbb C}_{1}\times_{{\mathbb C}_{0}} F \to F=F_{0}$, $p_{0}=f:{\mathbb F}_{0}=F\to {\mathbb C}_{0}$, $p_{1}=\pi_{1}:F_{1}={\mathbb C}_{1}\times_{{\mathbb C}_{0}} F \to {\mathbb C}_{1}$. The discrete opfibration corresponding to a  diagram $F\in [{\mathbb C}, {\cal E}]$ will also be denoted by $\pi^{opf}_{F}:{\int^{opf} F} \to {\mathbb C}$. 

For internal diagrams $G\in [{\mathbb C}^{\textrm{op}}, {\cal E}]$ it is also natural to consider the \emph{discrete fibration} corresponding to $G$, i.e. the opposite functor ${\pi^{opf}_{G}}^{\textrm{op}}:{{\int^{opf} G}}^{\textrm{op}} \to {{\mathbb C}^{\textrm{op}}}^{\textrm{op}}={\mathbb C}$. We shall denote this functor by ${\pi^{f}_{G}}:{{\int^{f} G}} \to {\mathbb C}$.

For any internal category $\mathbb C$ in a topos $\cal E$, the category $[{\mathbb C}, {\cal E}]$ of diagrams of shape $\mathbb{C}$ in $\cal E$ is equivalent to the category $\textbf{DOpf}\slash {\mathbb C}$ of \emph{discrete opfibrations} over $\mathbb C$ (cf. Proposition B2.5.3 \cite{El}); we shall denote this equivalence by 
\[
\tau^{\mathbb C}_{\cal E}:[{\mathbb C}, {\cal E}]\to \textbf{DOpf}\slash {\mathbb C}. 
\] 

A diagram of shape $\mathbb{C}$ in $\cal E$ lies in the subcategory $\bf{Tors}(\mathbb{C},\cal{E})$ of $[{\mathbb C}, {\cal E}]$ if and only if the domain of the corresponding discrete opfibration is a filtered internal category in $\cal E$.   

Any internal functor $H:{\mathbb C} \to {\mathbb D}$ between internal categories $\mathbb C$ and $\mathbb D$ in a topos $\cal E$ induces a functor $[{\mathbb D}, {\cal E}]\to [{\mathbb C}, {\cal E}]$, denoted $D \to D\circ H$, which corresponds, at the level of $\cal E$-indexed categories, to the composition functor with the indexed functor corresponding to $H$, and, at the level of discrete opfibrations associated to the internal diagrams, as the pullback of them along the functor $H$. The latter pullbacks in the category of internal categories in $\cal E$ are computed `pointwise' as pullbacks in $\cal E$, and they are preserved by the dualizing functor ${\mathbb C}\to {\mathbb C}^{\textrm{op}}$.  

Let us recall from \cite{PS} the notion of (indexed) colimit of a $\cal E$-indexed functor, where $\cal E$ is a cartesian category. We shall denote by $!_{I}$ the unique arrow from an object $I$ of $\cal E$ to the terminal object $1$ of $\cal E$. Given a $\cal S$-indexed functor $\Gamma:\underline{X}\to \underline{A}$ and an object $A$ of $A_{1}$, we denote by $\Delta^{\underline{X}}_{\underline{A}}A$ the $\cal E$-indexed functor $\underline{X}\to \underline{\cal A}$ assigning to any $E\in {\cal E}$ the constant functor on $X_{E}$ with value ${\cal A}_{!_{I}}(A)$.  

Let $\underline{X}$ be any indexed category and $\Gamma:\underline{X}\to \underline{A}$ any indexed functor. An indexed cocone $\mu:\Gamma \to A$ consists of an object $A$ in $A_{1}$ together with an indexed natural transformation $\mu: \Gamma \to {\Delta}^{\underline{X}}_{\underline{A}}A$, i.e. for each $I$ in $\cal E$ an ordinary cocone $\mu_{I}: \Gamma_{I} \to \Delta ({\cal A}_{!_{I}}(A))$ such that for each arrow $\alpha:J \to I $ in $\cal E$, $\alpha^{\ast} \circ \mu_{I}=\mu_{J}\circ \alpha^{\ast}$. If $\mu: \Gamma \to {\Delta}^{\underline{X}}_{\underline{A}}A$ is a universal such cone we say that it is a colimit cone over the indexed functor $\Gamma$ with vertex $A$. If, furthermore, for any object $I$ of $\cal E$, the localization $\mu\slash I:\Gamma\slash I \to ({\Delta}^{\underline{X}}_{\underline{A}}A)\slash I= {\Delta}^{\underline{X}\slash I}_{\underline{A}\slash I}({\cal A}_{!_{I}}(A))$ of $\mu$ at $I$ is a colimit cone we say that $\mu$ is the indexed colimit of $\Gamma\slash I$. 

We can describe this notion more explicitly in the particular case of $\cal E$-indexed functors with values in the $\cal E$-indexed category $\underline{{\cal E}}_{\cal E}$. Let us suppose that $\underline{\cal A}_{\cal E}$ is a $\cal E$-indexed category and $D:\underline{\cal A}_{\cal E} \to \underline{{\cal E}}_{\cal E}$ is an indexed functor. A cocone $\mu$ over $D$ consists of an object $U$ of $\cal E$ and, for each object $E$ of $\cal E$, of a cocone $\mu_{E}:D_{E}\to \Delta (!_{E}^{\ast}(U))$ over the diagram $D_{E}:{\cal A}_{E}\to {\cal E}\slash E$ such that for any arrow $\alpha:E' \to E$ in $\cal E$ we have $\alpha^{\ast}(\mu_{E}(c))=\mu_{E'}({\cal A}_{\alpha}(c))$ for all $c\in {\cal A}_{E}$, that is $\alpha^{\ast}\mu_{E}=\mu_{E'}{\cal A}_{\alpha}$ as arrows $D_{E'}({\cal A}_{\alpha}(c))= \alpha^{\ast}(D_{E}(c)) \to !_{E'}(U)=\alpha^{\ast}(!_{E}^{\ast}(U))$ in ${\cal E}\slash E'$, where $\alpha^{\ast}:{\cal E}\slash E \to {\cal E}\slash E'$ is the pullback functor (notice that $D_{E'} \circ {\cal A}_{\alpha} = \alpha^{\ast} \circ D_{E}$ since $D$ is a $\cal E$-indexed functor).   

Note that if $\underline{\cal A}_{\cal E}$ is the $\cal E$-externalization of an internal category $\mathbb C$ in $\cal E$ and $D:\underline{\cal A}_{\cal E} \to \underline{{\cal E}}_{\cal E}$ is an indexed functor corresponding to an internal diagram $D\in [{\mathbb C}, {\cal E}]$, the discrete opfibration associated to the internal diagram in $[!_{E}^{\ast}({\mathbb C}), {\cal E}\slash E]$ corresponding to the localization $D\slash E:\underline{\cal A}_{\cal E}\slash E \to \underline{{\cal E}}_{\cal E}\slash E\cong \underline{{\cal E}\slash E}_{{\cal E}\slash E}$ is the image of the discrete opfibration associated to $D$ under the pullback functor $!_{E}^{\ast}:{\cal E}\to {\cal E}\slash E$ along the unique arrow $!_{E}:E\to 1_{\cal E}$.

The colimit $colim_{\cal E}(D)$ in $\cal E$ of an internal diagram $D\in [{\mathbb C}, {\cal E}]$, where $\mathbb C$ is an internal category in $\cal E$, is defined to be the coequalizer of the two arrows $d_{0}^{\int^{opf} D}, d_{1}^{ \int^{opf} D}$ (recall that $\int^{opf} D$ is the domain of the discrete opfibration over $\mathbb C$ corresponding to the diagram $D$). 

For any internal diagram $G\in [{\mathbb F}, {\cal E}]$, its colimit $colim_{\cal E}(G)$ is isomorphic to the $\cal E$-indexed colimit $\underline{colim}_{{\cal E}}(\underline{G}_{\cal E})$ of the $\cal E$-indexed functor $\underline{G}_{\cal E}:\underline{{\mathbb F}}_{\cal E} \to \underline{\cal{E}}_{\cal E}$ corresponding to it under the equivalence $[{\mathbb F}, {\cal E}]\simeq [\underline{{\mathbb F}}_{\cal E}, \underline{\cal{E}}_{\cal E}]_{\cal E}$. 

If $D\in [{\mathbb C}, {\cal E}]$ is an internal diagram in $\cal E$ with colimit $colim_{\cal E}(D)=coeq(d_{1}^{{\int^{opf} D}}, d_{1}^{{\int^{opf} D}})$, the colimiting $\cal E$-indexed cocone $(colim_{\cal E}(D), \mu)$ of the corresponding $\cal E$-indexed functor $\overline{D}_{\cal E}:\underline{{\mathbb C}}_{\cal E} \to \underline{\cal E}_{\cal E}$ can be described as follows. Let us denote by $c$ the canonical coequalizer arrow $({\int^{opf} D})_{0}\to coeq(d_{0}^{{\int^{opf} D}}, d_{1}^{{\int^{opf} D}})$ in $\cal E$. For any object $E$ of $\cal E$ and any object $x:E\to {\mathbb C}_{0}$ of the category ${{\mathbb C}_{\cal E}}E$, the arrow $\mu_{E}(x):D_{E}(x)=r_{x}\to !_{E}^{\ast}(colim_{\cal E}(D))=coeq(d_{0}^{{\int^{opf} D}}, d_{1}^{{\int^{opf} D}}) \times E$ is equal to $<c\circ z_{x}, r_{x}>$, where the arrow $z_{x}$ is defined by the following pullback diagram:
\[  
\xymatrix {
R_{x} \ar[d]^{r_{x}} \ar[r]^{z_{x}} & ({\int^{opf} D})_{0} \ar[d]^{{(\pi^{opf}_{D})}_{0}}\\
E  \ar[r]_{x} & {\mathbb C}_{0}.} 
\]

\subsection{$\cal E$-filtered indexed categories}\label{filt} 

The following definition will be important in the sequel.

\begin{definition}\label{Efiltered}
Let $\cal E$ be a Grothendieck topos and $\underline{{\cal A}}$ be a $\cal E$-indexed category. We say that $\underline{{\cal A}}$ is \emph{${\cal E}$-filtered} if the following conditions are satisfied:
\begin{enumerate}[(a)]

\item For any object $E$ of $\cal E$ there exists an epimorphic family $\{e_{i}:E_{i} \to E \textrm{ | } i\in I\}$ in $\cal E$ and for each $i\in I$ an object $b_{i}$ of the category ${\cal A}_{E_{i}}$;

\item For any $E\in {\cal E}$ and any objects $a$ and $b$ of the category ${\cal A}_{E}$ there exists an epimorphic family $\{e_{i}:E_{i} \to E \textrm{ | } i\in I\}$ in $\cal E$ and for each $i\in I$ an object $c_{i}$ of the category ${\cal A}_{E_{i}}$ and arrows $f_{i}:{\cal A}_{e_{i}}(a)\to c_{i}$ and $g_{i}:{\cal A}_{e_{i}}(b)\to c_{i}$ in the category ${\cal A}_{E_{i}}$;  

\item For any object $E$ of $\cal E$ and any two parallel arrows $u, v:a\to b$ in the category ${\cal A}_{E}$ there exists an epimorphic family $\{e_{i}:E_{i} \to E \textrm{ | } i\in I\}$ in $\cal E$ and for each $i\in I$ an object $c_{i}$ of the category ${\cal A}_{E_{i}}$ and an arrow $w_{i}:{\cal A}_{e_{i}}(b)\to c_{i}$ in ${\cal A}_{E_{i}}$ such that $w_{i}\circ {\cal A}_{e_{i}}(u)=w_{i}\circ {\cal A}_{e_{i}}(v)$.    
\end{enumerate} 
\end{definition}

The $\cal E$-externalization of any internal filtered category in $\cal E$ is $\cal E$-filtered, but it is not true that if the externalization of an internal category $\mathbb C$ in $\cal E$ is $\cal E$-filtered then $\mathbb C$ is filtered as an internal category in $\cal E$.

A standard example of indexed $\cal E$-filtered categories is provided by indexed categories of elements of flat functors with values in $\cal E$, in the sense of the following definition.

\begin{definition}
Let $P:{\cal C}^{\textrm{op}}\to {\cal E}$ be a functor. The ${\cal E}$-\emph{indexed category of elements} $\underline{\int P}_{\cal E}$ of $P$ assigns to any object $E$ of $\cal E$ the category $\underline{\int P}_{E}$ whose objects are the pairs $(c, x)$ where $c$ is an object of $\cal C$ and $x:E\to F(c)$ is a an arrow in $\cal E$, and whose arrows $(c, x)\to (d, y)$ are the arrows $f:c\to d$ in $\cal C$ such that $F(f)\circ y=x$, and to any arrow $e:E'\to E$ in $\cal E$ the functor $\underline{\int P}_{e}:\underline{\int P}_{E} \to \underline{\int P}_{E'}$ sending any object $(c, x)$ of $\underline{\int P}_{E}$ to the object $(c, x\circ e)$ of $\underline{\int P}_{E'}$ and acting accordingly on the arrows.   
\end{definition}

\begin{proposition}
Let $F:{\cal C}^{\textrm{op}}\to {\cal E}$ be a flat functor. Then the indexed category $\underline{\int F}_{\cal E}$ is $\cal E$-filtered. 
\end{proposition}

\begin{proofs}
Straightforward from the characterization of flat functors as filtering functors given in chapter VII of \cite{MM}.
\end{proofs}

\subsection{Indexation of internal diagrams}

Given an internal diagram $D\in [{\mathbb C}, {\cal E}]$, the corresponding $\cal E$-indexed functor $\overline{D}_{\cal E}:\underline{{\mathbb C}}_{\cal E} \to \underline{\cal E}_{\cal E}$ can be described as follows. For any object $E$ of $\cal E$, $D_{E}:{\mathbb C}_{E}\to {\cal E}\slash E$ sends any object $x:E\to {\mathbb C}_{0}$ of ${\mathbb C}_{E}$ to the object $r_{x}:R_{x}\to E$ of ${\cal E}\slash E$ obtained by pulling $({\pi^{opf}_{F}})_{0}:({{\int^{opf} F}})_{0} \to {\mathbb C}_{0}$ back along $x$, and any arrow $h:E \to {\mathbb C}_{1}$ of ${\mathbb C}_{E}$ from $x:E\to {\mathbb C}_{0}$ to $x':E\to {\mathbb C}_{0}$ to the arrow $D_{E}(h):r_{x} \to r_{x'}$ in ${\cal E}\slash E$ defined as follows. Consider the pullback squares  
\[  
\xymatrix {
S \ar[d]^{s} \ar[r]^{u} & ({\int^{opf} D})_{1} \ar[d]^{({\pi^{opf}_{D}})_{1}}\\
E  \ar[r]_{h} & {\mathbb C}_{1}} 
\]
and   
\[  
\xymatrix {
R_{x'} \ar[d]^{r_{x'}} \ar[r]^{z_{x'}} & ({\int^{opf} D})_{0} \ar[d]^{({\pi^{opf}_{D}})_{0}}\\
E  \ar[r]_{x'} & {\mathbb C}_{0}.} 
\]

Since $d_{1}^{\mathbb C}\circ h=x'$ we have that $({\pi^{opf}_{D}})_{0} \circ  d_{1}^{\int^{opf} D} \circ u =x' \circ s$ and hence by the universal property of the first pullback square, there exists a unique arrow $\beta:S\to R_{x'}$ such that $r_{x'}\circ \beta=s$ and $z_{x'}\circ \beta= d_{1}^{\int^{opf} D} \circ u$. 

As $\pi^{opf}_{F}:{\int^{opf} F} \to {\mathbb C}$ is a discrete opfibration, the diagram
\[  
\xymatrix {
({\int^{opf} D})_{1} \ar[d]^{({\pi^{opf}_{D}})_{1}} \ar[r]^{d_{0}^{{\int^{opf} D}}} & ({\int^{opf} D})_{0} \ar[d]^{({\pi^{opf}_{D}})_{0}}\\
{\mathbb C}_{1}  \ar[r]_{d_{0}^{\mathbb C}} & {\mathbb C}_{0}} 
\] 
is a pullback. `Composing' it with the first pullback square thus yields a pullback square

\[  
\xymatrix {
S \ar[d]^{s} \ar[rr]^{d_{0}^{{\int^{opf} D}} \circ u} & & ({\int^{opf} D})_{0} \ar[d]^{({\pi^{opf}_{D}})_{0}}\\
E  \ar[rr]_{d_{0}^{\mathbb C} \circ h} & & {\mathbb C}_{0}.} 
\]

Now, since $({\pi^{opf}_{D}})_{0}\circ z_{x}=x\circ r_{x}=d_{0}^{\mathbb C} \circ h \circ r_{x}$, the universal property of this pullback square provides a unique arrow $\gamma:R_{x}\to S$ such that $s\circ \gamma=r_{x}$ and $d_{0}^{{\int^{opf} D}} \circ u \circ \gamma=z_{x}$. We define $D_{E}(h):r_{x}\to r_{x'}$ in ${\cal E}\slash E$ to be equal to the composite arrow $\beta \circ \gamma:R_{x} \to R_{x'}$.

\subsection{Colimits and tensor products}\label{seccol}

For any internal category ${\mathbb F}$ in $\cal E$, we denote by $coeq({\mathbb F})$ the coequalizer of the two arrows $d_{0}^{\mathbb F}$ and $d_{1}^{\mathbb F}$. 

From the above discussion it follows that for any $D, P \in [{\mathbb C}, {\cal E}]$, $colim_{\cal E}(P \circ \pi_{D}^{opf})\cong colim_{\cal E}(D \circ \pi_{P}^{opf})$. Indeed, if we consider the pullback square
\[  
\xymatrix {
{\mathbb R} \ar[d]^{q^{D}} \ar[r]^{q^{P}} & {\int^{opf} P} \ar[d]^{\pi^{opf}_{P}}\\
{\int^{opf} D}  \ar[r]_{\pi_{D}^{opf}} & {\mathbb C}} 
\]
in the category ${\textbf{cat}}({\cal E})$ of internal categories in $\cal E$, we have that $colim_{\cal E}(P \circ \pi_{D}^{opf})\cong coeq(dom(q^{D}))=coeq({\mathbb R})=coeq(dom(q^{P}))\cong colim_{\cal E}(D \circ \pi^{opf}_{P})$. 

Similarly, by exploiting the fact that the operation ${\mathbb F}\to coeq({\mathbb F})$ on internal categories in $\cal E$ is invariant under the dualization functor ${\mathbb F} \to {\mathbb F}^{\textrm{op}}$, we obtain another commutation result, which we shall use in the sequel: for any internal diagram $F\in [{\mathbb C}^{\textrm{op}}, {\cal E}]$ and any internal diagram $P\in [{\mathbb C}, {\cal E}]$, we have a natural isomorphism $colim_{\cal E}(P \circ \pi_{F}^{f})\cong colim_{\cal E}(F \circ {\pi_{P}^{f}})$. To prove this, we consider the following pullback squares:
\[  
\xymatrix {
{\mathbb R} \ar[d]^{q^{F}} \ar[r]^{q^{P}} & {\int^{opf} P} \ar[d]^{\pi^{opf}_{P}} & & & {\mathbb S} \ar[d]^{r^{P}} \ar[r]^{r^{F}} & {\int^{opf} F} \ar[d]^{\pi^{opf}_{F}} \\
{\int^{f} F}  \ar[r]_{\pi_{F}^{f}} & {\mathbb C} & & & {\int^{f} P} \ar[r]_{\pi_{P}^{f}} & {\mathbb C}^{\textrm{op}}.} 
\]

We have that 
\[
colim_{\cal E}(P \circ \pi_{F}^{f})\cong coeq(dom(q^{F}))=coeq({\mathbb R}),
\] 
while
\[ 
colim_{\cal E}(F \circ \pi_{P^{f}})=coeq(dom(r^{P}))=coeq({\mathbb S}).
\]
But the fact that the dualization functor preserves pullbacks in ${\textbf{cat}}({\cal E})$ implies that ${\mathbb S}\simeq {\mathbb R}^{\textrm{op}}$, whence $coeq({\mathbb R})\cong coeq({\mathbb S})$, as required.  

Summarizing, we have the following 

\begin{proposition}\label{commut}
Let ${\mathbb C}$ be an internal category in a topos $\cal E$, $P$ and $D$ internal diagrams in $[{\mathbb C}, {\cal E}]$ and $F$ an internal diagram in $[{\mathbb C}^{\textrm{op}}, {\cal E}]$. Then we have natural isomorphisms
\begin{enumerate}[(i)]
\item $colim_{\cal E}(P \circ \pi_{D}^{opf})\cong colim_{\cal E}(D \circ \pi_{P}^{opf})$; 

\item $colim_{\cal E}(P \circ \pi_{F}^{f})\cong colim_{\cal E}(F \circ {\pi_{P}^{f}})$.
\end{enumerate} 
\end{proposition}\qed
   
Let us now proceed to applying this proposition in the context of a functor $F:{\cal C}^{\textrm{op}} \to {\cal E}$, where $\cal C$ is a small category and $\cal E$ is a Grothendieck topos, and a functor $P:{\cal C}\to \Set$. To this end, we explicitly describe the discrete opfibration corresponding to the diagram of shape $\mathbb C$ in $\cal E$, where $\mathbb C$ is the internalization of $\cal C$ in $\cal E$, associated to a functor $G:{\cal C} \to {\cal E}$. We denote by $Ob({\cal C})$ the set of objects of $\cal C$ and by $Arr({\cal C})$ the set of arrows of $\cal C$. 

We have that $F=\mathbin{\mathop{\textrm{ $\coprod$}}\limits_{c\in {\cal C}}}G(c)$ and that $f:F=\mathbin{\mathop{\textrm{ $\coprod$}}\limits_{c\in {\cal C}}}G(c)\to {\mathbb C}_{0}=\mathbin{\mathop{\textrm{$\coprod$}}\limits_{c\in {\cal C}}}1_{\cal E}$ is equal to $\mathbin{\mathop{\textrm{ $\coprod$}}\limits_{c\in {\cal C}}}!_{G(c)}$, where $!_{G(c)}$ is the unique arrow $G(c)\to 1_{\cal E}$ in $\cal E$ (for any $c\in {\cal C}$).  

Let $J_{f}:G(dom(f)) \to \mathbin{\mathop{\textrm{ $\coprod$}}\limits_{f\in Arr({\cal C})}}G(dom(f))$ (resp. $\mu_{c}:G(c) \to \mathbin{\mathop{\textrm{ $\coprod$}}\limits_{c\in Ob({\cal C})}}G(c)$, $\kappa_{f}: 1_{\cal E} \to \mathbin{\mathop{\textrm{ $\coprod$}}\limits_{f\in Arr({\cal C})}}1_{\cal E}$, $\lambda_{c}: 1_{\cal E} \to \mathbin{\mathop{\textrm{ $\coprod$}}\limits_{c\in Ob({\cal C})}}1_{\cal E}$) be the canonical coproduct arrows.

First, let us show that the diagram
\[  
\xymatrix {
\mathbin{\mathop{\textrm{ $\coprod$}}\limits_{f\in Arr({\cal C})}}G(dom(f)) \ar[r]^{d_{0}^{G}} \ar[d]^{\mathbin{\mathop{\textrm{ $\coprod$}}\limits_{f\in Arr({\cal C})}}!_{G(dom(f))}} & \mathbin{\mathop{\textrm{ $\coprod$}}\limits_{c\in Ob({\cal C})}}G(c) \ar[d]^{\mathbin{\mathop{\textrm{ $\coprod$}}\limits_{c\in Ob({\cal C})}}!_{G(c)}} \\
\mathbin{\mathop{\textrm{ $\coprod$}}\limits_{f\in Arr({\cal C})}}1_{\cal E}  \ar[r]_{d_{0}^{\mathbb C}} & \mathbin{\mathop{\textrm{ $\coprod$}}\limits_{c\in Ob({\cal C})}}1_{\cal E}.} 
\]
is a pullback, where the arrow $d_{0}^{G}$ is defined by setting $d_{0}^{G}\circ J_{f}=\mu_{dom(f)}$ (for any $f\in Arr({\cal C})$). We have to prove that, for any object $E$ of $\cal E$ and arrows $\alpha:E \to \mathbin{\mathop{\textrm{ $\coprod$}}\limits_{f\in Arr({\cal C})}}1_{\cal E}$ and $\beta:E \to \mathbin{\mathop{\textrm{ $\coprod$}}\limits_{c\in Ob({\cal C})}}G(c)$ such that $\mathbin{\mathop{\textrm{ $\coprod$}}\limits_{c\in Ob({\cal C})}}!_{G(c)} \circ \beta= d_{0}^{\mathbb C} \circ \alpha$, there exists a unique arrow $\gamma:E \to \mathbin{\mathop{\textrm{ $\coprod$}}\limits_{f\in Arr({\cal C})}}!_{G(dom(f))}$ such that $\alpha=\mathbin{\mathop{\textrm{ $\coprod$}}\limits_{f\in Arr({\cal C})}}G(dom(f)) \circ \gamma$ and $\beta=d_{0}^{G} \circ \gamma$. To this end, consider, for any $c\in Ob({\cal C})$ and $f\in Arr({\cal C})$, the commutative diagram
\[  
\xymatrix {
E_{(f, c)} \ar[d]^{p_{c}} \ar[r]^{q_{f}} & E_{f} \ar[d]^{y_{f}} \ar[r]^{\alpha_{f}} & 1_{{\cal E}} \ar[d]^{\kappa_{f}} \ar[dr]^{id} & \\
E_{c} \ar[d]^{\beta_{c}} \ar[r]^{z_{c}} & E \ar[r]^{\alpha} \ar[d]^{\beta} & \mathbin{\mathop{\textrm{ $\coprod$}}\limits_{f\in Arr({\cal C})}} 1_{{\cal E}} \ar[d]^{d_{0}^{\mathbb C}} & 1_{{\cal E}} \ar[dl]^{\lambda_{dom(f)}} \\
G(c) \ar[dr]_{!_{G(c)}} \ar[r]^{\kappa_{c}} & \mathbin{\mathop{\textrm{ $\coprod$}}\limits_{c\in Ob({\cal C})}}G(c) \ar[r]^{\mathbin{\mathop{\textrm{ $\coprod$}}\limits_{c\in Ob({\cal C})}}!_{G(c)}} & \mathbin{\mathop{\textrm{ $\coprod$}}\limits_{c\in Ob({\cal C})}}1_{\cal E}   \\
& 1_{\cal E} \ar[ur]_{\lambda_{c}}, } 
\]
where all the squares except for the lower-right one are pullbacks. The commutativity of the diagram, combined with the fact that distinct coproduct arrows are disjoint from each other, immediately implies that for any pair $(f, c)$ such that $c\neq dom(f)$, we have $E_{(f, c)}\cong 0_{\cal E}$. On the other hand, the stability of coproducts under pullbacks implies that $E\cong \mathbin{\mathop{\textrm{ $\coprod$}}\limits_{c\in Ob({\cal C})}}E_{c}$, and $E\cong \mathbin{\mathop{\textrm{ $\coprod$}}\limits_{f\in Arr({\cal C})}}E_{f}$; whence $E\cong \mathbin{\mathop{\textrm{ $\coprod$}}\limits_{(f, c)\in Arr({\cal C}) \times Ob({\cal C})}}E_{(f, c)}\cong \mathbin{\mathop{\textrm{ $\coprod$}}\limits_{(f, c) \textrm{ | } dom(f)=c}}E_{(f, c)}$, with canonical coproduct arrows $\xi_{(f,c)}=y_{f}\circ q_{f}=z_{c}\circ p_{c}:E_{(f, c)}\to E$. We define, for each pair $(f, c)$ such that $c=dom(f)$, the arrow $\gamma_{(f,c)}:E_{(f, c)} \to \mathbin{\mathop{\textrm{ $\coprod$}}\limits_{f\in Arr({\cal C})}}G(dom(f))$ as the composite $J_{f}\circ \beta_{c}\circ p_{c}$, and set $\gamma$ equal to the arrow $\mathbin{\mathop{\textrm{ $\coprod$}}\limits_{(f, c) \textrm{ | } dom(f)=c}}\gamma_{(f, c)}: E \to \mathbin{\mathop{\textrm{ $\coprod$}}\limits_{f\in Arr({\cal C})}}G(dom(f))$. We have to verify that $\alpha= \mathbin{\mathop{\textrm{ $\coprod$}}\limits_{f\in Arr({\cal C})}}!_{G(dom(f))} \circ \gamma$ and $\beta=d_{0}^{G} \circ \gamma$ or, equivalently, that for any pair $(f, c)$ such that $dom(f)=c$, we have: 
\begin{enumerate}[(1)]
\item $\alpha \circ \xi_{(f, c)}= \mathbin{\mathop{\textrm{ $\coprod$}}\limits_{f\in Arr({\cal C})}}!_{G(dom(f))} \circ \gamma_{(f, c)}$ and

\item $\beta \circ \xi_{(f, c)}=d_{0}^{G} \circ \gamma_{(f, c)}$.
\end{enumerate}  

To prove $(1)$, we preliminarily show that for any pair $(f, c)$ such that $c=dom(f)$, we have $\alpha_{f}\circ q_{f}=!_{G(c)}\circ \beta_{c}\circ p_{c}$. Since the arrow $\lambda_{c}=\lambda_{dom(f)}:1_{\cal E}\to \mathbin{\mathop{\textrm{ $\coprod$}}\limits_{c\in Ob({\cal C})}}1_{\cal E}$ is monic, it is equivalent to prove that $\lambda_{dom(f)}\circ \alpha_{f}\circ q_{f}=\lambda_{c}\circ !_{G(c)}\circ \beta_{c}\circ p_{c}$. But the commutativity of the above diagram yields $\lambda_{dom(f)}\circ \alpha_{f}\circ q_{f}=d_{0}^{\mathbb C} \circ \kappa_{f} \circ \alpha_{f}\circ q_{f}=d_{0}^{\mathbb C} \circ \alpha \circ y_{f} \circ q_{f}=\mathbin{\mathop{\textrm{ $\coprod$}}\limits_{c\in Ob({\cal C})}}!_{G(c)} \circ \beta \circ y_{f} \circ q_{f}=\mathbin{\mathop{\textrm{ $\coprod$}}\limits_{c\in Ob({\cal C})}}!_{G(c)} \circ \beta \circ z_{c} \circ p_{c}=\mathbin{\mathop{\textrm{ $\coprod$}}\limits_{c\in Ob({\cal C})}}!_{G(c)} \circ \mu_{c} \circ \beta_{c}\circ p_{c}=\lambda_{c}\circ !_{G(c)} \circ \beta_{c} \circ p_{c}$, as required.

We thus have $\alpha\circ \xi_{(f, c)}=\alpha \circ y_{f}\circ q_{f}=\kappa_{f}\circ \alpha_{f}\circ q_{f}=\kappa_{f}\circ !_{G(c)} \circ \beta_{c}\circ p_{c}= \mathbin{\mathop{\textrm{ $\coprod$}}\limits_{f\in Arr({\cal C})}}!_{G(dom(f))} \circ \gamma_{(f, c)}$. This proves condition $(1)$. 

Further, $\beta \circ \xi_{(f, c)}= \beta \circ z_{c} \circ p_{c}=\mu_{c} \circ \beta_{c}\circ p_{c}=\mu_{dom(f)}\circ \beta_{c}\circ p_{c}=d_{0}^{G}\circ J_{f}\circ \beta_{c}\circ p_{c}=d_{0}^{G}\circ \gamma_{(f, c)}$. This proves condition $(2)$.

Let us define the arrow
\[
d_{1}^{G}:\mathbin{\mathop{\textrm{ $\coprod$}}\limits_{f\in Arr({\cal C})}}G(dom(f)) \to \mathbin{\mathop{\textrm{ $\coprod$}}\limits_{c\in Ob({\cal C})}}G(c)
\]
by setting, for each $f\in Arr({\cal C})$, $d_{1}^{G}\circ J_{f}=\mu_{cod(f)}\circ G(f)$.

The discrete opfibration $p:{\mathbb F}\to {\mathbb C}$ corresponding to $G$ can be described as follows: ${\mathbb F}_{0}=\mathbin{\mathop{\textrm{ $\coprod$}}\limits_{c\in Ob({\cal C})}}G(c)$, ${\mathbb F}_{1}=\mathbin{\mathop{\textrm{ $\coprod$}}\limits_{f\in Arr({\cal C})}}G(dom(f))$, the domain and codomain arrows $d_{0}^{\mathbb F}, d_{1}^{\mathbb F}:{\mathbb F}_{1}\to {\mathbb F}_{0}$ are respectively equal to $d_{0}^{G}$ and to $d_{1}^{G}$, $p_{0}:{\mathbb F}_{0}\to {\mathbb C}_{0}$ is equal to $\mathbin{\mathop{\textrm{ $\coprod$}}\limits_{c\in Ob({\cal C})}}!_{G(c)}$ and $p_{1}:{\mathbb F}_{1}\to {\mathbb C}_{1}$ is equal to $\mathbin{\mathop{\textrm{ $\coprod$}}\limits_{f\in Arr({\cal C})}}!_{G(dom(f))}$. The composition law in the internal category $\mathbb F$ is defined in the obvious way. 

We leave to the reader the straightforward task of verifying that this is indeed the discrete opfibration corresponding to the functor $G$ via the composite of the equivalence $[{\cal C}, {\cal E}]\simeq [{\mathbb C}, {\cal E}]$ with the equivalence 
\[
\tau^{\mathbb C}_{\cal E}:[{\mathbb C}, {\cal E}]\to \textbf{DOpf}\slash {\mathbb C}. 
\] 

Recalling that, for any functor $F:{\mathbb C}^{\textrm{op}}\to {\cal E}$, the discrete fibration $\pi_{F}^{f}:{\int^{f} F}\to {\mathbb C}$ associated to it is equal to ${\pi_{F}^{opf}}^{\textrm{op}}:{{\int^{opf} F}^{\textrm{op}}} \to {{\mathbb C}^{\textrm{op}}}^{\textrm{op}}={\mathbb C}$, we deduce the following explicit description of the discrete fibration $\pi^{f}_{F}:{\int^{f} F} \to {\mathbb C}$ associated to a functor $F:{\cal C}^{\textrm{op}}\to {\cal E}$: ${({\int^{f} F})}_{0}=\mathbin{\mathop{\textrm{ $\coprod$}}\limits_{c\in Ob({\cal C})}}F(c)$, ${({\int^{f} F})}_{1}=\mathbin{\mathop{\textrm{ $\coprod$}}\limits_{f\in Arr({\cal C})}}F(cod(f))$, the domain and codomain arrows $d_{0}^{F}, d_{1}^{F}:{({\int^{f} F})}_{1}\to {({\int^{f} F})}_{0}$ are defined by the conditions $d_{0}^{F}\circ J_{f}=\mu_{cod(f)}$ and $d_{1}^{F}\circ J_{f}=\mu_{dom(f)}\circ F(f)$ for all $f\in Arr({{\cal C}})$ (where $\mu_{c}:F(c) \to \mathbin{\mathop{\textrm{ $\coprod$}}\limits_{c\in Ob({\cal C})}}F(c)$ and $J_{f}:F(cod(f)) \to \mathbin{\mathop{\textrm{ $\coprod$}}\limits_{f\in Arr({\cal C})}}F(cod(f))$ are the canonical coproduct arrows), $({\pi_{F}^{f}})_{0}:{({\int^{f} F})}_{0} \to {\mathbb C}_{0}$ is equal to $\mathbin{\mathop{\textrm{ $\coprod$}}\limits_{c\in Ob({\cal C})}}!_{F(c)}$ and $({\pi_{F}^{f}})_{1}:{({\int^{f} F})}_{1}\to {\mathbb C}_{1}$ is equal to $\mathbin{\mathop{\textrm{ $\coprod$}}\limits_{f\in Arr({\cal C})}}!_{F(cod(f))}$. The composition law in the internal category ${\int^{f} F}$ is defined in the obvious way.

\begin{theorem}\label{tensorpr}
Let $F:{\cal C}^{\textrm{op}}\to {\cal E}$ be a functor from the opposite of a small category $\cal C$ to a Grothendieck topos $\cal E$ and $P:{\cal C}\to \Set$ be a functor. Then the following three objects are naturally isomorphic:

\begin{enumerate}[(i)]
\item $colim(F\circ \pi^{f}_{P})$

\item $colim_{\cal E}(\overline{F}\circ \pi_{P_{\cal E}}^{f})\cong colim_{\cal E}(P_{\cal E}\circ \pi_{\overline{F}}^{f})$ (cf. Proposition \ref{commut})

\item $\underline{colim}_{\cal E}(\underline{P_{\cal E}}_{\cal E}\circ \underline{\pi_{\overline{F}}^{f}}_{\cal E})$,
\end{enumerate}
where $P_{\cal E}$ is the internal diagram in $[{\mathbb C}, {\cal E}]$ given by $\overline{\gamma_{\cal E}^{\ast} \circ P}$.
\end{theorem}

\begin{proofs}
The isomorphism between $colim_{\cal E}(P_{\cal E}\circ \pi_{\overline{F}}^{f})$ and $\underline{colim}_{\cal E}(\underline{P_{\cal E}}_{\cal E}\circ \underline{\pi_{\overline{F}}^{f}}_{\cal E})$ follows from the general fact that for any internal diagram $G\in [{\mathbb D}, {\cal E}]$ its colimit $colim_{\cal E}(G)$ is isomorphic to the $\cal E$-indexed colimit $\underline{colim}_{{\cal E}}(\underline{G}_{\cal E})$. It thus remains to prove the isomorphism between $colim(F\circ \pi^{f}_{P})$ and $colim_{\cal E}(\overline{F}\circ \pi_{P_{\cal E}}^{f})$. To this end, we recall the following three general facts:
\begin{enumerate}[(1)]

\item For any functor $G:{\cal D}\to {\cal E}$, its colimit $colim(G)$ is naturally isomorphic to the colimit $colim_{\cal E}(\overline{G})$;  

\item For any functor $H:{\cal D}\to {\cal C}$ between small categories $\cal C$ and $\cal D$ and any functor $M:{\cal C}\to {\cal E}$, we have a natural isomorphism $\overline{F\circ H}\cong \overline{F}\circ \gamma_{\cal E}^{\ast}(H)$; 

\item For any geometric morphism $f:{\cal F}\to {\cal E}$, the diagram
\[  
\xymatrix {
[{\mathbb C}, {\cal E}] \ar[d]^{f^{\ast}(-)} \ar[r]^{\tau^{\mathbb C}_{\cal E}   } & \textbf{DOpf}\slash {\mathbb C}  \ar[d]^{f^{\ast}(-)} \\
[f^{\ast}({\mathbb C}), {\cal E}]  \ar[r]^{\tau^{f^{\ast}({\mathbb C})}_{\cal F}  } & \textbf{DOpf}\slash f^{\ast}({\mathbb C})} 
\]
commutes.
\end{enumerate} 

We therefore have that
\[
colim(F\circ \pi^{f}_{P})\cong colim_{\cal E}(\overline{F\circ \pi_{P}^{f}})\cong colim_{\cal E}(\overline{F}\circ \gamma_{\cal E}^{\ast}(\pi^{f}_{P}))\cong colim_{\cal E}(\overline{F}\circ \pi_{P_{\cal E}}^{f}),
\]
where the first isomorphism follows from $(1)$, the second from $(2)$ and the third from $(3)$. 
\end{proofs}

We shall indicate the three isomorphic objects of the theorem by the symbol $F\otimes_{\cal C} P$. 

The following lemma is essentially contained in the proof of Giraud's theorem (cf. \cite{Giraud} or, for instance, the Appendix of \cite{MM}), but we were not able to find its exact statement in the literature; we thus provide a proof of it for the reader's convenience.   

\begin{lemma}
Let $\cal E$ be a Grothendieck topos and $\{e_{i}:E_{i}\to E \textrm{ | } i\in I\}$ an epimorphic family in $\cal E$. Then the arrow $\mathbin{\mathop{\textrm{ $\coprod$}}\limits_{i\in I}e_{i}}:\mathbin{\mathop{\textrm{ $\coprod$}}\limits_{i\in I}}E_{i} \to E$ yields an isomorphism $(\mathbin{\mathop{\textrm{ $\coprod$}}\limits_{i\in I}}E_{i})\slash R\cong E$, where $R$ is the equivalence relation in $\cal E$ on the object $\mathbin{\mathop{\textrm{ $\coprod$}}\limits_{i\in I}}E_{i}$ given by the subobject $\mathbin{\mathop{\textrm{ $\coprod$}}\limits_{(i, j)\in I\times I}}E_{i}\times_{E} E_{j}\mono \mathbin{\mathop{\textrm{ $\coprod$}}\limits_{(i, j)\in I\times I}}E_{i}\times E_{j}\cong \mathbin{\mathop{\textrm{ $\coprod$}}\limits_{i\in I}}E_{i} \times \mathbin{\mathop{\textrm{ $\coprod$}}\limits_{j\in I}}E_{j}$; in particular, the arrow $\mathbin{\mathop{\textrm{ $\coprod$}}\limits_{i\in I}}e_{i}: \mathbin{\mathop{\textrm{ $\coprod$}}\limits_{i\in I}}E_{i} \to E$ is the coequalizer in $\cal E$ of the two canonical arrows $\mathbin{\mathop{\textrm{ $\coprod$}}\limits_{(i, j)\in I\times I}}E_{i}\times_{E} E_{j} \to \mathbin{\mathop{\textrm{ $\coprod$}}\limits_{i\in I}}E_{i}$.    
\end{lemma}

\begin{proofs}
Let $R$ be the kernel pair of the epimorphism $\mathbin{\mathop{\textrm{ $\coprod$}}\limits_{i\in I}e_{i}}$, that is the pullback of this arrow along itself; then, by the well-known exactness properties of Grothendieck toposes, $R$ is an equivalence relation on $\mathbin{\mathop{\textrm{ $\coprod$}}\limits_{i\in I}E_{i}}$ such that the coequalizer in $\cal E$ of the two associated projections is isomorphic to $q$. Now, the fact that pullbacks preserve coproducts in a Grothendieck topos implies that $R$ is isomorphic to the subobject $\mathbin{\mathop{\textrm{ $\coprod$}}\limits_{(i, j)\in I\times I}}E_{i}\times_{E} E_{j}\mono \mathbin{\mathop{\textrm{ $\coprod$}}\limits_{(i, j)\in I\times I}}E_{i}\times E_{j}\cong \mathbin{\mathop{\textrm{ $\coprod$}}\limits_{i\in I}}E_{i} \times \mathbin{\mathop{\textrm{ $\coprod$}}\limits_{j\in I}}E_{j}$, as required.
\end{proofs}

\begin{remark}
The lemma admits an obvious generalization to arbitrary separating sets for the topos $\cal E$ (cf. the Appendix of \cite{MM}).
\end{remark}

\begin{corollary}\label{cordesc}
Let $a:A\to E$ and $l:L\to E$ be objects of the topos ${\cal E}\slash E$, and $\{e_{i}:E_{i}\to E \textrm{ | } i\in I\}$ an epimorphic family in $\cal E$. Then a family of arrows $\{f_{i}:e_{i}^{\ast}(a)\to e_{i}^{\ast}(l) \textrm{ | } i\in I\}$ in the toposes ${\cal E}\slash E_{i}$ defines a (unique) arrow $f:a\to l$ in ${\cal E}\slash E$ such that $e_{i}^{\ast}(f)=f_{i}$ for all $i\in I$ if and only if for every $i,j \in I$, $q_{i}^{\ast}(f_{i})=q_{j}^{\ast}(f_{j})$, where the arrows $q_{i}$ and $q_{j}$ are defined by the following pullback square:
\[  
\xymatrix {
E_{i, j} \ar[d]^{q_{j}} \ar[rd]^{e_{i,j}} \ar[r]^{q_{i}} & E_{i} \ar[d]^{e_{i}} \\
E_{j} \ar[r]^{e_{j}} & E.} 
\]
\end{corollary}\qed

\subsection{$\cal E$-final subcategories}\label{finsub}

The following definition will be important in the sequel. We shall borrow the notation from section \ref{seccol}.

\begin{definition}\label{fulldefinition}
Let ${\underline{\cal A}}_{\cal E}$ be a $\cal E$-indexed category and $i:{\underline{\cal B}}_{\cal E} \to {\underline{\cal A}}_{\cal E}$ be a $\cal E$-indexed functor. 

\begin{enumerate}[(a)]

\item We say that $i$ is \emph{$\cal E$-final} if for every $E\in {\cal E}$ and $x\in {\cal A}_{E}$ there exists a non-empty set ${\cal E}_{E}^{x}$ of triplets $\{(e_{i}, b_{i}, f_{i}) \textrm{ | } i\in I\}$ such that the family $\{e_{i}:E_{i} \to E \textrm{ | } i\in I\}$ is epimorphic, $b_{i}$ is an object of ${\cal B}_{E_{i}}$ and $f_{i}:{\cal A}_{e_{i}}(x) \to i_{E_{i}}(b_{i})$ is an arrow in ${\cal A}_{E_{i}}$ (for each $i\in I$) with the property that for any triplets $\{(e_{i}, b_{i}, f_{i}) \textrm{ | } i\in I\}$ and $\{(e_{j}', c_{j}, f_{j}') \textrm{ | } j\in J\}$ in ${\cal E}_{E}^{x}$ there exists an epimorphic family $\{g^{i, j}_{k}:E^{i, j}_{k}\to E_{i, j} \textrm{ | } k\in K_{i, j}\}$ and for each $k\in K_{i, j}$ an object $d^{i, j}_{k}$ of ${\cal B}_{E^{i, j}_{k}}$ and arrows $r_{k}^{i, j}:{\cal B}_{q_{i}\circ g^{i, j}_{k}}(b_{i})\to d^{i, j}_{k}$ and $s_{k}^{i, j}:{\cal B}_{q_{j}\circ g^{i, j}_{k}}(c_{j})\to d^{i, j}_{k}$ in ${\cal B}_{E^{i, j}_{k}}$ such that $i_{E^{i, j}_{k}}(r_{k}^{i, j}) \circ {\cal A}_{q_{i}}(f_{i})=i_{E^{i, j}_{k}}(s_{k}^{i, j}) \circ {\cal A}_{q_{j}}(f_{j}')$. 

\item We say that a $\cal E$-indexed subcategory ${\underline{\cal B}}_{\cal E}$ of a $\cal E$-indexed category ${\underline{\cal A}}_{\cal E}$ is a \emph{$\cal E$-final subcategory} of ${\underline{\cal A}}_{\cal E}$ if the canonical $\cal E$-indexed inclusion functor $i:{\underline{\cal B}}_{\cal E} \hookrightarrow {\underline{\cal A}}_{\cal E}$ is $\cal E$-final, in other words if for every $E\in {\cal E}$ and $x\in {\cal A}_{E}$ there exists a non-empty set ${\cal E}_{E}^{x}$ of triplets $\{(e_{i}, b_{i}, f_{i}) \textrm{ | } i\in I\}$ such that the family $\{e_{i}:E_{i} \to E \textrm{ | } i\in I\}$ is epimorphic, $b_{i}$ is an object of ${\cal B}_{E_{i}}$ and $f_{i}:{\cal A}_{e_{i}}(x) \to b_{i}$ is an arrow in ${\cal A}_{E_{i}}$ (for each $i\in I$) with the property that for any triplets $\{(e_{i}, b_{i}, f_{i}) \textrm{ | } i\in I\}$ and $\{(e_{j}', c_{j}, f_{j}') \textrm{ | } j\in J\}$ in ${\cal E}_{E}^{x}$ there exists an epimorphic family $\{g^{i, j}_{k}:E^{i, j}_{k}\to E_{i, j} \textrm{ | } k\in K_{i, j}\}$ and for each $k\in K_{i, j}$ an object $d^{i, j}_{k}$ of ${\cal B}_{E^{i, j}_{k}}$ and arrows $r_{k}^{i, j}:{\cal B}_{q_{i}\circ g^{i, j}_{k}}(b_{i})\to d^{i, j}_{k}$ and $s_{k}^{i, j}:{\cal B}_{q_{j}\circ g^{i, j}_{k}}(c_{j})\to d^{i, j}_{k}$ in ${\cal B}_{E^{i, j}_{k}}$ such that $r_{k}^{i, j} \circ {\cal A}_{q_{i}}(f_{i})=s_{k}^{i, j} \circ {\cal A}_{q_{j}}(f_{j}')$. 

\item We say that a $\cal E$-indexed subcategory ${\underline{\cal B}}_{\cal E}$ of a $\cal E$-indexed category ${\underline{\cal A}}_{\cal E}$ is a \emph{$\cal E$-strictly final subcategory} of ${\underline{\cal A}}_{\cal E}$ if for any arrow $f:a\to b$ in ${\cal A}_{E}$ there exists an epimorphic family $\{e_{i}:E_{i}\to E \textrm{ | } i\in I\}$ in $\cal E$ such that for any $i\in I$, the arrow ${\cal A}_{e_{i}}(f):{\cal A}_{e_{i}}(a)\to {\cal A}_{e_{i}}(b)$ lies in ${\cal B}_{E_{i}}$.      

\item We say that a $\cal E$-indexed subcategory ${\underline{\cal B}}_{\cal E}$ of a $\cal E$-indexed category ${\underline{\cal A}}_{\cal E}$ is \emph{$\cal E$-full} if for every arrow $f:b\to b'$ in ${\cal A}_{E}$, where $b$ and $b'$ are objects of ${\cal B}_{E}$, there exists an epimorphic family $\{e_{i}:E_{i}\to E \textrm{ | } i\in I\}$ such that for any $i\in I$ the arrow ${\cal A}_{e_{i}}(f)$ lies in ${\cal B}_{E_{i}}$.
\end{enumerate}

\end{definition}

\begin{remarks}\label{fullfinal}

\begin{enumerate}[(a)]
\item Let ${\underline{\cal B}}_{\cal E}$ be a $\cal E$-indexed subcategory of ${\underline{\cal A}}_{\cal E}$ with the property that for every object $E$ of $\cal E$ and any object $a\in {\cal A}_{E}$ $\alpha:E'\to E$, there exists an epimorphic family $\{e_{i}:E_{i}\to E \textrm{ | } i\in I\}$ in $\cal E$ such that for any $i\in I$, the object ${\cal A}_{e_{i}}(a)$ lies in ${\cal B}_{E_{i}}$ and for any arrow $f:a\to b$ in ${\cal A}_{E}$ where $a, b\in {\cal B}_{E}$, there exists an epimorphic family $\{e_{i}:E_{i}\to E \textrm{ | } i\in I\}$ in $\cal E$ such that for any $i\in I$, the arrow ${\cal A}_{e_{i}}(f):{\cal A}_{e_{i}}(a)\to {\cal A}_{e_{i}}(b)$ lies in ${\cal B}_{E_{i}}$. Then ${\underline{\cal B}}_{\cal E}$ is a $\cal E$-strictly final subcategory of ${\underline{\cal A}}_{\cal E}$. 

\item Every $\cal E$-strictly final subcategory is a $\cal E$-final subcategory.  

\item If $i$ is a $\cal E$-full embedding of a $\cal E$-indexed category ${\underline{\cal B}}_{\cal E}$ into a $\cal E$-filtered $\cal E$-indexed category ${\underline{\cal A}}_{\cal E}$ then $i$ is $\cal E$-final if and only if for every $E\in {\cal E}$ and $x\in {\cal A}_{E}$ there exists a non-empty epimorphic family $\{e_{i}:E_{i} \to E \textrm{ | } i\in I\}$ in $\cal E$ and for each $i\in I$ an object $b_{i}$ of ${\cal B}_{E_{i}}$ and an arrow $f_{i}:{\cal A}_{e_{i}}(x) \to b_{i}$ of ${\cal A}_{E_{i}}$.
\end{enumerate}
\end{remarks}

\begin{proposition}\label{finalfiltered}
Let ${\underline{\cal B}}_{\cal E}$ be a $\cal E$-full $\cal E$-indexed subcategory of a $\cal E$-indexed category ${\underline{\cal A}}_{\cal E}$. Then ${\underline{\cal B}}_{\cal E}$ is $\cal E$-filtered if and only if ${\underline{\cal A}}_{\cal E}$ is $\cal E$-filtered.     
\end{proposition}

\begin{proofs}
The proof is entirely analogous to the classical one and left to the reader.  
\end{proofs}

The following result represents a $\cal E$-indexed version of the classical theorem formalizing the behaviour of colimits with respect to final subcategories.

\begin{theorem}\label{fincol}
Let $i:{\underline{\cal B}}_{\cal E}\hookrightarrow {\underline{\cal A}}_{\cal E}$ be a $\cal E$-final functor and $D:{\underline{\cal A}}_{\cal E} \to \underline{\cal E}_{\cal E}$ a $\cal E$-indexed functor. Then $D$ admits a colimit (resp. a $\cal E$-indexed colimit) if and only if $D\circ i$ admits a colimit (resp. a $\cal E$-indexed colimit), and the two colimits are equal.
\end{theorem}

\begin{proofs}
First, let us show that any cocone $\lambda$ over $D\circ i$ with vertex $V$ can be (uniquely) extended to a cocone $\tilde{\lambda}$ over $D$. For any object $a$ of ${\cal A}_{E}$, we have to define an arrow $\tilde{\lambda}_{E}(a):D_{E}(a)\to !_{E}^{\ast}(V)$ in ${\cal E}\slash E$. By our hypotheses, there exists an epimorphic family ${\cal E}=\{e_{i}:E_{i}\to E \textrm{ | } i\in I\}$ in $\cal E$ and a family of arrows $\{f_{i}:{\cal A}_{e_{i}}(a) \to i_{E_{i}}(b_{i}) \textrm{ | } i\in I\}$, where $b_{i}$ lies in ${\cal B}_{E_{i}}$. Consider, for each $i\in I$, the arrow $\lambda_{E_{i}}(b_{i}) \circ D_{E_{i}}(f_{i}):D_{E_{i}}({\cal A}_{e_{i}}(a))=e_{i}^{\ast}(D_{E}(a)) \to !_{E_{i}}(V)=e_{i}^{\ast}(!_{E}^{\ast}(V))$. By applying the condition in the definition of $\cal E$-final functor to the pair $((e_{i}, f_{i}, b_{i}), (e_{i}, f_{i}, b_{i}))$ and exploiting the fact that $\lambda$ is a cocone over $D\circ i$, we obtain that ${q_{i}}^{\ast}(\lambda_{E_{i}}(b_{i}) \circ D_{E_{i}}(f_{i}))={q_{j}}^{\ast}(\lambda_{E_{j}}(b_{j}) \circ D_{E_{j}}(f_{j}))$ (here we use the notation of Corollary \ref{cordesc}). By Corollary \ref{cordesc}, the family of arrows $\lambda_{E_{i}}({\cal A}_{e_{i}}(a))$ (for $i\in I$) thus induces a unique arrow $u_{{\cal E}}=D_{E}(a)\to !_{E}^{\ast}(V)$. To be able to set $\tilde{\lambda}_{E}(a)$ equal to this arrow, we have to show that such definition does not depend on the choice of the family $\{(e_{i}, f_{i})\}$. But this follows similarly as above, by applying the condition in the definition of $\cal E$-final functor and invoking the fact that $\lambda$ is a cocone over $D\circ i$. 

Notice that if $a\in {\cal B}_{E}$ then $\tilde{\lambda}_{E}(a)=\lambda_{E}(a)$, and that for any cocone $\xi$ over $D$, $\xi=\tilde{\xi \circ i}$.
      
Now that we have showed that for any $a\in {\cal A}_{E}$ the definition of the arrow $\tilde{\lambda}_{E}(a):D_{E}(a)\to !_{E}^{\ast}(V)$ in the topos ${\cal E}\slash E$ is well-posed, it remains to prove that the assignment $a \to \tilde{\lambda}_{E}(a)$ defines a cocone over the diagram $D_{E}$ with vertex $!_{E}^{\ast}(V)$, i.e. that for any arrow $f:a\to b$ in ${\cal A}_{E}$, $\tilde{\lambda}_{E}(b)\circ D_{E}(f)=\tilde{\lambda}_{E}(a)$ as arrows in ${\cal E}\slash E$. Further, we have to show that the assignment $E \to \tilde{\lambda}_{E}$ defines a $\cal E$-indexed cone on the $\cal E$-indexed functor $D$, i.e. that for any arrow $\alpha:E' \to E$ in $\cal E$ and any $a\in {\cal A}_{E}$ we have $\alpha^{\ast}(\tilde{\lambda}_{E}(a))=\tilde{\lambda}_{E'}({\cal A}_{\alpha}(a))$. This can be easily deduced from the fact that $\lambda$ is an indexed cocone over $D\circ i$. 

The second part of the theorem, for indexed colimits, follows from the first part by noticing that any localization of an indexed $\cal E$-final functor with respect to a localization functor ${\cal E}\to {\cal E}\slash E$ is a ${\cal E}\slash E$-final functor. So it will be sufficient to show the first part. 

To complete the proof of the theorem, it remains to show that for any cocone $(U, \mu)$ over $D$, $(U, \mu)$ is colimiting over $D$ if and only if $(U, \mu\circ i)$ is colimiting over $D$.

Suppose that $D$ has a colimiting cocone $\mu$ with vertex $U$. We want to prove that $\mu\circ i$ is a colimiting cocone for the functor $D\circ i$.  
Let $\lambda$ be a cocone over the diagram $D\circ i$ with vertex $V$; then, as we have just proved, $\tilde{\lambda}$ is a cocone over $D$ with vertex $V$; therefore $(U, \mu)$ factors through $(V, \tilde{\lambda})$, by an arrow $z:U\to V$ in $\cal E$. Clearly, $z$ yields in particular a factorization of $(U, \mu\circ i)$ across $(V, \lambda)$. The uniqueness of the factorization of $(U, \mu\circ i)$ across $(V, \lambda)$ follows from the fact that for any two cocones $\xi$ and $\chi$ over $D\circ i$, if $\xi$ factors through $\chi$ by an arrow $w$ then $\tilde{\xi}$ factors through $\tilde{\chi}$ by the same arrow. This proves that if $(U, \mu)$ is a colimiting cocone over the diagram $D$ then $(U, \mu\circ i)$ is a colimiting cocone over the diagram $D\circ i$.

Conversely, suppose that $(U, \mu \circ i)$ is a colimiting cocone over the diagram $D\circ i$.  $(Z, \tilde{\chi})$ is a colimiting cocone over the diagram $D$. For any cocone $(Z, \chi)$ over $D$, we have that $(Z, \chi \circ i)$ is a cocone over $D\circ i$; therefore $(U, \mu\circ i)$ factors uniquely through $(Z, \chi \circ i)$ or equivalently (by similar arguments as above), $(Z, \tilde{\chi})$ factors uniquely through $(U, \mu)$. Therefore $(U, \mu)$ is a colimiting cocone over $D$, as required.
\end{proofs}

Let $\underline{\Set}_{\cal E}$ be the $\cal E$-indexed category given by: ${\Set}_{E}=\Set$ for all $E\in {\cal E}$ and ${\Set}_{\alpha}=1_{\Set}$ (where $1_{\Set}$ is the identical functor on $\Set$) for all arrows $\alpha$ in $\cal E$. There is a $\cal E$-indexed functor 
\[
\gamma_{\cal E}:\underline{\Set}_{\cal E} \to \underline{\cal E}_{{\cal E}}
\]
defined by: ${\gamma_{\cal E}}_{E}=\gamma_{{\cal E}\slash E}^{\ast}:\Set \to {\cal E}\slash E$ (for any object $E$ of $\cal E$). 

We shall consider $\cal E$-indexed functors obtained by composing the $\cal E$-indexed functor $\gamma_{\cal E}:\underline{\Set}_{\cal E} \to \underline{\cal E}_{{\cal E}}$ with a $\cal E$-indexed functor $D:\underline{\cal A}_{\cal E}\to \underline{\Set}_{\cal E}$. 

Notice that a $\cal E$-indexed functor $D:\underline{\cal A}_{\cal E}\to \underline{\Set}_{\cal E}$ consists of a functor $D_{E}:{\cal A}_{E}\to \Set$ for each object $E$ of $\cal E$ such that for any arrow $\alpha:E'\to E$ in $\cal E$, $D_{E'}\circ {\cal A}_{\alpha}=D_{E}$.

Let us give an explicit description of the cocones over the $\cal E$-indexed functor $\gamma_{\cal E}\circ D$. Specializing the general definition, we obtain that a cocone $(U, \mu)$ over $\gamma_{\cal E}\circ D$ consists of an object $U$ of $\cal E$ and of an arrow $\mu_{E}(a):D_{E}(a)\to !_{E}^{\ast}(U)$ in ${\cal E}\slash E$ (for each objects $E$ of $\cal E$ and $a$ of ${\cal A}_{E}$) such that for any arrow $f:a\to b$ in ${\cal A}_{E}$, $\mu_{E'}(b)\circ \gamma_{{\cal E}\slash E}^{\ast}(D_{E}(f))=\mu_{E}(a)$ and for any arrow $\alpha:E'\to E$ in $\cal E$, $\alpha^{\ast}(\mu_{E}(a))=\mu_{E'}({\cal A}_{\alpha}(a))$.

By using the well-known adjunction between $\gamma_{\cal E}^{\ast}$ and the global sections functor ${\cal E}\to \Set$, we can alternatively present the above set of data as follows: a cocone $(U, \mu)$ over $\gamma_{\cal E}\circ D$ consists of an object $U$ of $\cal E$ and of a function $\mu_{E}(a):D_{E}(a)\to Hom_{\cal E}(E, U)$ (for each objects $E$ of $\cal E$ and $a$ of ${\cal A}_{E}$) such that for any arrow $f:a\to b$ in ${\cal A}_{E}$, $\mu_{E'}(b)\circ D_{E}(f)=\mu_{E}(a)$ and for any arrow $\alpha:E'\to E$ in $\cal E$, $Hom_{\cal E}(\alpha, U)\circ \mu_{E}(a)   =\mu_{E'}({\cal A}_{\alpha}(a))$ (notice that the domains of these two arrows are the same since $D_{E}(a)=D_{E'}({\cal A}_{\alpha}(a))$).

\begin{theorem}\label{finalfunc}
Let $F:{\cal C}^{\textrm{op}}\to {\cal E}$ be a functor. Then the $\cal E$-indexed subcategory $\underline{\int F}_{\cal E}$ of $\underline{\int F}^{f}_{\cal E}$ is strictly final. 
\end{theorem}

\begin{proofs}
By definition of ${\underline{{\int F}^{f}}_{\cal E}}$, for any object $E$ of $\cal E$ the category ${{\int F}^{f}}_{E}$ has as objects the arrows $x:E\to \mathbin{\mathop{\textrm{ $\coprod$}}\limits_{c\in Ob({\cal C})}}F(c)$ and as arrows $x\to x'$ the arrows $z:E\to \mathbin{\mathop{\textrm{ $\coprod$}}\limits_{f\in Arr({\cal C})}}F(cod(f))$ such that $d_{0}^{F}\circ z=x$ and $d_{1}^{F}\circ z=x'$. For any $E\in {\cal E}$, the category ${{\int F}}_{E}$ thus identifies as a subcategory of the category ${{\int F}^{f}}_{E}$, through the assignment sending any pair $(c, x)$, where $x:E\to F(c)$, to the arrow $\mu_{c}\circ x: E\to \mathbin{\mathop{\textrm{ $\coprod$}}\limits_{c\in Ob({\cal C})}}F(c)$ and any arrow $(c, x)\to (c', x')$ to the arrow $J_{f}\circ x':E\to \mathbin{\mathop{\textrm{ $\coprod$}}\limits_{f\in Arr({\cal C})}}F(cod(f))$ (here we use the notations of section \ref{seccol}). It is clear that these assignments make ${\underline{{\int F}}_{\cal E}}$ into a $\cal E$-indexed subcategory of ${\underline{{\int F}^{f}}_{\cal E}}$. 

To prove our thesis, we shall apply the criterion of Remark \ref{fullfinal}(a). To show that the $\cal E$-indexed subcategory $\underline{\int F}_{\cal E}$ satisfies the first condition in the remark, we observe that for any object $x:E\to \mathbin{\mathop{\textrm{ $\coprod$}}\limits_{c\in Ob({\cal C})}}F(c)$ of the category ${{\int F}^{f}}_{E}$, if we consider the pullbacks
\[  
\xymatrix {
E_{c} \ar[d]^{e_{c}} \ar[rr]^{x_{c}} & &  F(c) \ar[d]^{\mu_{c}} \\
E \ar[rr]^{x} & & \mathbin{\mathop{\textrm{ $\coprod$}}\limits_{c\in Ob({\cal C})}}F(c)} 
\]
of $x$ along the coproduct arrows $\mu_{c}:F(c)\to \mathbin{\mathop{\textrm{ $\coprod$}}\limits_{c\in Ob({\cal C})}}F(c)$, we obtain an epimorphic family $\{e_{c}:E_{c}\to E \textrm{ | } c\in Ob({\cal C})\}$ such that ${\int F}^{f}_{e_{c}}(x)$ lies in the subcategory ${{\int F}}_{E_{c}}$.      

It remains to prove that the second condition of Remark \ref{fullfinal}(a) is satisfied. Let us suppose that $(c, x)$ and $(c', x')$ are objects of the subcategory ${{\int F}}_{E}$, and suppose that $\alpha:E \to \mathbin{\mathop{\textrm{ $\coprod$}}\limits_{f\in Arr({\cal C})}}F(cod(f))$ is an arrow such that $d_{0}^{F}\circ \alpha= \mu_{c}\circ x$ and $d_{1}^{F}\circ \alpha=\mu_{c'}\circ x'$:
\[  
\xymatrix {
E \ar[d]^{x} \ar[rr]^{\alpha} & & \mathbin{\mathop{\textrm{ $\coprod$}}\limits_{f\in Arr({\cal C})}}F(cod(f)) \ar[d]^{d_{0}^{F}} & &  E \ar[d]^{x'} \ar[rr]^{\alpha} & &  \mathbin{\mathop{\textrm{ $\coprod$}}\limits_{f\in Arr({\cal C})}}F(cod(f))   \ar[d]^{d_{1}^{F}} \\
F(c)  \ar[rr]_{\mu_{c}} & & \mathbin{\mathop{\textrm{ $\coprod$}}\limits_{c\in Ob({\cal C})}}F(c) & &  F(c') \ar[rr]_{\mu_{c'}} & & \mathbin{\mathop{\textrm{ $\coprod$}}\limits_{c\in Ob({\cal C})}}F(c).} 
\]
Consider, for each $f\in Arr({\cal C})$, the following pullback square:   
\[  
\xymatrix {
E_{f} \ar[d]^{e_{f}} \ar[rr]^{\alpha_{f}} & &  F(cod(f)) \ar[d]^{J_{f}} \\
E \ar[rr]^{\alpha} & & \mathbin{\mathop{\textrm{ $\coprod$}}\limits_{f\in Arr({\cal C})}}F(cod(f)).} 
\]

The commutativity of the three diagrams above, combined with the fact that distinct coproduct arrows are disjoint, implies that for any arrow $f$ of $\cal C$ such that $dom(f)\neq c$ or $cod(f)\neq c'$ we have $E_{f}\cong 0_{\cal E}$. We can therefore restrict our attention to the arrows $f$ such that $dom(f)=c$ and $cod(f)=c'$. For any such arrow $f$, we have $\alpha_{f}=x'\circ e_{f}$ (equivalently, since $\mu_{c'}=\mu_{cod(f)}$ is monic, $\mu_{c'}\circ \alpha_{f}=\mu_{c'}\circ x'\circ e_{f}$). Indeed, $\mu_{c'}\circ \alpha_{f}=d_{1}^{F}\circ J_{f}\circ \alpha_{f}=d_{1}^{F}\circ \alpha \circ e_{f}=\mu_{c'}\circ x' \circ e_{f}$. Therefore $\alpha \circ e_{f}=J_{f}\circ x'\circ e_{f}$; indeed, $\alpha \circ e_{f}=J_{f}\circ \alpha_{f}=J_{f}\circ x'\circ e_{f}$. Lastly, we observe that the arrow $f$ defines an arrow $(c, x\circ e_{f})\to (c', x'\circ e_{f})$ in the subcategory $({{\int F}})_{E_{f}}$, i.e. $F(f)\circ x'\circ e_{f}=x\circ e_{f}$. Indeed, the arrow $\mu_{dom(f)}=\mu_{c'}$ is monic and we have $\mu_{dom(f)}\circ F(f)\circ x'\circ e_{f}= d_{0}^{F}\circ J_{f}\circ x' \circ e_{f}=J_{f}\circ \alpha \circ e_{f}=\mu_{c}\circ x \circ e_{f}$.

The epimorphic family $\{e_{f}:E_{f}\to E \textrm{ | } dom(f)=c \textrm{ and } cod(f)=c'\}$ thus satisfies the property that the arrow $({{\int F}^{f}_{\cal E}})_{e_{f}}(\alpha)$ lies in the subcategory $({{\int F}})_{E_{f}}$. Our proof is therefore complete.
\end{proofs}

\begin{theorem}\label{altdes}
Let $\cal C$ be a small category, $\cal E$ a Grothendieck topos, $F:{\cal C}^{\textrm{op}}\to {\cal E}$ and $P:{\cal C}\to \Set$ functors and $P_{\cal E}$ the internal diagram in $[{\mathbb C}, {\cal E}]$ given by $\overline{\gamma_{\cal E}^{\ast} \circ P}$. Then the restriction of the functor $(\underline{P_{\cal E}}_{\cal E}\circ \underline{\pi_{\overline{F}}^{f}}_{\cal E})\cong {\underline{{P_{\cal E}} \circ \pi_{\overline{F}}^{f}}}_{\cal E}$ to the $\cal E$-indexed subcategory $\underline{\int F}_{\cal E}$ of $\underline{\int F}^{f}_{\cal E}$ is naturally isomorphic to $\gamma_{\cal E}\circ z_{\cal E}$, where $z_{\cal E}:\underline{\int F}_{\cal E} \to \underline{\Set}_{\cal E}$ is the $\cal E$-indexed functor defined by: ${z_{\cal E}}_{E}((c, x))=P(c)$ and ${z_{\cal E}}_{E}(f)=P(f)$ (for any $E\in {\cal E}$, object $(c, x)$ and arrow $f$ in the category ${\int F}_{E}$).   
\end{theorem}

\begin{proofs}
We shall exhibit an isomorphism
\[
{{\underline{{P_{\cal E}} \circ \pi_{\overline{F}}^{f}}}_{\cal E}}_{E}((c, x))\cong ({\gamma_{\cal E}\circ z_{\cal E}})_{E}((c, x))
,\]
natural in $E\in {\cal E}$ and $(c, x)\in {\int F}_{E}$.  

Consider the following pullback diagram in $\textbf{cat}({\cal E})$:

\[  
\xymatrix {
\overline{F} \otimes_{\mathbb C} P \ar[d]^{t_{\overline{F}}} \ar[rr]^{t_{P}} & & \int^{opf} P_{\cal E} \ar[d]^{\pi^{opf}_{P_{\cal E}}} \\
\int^{f} \overline{F} \ar[rr]^{\pi_{\overline{F}}} & & {\mathbb C}.} 
\]

Since coproducts are stable under pullback in a topos, we have  that: 
 \[
(\overline{F} \otimes_{\mathbb C} P)_{0}={\mathbin{\mathop{\textrm{ $\coprod$}}\limits_{c\in Ob({\cal C})}}F(c) \times \gamma_{\cal E}^{\ast}(P(c))},
\]
\[
(\overline{F} \otimes_{\mathbb C} P)_{1}={\mathbin{\mathop{\textrm{ $\coprod$}}\limits_{f\in Arr({\cal C})}}F(cod(f)) \times \gamma_{\cal E}^{\ast}(P(dom(f)))},
\]
\[
{t_{F}}_{0}:(F \otimes_{\mathbb C} P)_{0}={\mathbin{\mathop{\textrm{ $\coprod$}}\limits_{c\in Ob({\cal C})}}F(c) \times \gamma_{\cal E}^{\ast}(P(c))} \to ({\int^{f} F})_{0}=\mathbin{\mathop{\textrm{ $\coprod$}}\limits_{c\in Ob({\cal C})}}F(c)
\]
is equal to the arrow $\mathbin{\mathop{\textrm{ $\coprod$}}\limits_{c\in Ob({\cal C})}}\pi_{F(c)}$, where $\pi_{F(c)}$ is the canonical projection $F(c)\times \gamma_{\cal E}^{\ast}(P(c))\to F(c)$,
\[
{t_{\overline{F}}}_{1}:(\overline{F} \otimes_{\mathbb C} P)_{1}={\mathbin{\mathop{\textrm{ $\coprod$}}\limits_{f\in Arr({\cal C})}}F(cod(f)) \times \gamma_{\cal E}^{\ast}(P(dom(f)))} \to ({\int^{f} F})_{1}=\mathbin{\mathop{\textrm{ $\coprod$}}\limits_{f\in Arr({\cal C})}}F(cod(f))
\]
is equal to the arrow $\mathbin{\mathop{\textrm{ $\coprod$}}\limits_{f\in Arr({\cal C})}}\pi'_{F(cod(f))}$, where $\pi'_{F(cod(f))}$ is the canonical projection $F(cod(f))\times \gamma_{\cal E}^{\ast}(P(dom(f)))\to F(cod(f))$ and the domain and codomain arrows 
\[
d_{0}^{\overline{F} \otimes_{\mathbb C} P}, d_{1}^{\overline{F} \otimes_{\mathbb C} P}:(\overline{F} \otimes_{\mathbb C} P)_{1} \to (\overline{F} \otimes_{\mathbb C} P)_{0}
\] 
are defined by the following conditions: 
\[
d_{0}^{\overline{F} \otimes_{\mathbb C} P} \circ W_{f}=Z_{cod(f)}\circ <\pi'_{F(cod(f))}, \gamma_{\cal E}^{\ast}(P(f))\circ \mu_{\gamma_{\cal E}^{\ast}(P(dom(f)))}>
\]
and 
\[
d_{1}^{\overline{F} \otimes_{\mathbb C} P} \circ W_{f}=Z_{dom(f)}\circ <F(f) \circ\pi'_{F(cod(f))}, \mu_{\gamma_{\cal E}^{\ast}(P(dom(f)))}>,
\]
where $\mu_{\gamma_{\cal E}^{\ast}(P(dom(f)))}$ is the canonical projection 
\[
F(cod(f))\times \gamma_{\cal E}^{\ast}(P(dom(f))) \to \gamma_{\cal E}^{\ast}(P(dom(f)))
\]
and 
\[
W_{f}:F(cod(f)) \times \gamma_{\cal E}^{\ast}(P(dom(f))) \to {\mathbin{\mathop{\textrm{ $\coprod$}}\limits_{f\in Arr({\cal C})}}F(cod(f)) \times \gamma_{\cal E}^{\ast}(P(dom(f)))},
\]
\[
Z_{c}:F(c) \times \gamma_{\cal E}^{\ast}(P(c)) \to {\mathbin{\mathop{\textrm{ $\coprod$}}\limits_{c\in Ob({\cal C})}}F(c) \times \gamma_{\cal E}^{\ast}(P(c))}
\]
are the canonical coproduct arrows (respectively for $f\in Arr({\cal C})$ and $c\in Ob({\cal C})$).

Now, the internal diagram $P_{\cal E}\circ \pi^{f}_{\overline{F}}$ corresponds precisely to the discrete opfibration $t_{\overline{F}}$ and hence the $\cal E$-indexed functor ${\underline{{P_{\cal E}} \circ \pi_{\overline{F}}^{f}}}_{\cal E}$ sends any generalized element $x:E \to \mathbin{\mathop{\textrm{ $\coprod$}}\limits_{c\in Ob({\cal C})}}F(c)$ to the object of ${\cal E}\slash E$ given by the pullback of ${t_{\overline{F}}}_{0}$ along it. In particular, since in a topos any diagram of the form 
\[  
\xymatrix {
A_{i} \ar[d]^{f_{i}} \ar[rr] & & {\mathbin{\mathop{\textrm{ $\coprod$}}\limits_{i\in I}}A_{i}} \ar[d]^{{\mathbin{\mathop{\textrm{ $\coprod$}}\limits_{i\in I}}f_{i}}} \\
B_{i} \ar[rr] & & {\mathbin{\mathop{\textrm{ $\coprod$}}\limits_{i\in I}}B_{i}},}  
\]
where the horizontal arrows are the canonical coproduct arrows, is a pullback and products commute with pullbacks, the functor $({{P_{\cal E}} \circ \pi_{\overline{F}}^{f}})_{E}$ sends any object $(c, x)$ of the category ${\int F}_{E}$ to the object given by the canonical projection $E\times \gamma_{\cal E}^{\ast}(P(c)) \to E$. But this object is canonically isomorphic to $\gamma_{{\cal E}\slash E}^{\ast}(P(c))$ in ${\cal E}\slash E$, which is the value at $(c, x)$ of the functor $(\gamma_{\cal E} \circ z_{\cal E})_{E}$. Our proof is therefore complete.  
\end{proofs}

Let us now explicitly describe the $\cal E$-indexed colimiting cocone on the restriction to the $\cal E$-indexed subcategory $\underline{{\int \overline{F}}}_{\cal E}$ of the $\cal E$-indexed diagram ${\underline{{P_{\cal E}} \circ \pi_{\overline{F}}^{f}}}_{\cal E}$ (cf. section \ref{backgr} for the background). The vertex of this colimiting cocone is the codomain of the coequalizer $w:(\overline{F} \otimes_{\mathbb C} P)_{0} \to \underline{colim}_{\cal E}(\underline{P_{\cal E}}_{\cal E}\circ \underline{\pi_{\overline{F}}^{f}}_{\cal E})$ of the pair of arrows $d_{0}^{\overline{F} \otimes_{\mathbb C} P}, d_{1}^{\overline{F} \otimes_{\mathbb C} P}:(\overline{F} \otimes_{\mathbb C} P)_{1} \to (\overline{F} \otimes_{\mathbb C} P)_{0}$. For any object $(c, x)$ of the category ${\int F}_{E}$, the colimit arrow $({{P_{\cal E}} \circ \pi_{\overline{F}}^{f}})_{E})((c, x))\cong \gamma_{{\cal E}\slash E}^{\ast}(P(c)) \to \underline{colim}_{\cal E}(\underline{P_{\cal E}}_{\cal E}\circ \underline{\pi_{\overline{F}}^{f}}_{\cal E})$ is given by the composite $!_{E}^{\ast}(w) \circ h_{(c, x)}$, where $h_{(c, x)}$ is the arrow $!_{E}^{\ast}(Z_{c})\circ <x\times 1_{\gamma_{{\cal E}}^{\ast}(P(c))}, \pi_{E}>:\gamma_{{\cal E}\slash E}^{\ast}(P(c))\cong E\times \gamma_{{\cal E}}^{\ast}(P(c)) \to {\mathbin{\mathop{\textrm{ $\coprod$}}\limits_{c\in Ob({\cal C})}}F(c) \times \gamma_{\cal E}^{\ast}(P(c))}\times E \cong !_{E}^{\ast}({\mathbin{\mathop{\textrm{ $\coprod$}}\limits_{c\in Ob({\cal C})}}F(c) \times \gamma_{\cal E}^{\ast}(P(c))})$; in other words, the function $\xi_{(c, x)}: P(c)\to Hom_{\cal E}(E, \underline{colim}_{\cal E}(\underline{P_{\cal E}}_{\cal E}\circ \underline{\pi_{\overline{F}}^{f}}_{\cal E}))$ corresponding to it assigns to every element $a\in P(c)$ the arrow $w\circ Z_{c}\circ <x, r_{a}\circ !_{E}>$, where $r_{a}:1_{\cal E} \to \gamma_{\cal E}^{\ast}(P(c))\cong {\mathbin{\mathop{\textrm{ $\coprod$}}\limits_{a\in P(c)}}1_{\cal E}}$ is the coproduct arrow corresponding to the element $a$ and $!_{E}$ is the unique arrow $E\to 1_{\cal E}$ in $\cal E$.

It is immediate to verify that, under the isomorphism 
\[
\underline{colim}_{\cal E}(\underline{P_{\cal E}}_{\cal E}\circ \underline{\pi_{\overline{F}}^{f}}_{\cal E})\cong colim(F\circ \pi^{f}_{P})
\]
of Theorem \ref{tensorpr}, the functions $\xi_{(c, x)}:P(c)\to \underline{colim}_{\cal E}(\underline{P_{\cal E}}_{\cal E}\circ \underline{\pi_{\overline{F}}^{f}}_{\cal E})\cong  colim(F\circ \pi^{f}_{P})$ admit the following description in terms of the colimiting arrows $\chi_{(c, a)}:F(c)\to colim(F\circ \pi^{f}_{P})$: for any $a\in P(c)$, $\xi_{(c, x)}(a)=\chi_{(c, a)}\circ x$.

\subsection{A characterization of $\cal E$-indexed colimits}\label{charindcol}

In this section we establish necessary and sufficient conditions for a $\cal E$-indexed cocone on a small diagram in $\cal E$ to be a ($\cal E$-indexed) colimit cone. 

Before stating and proving the main theorem of this section, we need to introduce some relevant definitions and a few technical lemmas.

Given a $\cal E$-indexed functor $D:\underline{\cal A}_{\cal E} \to \underline{{\cal E}}_{\cal E}$ and an object $R$ of the topos ${\cal E}$, we denote by ${\cal I}_{R}^{D}$ the set of pairs of the form $(x, y)$, where $x$ is an object of the category ${\cal A}_{R}$ and $y$ is an arrow $1_{{\cal E}\slash R}\to D_{R}(x)$ in the topos ${\cal E}\slash R$. On such a set we consider the equivalence relation ${\cal R}_{R}^{D}$ generated by the pairs of the form $((x, y), (x', y'))$, where there exists an arrow $f:x \to x'$ in the category ${\cal A}_{R}$ such that $y'=D_{R}(f)\circ y$. 

The following remarks are useful in connection with the application of the localization technique.

\begin{remarks}\label{rem3}
\begin{enumerate}[(a)]

\item For any $E\in {\cal E}$, we have a natural bijection ${\cal I}^{D}_{E}\cong {\cal I}^{D\slash E}_{1_{{\cal E}\slash E}}$, and the relation ${\cal R}^{D}_{E}$ corresponds to the relation ${\cal R}^{D\slash E}_{1_{{\cal E}\slash E}}$ under this bijection. 

\item For any arrow $f:F\to E$ and any object $x\in {\cal A}_{E}$, $D_{F}({\cal A}_{f}(x))=f^{\ast}(D_{E}(x))$ (this follows from the fact that $D$ is a $\cal E$-indexed functor). We shall denote the canonical arrow $dom(D_{F}({\cal A}_{f}(x))) \to dom(D_{E}(x))$ by the symbol $r_{f}$. 

\item For any object $x$ of ${\cal A}_{E}$ and any arrow $y:(f:F\to E)\to D_{E}(x)$ in the topos ${\cal E}\slash E$, there exists a unique arrow $y_{f}:(1_{F}:F\to E) \to D_{F}({\cal A}_{f}(x))$ in the topos ${\cal E}\slash F$ such that $r_{f}\circ y_{f}=y$ (where $r_{f}$ is considered here as an arrow in ${\cal E}\slash E$ in the obvious way);

\item For any $\cal E$-indexed cocone $\mu$ over $D$ with vertex $U$, any object $x$ of ${\cal A}_{E}$ and any arrow $y:(f:F\to E)\to D_{E}(x)$ in the topos ${\cal E}\slash E$, we have $\mu_{E}(x)\circ y=(1_{U}\times f)\circ \mu_{F}({\cal A}_{f}(x))\circ y_{f}$. Moreover, we have pullback squares  
\[  
\xymatrix {
D_{F}({\cal A}_{f}(x)) \ar[d]^{\mu_{F}({\cal A}_{f}(x))} \ar[r]^{r_{f}} & D_{E}(x) \ar[d]^{\mu_{E}(x)}\\
U\times F \ar[d]^{\pi_{F}}  \ar[r]_{1_{U}\times f} & U\times E \ar[d]^{\pi_{E}} \\
F \ar[r]^{f} & E} 
\]
in $\cal E$, where $\pi_{E}$ and $\pi_{F}$ are the obvious canonical projections. In particular, for any other pair $(x', y')$, $\mu_{E}(x)\circ y=\mu_{E}(x')\circ y'$ if and only if $\mu_{F}({\cal A}_{f}(x))\circ y_{f}=\mu_{F}({\cal A}_{f}(x'))\circ y_{f}$; 

\item If $\underline{\cal A}_{\cal E}$ is the $\cal E$-externalization of an internal category $\mathbb C$ in $\cal E$, $x$ is a generalized element $E\to {\mathbb C}_{0}$ and $y$ is an arrow $(f:F\to E) \to D_{E}(x)$ in ${\cal E}\slash E$, there exists a unique arrow $y_{f}:(1_{F}:F\to F) \to D_{F}(x\circ f)$ in ${\cal E}\slash F$ such that $z_{x\circ f}\circ y_{f}=z_{x}\circ y$ (this arrow is provided by the universal property of the pullback square defining $D_{F}(x\circ f)$), which coincides with the arrow $y_{f}$ defined in Remark \ref{rem3}(c);    

\item If a pair $((x, y), (x', y'))$ belongs to ${\cal R}_{R}^{D}$ then for any indexed cocone $\lambda$ over $D$ with vertex $V$, $\lambda_{R}(x)\circ y=\lambda_{R}(x')\circ y'$. Indeed, for any pair $((x, y), (x', y'))$ with the property that there exists an arrow $f:x \to x'$ in the category ${\cal A}_{R}$ such that $y'=D_{E}(f)\circ y$, we have that $\lambda_{R}(x')\circ y'=\lambda_{R}(x')\circ D_{E}(f)\circ y=\lambda_{R}(x)\circ y$ (since $\lambda_{R}$ is a cocone).   
\end{enumerate}
\end{remarks}

\begin{lemma}\label{lifting1}
Let $D\in [{\mathbb C}, {\cal E}]$ be an internal diagram, $x:1\to {\mathbb C}_{0}$ and $f:1\to {\mathbb C}_{1}$ generalized elements such that $d_{0}^{\mathbb C}\circ f=x$ and $m:F\to D_{1}(x)$ an arrow in $\cal E$. Then the unique arrow $\chi_{m, f}:F\to (\int^{opf} D)_{1}$ in $\cal E$ such that $d_{0}^{{\int^{opf} D}}\circ \chi_{m, f}=z_{x}\circ m$ and $(\pi^{opf}_{D})_{1}\circ \chi_{m, f}=f\circ !_{F}$ (provided by the universal property of the pullback square corresponding to the domains of the discrete opfibration associated to $D$) satisfies the property that $z_{x}\circ m=d_{0}^{{\int^{opf} D}}\circ \chi_{m, f}$ (by definition) and $z_{x'}\circ D_{1}(f)\circ m=d_{1}^{{\int^{opf} D}}\circ \chi_{m, f}$, where $x'$ is equal to $d_{1}^{\mathbb C}\circ f$ and $d_{0}^{{\int^{opf} D}}, d_{1}^{{\int^{opf} D}}$ are respectively the domain and codomain arrows $(\int^{opf} D)_{1} \to (\int^{opf} D)_{0}$. 
\end{lemma}

\begin{proofs}
By definition of $D_{1}(f)$, we have that $z_{x'}\circ D_{1}(f)=d_{1}^{{\int^{opf} D}}\circ r \circ \xi$, where $r$ is defined by the pullback square
\[  
\xymatrix {
H \ar[d]^{!_{H}} \ar[r]^{r} & ({\int^{opf} D})_{1} \ar[d]^{{\pi^{opf}_{D}}_{1}}\\
1_{\cal E}  \ar[r]_{f} & {\mathbb C}_{1}} 
\]
and $\xi$ is the unique arrow $D_{1}(x)\to H$ such that $z_{x}=d_{0}^{{\int^{opf} D}}\circ r\circ \xi$, provided by the universal property of the pullback square  

\[  
\xymatrix {
H \ar[d]^{s} \ar[r]^{d_{0}^{{\int^{opf} D}} \circ r} & ({\int^{opf} D})_{0} \ar[d]^{({\pi^{opf}_{D}})_{0}}\\
1_{\cal E}  \ar[r]_{d_{0}^{\mathbb C} \circ f} & {\mathbb C}_{0}.} 
\]
 
It thus remains to prove that $\chi_{m, f}=r\circ \xi \circ m$. But this immediately follows from the universal property of the pullback 

\[  
\xymatrix {
{\int^{opf} D}_{1} \ar[d]^{{\pi^{opf}_{D}}_{1}} \ar[r]^{d_{0}^{{\int^{opf} D}}} & ({\int^{opf} D})_{0} \ar[d]^{{\pi^{opf}_{D}}_{0}}\\
{\mathbb C}_{1}  \ar[r]_{d_{0}^{\mathbb C}} & {\mathbb C}_{0},} 
\]
since $d_{0}^{{\int^{opf} D}}\circ \chi_{m, f}=d_{0}^{{\int^{opf} D}}\circ (r\circ \xi \circ m)$ and ${\pi^{opf}_{D}}_{1} \circ \chi_{m, f}={\pi^{opf}_{D}}_{1}\circ (r\circ \xi \circ m)$.

\end{proofs}

\begin{lemma}\label{downl}
Let $D$ be an internal diagram in $[{\mathbb C}, {\cal E}]$ and $\chi$ be an arrow $1_{\cal E}\to (\int^{opf}D)_{1}$. Let $x=(\pi^{opf}_{D})_{0}\circ d_{0}^{{\int^{opf} D}}\circ \chi$, $x'=(\pi^{opf}_{D})_{0} \circ d_{1}^{{\int^{opf} D}}\circ \chi$ and $f=(\pi^{opf}_{D})_{1}\circ \chi$. Then $d_{0}^{\mathbb C}\circ f=x$, $d_{1}^{\mathbb C}\circ f=x'$ and $D_{1}(f)\circ \chi_{0}=\chi_{1}$, where $\chi_{0}$ is the unique arrow $1_{\cal E}\to D_{1}(x)$ such that $z_{x}\circ \chi_{0}=d_{0}^{{\int^{opf} D}}\circ \chi$ (which exists by the universal property of the pullback square defining $D_{1}(x)$) and $\chi_{1}$ is the unique arrow $1_{\cal E}\to D_{1}(x')$ such that $z_{x'}\circ \chi_{1}=d_{1}^{{\int^{opf} D}}\circ \chi$ (which exists by the universal property of the pullback square defining $D_{1}(x')$.
\end{lemma}

\begin{proofs}
The first two identities follow straightforwardly from the definition. It remains to prove that $D_{1}(f)\circ \chi_{0}=\chi_{1}$. Let us refer to the proof of Lemma \ref{lifting1} for notation. By the universal property of the pullback square defining $D_{1}(x')$, it is equivalent to verify that $z_{x'}\circ D_{1}(f)\circ \chi_{0}=z_{x'}\circ \chi_{1}$. But $z_{x'}\circ D_{1}(f)\circ \chi_{0}= d_{1}^{{\int^{opf} D}}\circ r\circ \xi \circ \chi_{0}$, while $z_{x'}\circ \chi_{1}=d_{1}^{{\int^{opf} D}}\circ \chi$; so to prove the desired equality it suffices to show that $r\circ \xi \circ \chi_{0}=\chi$. This follows from the fact that $d_{0}^{{\int^{opf} D}}\circ r\circ \xi \circ \chi_{0}=z_{x}\circ \chi_{0}=d_{0}^{{\int^{opf} D}}\circ \chi$ by virtue of the universal property of the pullback square given by the domain of the discrete opfibration associated to $D$. 
\end{proofs}

\begin{proposition}\label{lifting2}
Let $D\in [{\mathbb C}, {\cal E}]$ be an internal diagram, $x,x':1\to {\mathbb C}_{0}$ generalized elements and $y:1_{\cal E} \to D_{1}(x)$, $y':1_{\cal E} \to D_{1}(x')$ arrows in $\cal E$. Then $(z_{x}\circ y, z_{x'}\circ y')$ belongs to the equivalence relation on $Hom_{\cal E}(1_{\cal E}, (\int^{opf} D)_{0})$ generated by the pairs of the form $(d_{0}^{{\int^{opf} D}}\circ a, d_{1}^{{\int^{opf} D}}\circ a)$ for some arrow $a:1_{\cal E} \to (\int^{opf} D)_{1}$ if and only if $((x, y), (x', y'))$ belongs to ${\cal R}^{D}_{1_{\cal E}}$. 
\end{proposition}

\begin{proofs}
Let us first prove the `if' direction. We have that $((x, y), (x', y'))$ belongs to ${\cal R}^{D}_{1_{\cal E}}$ if and only if there exists a finite sequence 
\[
(x_{0}, y_{0})=(x, y), \ldots,  (x_{n}, y_{n})=(x', y')
\]
such that for any $i\in \{0, \ldots, n-1\}$

\begin{enumerate}[(1)]
\item either there exists $f_{i}:1\to {\mathbb C}_{1}$ such that $d_{0}^{\mathbb C}\circ f_{i}=x_{i}$, $d_{1}^{\mathbb C}\circ f_{i}=x_{i+1}$ and $D_{1}(f_{i})\circ y_{i}=y_{i+1}$ or

\item there exists $g_{i}:1\to {\mathbb C}_{1}$ such that $d_{0}^{\mathbb C}\circ g_{i}=x_{i+1}$, $d_{1}^{\mathbb C}\circ g_{i}=x_{i}$ and $D_{1}(g_{i})\circ y_{i+1}=y_{i}$.
\end{enumerate} 

In case $(1)$, Lemma \ref{lifting1} implies (by taking $m$ to be $y_{i}$ and $f$ to be $f_{i}$) that $z_{x_{i}}\circ y_{i}=d_{0}^{{\int^{opf} D}}\circ \chi_{y_{i}, f_{i}}$ and $z_{x_{i+1}}\circ y_{i+1}=d_{1}^{{\int^{opf} D}}\circ \chi_{y_{i}, f_{i}}$ (here the notation is that of the lemma). In case $(2)$, Lemma \ref{lifting1} implies (by taking $m$ to be $y_{i+1}$ and $f$ to be $g_{i}$) that $z_{x_{i+1}}\circ y_{i+1}=d_{0}^{{\int^{opf} D}}\circ \chi_{y_{i+1}, g_{i}}$ and $z_{x_{i}}\circ y_{i}=d_{1}^{{\int^{opf} D}}\circ \chi_{y_{i+1}, g_{i}}$.

From this it clearly follows that $(z_{x}\circ y, z_{x'}\circ y')$ belongs to the equivalence relation $T$ on $Hom_{\cal E}(1_{\cal E}, (\int^{opf} D)_{0})$ generated by the pairs of the form $(d_{0}^{{\int^{opf} D}}\circ a, d_{1}^{{\int^{opf} D}}\circ a)$ for some arrow $a:1_{\cal E} \to (\int^{opf} D)_{1}$, as required.

Let us now prove the `only if' direction. If $(z_{x}\circ y, z_{x'}\circ y')$ belongs to the equivalence relation $T$ defined above, there exists a finite sequence $\chi_{0}, \ldots, \chi_{n}$ of arrows $1_{\cal E}\to (\int^{opf} D)_{1}$ and of arrows $e_{i}$ (which are either $d_{0}^{{\int^{opf} D}}$ or $d_{1}^{{\int^{opf} D}}$ - we denote by $e_{i}^{\textrm{op}}$ the arrow $d_{1}^{{\int^{opf} D}}$ if $e_{i}$ is $d_{0}^{{\int^{opf} D}}$ and the arrow $d_{0}^{{\int^{opf} D}}$ if $e_{i}$ is $d_{1}^{{\int^{opf} D}}$) such that $e_{0}\circ \chi_{0}=z_{x}\circ y$, $e_{i+1}\circ \chi_{i+1}=e_{i}^{\textrm{op}}\circ \chi_{i}$ for all $i\in \{0, \ldots, n-1\}$ and $e_{n}\circ \chi_{n}=z_{x'}\circ y'$. 

To deduce our thesis, it suffices to apply Lemma \ref{downl} noticing that if $e_{i}\circ \chi=z_{x''}\circ y''$ then $\chi_{0}=y''$ if $e_{i}$ is $d_{0}^{{\int^{opf} D}}$ and $\chi_{1}=y''$ if $e_{i}$ is $d_{1}^{{\int^{opf} D}}$ (where $\chi_{0}$ and $\chi_{1}$ are the arrows defined in the statement of Lemma \ref{downl}) since $z_{x''}\circ \chi_{0}=d_{0}^{{\int^{opf} D}}\circ \chi=z_{x''}\circ y''$ (if $e_{i}$ is $d_{0}^{{\int^{opf} D}}$) and $z_{x''}\circ \chi_{1}=d_{1}^{{\int^{opf} D}}\circ \chi=z_{x''}\circ y''$ (if $e_{i}$ is $d_{1}^{{\int^{opf} D}}$).     
\end{proofs}

\begin{theorem}\label{charact}
Let $D:\underline{\cal A}_{\cal E} \to \underline{{\cal E}}_{\cal E}$ be a $\cal E$-indexed functor, where $\underline{\cal A}_{\cal E}$ is equivalent to the $\cal E$-externalization of an internal category in $\cal E$. Then a cocone $\mu$ over $D$ with vertex $U$ is an indexed colimiting cocone for $D$ if and only if the following conditions are satisfied:
\begin{enumerate}[(i)]

\item For any object $F$ of $\cal E$ and arrow $h:F\to U$ in the topos $\cal E$, there exists an epimorphic family $\{f_{i}:F_{i}\to F \textrm{ | } i\in I\}$ in $\cal E$ and for each $i\in I$ an object $x_{i}\in {\cal A}_{F_{i}}$ and an arrow $\alpha_{i}:1_{{\cal E}\slash F_{i}} \to D_{F_{i}}(x_{i})$ in the topos ${\cal E}\slash F_{i}$ such that $<h\circ f_{i}, 1_{F_{i}}>=\mu_{F_{i}}(x_{i})\circ \alpha_{i}$ as arrows $1_{{\cal E}\slash F_{i}} \to F_{i}^{\ast}(U)$ in ${\cal E}\slash F_{i}$;

\item For any pairs $(x, y)$ and $(x', y')$, where $x$ and $x'$ are objects of ${\cal A}_{E}$, $y$ is an arrow $(f:F\to E) \to D_{E}(x)$ in ${\cal E}\slash E$ and $y'$ is an arrow $(f:F\to E) \to D_{E}(x')$ in ${\cal E}\slash E$, we have $\mu_{E}(x)\circ y=\mu_{E}(x')\circ y'$ if and only if  there exists an epimorphic family $\{f_{i}:F_{i}\to F \textrm{ | } i\in I\}$ in $\cal E$ such that the pair $(({\cal A}_{f\circ f_{i}}(x), f_{i}^{\ast}(y)), ({\cal A}_{f\circ f_{i}}(x'), f_{i}^{\ast}(y'))$ belongs to the relation ${\cal R}_{F_{i}}^{D}$.
\end{enumerate}
\end{theorem}

\begin{proofs}
First, let us prove that, under the assumption that $\underline{\cal A}_{\cal E}$ is the $\cal E$-externalization of an internal category $\mathbb C$ in $\cal E$, the colimit cone for $D$ satisfies the two conditions of the theorem. 

Under this assumption, condition $(i)$, applied to the colimiting cocone $\mu$ for $D$, rewrites as follows: for any object $F$ of $\cal E$ and arrow $h:F\to colim(D)$ in ${\cal E}$, there exists an epimorphic family $\{f_{i}:F_{i}\to F \textrm{ | } i\in I\}$ in $\cal E$ and for each $i\in I$ a generalized element $x_{i}:F_{i}\to {\mathbb C}_{0}$ and an arrow $\alpha_{i}:1_{{\cal E}\slash F_{i}} \to D_{F_{i}}(x_{i})$ in the topos ${\cal E}\slash F_{i}$ such that $<h\circ f_{i}, 1_{F_{i}}>= c\circ z_{x_{i}} \circ \alpha_{i}$ as arrows $(f_{i}:1_{{\cal E}\slash F_{i}} \to F_{i}^{\ast}(colim(D))$ in ${\cal E}\slash F_{i}$ (where the notation is that of section \ref{backgr}).   

Consider the following pullback square:
\[  
\xymatrix {
F' \ar[d]^{h'} \ar[rr]^{f'} & & F \ar[d]^{h} \\
({\int^{opf}D})_{0} \ar[rr]^{c} & & colim(D).}  
\] 
Since the arrow $c$ is an epimorphism, the arrow $f':F' \to F$ is an epimorphism. This arrow will form, by itself, the single element of an epimorphic family satisfying condition $(i)$. We set $x'$ equal to the composite $({\pi^{opf}_{D}})_{0} \circ h'$ (recall that $({\pi_{D}^{opf}})_{0}:({\int^{opf}D})_{0} \to {\mathbb C}_{0}$ is the object component of the discrete opfibration associated to $D$). Consider the following pullback square defining $D_{F'}(x')$:
\[  
\xymatrix {
D_{F'}(x') \ar[d]^{v} \ar[rr]^{u} & & (\int^{opf} D)_{0} \ar[d]^{(\pi^{opf}_{D})_{0}} \\
F' \ar[rr]^{x'} & & {\mathbb C}_{0}.}  
\] 
By the universal property of this pullback square, there exists a unique arrow $\alpha_{x'}:F'\to D_{F'}(x')$ in $\cal E$ such that $u\circ \alpha_{x'}=h'$ and $v\circ \alpha_{x'}=1_{F'}$. So $\alpha_{x'}$ is an arrow $1_{{\cal E}\slash F'}\to D_{F}(x')$ in the topos ${\cal E}\slash F'$, and $<h\circ f', 1_{F'}>=<c\circ u\circ \alpha_{x'}, v\circ \alpha_{x'}>$, as required.

Under the assumption that $\underline{\cal A}_{\cal E}$ is the $\cal E$-externalization of an internal category $\mathbb C$ in $\cal E$, condition $(ii)$, applied to the colimiting cocone $\mu$ for $D$, rewrites as follows: for any pairs $(x, y)$ and $(x', y')$, where $x$ and $x'$ is a generalized elements $E \to {\mathbb C}_{0}$, $y$ is an arrow $(f:F\to E) \to D_{E}(x)$ in ${\cal E}\slash E$ and $y'$ is an arrow $(f:F\to E) \to D_{E}(x')$ in ${\cal E}\slash E$, we have $c\circ z_{x}\circ y=c\circ z_{x'}\circ y'$ if and only if there exists an epimorphic family $\{f_{i}:F_{i}\to F \textrm{ | } i\in I\}$ in $\cal E$ such that the pair $((x\circ f \circ f_{i}, y\circ f_{i}), (x'\circ f\circ f_{i}, y'\circ f_{i}))$ belongs to the relation ${\cal R}_{F_{i}}^{D}$.   

Thanks to the localization technique, we can suppose $E=1_{\cal E}$ without loss of generality. Indeed, condition $(ii)$ for the diagram $D$, the cocone $\mu$ and the object $E$ is equivalent to condition $(ii)$ for the diagram $D\slash E$, the cocone $\mu\slash E$ and the object $1_{{\cal E}\slash E}$, and $\mu$ is, by our assumption, an indexed colimiting cocone and hence stable under localization.

By Lemma \ref{lemma2}, we have that $c\circ z_{x}\circ y=c\circ z_{x'}\circ y'$ if and only if there exists an epimorphic family $\{f_{i}: F_{i} \to F \textrm{ | } i\in I \}$ in $\cal E$ such that for each $i\in I$ $(z_{x}\circ y \circ f_{i}, z_{x'}\circ y' \circ f_{i})$ belongs to the equivalence relation on the set $Hom_{{\cal E}}(F_{i}, (\int^{opf} D)_{0})$ generated by the set of pairs of the form $(d_{0}^{{\int^{opf} D}}\circ a, d_{1}^{{\int^{opf} D}}\circ a)$, for a generalized element $a:F_{i}\to (\int^{opf} D)_{0}$ in $\cal E$. 

By Remark \ref{rem3}(e), we have that $(z_{x}\circ y \circ f_{i}, z_{x'}\circ y' \circ f_{i})=(z_{x\circ f\circ f_{i}} \circ (y\circ f_{i}), z_{x'\circ f\circ f_{i}} \circ (y'\circ f_{i}))$. Our thesis then follows from Proposition \ref{lifting2}, applied to the toposes ${\cal E}\slash F_{i}$ and the pairs $((x\circ f \circ f_{i}, y\circ f_{i}), (x'\circ f\circ f_{i}, y'\circ f_{i}))$.
 
Let us now prove that the conditions of the theorem are sufficient for $\mu$ to be a colimiting cocone. For this part of the theorem, we shall not need to assume the $\cal E$-indexed category $\underline{\cal A}_{\cal E}$ to be small.

Since conditions $(i)$ and $(ii)$ are both stable under localization, it suffices to prove that $\mu$ is a universal colimiting cocone over the diagram $D$. Suppose that $\lambda$ is an indexed cocone over $D$ with vertex $V$. We have to prove that there exists a unique arrow $l:U\to V$ in $\cal E$ such that for any $E\in {\cal E}$ and object $x$ of ${\cal A}_{E}$, $(l\times 1_{E})\circ \mu_{E}(x)=\lambda_{E}(x)$. To define the arrow $l$ we shall define a function $L:Hom_{\cal E}(E, U)\to Hom_{\cal E}(E, V)$ natural in $E\in {\cal E}$. Given $h\in Hom_{\cal E}(E, U)$, by condition $(i)$ in the statement of the theorem, there exists an epimorphic family ${\cal E}=\{e_{i}:E_{i}\to E \textrm{ | } i\in I\}$ in $\cal E$ and for each $i\in I$ an object $x_{i}$ of ${\cal A}_{E_{i}}$ and an arrow $\alpha_{i}:1_{{\cal E}\slash E_{i}} \to D_{E_{i}}(x_{i})$ in ${\cal E}\slash E_{i}$ such that $<h\circ e_{i}, 1_{E_{i}}>=\mu_{E_{i}}(x_{i})\circ \alpha_{i}$. The family of arrows $\{\lambda_{E_{i}}(x_{i})\circ \alpha_{i} \textrm{ | } i\in I\}$ satisfies the hypotheses of Corollary \ref{cordesc}. Indeed, using the notation of the corollary, we have that $q_{i}^{\ast}(\lambda_{E_{i}}(x_{i})\circ \alpha_{i})=q_{j}^{\ast}(\lambda_{E_{j}}(x_{j})\circ \alpha_{j})$ for all $i, j\in I$; this can be easily proved by using condition $(ii)$ of the theorem, the fact that the same identity holds for the cocone $\mu$ (in place of $\lambda$) and Remark \ref{rem3}(d). Therefore there exists a unique arrow $L^{\cal E}_{h}:E\to V$ such that $<L^{\cal E}_{h}\circ e_{i}, 1_{E_{i}}>=\lambda_{E_{i}}(x_{i})\circ \alpha_{i}$ for all $i\in I$. In order to be able to define $L(h)$ as equal to $L^{\cal E}_{h}$ we need to check that for any other epimorphic family ${\cal E}'=\{e_{j}':E_{j}'\to E \textrm{ | } j\in J\}$ in $\cal E$, we have $L^{\cal E}_{h}=L^{{\cal E}'}_{h}$. To this end, consider the fibered product of the epimorphic families ${\cal E}$ and ${\cal E}'$, that is the family of arrows $p_{i, j}=e_{i}\circ f_{i}=e_{j}'\circ f_{j}':F_{i, j}\to E$ (for $i, j\in I\times J$), where the arrows $f_{i}$ and $f_{j}'$ are defined by the following pullback square:
\[  
\xymatrix {
F_{i, j} \ar[d]^{f_{j}'} \ar[r]^{f_{i}} & E_{i} \ar[d]^{e_{i}} \\
E_{j} \ar[r]^{e_{j}'} & E.} 
\]  
It clearly suffices to verify that for any $i, j\in I$, we have $L^{\cal E}_{h}\circ p_{i, j}=L^{{\cal E}'}_{h}\circ p_{i, j}$. Now, we have that $<L^{{\cal E}}_{h}\circ e_{i}, 1_{E_{i}}>=\lambda_{E_{i}}(x_{i})\circ \alpha_{i}$, while $<L^{{\cal E}'}_{h}\circ e_{j}', 1_{E_{j}'}>=\lambda_{E_{j}'}(x_{j}')\circ \beta_{j}$. Applying respectively $f_{i}^{\ast}$ and $f_{j}'^{\ast}$ to these two identities, we obtain that $<L^{{\cal E}}_{h}\circ p_{i, j}, 1_{F_{i, j}}>=\lambda_{F_{i, j}}(x_{i, j})\circ f_{i}^{\ast}(\alpha_{i})$ and $<L^{{\cal E}'}_{h}\circ p_{i, j}, 1_{F_{i, j}}>=\lambda_{F_{i, j}}(x_{i, j}')\circ f_{j}'^{\ast}(\beta_{j})$, where $x_{i, j}={\cal A}_{f_{i}}(x_{i})$ and $x_{i, j}'={\cal A}_{f_{j}'}(x_{j}')$, from which it follows that the condition $L^{\cal E}_{h}\circ p_{i, j}=L^{{\cal E}'}_{h}\circ p_{i, j}$ is equivalent to the condition           
$\lambda_{F_{i, j}}(x_{i, j})\circ f_{i}^{\ast}(\alpha_{i})=\lambda_{F_{i, j}}(x_{i, j}')\circ f_{j}'^{\ast}(\beta_{j})$. This condition can be proved by using condition $(ii)$ of the theorem, the fact that the same identity holds with the cocone $\mu$ in place of $\lambda$ and Remark \ref{rem3}(d). More specifically, the identities $<h\circ e_{i}, 1_{E_{i}}>=\mu_{E_{i}}(x_{i})\circ \alpha_{i}$ and $<h\circ e_{j}', 1_{E_{j}'}>=\mu_{E_{j}'}(x_{j}')\circ \beta_{i}$ imply, applying respectively $f_{i}^{\ast}$ and $f_{j}'^{\ast}$ to them, the identity $\mu_{F_{i, j}}(x_{i, j})\circ f_{i}^{\ast}(\alpha_{i})=\mu_{F_{i, j}}(x_{i, j}')\circ f_{j}'^{\ast}(\beta_{j})$. From condition $(ii)$ of the theorem, it thus follows that for any $i, j \in I\times J$, there exists an epimorphic family $\{g_{k}^{i, j}:G_{k}^{i, j} \to F_{i, j} \textrm{ | } k\in K_{i,j}\}$ such that the pair $(({\cal A}_{g_{k}^{i, j}(x_{i, j})}, f_{i}^{\ast}(\alpha_{i})\circ g_{k}^{i, j}), ({\cal A}_{g_{k}^{i, j}(x_{i, j}')}, f_{j}'^{\ast}(\beta_{j})\circ g_{k}^{i, j}))$ belongs to the relation ${\cal R}_{G_{k}^{i, j}}^{D}$. This in turn implies (by Remark \ref{rem3}(d)) that $\lambda_{G_{k}^{i,j}}({\cal A}_{g_{k}^{i, j}}(x_{i, j})) \circ (f_{i}^{\ast}(\alpha_{i})\circ g_{k}^{i, j})= \lambda_{G_{k}^{i,j}}({\cal A}_{g_{k}^{i, j}}(x_{i, j}')) \circ (f_{j}'^{\ast}(\beta_{j})\circ g_{k}^{i, j})$. But $\lambda_{G_{k}^{i,j}}({\cal A}_{g_{k}^{i, j}}(x_{i, j})) \circ (f_{i}^{\ast}(\alpha_{i})\circ g_{k}^{i, j})= (V\times g_{k}^{i, j})^{\ast}(\lambda_{F_{i, j}}(x_{i, j})\circ f_{i}^{\ast}(\alpha_{i}))$ and $\lambda_{G_{k}^{i,j}}({\cal A}_{g_{k}^{i, j}}(x_{i, j}')) \circ (f_{j}'^{\ast}(\beta_{j})\circ g_{k}^{i, j})= (V\times g_{k}^{i, j})^{\ast}(\lambda_{F_{i, j}}(x_{i, j}')\circ f_{j}'^{\ast}(\beta_{j}))$, whence $(V\times g_{k}^{i, j})^{\ast}(\lambda_{F_{i, j}}(x_{i, j})\circ f_{i}^{\ast}(\alpha_{i}))=(V\times g_{k}^{i, j})^{\ast}(\lambda_{F_{i, j}}(x_{i, j}')\circ f_{j}'^{\ast}(\beta_{j}))$, equivalently (since the $g_{k}^{i, j}$ are jointly epimorphic) $\lambda_{F_{i, j}}(x_{i, j})\circ f_{i}^{\ast}(\alpha_{i})=\lambda_{F_{i, j}}(x_{i, j}')\circ f_{j}'^{\ast}(\beta_{j})$, as required.         

It is clear that the assignment $h \to L(h)$ defined above is natural in $E\in {\cal E}$; it therefore remains to show that the resulting arrow $l:U\to V$ satisfies the required property that for any object $E\in {\cal E}$ and any object $x$ of the category ${\cal A}_{E}$, $(l\times 1_{E})\circ \mu_{E}(x)=\lambda_{E}(x)$. Take $s$ to be the canonical arrow $s:D_{E}(x)\to E$ in $\cal E$ and set $x'={\cal A}_{s}(x)$, $E'=D_{E}(x)$; then $x\in {\cal A}_{E'}$ and, considered the pullback square
\[  
\xymatrix {
D_{E'}(x') \ar[d]^{t} \ar[r]^{r} & D_{E}(x) \ar[d]^{s} \\
E' \ar[r]^{s} & E,} 
\]
the unique arrow $\alpha:E'\to D_{E'}(x')$ such that $r\circ \alpha=1_{E'}$ and $t\circ \alpha=1_{E'}$ satisfies the property that $<h, 1_{E'}>=\mu_{E'}(x')\circ \alpha$. So the epimorphic family $\{1_{E'}:E'\to E'\}$ satisfies condition $(i)$ of the theorem with respect to the arrow $h=\pi_{U}\circ \mu_{E}(x)$, where $\pi_{U}:U\times E \to U$ is the canonical projection, and hence $l\circ h=L(h)=\lambda_{E'}(x')\circ \alpha$; therefore, $(l\times 1_{E})\circ \mu_{E}(x)=\lambda_{E}(x)$, as required.  
 
The proof of the theorem is now complete.
\end{proofs}

\begin{remark}
\begin{enumerate}[(a)]

\item The proof of the theorem shows that the sufficiency of the conditions of the theorem holds more in general for any (i.e., not necessarily small) $\cal E$-indexed category $\underline{\cal A}_{\cal E}$;

\item The conditions in the statement of the theorem are both stable under localization; that is, if the cone $\mu$ over the diagram $D$ satisfies them then for any $E\in {\cal E}$ the cone $\mu\slash E$ over the diagram $D\slash E$ does. 

\end{enumerate}

\end{remark}

\begin{proposition}
Let ${\underline{\cal B}}_{\cal E}$ be a $\cal E$-final subcategory of an indexed category ${\underline{\cal A}}_{\cal E}$, $D:{\underline{\cal A}}_{\cal E}\to {\underline{\cal E}}_{\cal E}$ a $\cal E$-indexed functor and $\mu$ a $\cal E$-indexed cocone over $D$. Then $\mu$ satisfies the conditions of Theorem \ref{charact} with respect to the diagram $D$ if and only if $\mu \circ i$ satisfies them with respect to the diagram $D\circ i$, where $i$ is the canonical embedding ${\underline{\cal B}}_{\cal E} \hookrightarrow {\underline{\cal A}}_{\cal E}$.
\end{proposition}

\begin{proofs}
Suppose that $\mu$ satisfies the conditions of Theorem \ref{charact} with respect to the diagram $D$. The fact that $\mu \circ i$ satisfies condition $(i)$ of the theorem with respect to the diagram $D\circ i$ immediately follows from the fact that $\mu$ does with respect to the diagram $D$, by using the fact that for each object $x_{i}\in {\cal A}_{F_{i}}$ there exists an epimorphic family $\{g_{k}^{i}:G^{i}_{k} \to F_{i} \textrm{ | } k\in K_{i}\}$ such that ${\cal A}_{g_{k}^{i}}(x_{i})$ lies in ${\cal B}_{G^{i}_{k}}$ and for any arrow $\alpha_{i}:1_{{\cal E}\slash F_{i}} \to D_{F_{i}}(x_{i})$ in the topos ${\cal E}\slash F_{i}$, its image under the pullback functor ${g_{k}^{i}}^{\ast}$ is an arrow $1_{{\cal E}\slash G^{i}_{k}} \to D_{G^{i}_{k}}({\cal A}_{g^{i}_{k}}(x_{i}))$ in the topos ${\cal E}\slash G^{i}_{k}$. Let us now show that $\mu \circ i$ satisfies condition $(ii)$ with respect to the diagram $D\circ i$.

Let us first establish the following fact $(\ast)$: for any pairs $(x, y)$ and $(x', y')$ in ${\cal I}_{R}^{D}$ (see the beginning of this section for the notation), if $((x, y), (x', y'))\in {\cal R}_{R}^{D}$, there exists an epimorphic family $\{r_{i}:R_{i}\to R \textrm{ | } i\in I\}$ in $\cal E$ such that $(({\cal A}_{r_{i}}(x), r_{i}^{\ast}(y)), ({\cal A}_{r_{i}}(x'), r_{i}^{\ast}(y')))\in {\cal R}_{R_{i}}^{D\circ i}$.  
 
If  $((x, y), (x', y'))\in {\cal R}_{R}^{D}$, there exists a finite sequence
\[
(x, y)=(x_{0}, y_{0}), \ldots, (x_{n}, y_{n})=(x', y')
\]
of pairs in ${\cal I}_{R}^{D}$ such that for any $j\in \{0, \ldots, n-1\}$, either there exists an arrow $f_{j}:x_{j}\to x_{j+1}$ in ${\cal A}_{R}$ such that $y_{j+1}=D_{R}(f_{j})\circ y_{j}$ or there exists an arrow $f_{j}:x_{j+1}\to x_{j}$ in ${\cal A}_{R}$ such that $y_{j}=D_{R}(f_{j})\circ y_{j+1}$. By applying the definition of $\cal E$-final subcategory a finite number of times (using the fact that the fibered product of epimorphic families is again an epimorphic family), we can find an epimorphic family $\{r_{i}:R_{i}\to R \textrm{ | } i\in I\}$ in $\cal E$ such that for every $i\in I$ and any $j\in \{0, \ldots, n-1\}$, the arrow ${\cal A}_{r_{i}}(f_{j})$ lies in ${\cal B}_{R_{i}}$. From this it immediately follows that $(({\cal A}_{r_{i}}(x), r_{i}^{\ast}(y)), ({\cal A}_{r_{i}}(x'), r_{i}^{\ast}(y')))\in {\cal R}_{R_{i}}^{D\circ i}$, as required.  
       
Using fact $(\ast)$, the proof of the fact that $\mu \circ i$ satisfies condition $(ii)$ of the theorem with respect to the diagram $D\circ i$ follows straightforwardly from the fact that $\mu$ does with respect to the diagram $D$. 

Conversely, let us suppose that $\mu \circ i$ satisfies the conditions of the theorem with respect to $D\circ i$ and deduce that $\mu$ does with respect to $D$. The fact that the validity of condition $(i)$ for $\mu \circ i$ with respect to $D\circ i$ implies the validity of condition $(i)$ for $\mu$ with respect to $D$ is obvious. Concerning condition $(ii)$ for $(\mu, D)$, this can be deduced from condition $(ii)$ for $(\mu\circ i, D\circ i)$ by using the fact that, given $(x, y)$ and $(x', y')$ as in the statement of the condition for $(\mu, D)$, there exists an epimorphic family 
$\{e_{j}:E_{j}\to E \textrm{ | } j\in J\}$ such that ${\cal A}_{e_{j}}(x)$ and ${\cal A}_{e_{j}}(x')$ lie in ${\cal B}_{E_{j}}$ for all $j\in J$.

The proof of the proposition is now complete.
\end{proofs}

\begin{remark}\label{remprop}
The proposition shows in particular that (the necessity of the conditions of) Theorem \ref{charact} not only holds for $\cal E$-indexed  categories $\underline{\cal A}_{\cal E}$ which are $\cal E$-small (i.e., which are equivalent to the $\cal E$-externalization of an internal category in $\cal E$.), but for any $\cal E$-indexed category which is a $\cal E$-final subcategory of a $\cal E$-small $\cal E$-indexed category (for instance, to the $\cal E$-indexed categories of the form $\underline{\int F}_{\cal E}$, for a functor $F:{\cal C}^{\textrm{op}}\to {\cal E}$ - cf. Theorem \ref{finalfunc}). 
\end{remark}

\begin{proposition}
Let $D:\underline{\cal A}_{\cal E} \to \underline{{\cal E}}_{\cal E}$ an indexed functor, where $\underline{\cal A}_{\cal E}$ is a $\cal E$-filtered category, $R$ an object of $\cal E$, $(x, y), (x', y')$ pairs in ${\cal I}_{R}^{D}$. Then there exists an epimorphic family $\{r_{i}:R_{i} \to R \textrm{ | } i\in I\}$ in $\cal E$ such that for all $i\in I$, $(({\cal A}_{r_{i}}(x), r_{i}^{\ast}(y)), ({\cal A}_{r_{i}}(x'), r_{i}^{\ast}(y')))\in {\cal R}^{D}_{R_{i}}$ if and only if there exists an epimorphic family $\{e_{j}:E_{j} \to R \textrm{ | } j\in J\}$ in $\cal E$ with the property that for any $j\in J$ there exist arrows $f_{j}:{\cal A}_{e_{j}}(x)\to z$ and $g_{j}:{\cal A}_{e_{j}}(x')\to z$ in the category ${\cal A}_{E_{j}}$ such that $D_{E_{j}}(f_{j})\circ y=D_{E_{j}}(g_{j})\circ y'$. 
\end{proposition}

\begin{proofs}
The `if' direction is obvious, so it remains to prove the `only if' one. It clearly suffices to prove, by induction on the length of a sequence $(x_{1}, y_{1}), \ldots, (x_{n}, y_{n})$ of pairs in ${\cal I}_{R}^{D}$ with the property that for any $i\in \{1, \ldots, n-1\}$, either there exists an arrow $f_{j}:x_{j}\to x_{j+1}$ in ${\cal A}_{R}$ such that $y_{j+1}=D_{R}(f_{j})\circ y_{j}$ or there exists an arrow $f_{j}:x_{j+1}\to x_{j}$ in ${\cal A}_{R}$ such that $y_{j}=D_{R}(f_{j})\circ y_{j+1}$, that there exists an epimorphic family $\{r_{i}:R_{i} \to R \textrm{ | } i\in I\}$ in $\cal E$ such that for all $i\in I$, $(({\cal A}_{r_{i}}(x_{1}), r_{i}^{\ast}(y_{1})), ({\cal A}_{r_{i}}(x_{n}), r_{i}^{\ast}(y_{n})))\in {\cal R}^{D}_{R_{i}}$. The case $n=1$ is obvious. Suppose now that the claim is valid for $n$ and prove that it holds for $n+1$. There exists an epimorphic family $\{e_{j}:E_{j} \to R \textrm{ | } j\in J\}$ in $\cal E$ such that for any $j\in J$ there exist arrows $f_{j}:{\cal A}_{e_{j}}(x_{1})\to z$ and $g_{j}:{\cal A}_{e_{j}}(x_{n})\to z$ in the category ${\cal A}_{E_{j}}$ such that $D_{E_{j}}(f_{j})\circ y_{n}=D_{E_{j}}(g_{j})\circ y'_{n}$. On the other hand, there exists either an arrow $g:x_{n}\to x_{n+1}$ in ${\cal A}_{R}$ such that $D_{R}(g)\circ y_{n}=y_{n+1}$ or an arrow $h:x_{n+1}\to x_{n}$ in ${\cal A}_{R}$ such that $D_{R}(g)\circ y_{n+1}=y_{n}$. In the latter case, the thesis follows straightforwardly; so we can concentrate on the case when there exists an arrow $g:x_{n}\to x_{n+1}$ in ${\cal A}_{R}$ such that $D_{R}(g)\circ y_{n}=y_{n+1}$.    
 
By using the definition of $\cal E$-final subcategory, we obtain an epimorphic family $\{e_{j}:E_{j}\to R \textrm{ | } j\in J\}$ in $\cal E$ and for each $j\in J$ two arrows $f_{j}:{\cal A}_{e_{j}}(x_{1})\to z$ and $g_{j}:{\cal A}_{e_{j}}(x_{n})\to z$ such that $D_{E_{j}}(f_{j})\circ e_{j}^{\ast}(y_{1})=D_{E_{j}}(g_{j})\circ e_{j}^{\ast}(y_{n})$. Similarly, we obtain, for each $j\in J$, an epimorphic family $\{f_{k}^{j}:F_{k}^{j} \to E_{j} \textrm{ | } k\in K_{j}\}$ in $\cal E$ and for each $k\in K_{j}$ arrows $m_{k}^{j}:{\cal A}_{f_{k}^{j}}(z)\to w_{k}^{j}$ and $n_{k}^{j}:{\cal A}_{f_{k}^{j}}({\cal A}_{e_{j}}(x_{n+1}))\to w_{k}^{j}$ such that $D_{F_{k}^{j}}(m_{k}^{j})\circ D_{F_{k}^{j}}({\cal A}_{f_{k}^{j}}(g_{j}))\circ (f_{k}^{j}\circ e_{j})^{\ast}(y_{n})=D_{F_{k}^{j}}(n_{k}^{j})\circ D_{F_{k}^{j}}({\cal A}_{f_{k}^{j}}(g))\circ (f_{k}^{j}\circ e_{j})^{\ast}(y_{n+1})$. Now, by definition of $\cal E$-final subcategory, for each $j\in J$ and $k\in K_{j}$ there exists an epimorphic family $\{g_{l}^{k,j}:G_{l}^{k,j} \to F_{k}^{j} \textrm{ | } l\in L_{k, j}\}$ and for each $l\in L_{k, j}$ an arrow $p_{l}^{k, l}:{\cal A}_{g_{l}^{k, j}}(w_{k}^{j})\to v_{l}^{k, j}$ such that $p_{l}^{k, l}\circ m_{k}^{j}\circ {\cal A}_{f_{k}^{j}}(g_{j})=p_{l}^{k, j}\circ n_{k}^{j}\circ {\cal A}_{f_{k}^{j}}({\cal A}_{e_{j}}(g))$.
For any $j\in J, k\in K_{j}$ and $l\in L_{k, j}$, set $a_{j,k, l}:{\cal A}_{e_{j}\circ f_{k}^{j}\circ g_{l}^{k, j}}(x_{1})\to v_{l}^{k, j}$ equal to $p_{l}^{k, j}\circ {\cal A}_{g_{l}^{j, k}}(m_{k}^{j})\circ {\cal A}_{g_{l}^{k, j} \circ f_{k}^{j}}(f_{j})$ and $b_{j,k, l}:{\cal A}_{e_{j}\circ f_{k}^{j}\circ g_{l}^{k, j}}(x_{n})\to v_{l}^{k, j}$ equal to $p_{l}^{k, j}\circ {\cal A}_{g_{l}^{j, k}}(n_{k}^{j})\circ {\cal A}_{g_{l}^{k, j} \circ f_{k}^{j} \circ e_{j}}(g)$. It is readily seen that $D_{G_{l}^{k, j}}(a_{j, k, l})\circ (e_{j}\circ f_{k}^{j}\circ g_{l}^{k, j})^{\ast}(y_{1})=D_{G_{l}^{k, j}}(b_{j, k, l})\circ (e_{j}\circ f_{k}^{j}\circ g_{l}^{k, j})^{\ast}(y_{n+1})$, whence the required condition is satisfied (by taking as epimorphic family $\{ e_{j}\circ f_{k}^{j}\circ g_{l}^{k, j} \textrm{ | } j\in J, k\in K_{j}, l\in L_{k, j}\}$).   
\end{proofs}

Combining the proposition with Theorem \ref{charact}, we immediately obtain the following result.

\begin{corollary}\label{corchar}
Let $D:\underline{\cal A}_{\cal E} \to \underline{{\cal E}}_{\cal E}$ an indexed functor, where $\underline{\cal A}_{\cal E}$ is a $\cal E$-filtered $\cal E$-final subcategory of a small $\cal E$-indexed category. Then a cocone $\mu$ over $D$ with vertex $U$ is a indexed colimiting cocone for $D$ if and only if the following conditions are satisfied:
\begin{enumerate}[(i)]

\item For any object $E$ of $\cal E$ and arrow $h:F \to U$ in the topos ${\cal E}$, there exists an epimorphic family $\{f_{i}:F_{i}\to F \textrm{ | } i\in I\}$ in $\cal E$) and for each $i\in I$ an object $x_{i}\in {\cal A}_{F_{i}}$ and an arrow $\alpha_{i}:1_{{\cal E}\slash F_{i}} \to D_{F_{i}}(x_{i})$ in the topos ${\cal E}\slash F_{i}$ such that $<h\circ f_{i}, 1_{F_{i}}>=\mu_{F_{i}}(x_{i})\circ \alpha_{i}$ as arrows $1_{{\cal E}\slash F_{i}} \to F_{i}^{\ast}(U)$ in ${\cal E}\slash F_{i}$;

\item For any pairs $(x, y)$ and $(x', y')$, where $x$ and $x'$ are objects of ${\cal A}_{E}$, $y$ is an arrow $(f:F\to E) \to D_{E}(x)$ in ${\cal E}\slash E$ and $y'$ is an arrow $(f:F\to E) \to D_{E}(x')$ in ${\cal E}\slash E$, we have that $\mu_{E}(x)\circ y=\mu_{E}(x')\circ y'$ if and only if there exists an epimorphic family $\{f_{i}:F_{i}\to F \textrm{ | } i\in I\}$ in $\cal E$ and for each $i\in I$ arrows $g_{i}:{\cal A}_{f\circ f_{i}}(x) \to z_{i}$ and $h_{i}:{\cal A}_{f\circ f_{i}}(x') \to z_{i}$ in the category ${\cal A}_{F_{i}}$ such that $D_{F_{i}}(g_{i})\circ f_{i}^{\ast}(y)=D_{F_{i}}(h_{i})\circ f_{i}^{\ast}(y')$.  
\end{enumerate}

Moreover, if $\mu$ is colimiting for $D$ then the following `joint embedding property' holds: for any pairs $(x, y)$ and $(x', y')$, where $x$ and $x'$ are objects of ${\cal A}_{E}$, $y$ is an arrow $(f:F\to E) \to D_{E}(x)$ in ${\cal E}\slash E$ and $y'$ is an arrow $(f:F\to E) \to D_{E}(x')$ in ${\cal E}\slash E$, there exists an epimorphic family $\{e_{i}:F_{i}\to F \textrm{ | } i\in I\}$ in $\cal E$ and for each $i\in I$ arrows $g_{i}:{\cal A}_{f\circ e_{i}}(x) \to z_{i}$ and $h_{i}:{\cal A}_{f\circ e_{i}}(x') \to z_{i}$ in the category ${\cal A}_{F_{i}}$ such that $D_{F_{i}}(g_{i})\circ f_{i}^{\ast}(y)=D_{F_{i}}(h_{i})\circ f_{i}^{\ast}(y')$ and the  diagram 

\[  
\xymatrix {
D_{F_{i}}({\cal A}_{f\circ e_{i}}(x)) \ar[dr]_{D_{F_{i}}(g_{i})} \ar[rrrd]^{\xi_{(F_{i}, {\cal A}_{f\circ e_{i}}(x))}} & &  & \\
& D_{F_{i}}(z_{i}) \ar[rr]^{\xi_{(E_{i}, z_{i})}} & &  colim(D)\times F_{i} \\
D_{F_{i}}({\cal A}_{f\circ e_{i}}(x')) \ar[ur]^{D_{F_{i}}(h_{i})} \ar[rrru]_{\xi_{(F_{i}, {\cal A}_{f\circ e_{i}}(x'))}} & & &,} 
\]
where $\xi_{(E, x)}:D_{E}(x)\to colim(D)\times E$ are the colimit arrows, commutes. 
\end{corollary}

\begin{remark}\label{wlof}
\begin{enumerate}[(a)]
\item In the statement of the corollary, we can suppose without loss of generality the object $f:F\to E$ of the topos  ${\cal E}\slash E$ to be equal to the terminal $1_{E}:E\to E$. Indeed, by Remark \ref{rem3}, for any object $x$ of ${\cal A}_{E}$ and any arrow $y:(f:F\to E)\to D_{E}(x)$ in the topos ${\cal E}\slash E$, $\mu_{E}(x)\circ y=\mu_{E}(x')\circ y'$ if and only if $\mu_{F}({\cal A}_{f}(x))\circ y_{f}=\mu_{F}({\cal A}_{f}(x'))\circ y_{f}$.

\item For any pairs $(x, y)$ and $(x', y')$, where $x$ and $x'$ are objects of ${\cal A}_{E}$, $y$ is an arrow $(f:F\to E) \to D_{E}(x)$ in ${\cal E}\slash E$, $y'$ is an arrow $(f:F\to E) \to D_{E}(x')$ in ${\cal E}\slash E$, and any epimorphic family $\{g_{j}:G_{j}\to F \textrm{ | } j\in J\}$ in $\cal E$, condition $(ii)$ of Corollary \ref{corchar} is satisfied by the pair $((x, y), (x', y'))$ if and only if it is satisfied by the pair $((x, y\circ g_{j}), (x', y'\circ g_{j}))$ for all $j\in J$.
\end{enumerate}
\end{remark}

Let us now apply Corollary \ref{corchar} in the context of a flat functor $F:{\cal C}^{\textrm{op}}\to {\cal E}$ and a functor $P:{\cal C}\to \Set$, where $\cal C$ is a small category and $\cal E$ is a Grothendieck topos. Consider the restriction of the functor $(\underline{P_{\cal E}}_{\cal E}\circ \underline{\pi_{\overline{F}}^{f}}_{\cal E})\cong {\underline{{P_{\cal E}} \circ \pi_{\overline{F}}^{f}}}_{\cal E}$ (where $P_{\cal E}$ is the internal diagram $[{\mathbb C}, {\cal E}]$ given by $\overline{\gamma_{\cal E}^{\ast} \circ P}$) to the $\cal E$-indexed subcategory $\underline{\int F}_{\cal E}$ of $\underline{\int F}^{f}_{\cal E}$. Recall that, by Theorem \ref{altdes}, the $\cal E$-indexed functor $(\underline{P_{\cal E}}_{\cal E}\circ \underline{\pi_{\overline{F}}^{f}}_{\cal E})$ is naturally isomorphic to the composite functor $\gamma_{\cal E}\circ z_{\cal E}$, where $z_{\cal E}:\underline{\int F}_{\cal E} \to \underline{\Set}_{\cal E}$ is the $\cal E$-indexed functor defined by: ${z_{\cal E}}_{E}((c, x))=P(c)$ and $z_{\cal E}E(f)=P(f)$ (for any $E\in {\cal E}$, object $(c, x)$ and arrow $f$ in the category ${\int F}_{E}$).
 
As observed in section \ref{finsub}, a $\cal E$-indexed cocone $\mu$ over this diagram with vertex $U$ can be identified with a family of functions $\mu_{(c, x)}:P(c)\to Hom_{\cal E}(E, U)$ indexed by the pairs $(c, x)$ consisting of an object $c$ of $\cal C$ and a generalized element $x:E\to F(c)$ satisfying the following properties:

\begin{enumerate}[(i)]
\item For any generalized element $x:E\to F(c)$ and any arrow $f:d\to c$ in $\cal C$, $\mu_{(c, x)} \circ P(f)=\mu_{(c, F(f)\circ x)}$;

\item For any generalized element $x:E\to F(c)$ and any arrow $e:E' \to E$, $\mu_{(c, x\circ e)}=Hom_{\cal E}(e, U) \circ \mu_{(c, x)}$.
\end{enumerate}  

The following theorem characterizes the cocones $\mu$ which are colimiting.

\begin{theorem}\label{thmres}
Let $\cal C$ be a small category, $\cal E$ a Grothendieck topos, $P:{\cal C}\to \Set$ a functor and $F:{\cal C}^{\textrm{op}}\to {\cal E}$ a flat functor. Then a $\cal E$-indexed cocone $\mu=\{\mu_{(c, x)}:P(c)\to Hom_{\cal E}(E, U) \textrm{ | } c\in {\cal C}, \textrm{ } x:E\to F(c) \textrm{ in } {\cal E} \}$ with vertex $U$ over the diagram given by the restriction of the functor $(\underline{P_{\cal E}}_{\cal E}\circ \underline{\pi_{\overline{F}}^{f}}_{\cal E})\cong {\underline{{P_{\cal E}} \circ \pi_{\overline{F}}^{f}}}_{\cal E}$ (where $P_{\cal E}$ is the internal diagram $[{\mathbb C}, {\cal E}]$ given by $\overline{\gamma_{\cal E}^{\ast} \circ P}$) to the $\cal E$-indexed subcategory $\underline{\int F}_{\cal E}$ is colimiting if and only if the following conditions are satisfied:  

\begin{enumerate}[(i)] 
\item For any generalized element $x:E\to U$ there exists an epimorphic family $\{e_{i}:E_{i}\to E \textrm{ | } i\in I\}$ in $\cal E$, for each index $i\in I$ an object $c_{i}$ in $\cal C$, a generalized element $x_{i}:E_{i}\to F(c_{i})$ and an element $y_{i}\in P(c_{i})$ such that $\mu_{(c_{i}, x_{i})}(y_{i})=x\circ e_{i}$; 

\item For any pairs $(a, x)$ and $(b, x')$, where $a$ and $b$ are objects of $\cal C$ and $x:E\to F(a)$, $x':E\to F(b)$ are generalized elements, and any elements $y\in P(a)$ and $y'\in P(b)$, $\mu_{(a, x)}(y)=\mu_{(b, x')}(y')$ if and only if there exists an epimorphic family $\{e_{i}:E_{i}\to E \textrm{ | } i\in I\}$ in $\cal E$ in $\cal E$ and for each index $i\in I$ an object $c_{i}$ of $\cal C$, arrows $f_{i}:a\to c_{i}$, $g_{i}:b\to c_{i}$ in $\cal C$ and a generalized element $x_{i}:E_{i}\to F(c_{i})$ such that $<x, x'>\circ e_{i}=<F(f_{i}), F(g_{i})>\circ x_{i}$ for all $i\in I$ and $P(f_{i})(y)=P(g_{i})(y')$.
\end{enumerate}

Moreover, if $\mu$ is colimiting the following `joint embedding property' holds: for any pairs $(a, x)$ and $(b, x')$, where $a$ and $b$ are objects of $\cal C$ and $x:E\to F(a)$, $x':E\to F(b)$ are generalized elements, there exists an epimorphic family $\{e_{i}:E_{i}\to E \textrm{ | } i\in I\}$ in $\cal E$ and for each index $i\in I$ an object $c_{i}$ of $\cal C$, arrows $f_{i}:a\to c_{i}$, $g_{i}:b\to c_{i}$ in $\cal C$ and a generalized element $x_{i}:E_{i}\to F(c_{i})$ such that $<x, x'>\circ e_{i}=<F(f_{i}), F(g_{i})>\circ x_{i}$ for all $i\in I$ and (by Lemma \ref{lemmaxi}) the following diagram commutes:

\[  
\xymatrix {
P(a) \ar[dr]^{P(f)} \ar[rr]^{\mu_{(a, x)}} & & Hom_{\cal E}(E, U) \ar[dr]^{Hom_{\cal E}(e_{i}, M)} & \\
& P(c_{i}) \ar[rr]^{\mu_{(c_{i}, x_{i})}} & &  Hom_{\cal E}(E_{i}, U) \\
P(b) \ar[ur]_{P(g)} \ar[rr]^{\mu_{(b, x')}} & & Hom_{\cal E}(E, U) \ar[ur]_{Hom_{\cal E}(e_{i}, U)} &} 
\]

In fact, the epimorphic family $\{e_{i}:E_{i}\to E \textrm{ | } i\in I\}$ can be taken to be the pullback of the family of arrows $<F(f), F(g)>:F(c)\to F(a)\times F(b)$ (for all spans $(f:a\to c,\textrm{ }g:b\to c)$ in the category $\cal C$) along the arrow $<x,x'>:E\to F(a)\times F(b)$.  
\end{theorem} 

\begin{proofs}
The theorem can be deduced from Corollary \ref{corchar} by reasoning as follows. Concerning condition $(i)$, its equivalence with condition $(i)$ of Corollary \ref{corchar} follows from Remark \ref{wlof}(a) and the fact that for any object $(z, c)$ of the category ${\int F}_{E}$ and any arrow $\alpha:1_{{\cal E}\slash E} \to D_{E}((z, c))$ in the topos ${\cal E}\slash E$, denoting by $!_{f_{r}}:{f_{r}:F_{r} \to E} \to 1_{{\cal E}\slash E}$ the pullback of the coproduct arrow $u_{r}:1_{{\cal E}\slash E} \to \gamma_{{\cal E}\slash E}^{\ast}(P(c))\cong D_{E}((z, c))$ in ${\cal E}\slash E$ (for each $r\in P(c)$), the arrows $!_{f_{r}}$ define, for $r\in P(c)$, an epimorphic family in ${\cal E}\slash E$ on the object $1_{{\cal E}\slash E}$; indeed, this shows that we can suppose without loss of generality $\alpha$ to factor through one of the coproduct arrows $u_{r}$ and hence to correspond to an element $y\in P(c)$. 

Concerning condition $(ii)$, by Remark \ref{wlof} we can suppose without loss of generality that the pairs $(x, y)$ and $(x', y')$, where $x$ (resp. $x'$) is an object of ${\int F}_{E}$, $y$ is an arrow $(f:F\to E) \to D_{E}(x)$ in ${\cal E}\slash E$ and $y'$ is an arrow $(f:F\to E) \to D_{E}(x')$ in ${\cal E}\slash E$, have respectively the form $((x, c), y)$, where  $y=u_{r}\circ !_{{f_{r}:F_{r} \to E}}$ and $((x', c'), y')$, where $y'=u_{r'}\circ !_{f_{r'}}$ (with the above notation), for some $r, r'\in P(c)$. Thus we have that $\mu_{E}(x)\circ y=\mu_{E}(x')\circ y'$ if and only if $(Hom_{\cal E}(f_{r}, U)\circ \mu_{(x, c)})(r)=(Hom_{\cal E}(f_{r'}, U)\circ \mu_{(x', c')})(r')$; but $(Hom_{\cal E}(f_{r}, U)\circ \mu_{(x, c)})(r)=\mu_{(c, x\circ f_{r})}(r)$ and $(Hom_{\cal E}(f_{r'}, U)\circ \mu_{(x', c')})(r')=\mu_{(c', x'\circ f_{r'})}(r')$, whence $(Hom_{\cal E}(f_{r}, U)\circ \mu_{(x, c)})(r)=(Hom_{\cal E}(f_{r'}, U)\circ \mu_{(x', c')})(r')$ if and only if $\mu_{(c, x\circ f_{r})}(r)=\mu_{(c', x\circ f_{r'})}(r')$.        
\end{proofs}

\subsection{Explicit calculation of set-indexed colimits}\label{expcolim}

Recall from \cite{MM} (cf. p. 355) that the colimit of a functor $H:{\cal I}\to {\cal E}$ with values in a Grothendieck topos $\cal E$ can always be realized as the coequalizer $q: \mathbin{\mathop{\textrm{ $\coprod$}}\limits_{i\in I}}H(i) \to colim(H)$ of the pair of arrows 
\[
a,b:\mathbin{\mathop{\textrm{ $\coprod$}}\limits_{u:i\to j\textrm{ in } {\cal I}}}H(i) \to \mathbin{\mathop{\textrm{ $\coprod$}}\limits_{i\in I}}H(i)
\]
defined by: 
\[
a\circ \lambda_{u}=\kappa_{i}
\]
and
\[
b\circ \lambda_{u}=\kappa_{j}\circ H(u)
\]
(for every arrow $u:i\to j$ in $\cal I$), where $\lambda_{u}:H(dom(u))\to \mathbin{\mathop{\textrm{ $\coprod$}}\limits_{u:i\to j\textrm{ in } {\cal I}}}H(i)$ (for any $u:i\to j$ in $\cal I$) and $\kappa_{i}:H(i)\to \mathbin{\mathop{\textrm{ $\coprod$}}\limits_{i\in I}}H(i)$ (for any $i\in {\cal I}$) are the canonical coproduct arrows. 

Given two arrows $s,t:E\to \mathbin{\mathop{\textrm{ $\coprod$}}\limits_{i\in I}}H(i)$, for every objects $i, j\in {\cal I}$, we can consider the following pullback diagrams:

\[  
\xymatrix {
E^{s}_{i} \ar[d]^{p^{s}_{i}} \ar[rr]^{a^{s}_{i}} & & H(i) \ar[d]^{\kappa_{i}} & E^{t}_{j} \ar[d]^{p_{j}^{t}} \ar[rr]^{a_{j}^{t}} & & H(j) \ar[d]^{\kappa_{j}} \\
E \ar[rr]_{s} & & \mathbin{\mathop{\textrm{ $\coprod$}}\limits_{i\in I}}H(i) & E \ar[rr]_{t} & & \mathbin{\mathop{\textrm{ $\coprod$}}\limits_{i\in I}}H(i) }
\]
Let us also consider, for each $i, j\in {\cal J}$, the pullback square
\[  
\xymatrix {
E^{s,t}_{i,j} \ar[d]^{q^{t}_{i,j}} \ar[r]^{q^{s}_{i,j}} & E_{i}^{s} \ar[d]^{p^{s}_{i}}\\
E_{j}^{t}  \ar[r]_{p^{t}_{j}} & E} 
\]  
and denote by $r^{s,t}_{i,j}:E^{s,t}_{i,j} \to E$ the arrow $p^{s}_{i} \circ q^{s}_{i,j}= p^{t}_{j} \circ q^{t}_{i,j}$. Notice that the family of arrows $\{r^{s,t}_{i,j} \textrm{ | } i,j\in {\cal I}\}$ is epimorphic.

We shall use these notations throughout the section.

The following lemma provides a characterization of the coequalizer of the functor $H$, which will be useful to us later on.

\begin{lemma}\label{reduction}
Let $p:\mathbin{\mathop{\textrm{ $\coprod$}}\limits_{i\in I}}H(i) \to A$ be an epimorphism. Then $p$ is isomorphic to the canonical map $q: \mathbin{\mathop{\textrm{ $\coprod$}}\limits_{i\in I}}H(i) \to colim(H)$ if and only if for any object $E$ of $\cal E$, any objects $i, j\in {\cal I}$ and any generalized elements $z:E\to H(i)$ and $w:E\to H(j)$, $(\kappa_{i}\circ z,\kappa_{j}\circ w)\in R_{E}$ if and only if $p\circ z=p\circ w$.
\end{lemma}

\begin{proofs}
As in a topos every epimorphism is the coequalizer of its kernel pair, any two epimorphisms $\alpha:B\to C$ and $\beta:B\to D$ in $\cal E$ are isomorphic in $B\slash {\cal E}$ if and only if their kernel pairs are isomorphic as subobjects of $B\times B$. In particular, $p$ is isomorphic to $q$ if and only if the kernel pair of $p$ is isomorphic to $R$. 

The kernel pair $S$ of $p$ satisfies the property that for any generalized elements $s, t:E\to \mathbin{\mathop{\textrm{ $\coprod$}}\limits_{i\in I}}H(i)$, $(s,t)\in S_{E}$ if and only if $p\circ s=p\circ t$. It thus remains to prove that the following two conditions are equivalent:
\begin{enumerate}[(1)]
\item for any generalized elements $s, t:E\to \mathbin{\mathop{\textrm{ $\coprod$}}\limits_{i\in I}}H(i)$, $p\circ s=p\circ t$ if and only if $(s, t)\in R_{E}$;

\item for any object $E$ of $\cal E$, any objects $i, j\in {\cal I}$ and any generalized elements $z:E\to H(i)$ and $w:E\to H(j)$, $p\circ \kappa_{i}\circ z=p\circ \kappa_{j}\circ w$ if and only if $(\kappa_{i}\circ z,\kappa_{j}\circ w)\in R_{E}$. 
\end{enumerate}

It is clear that $(1)$ implies $(2)$. To prove the converse direction, we observe that, as we remarked above, $(s, t)\in R_{E} $ if and only if for every $i, j\in {\cal I}$ $(\kappa_{i}\circ a_{i}^{s} \circ q^{s}_{i,j}, \kappa_{j} \circ a_{j}^{t} \circ q^{t}_{i,j})\in R_{E^{s,t}_{i,j}}$; but the latter condition is equivalent, under assumption $(2)$, to the requirement that for every $i, j\in {\cal I}$, $p\circ \kappa_{i}\circ a_{i}^{s} \circ q^{s}_{i,j} = p\circ \kappa_{j} \circ a_{j}^{t} \circ q^{t}_{i,j}$, where $s\circ r^{s,t}_{i,j}=\kappa_{i}\circ a_{i}^{s} \circ q^{s}_{i,j}$ and $t\circ r^{s,t}_{i,j}=\kappa_{j} \circ a_{j}^{t} \circ q^{t}_{i,j}$, that is, as the family of arrows $r^{s,t}_{i,j}$ is epimorphic, to the condition $p\circ s=p\circ t$, as required.    
\end{proofs}

\begin{lemma}\label{lemma1}
Let $B$ be an object of a topos $\cal E$. Then an equivalence relation in $\cal E$ on $B$ can be identified with an assignment $E \to R_{E}$ of an equivalence relation on the set $Hom_{\cal E}(E, B)$ satisfying the following properties:
\begin{enumerate}[(i)]
\item For any arrow $f:E\to E'$ in $\cal E$, if $(h,k)\in R_{E'}$ then $(h\circ f, k\circ f)\in R_{E}$;

\item For any epimorphic family $\{e_{i}:E_{i}\to E \textrm{ | } i\in I\}$ in $\cal E$ and any arrows $f,g\in Hom_{\cal E}(E, B)$, $(f\circ e_{i}, g\circ e_{i})\in R_{E_{i}}$ for all $i\in I$ implies $(f,g)\in R_{E}$.
\end{enumerate}
\end{lemma}

\begin{proofs}
Clearly, for any object $A$ of $\cal C$, the subobjects of $A$ can be identified with the $J_{\cal E}$-closed sieves on $A$, or equivalently with the $c_{J_{\cal E}}$-closed subobjects of the representable $Hom_{\cal E}(- A)$, where $c_{J_{\cal E}}$ is the closure operation on subobjects corresponding to the canonical Grothendieck topology $J_{\cal E}$ on $\cal E$. Applying this remark to $A=B\times B$ and noticing that the concept of equivalence relation is cartesian and hence preserved and reflected by cartesian fully faithful functors, in particular by the Yoneda embedding ${\cal E}\to [{\cal E}^{\textrm{op}}, \Set]$, we immediately deduce the thesis.  
\end{proofs}

Lemma \ref{lemma1} implies that $(s, t)\in R_{E} $ if and only if for every $i, j\in {\cal I}$ $(\kappa_{i}\circ a_{i}^{s} \circ q^{s}_{i,j}, \kappa_{j} \circ a_{j}^{t} \circ q^{t}_{i,j})\in R_{E^{s,t}_{i,j}}$; indeed, $s\circ r^{s,t}_{i,j}=\kappa_{i}\circ a_{i}^{s} \circ q^{s}_{i,j}$ and $t\circ r^{s,t}_{i,j}=\kappa_{j} \circ a_{j}^{t} \circ q^{t}_{i,j}$. 

This shows that, in order to completely describe the relation $R$, is suffices to consider elements of the form $\kappa_{i}\circ z$ and $\kappa_{j}\circ w$, where $z:E\to H(i)$ and $w:E\to H(j)$, and characterize when they belong to $R_{E}$.

\begin{lemma}\label{lemma2}
Let $\cal E$ be a topos, $f, g:A\to B$ arrows in $\cal E$ with coequalizer $q:B\to C$ and $R\mono B\times B$ the kernel pair of $q$. Then for any object $E$ of $\cal E$ and any elements $h,k\in Hom_{{\cal E}}(E, B)$, $<h,k>:E\to B\times B$ factors through $R$ if and only if there exists an epimorphic family $\{e_{i}: E_{i} \to E \textrm{ | } i\in I \}$ in $\cal E$ such that for each $i\in I$ $(h\circ e_{i}, k\circ e_{i})$ belongs to the equivalence relation on the set $Hom_{{\cal E}}(E_{i}, B)$ generated by the set of pairs of the form $(f\circ a, g\circ a)$, for a generalized element $a:E_{i}\to A$ in $\cal E$. 
\end{lemma}

\begin{proofs}
For any object $E$ of $\cal E$, consider the equivalence relation $R^{f,g}_{E}$ on the set $Hom_{\cal E}(E, B)$ consisting of the pairs $(h,k)$ with the property that there exists an epimorphic family $\{e_{i}: E_{i} \to E \textrm{ | } i\in I \}$ in $\cal E$ such that for each $i\in I$ $(h\circ e_{i}, k\circ e_{i})$ belongs to the equivalence relation on the set $Hom_{{\cal E}}(E_{i}, B)$ generated by the set of pairs of the form $(f\circ a, g\circ a)$, for a generalized element $a:E_{i}\to A$ in $\cal E$. The assignment $E\to R^{f,g}_{E}$ clearly satisfies the conditions of Lemma \ref{lemma1} and hence defines an equivalence relation on $B$ in $\cal E$. To prove that $R=R^{f,g}$, we argue as follows. Denoting by $\pi_{1}, \pi_{2}:R^{f,g}\to B$ the canonical projections, we have that for any $E\in {\cal E}$, $Hom(E, q)$ coequalizes the arrows $Hom(E , \pi_{1})$ and $Hom(E ,\pi_{2})$; therefore $q$ coequalizes $\pi_{1}$ and $\pi_{2}$ in $\cal E$. To prove that $q$ is actually the coequalizer of $\pi_{1}$ and $\pi_{2}$ in $\cal E$, we observe that any arrow $p:B\to D$ such that $p\circ \pi_{1}=p\circ \pi_{2}$ satisfies $p\circ f=p\circ g$ and hence, by definition of $q$, there exists a unique factorization of $p$ through $q$. As $q$ is the coequalizer of $R^{f,g}$ and in a topos all equivalence relations are the kernel pairs of their coequalizers, we conclude that $R^{f,g}=R$, as required.  
\end{proofs}

\begin{remark}\label{equivgen}
Under the hypotheses of the lemma, $R$ can be characterized as the equivalence relation on $B$ generated by the arrows $f$ and $g$, i.e. as the smallest equivalence relation on $B$ containing the image of the arrow $<f,g>:A\to B\times B$. Indeed, clearly $R$ contains the image of this arrow, and if $T$ is an equivalence relation on $B$ containing the image of the arrow $<f,g>$ then the coequalizer $z$ of $T$ factors through $q$ and hence the kernel pair of $q$, namely $R$, is contained in the kernel pair of $z$, namely $T$.   
\end{remark}

Recall that the internal language $\Sigma_{\cal E}$ of a topos $\cal E$ is the first-order signature consisting of one sort $\name{E}$ for any object $E$ of $\cal E$, one function symbol $\name{f}:\name{A}\to \name{B}$ for any arrow $f:A\to B$ in $\cal E$ in $\cal E$ and a relation symbol $\name{R}\mono \name{A_{1}}, \ldots, \name{A_{n}}$ for each subobject $R\mono A_{1}\times \cdots \times A_{n}$ in $\cal E$. There is a `tautological' $\Sigma_{\cal E}$-structure ${\cal S}_{\cal E}$ in $\cal E$, obtained by interpreting each sort as the corresponding object, each function symbol as the corresponding arrow and each relation symbol as the corresponding subobject.

We shall use the notation $x^{\name{C}}$ to denote a variable $x$ of sort $C$, and will omit the superscript when it can be unambiguously inferred from the context.

Lemma \ref{lemma2} can be reformulated in logical terms as follows.

\begin{lemma}\label{lemma3}
Let $\cal E$ be a topos, $f, g:A\to B$ arrows in $\cal E$ with coequalizer $q:B\to C$ and $R$ the kernel pair of $q$. Let $(z^{B}, z'^{B})\in G_{f,g}$ be an abbreviation for the formula
\[
z=z' \vee (\exists x)(\name{f}(x)=z \wedge \name{g}(x)=z') \vee  (\exists x')(\name{f}(x')=z' \wedge \name{g}(x')=z)
\]
over $\Sigma_{\cal E}$. Then the geometric bi-sequent
\[
(\name{R}(y, y') \dashv \vdash_{y^{B}, y'^{B}} \mathbin{\mathop{\textrm{\huge $\vee$}}\limits_{n\in {\mathbb N}}(\exists z_{1}\ldots \exists z_{n}(z_{1}=x \wedge z_{n}=x' \wedge \mathbin{\mathop{\textrm{\huge $\wedge$}}\limits_{1\lt i\lt n}}(z_{i-1}, z_{i})\in G_{f,g}))}
\]
is valid in the $\Sigma_{\cal E}$-structure ${\cal S}_{\cal E}$. 
\end{lemma}

\begin{proofs}
The interpretation of the formula $(z^{B}, z'^{B})\in G_{f,g}$ is the symmetric and reflexive closure of the relation $S_{f,g}$ on $B$ given by the image of the arrow $<f,g>:A\to B\times B$, that is the union of the diagonal subobject $\Delta:B\to B\times B$ of $B\times B$, $S_{f,g}\mono B\times B$ and the subobject $\tau \circ S_{f,g}$ of $B\times B$ given by the composite of $S_{f,g}$ with the exchange isomorphism $\tau:B\times B \to B\times B$. The interpretation of the formula $\mathbin{\mathop{\textrm{\huge $\vee$}}\limits_{n\in {\mathbb N}}(\exists z_{1}\ldots \exists z_{n}(z_{1}=x \wedge z_{n}=x' \wedge \mathbin{\mathop{\textrm{\huge $\wedge$}}\limits_{1\lt i\lt n}}(z_{i-1}, z_{i})\in G_{f,g}))}$ thus coincides with the equivalence relation on $B$ generated by $S_{f,g}$, that is with $R$ (cf. Remark \ref{equivgen}).   
\end{proofs}

Given a functor $H:{\cal I}\to {\cal E}$ from a small category $\cal I$ to a Grothendieck topos $\cal E$ and an object $E$ of $\cal E$, we define ${\cal I}^{E}_{H}$ as the set of pairs $(i, x)$, where $i$ is an object of $\cal I$ and $x$ is a generalized element $E\to H(i)$. In the sequel we shall occasionally identify a generalized element $x:E\to H(i)$ with its composite $\kappa_{i}\circ x$ with the coproduct arrow $\kappa_{i}$.

The following proposition can be established by means of arguments similar to those employed in the proof of Lemma \ref{lemma2}.

\begin{proposition}
Let $H:{\cal I}\to {\cal E}$ be a functor from a small category $\cal I$ to a Grothendieck topos $\cal E$. Then the equivalence relation $R$ on the object $\mathbin{\mathop{\textrm{ $\coprod$}}\limits_{i\in I}}H(i)$ given by the kernel pair of the canonical arrow $\mathbin{\mathop{\textrm{ $\coprod$}}\limits_{i\in I}}H(i) \to colim(H)$ satisfies the following property: for any object $E$ of $\cal E$ and any pairs $(i, x)$ and $(i', x')$ in ${\cal I}^{E}_{H}$, we have that $(\kappa_{i}\circ x, \kappa_{i'}\circ x')\in R_{E}$ if and only if there exists an epimorphic family $\{e_{k}:E_{k}\to E \textrm{ | } k\in K\}$ in $\cal E$ such that for every $k\in K$, the pair $((i, x\circ e_{k}), (i', x'\circ e_{k}))$ belongs to the equivalence relation on the set ${\cal I}^{E_{k}}_{H}$ generated by the pairs of the form $((j, y), (j', H(f)\circ y))$, where $f:j\to j'$ is an arrow in $\cal I$ and $y$ is a generalized element $E_{k}\to H(j)$.    
\end{proposition}

\begin{proofs}
Let us consider the pair of arrows $a,b:\mathbin{\mathop{\textrm{ $\coprod$}}\limits_{u:i\to j\textrm{ in } {\cal I}}}H(i) \to \mathbin{\mathop{\textrm{ $\coprod$}}\limits_{i\in I}}H(i)$ defined at the beginning of section \ref{expcolim}. The canonical arrow $\mathbin{\mathop{\textrm{ $\coprod$}}\limits_{i\in I}}H(i) \to colim(H)$ is the coequalizer of $a$ and $b$.
 
For any object $E$ of $\cal E$, let $R_{E}$ be the relation on the set $Hom_{\cal E}(E, \mathbin{\mathop{\textrm{ $\coprod$}}\limits_{i\in I}}H(i))$ consisting of the pairs $(x,x')$ with the property that for any $i, j\in {\cal I}$ there exists an epimorphic family $\{e^{(i, j)}_{k}: E^{(i, j)}_{k} \to E^{x, x'}_{i, j}  \textrm{ | } k\in K_{(i, j)} \}$ in $\cal E$ such that for each $k\in K_{(i, j)}$, $(x \circ r^{x, x'}_{i, j}\circ e^{(i, j)}_{k}, x \circ r^{x, x'}_{i, j}\circ e^{(i, j)}_{k})$ belongs to the equivalence relation $T_{k}^{(i, j)}$ on the set ${\cal I}^{E^{(i, j)}_{k}}_{H}$ generated by the set of pairs of the form $((n, y), (m,H(z)\circ y))$ (for $n\in {\cal I}$, $y:E^{(i, j)}_{k} \to H(n)$ and $z:n\to m$ in $\cal I$). 

It is readily seen that the assignment $E \to R_{E}$ satisfies the hypotheses of Lemma \ref{lemma1} (notice that the fibered product of two epimorphic families is again an epimorphic family), whence it defines an equivalence relation $R$ in $\cal E$ on the object $\mathbin{\mathop{\textrm{ $\coprod$}}\limits_{i\in I}}H(i)$ such that for any generalized elements $x, x':E \to \mathbin{\mathop{\textrm{ $\coprod$}}\limits_{i\in I}}H(i)$, $<x, x'>$ factors through $R$ if and only if $(x, x')\in R_{E}$.

The arrow $<a, b>$ factors through $R$ since for every $i, j\in {\cal I}$, $(a \circ r^{s, t}_{i, j}, b \circ r^{s, t}_{i, j})$ belongs to the equivalence relation $T_{k}^{(i, j)}$. 

By Remark \ref{equivgen}, $R$ thus contains the kernel pair of $q$. The converse inclusion follows from the fact that for any $(x, x')\in R_{E}$, $q\circ x=q\circ x'$; indeed, if $\{e^{(i, j)}_{k}: E^{(i, j)}_{k} \to E^{x, x'}_{i, j} \textrm{ | } k\in K_{(i, j)} \}$ in $\cal E$ is an epimorphic family such that for each $k\in K_{(i, j)}$, $(x \circ r^{x, x'}_{i, j}\circ e^{(i, j)}_{k}, x \circ r^{x, x'}_{i, j}\circ e^{(i, j)}_{k})$ belongs to the equivalence relation $T_{k}^{(i, j)}$ and hence we have $q \circ x \circ r^{x, x'}_{i, j}\circ e^{(i, j)}_{k}=q \circ x' \circ r^{x, x'}_{i, j}\circ e^{(i, j)}_{k}$ for all $k\in K_{(i, j)}$, that is $q\circ x=q\circ x'$. 
\end{proofs}

\begin{remark}
The proposition could be alternatively deduced as a consequence of Lemma \ref{lemma2}, but the proof of such implication would be more involved than the direct one that we have given above.  
\end{remark}

In the case of filtered indexing categories, the description of colimits given above simplifies. We shall be concerned in particular with functors of the form $F\circ \pi_{P}$, where $F$ is a flat functor with values in a Grothendieck topos and $\pi_{P}$ is the fibration associated to a set-valued functor $P$.

Recall from chapter VII of \cite{MM} that a functor $F:{\cal C}^{\textrm{op}}\to {\cal E}$ is flat if and only if it is filtering, i.e. the following conditions are satisfied:

\begin{enumerate}[(i)]
\item For any object $E$ of $\cal E$ there exists some epimorphic family $\{e_{i}: E_{i} \to E \textrm{ | } i\in I \}$ in $\cal E$ and for each index $i$ an object $b_{i}$ of $\cal C$ and a 
generalized element $E_{i} \to F(b_{i})$ in $\cal E$; 
\item For any two objects $c$ and $d$ in $\cal C$ and any generalized element $(x,y): E \to F(c) \times F(d)$ in $\cal E$ there is an epimorphic family $\{e_{i}: E_{i} \to E \textrm{ | } i\in I \}$ in $\cal E$ and for each index $i$ an object $b_{i}$ of $\cal C$ with arrows $u_{i}: c \to b_{i}$ and $v_{i}: d \to b_{i}$ in $\cal C$ and a generalized element $z_{i}:E_{i} \to F(b_{i})$ in $\cal E$ such that $<F(u_{i}), F(v_{i})>\circ z_{i}=<x, y>\circ e_{i}$; 
\item For any two parallel arrows $u, v: d \to c$ in $\cal C$ and any generalized element $x: E \to F(c)$ in $\cal E$ for which $F(u)\circ x = F(v)\circ x$, there is an epimorphic family $\{e_{i}: E_{i} \to E \textrm{ | } i\in I \}$ in $\cal E$ and for each index $i$ an arrow $w_{i}:c \to b_{i}$ and a generalized element $y_{i}:E_{i} \to F(b_{i})$ such that $w_{i}\circ u=w_{i}\circ v$ and $F(w_{i})\circ y_{i}=x \circ e_{i}$.
\end{enumerate}

Notice that can suppose $E=1_{\cal E}$ in condition $(i)$ without loss of generality and if all the arrows in the category $\cal C$ are monic then condition $(iii)$ rewrites as follows: for any two parallel arrows $u, v: D \to C$ in $\cal C$, either $u=v$ or the equalizer of $F(u)$ and $F(v)$ is zero. 

\begin{proposition}\label{colim}
Let $\cal C$ be a small category, $\cal E$ a Grothendieck topos,  $P:{\cal C}\to \Set$ a functor and $F:{\cal C}^{\textrm{op}}\to {\cal E}$ a flat functor. Let $\pi_{P}:\int P \to {\cal C}^{\textrm{op}}$ be the canonical projection from the category of elements of $P$ to ${\cal C}^{\textrm{op}}$. Then
\[
colim(F\circ \pi_{P})\cong (\mathbin{\mathop{\textrm{ $\coprod$}}\limits_{(c, x)\in \int P}}F(c))\slash R,
\]
where $R$ is the equivalence relation in $\cal E$ defined by saying that, for any objects $(c, x)$ and $(d, y)$ of the category $\int P$, the geometric bi-sequent
\[
\begin{array}{lll}
( \name{R}(\name{\xi_{(c, x)}}(u), \name{\xi_{(d, y)}}(v)) \dashv\vdash_{u^{F(c)}, v^{F(d)}} & & \\ 
\mathbin{\mathop{\textrm{\huge $\vee$}}\limits_{\stackrel{c \stackrel{f}{\rightarrow} a \stackrel{g}{\leftarrow} d \textrm{|}}{{\scriptscriptstyle P(f)(x)=P(g)(y)}}}} (\exists \xi^{F(a)})(\name{F(f)}(\xi)=u \wedge \name{F(g)}(\xi)=v)) & &    
\end{array}
\]
is valid in the $\Sigma_{\cal E}$-structure ${\cal S}_{\cal E}$. 

In particular, for any objects $(c, x), (c, y)$ of the category $\int P$, the geometric bi-sequent
\[
( \name{R}(\name{\xi_{(c, x)}}(u), \name{\xi_{(c, y)}}(u)) \dashv\vdash_{u^{F(c)}} \mathbin{\mathop{\textrm{\huge $\vee$}}\limits_{\stackrel{c\stackrel{f}{\rightarrow} a \textrm{|}}{{\scriptscriptstyle P(f)(x)=P(f)(y) }}}} (\exists \xi^{F(a)})(\name{F(f)}(\xi)=u)) 
\]
is valid in the $\Sigma_{\cal E}$-structure ${\cal S}_{\cal E}$. 
\end{proposition}

\begin{proofs}
Let $\mathbin{\mathop{\textrm{ $\coprod$}}\limits_{(c, x)\in \int P}}F(c)$ be the coproduct in $\cal E$ of the functor $F\circ \pi_{P}$, with canonical coproduct arrows $\xi_{(c, x)}:F(c)\to \mathbin{\mathop{\textrm{ $\coprod$}}\limits_{(c, x)\in \int P}}F(c)$, and let $T$ be the equivalence relation on $\mathbin{\mathop{\textrm{ $\coprod$}}\limits_{(c, x)\in \int P}}F(c)$ corresponding to the quotient $\mathbin{\mathop{\textrm{ $\coprod$}}\limits_{(c, x)\in \int P}}F(c) \to colim(F\circ \pi_{P})$. Then we have $colim(F\circ \pi_{P})\cong (\mathbin{\mathop{\textrm{ $\coprod$}}\limits_{(c, x)\in \int P}}F(c))\slash R$, and by Lemma \ref{lemma2}, the relation $R$ is the equivalence relation on the coproduct $\mathbin{\mathop{\textrm{ $\coprod$}}\limits_{(c, x)\in \int P}} F(c)$ generated by the image of the arrow $<a,b>$ where $s,t:\mathbin{\mathop{\textrm{ $\coprod$}}\limits_{z:(c, x)\to (d, y)\textrm{ in } \int P }}F(c) \to \mathbin{\mathop{\textrm{ $\coprod$}}\limits_{(c, x)\in \int P}}F(c)$ are the arrows defined above. 

Let $\lambda_{z}:F(c)\to \mathbin{\mathop{\textrm{ $\coprod$}}\limits_{z:(c, x)\to (d, y)\textrm{ in } \int P }}F(c)$ be the canonical coproduct arrows (for each arrow $z:(c, x)\to (d, y)$ in the category $\int P$). Notice that for any arrow $z:(c, x)\to (d, y)$ in $\int P$ the sequent $\top \vdash_{u'^{F(c)}} (\name{a}(\name{\lambda_{z}}(u')), \name{b}(\name{\lambda_{z}}(u'))) \in E_{(c, x), (d, y)}$ is valid in ${\cal S}_{\cal E}$. 

The fact that the right-to-left implication of the bi-sequent in the statement of the proposition holds in ${\cal S}_{\cal E}$ is clear. 

For any objects $(c, x), (d, y)$ in $\int P$, let the expression $(u^{F(c)}, v^{F(d)})\in E_{(c, x), (d, y)}$ be an abbreviation of the formula $\mathbin{\mathop{\textrm{\huge $\vee$}}\limits_{\stackrel{c\stackrel{f}{\rightarrow} a \stackrel{g}{\leftarrow} d \textrm{ | }}{{\scriptscriptstyle P(f)(x)=P(g)(y)}}}} (\exists \xi^{F(a)})(\name{F(f)}(\xi)=u \wedge \name{F(g)}(\xi)=v))$.

The validity of the left-to-right implication will clearly follow from Lemma \ref{lemma3} once we have shown that for any objects $(c, x), (d, y), (e, z)$ in $\int P$ the sequent 
\[
(u, v)\in E_{(c, x), (d, y)} \wedge (v, w)\in E_{(d, y), (e, z)} \vdash_{u^{F(c)}, v^{F(d)}, z^{F(e)}}  (u, w)\in E_{(c, x), (e, z)}
\]
is valid in the $\Sigma_{\cal E}$-structure ${\cal S}_{\cal E}$.

For simplicity we shall give the proof in the case ${\cal E}=\Set$, but all of our the arguments are formalizable in the internal logic of the topos and hence are valid in general. In fact, the lift of the proof from the set-theoretic to the topos-theoretic setting is made possible by the following logical characterization of flatness for a functor; specifically, the fact that a functor $F:{\cal C}^{\textrm{op}}\to {\cal E}$ is flat can be expressed in terms of the validity of the following sequents in the structure ${\cal S}_{\cal E}$: 
\[
(\top \vdash_{[]} \mathbin{\mathop{\textrm{\huge $\vee$}}\limits_{c\in Ob({\cal C})}(\exists x^{F(c)})(x=x)} );
\]
\[
(\top \vdash_{x^{F(a)}, y^{F(b)}} \mathbin{\mathop{\textrm{\huge $\vee$}}\limits_{a\stackrel{f}{\rightarrow} c \stackrel{g}{\leftarrow} b}}(\exists z^{F(c)})(\name{F(f)}(z)=x \wedge \name{F(g)}(z)=y))
\]
for any objects $a$, $b$ of $\cal C$;
\[
(\name{F(f)}(x)=\name{F(g)}(x) \vdash_{x^{F(a)}} \mathbin{\mathop{\textrm{\huge $\vee$}}\limits_{h:a\to c \textrm{ | } h\circ f=h\circ g}}(\exists z^{F(c)})(\name{F(h)}(z)=x))
\]
for any pair of arrows $f,g:b\to a$ in $\cal C$ with common domain and codomain.

Now, given objects $(c, x), (d,y), (e, z)\in \int P$ and elements $u\in F(c), v\in F(d), w\in F(e)$, suppose that $(u, v)\in E_{(c, x), (d, y)}$ and $(v, w)\in E_{(d, y), (e, z)}$; we want to prove that $(u, w)\in E_{(c, x), (e, z)}$. 

As $(u, v)\in E_{(c, x), (d, y)}$, there exist arrows $f:c\to a$ and $g:d \to a$ in $\cal C$ and an element $\xi \in F(a)$ such that $F(f)(\xi)=u$, $F(g)(\xi)=v$, $P(f)(x)=P(g)(y)$. Similarly, as $(v, w)\in E_{(d, y), (e, z)}$, there exist arrows $h:d\to b$ and $k:e\to b$ in $\cal C$ and an element $\chi \in F(b)$ such that $F(h)(\chi)=v$, $F(k)(\chi)=w$ and $P(h)(y)=P(k)(z)$.

Using condition $(ii)$ in the definition of flat functor, we obtain the existence of an object $m$ of $\cal C$, arrows $r:a\to m$ and $t:b\to m$ in $\cal C$ and an element $\epsilon \in F(m)$ such that $F(r)(\epsilon)=\xi$ and $F(t)(\epsilon)=\chi$. Now, consider the arrows $r\circ g, t\circ h:d \to m$; we have that $F(r\circ g)(\epsilon)=F(t\circ h)(\epsilon)$ and hence by condition $(iii)$ in the definition of flat functor there exists an arrow $s:m\to n$ in $\cal C$ such that $s\circ r\circ g=s\circ t\circ h$ and an element $\alpha \in F(n)$ such that $F(s)(\alpha)=\epsilon$. 

Consider the arrows $s\circ r\circ f$ and $s\circ t\circ k$. They satisfy the necessary conditions to ensure that $(u, w)\in E_{(c, x), (e, z)}$, i.e. $F(s\circ r\circ f)(\alpha)=u$, $F(s\circ t \circ k)(\alpha)=w$, $P(s\circ r \circ f)(x)=P(s\circ t \circ k)(z)$. Indeed, $F(s\circ r\circ f)(\alpha)=F(f)(F(r)(F(s)(\alpha)))=F(f)(F(r)(\epsilon))=F(f)(\xi)=u$, $F(s\circ t\circ k)(\alpha)=F(k)(F(t)(F(s)(\alpha)))=F(k)(F(t)(\epsilon))=F(k)(\chi)=w$ and $P(s\circ r\circ f)(x)=P(s\circ r)(P(f)(x))=P(s\circ r)(P(g)(y))=P(s\circ r\circ g)(y)=P(s\circ t \circ h)(y)=P(s\circ t)(P(h)(y))=P(s\circ t)(P(k)(z))=P(s\circ t\circ k)(z)$, as required.         

The second part of the proposition follows from the first by an easy application of condition $(iii)$ in the definition of flat functor. 
\end{proofs}

The following proposition represents the translation of Proposition \ref{colim} in the categorical language of generalized elements.

\begin{proposition}\label{colimcat}
Let $\cal C$ be a small category, $\cal E$ a Grothendieck topos with a separating set $S$,  $P:{\cal C}\to \Set$ a functor and $F:{\cal C}^{\textrm{op}}\to {\cal E}$ a flat functor. Let $\pi_{P}:\int P \to {\cal C}^{\textrm{op}}$ be the canonical projection from the category of elements of $P$ to ${\cal C}^{\textrm{op}}$. Then the equivalence relation $R$ defined by the formula
\[
colim(F\circ \pi_{P})\cong (\mathbin{\mathop{\textrm{ $\coprod$}}\limits_{(c, x)\in \int P}}F(c))\slash R
\] 
admits the following characterization: for any objects $(c, x), (d, y)$ of $\int P$ and generalized elements $u:E\to F(c)$, $v:E\to F(d)$ (where $E\in S$), $(\xi_{(c, x)} \circ u, \xi_{(d, y)} \circ v)\in R_{E}$ if and only if there exists an epimorphic family $\{e_{i}:E_{i}\to E \textrm{ | } i\in I, E_{i}\in S\}$ and for each index $i\in I$ an object $a_{i}\in {\cal C}$, a generalized element $h_{i}:E_{i}\to F(b_{i})$ and two arrows $f_{i}:c\to a_{i}$ and $f_{i}':d\to a_{i}$ in $\cal C$ such that $P(f_{i})(x)=P(f_{i}')(y)$ and $<F(f_{i}), F(f'_{i})>\circ h_{i}=<u, v>\circ e_{i}$.  

In particular, for any objects $(c, x)$ and $(c, y)$ of $\int P$ and any generalized element $u:E\to F(c)$ (where $E\in S$), $(\xi_{(c, x)} \circ u, \xi_{(d, y)} \circ v)\in R_{E}$ if and only if there exists an epimorphic family $\{e_{i}:E_{i}\to E \textrm{ | } i\in I, E_{i}\in S\}$ and for each index $i\in I$ an object $a_{i}\in {\cal C}$, a generalized element $h_{i}:E_{i}\to F(a_{i})$ and an arrow $f_{i}:c\to a_{i}$ in $\cal C$ such that $P(f_{i})(x)=P(f_{i}')(y)$ and $F(f_{i}) \circ h_{i}=u \circ e_{i}$. 
\end{proposition}\qed

\section{Extensions of flat functors}\label{ext}

In this section, we investigate the operation on flat functors induced by a geometric morphism of toposes. This will be relevant for the characterization of the class of theories classified by a presheaf topos addressed in section 5.

\subsection{General extensions}\label{gene}

Let $({\cal C}, J)$ and $({\cal D}, K)$ be Grothendieck sites, and $u:\Sh({\cal D}, K)\to \Sh({\cal C}, J)$ a geometric morphism. This morphism induces, via Diaconescu's equivalence, a functor 
\[
\xi_{{\cal E}}:\textbf{Flat}_{K}({\cal D}, {\cal E})\to \textbf{Flat}_{J}({\cal C}, {\cal E}),
\]
where we write $\textbf{Flat}_{Z}({\cal R}, {\cal E})$ for the category of flat $Z$-continuous functors from $\cal R$ to $\cal E$.

This functor can be described explicitly as follows. For any flat $K$-continuous functor $F:{\cal D} \to {\cal E}$ with corresponding geometric morphism $f_{F}:{\cal E}\to \Sh({\cal D}, K)$, the flat functor $\tilde{F}:=\xi_{\cal E}(F):{\cal C} \to {\cal E}$ is given by $f_{F}^{\ast}\circ u^{\ast} \circ l_{\cal C}$, where $l_{\cal C}:{\cal C}\to \Sh({\cal C}, J)$ is the composite of the Yoneda embedding $y_{\cal C}:{\cal C} \to [{\cal C}^{\textrm{op}}, \Set]$ with the associated sheaf functor $[{\cal C}^{\textrm{op}}, \Set] \to \Sh({\cal C}, J)$. 

For any natural transformation $\alpha:F\to G$ between flat functors in $\textbf{Flat}_{K}({\cal D}, {\cal E})$, the corresponding natural transformation ${f_{\alpha}}^{\ast}:{f_{F}}^{\ast}\to {f_{G}}^{\ast}$, applied to the functor $u^{\ast} \circ l_{\cal C}$, gives rise to a natural transformation $\tilde{\alpha}:\tilde{F}\to \tilde{G}$ which is precisely $\xi_{\cal E}(\alpha)$.  

We can express $\tilde{F}$ directly in terms of $F$ by using a colimit construction, as follows. 

For any $c\in {\cal C}$, the $K$-sheaf $u^{\ast}(l_{\cal C}(c)):{\cal D}^{\textrm{op}}\to \Set$, considered as an object of $[{\cal D}^{\textrm{op}}, \Set]$, can be canonically expressed as a colimit of representables, indexed over its category of elements ${\cal A}_{c}$; specifically, $u^{\ast}(l_{\cal C}(c)) \cong colim(y_{\cal D} \circ \pi_{c})$ in $[{\cal D}^{\textrm{op}}, \Set]$, where $\pi_{c}:{\cal A}_{c}\to {\cal D}$ is the canonical projection functor. As the associated sheaf functor $a_{K}:[{\cal D}^{\textrm{op}}, \Set] \to \Sh({\cal D}, K)$ preserves colimits, the functor $u^{\ast}(l_{\cal C}(c))$ is isomorphic in $\Sh({\cal D}, K)$ to the colimit of the functor $l_{{\cal D}}\circ \pi_{c}$. Therefore $\tilde{F}(c)=(f_{F}^{\ast}\circ u^{\ast} \circ l_{\cal C})(c)=f_{F}^{\ast}(colim(l_{{\cal D}}\circ \pi_{c}))\cong colim(f_{F}^{\ast}\circ l_{{\cal D}}\circ \pi_{c})=colim(F\circ \pi_{c})$. 

For any $(d, z)\in {\cal A}_{c}$, we write $\kappa^{F}_{(d,z)}:F(d)\to \tilde{F}(c)$ for the canonical colimit arrow.

\begin{proposition}\label{alphaext}
\begin{enumerate}[(i)]
\item For any natural transformation $\alpha:F\to G$ between flat functors $F,G:{\cal D} \to {\cal E}$ in $\textbf{Flat}_{K}({\cal D}, {\cal E})$, the natural transformation $\tilde{\alpha}:\tilde{F}\to \tilde{G}$ is characterized by the following condition: for any object $(d, z)$ of the category ${\cal A}_{c}$ the diagram 
\[  
\xymatrix {
\tilde{F}(c)  \ar[r]^{\tilde{\alpha}(c)} & \tilde{G}(c) \\
F(d) \ar[u]^{\kappa^{F}_{(d,z)}} \ar[r]_{\alpha(d)} & G(d) \ar[u]^{\kappa^{G}_{(d,z)}}} 
\]
commutes. 

\item The functor $\xi_{{\cal E}}:\textbf{Flat}_{K}({\cal D}, {\cal E})\to \textbf{Flat}_{J}({\cal C}, {\cal E})$ is natural in $\cal E$; that is, for any geometric morphism $f:{\cal F}\to {\cal E}$, the diagram
\[  
\xymatrix {
\textbf{Flat}_{K}({\cal D}, {\cal E}) \ar[d]^{f^{\ast} \circ -} \ar[r]^{\xi_{{\cal E}}} & \textbf{Flat}_{J}({\cal C}, {\cal E}) \ar[d]^{f^{\ast} \circ -} \\
\textbf{Flat}({\cal D}, {\cal F}) \ar[r]^{\xi_{{\cal F}}} & \textbf{Flat}({\cal C}, {\cal F}) } 
\]
commutes. 

\item For any arrow $f:c\to c'$ in $\cal C$ and any object $(d, z)$ in ${\cal A}_{c}$, the diagram

\[  
\xymatrix {
\tilde{F}(c)  \ar[r]^{\tilde{F}(f)} & \tilde{F}(c') \\
F(d) \ar[u]^{\kappa^{F}_{(d,z)}} \ar[ur]_{\kappa^{F}_{(d,u^{\ast}(l_{{\cal C}}(f))(z))}} & } 
\]
commutes.

\end{enumerate}
\end{proposition}

\begin{proofs}
$(i)$ The given square commutes as it is the naturality square for the natural transformation ${f_{\alpha}}^{\ast}:{f_{F}}^{\ast}\to {f_{G}}^{\ast}$ with respect to the canonical colimit arrow $y_{\cal D}(d) \to u^{\ast}(l_{\cal C}(c))$ in $[{\cal D}^{\textrm{op}}, \Set]$ corresponding to the object $(d, z)$ of ${\cal A}_{c}$.  

$(ii)$ This follows as an immediate consequence of the fact that Diaconescu's equivalences are natural in $\cal E$.

$(iii)$ Given a natural transformation $\beta:P\to P'$ of presheaves $P$ and $P'$ in $[{\cal D}^{\textrm{op}}, \Set]$, denoting by $\kappa^{P}:\int P \to {\cal D}$ and $\kappa^{P'}:\int P' \to {\cal D}$ the canonical projection functors and by $\kappa^{P}_{(d, z)}:y_{\cal D}d \to P$, $\kappa^{P'}_{(d', z')}:y_{\cal D}d' \to P$ the canonical colimit arrows in $[{\cal D}^{\textrm{op}}, \Set]$ (for $(d, z)\in \int P$ and $(d', z')\in \int P'$), we clearly have that $\beta \circ \kappa^{P}_{(d, z)}=k^{P'}_{(d, \beta(d)(z))}$. To obtain our thesis it suffices to apply this result to the natural transformation $\beta:u^{\ast}(l_{{\cal C}}(f)):u^{\ast}(l_{{\cal C}}(c)) \to u^{\ast}(l_{{\cal C}}(c'))$; indeed, by applying the functor $f^{\ast}_{F}$ to the resulting equality we obtain precisely the identity in the statement of the proposition.
\end{proofs}

\subsection{Extensions along embeddings of categories}\label{seccov}

Let $\cal D$ be a subcategory of a small category $\cal C$. Then the inclusion functor $i:{\cal D}\hookrightarrow {\cal C}$ induces a geometric morphism $E(i):[{\cal D}, \Set]\to [{\cal C}, \Set]$ and hence, by Diaconescu's equivalence, a functor 
\[
\xi_{{\cal E}}:\textbf{Flat}({\cal D}^{\textrm{op}}, {\cal E})\to \textbf{Flat}({\cal C}^{\textrm{op}}, {\cal E}).
\]

We can describe $\tilde{F}$ directly in terms of $F$ as follows. For any $c\in {\cal C}$, the functor $E(i)^{\ast}(y_{\cal C}(c)):{\cal D}\to \Set$ coincides with the functor $Hom_{{\cal C}}(c, -):{\cal D}\to \Set$. Its category of elements ${\cal A}_{c}$ has as objects the pairs $(d, h)$ where $d$ is an object of $\cal D$ and $h$ is an arrow $h:c\to d$ in $\cal C$ and as arrows $(d, h)\to (d', h')$ the arrows $k:d'\to d$ in $\cal D$ such that $k\circ h'=h$ in $\cal C$. As $f_{F}^{\ast}$ preserves all small colimits and $Hom_{{\cal C}}(c, -)$ is the colimit in $[{\cal D}, \Set]$ of the composition of the canonical projection $\pi_{c}:{\cal A}_{c}\to {\cal D}^{\textrm{op}}$ with the Yoneda embedding ${\cal D}^{\textrm{op}}\to [{\cal D}, \Set]$, the object $\tilde{F}(c)$ can be identified with the colimit of the composite functor $F\circ \pi_{c}:{\cal A}_{c}\to {\cal E}$. For any object $d$ of $\cal D$, the pair $(d, 1_{d})$ defines an object of the category ${\cal A}_{d}$ and hence a canonical colimit arrow $\chi_{d}:F(d)\to \tilde{F}(d)$, also denoted by $\chi_{d}^{F}$. For any object $(d, h)$ of ${\cal A}_{c}$, the colimit arrow $\kappa_{(d, h)}:F(d)\to \tilde{F}(c)$ is equal to the composite $\tilde{F}(h)\circ \chi_{d}$. 

\begin{remarks}\label{cov}
\begin{enumerate}[(i)]

\item The functor $\xi_{{\cal E}}:\textbf{Flat}({\cal D}^{\textrm{op}}, {\cal E})\to \textbf{Flat}({\cal C}^{\textrm{op}}, {\cal E})$ coincides with the (restriction to the subcategories of flat functors) of the left Kan extension functor $[{\cal D}^{\textrm{op}}, {\cal E}]\to [{\cal C}^{\textrm{op}}, {\cal E}]$ along the inclusion ${\cal D}^{\textrm{op}}\hookrightarrow {\cal C}^{\textrm{op}}$ (cf. Remark \ref{remadj} below).

\item For any object $d\in {\cal D}$ the diagram
\[  
\xymatrix {
\tilde{F}(d)  \ar[r]^{\tilde{\alpha}(d)} & \tilde{G}(d) \\
F(d) \ar[u]^{\chi^{F}_{d}} \ar[r]_{\alpha(d)} & G(d) \ar[u]^{\chi^{G}_{d}}} 
\]
commutes (cf. Proposition \ref{alphaext}(a)). 

\item For any object $d\in {\cal D}$, the arrow $\chi_{d}:F(d)\to \tilde{F}(d)$ is monic. Indeed, it can be identified with the image of the canonical subfunctor inclusion $Hom_{\cal D}(d, -)\hookrightarrow Hom_{\cal C}(d, -)$ under the inverse image $f_{F}^{\ast}$ of the geometric morphism $f_{F}$. From this, in view of Remark \ref{cov}(a), it follows at once that the functor $\xi_{{\cal E}}:\textbf{Flat}({\cal D}^{\textrm{op}}, {\cal E})\to \textbf{Flat}({\cal C}^{\textrm{op}}, {\cal E})$ is faithful. 

\item For any flat functors $F,G:{\cal D}^{\textrm{op}}\to {\cal E}$ and natural transformation $\beta:\tilde{F}\to \tilde{G}$, $\beta=\tilde{\alpha}$ for some natural transformation $\alpha:F\to G$ if and only if for every $d\in {\cal D}$ the diagram 
\[  
\xymatrix {
\tilde{F}(d)  \ar[r]^{\beta(d)} & \tilde{G}(d) \\
F(d) \ar[u]^{\chi^{F}_{d}} \ar[r]_{\alpha(d)} & G(d) \ar[u]^{\chi^{F}_{d}}} 
\]
commutes. One direction follows from Remark \ref{cov}(a). To prove the other one we observe that $\beta(c)=\tilde{\alpha}(c)$ for all $c\in {\cal C}$ if and only if for every object $(a, z)\in {\cal A}_{c}$, $\beta(c)\circ \kappa_{(a, z)}=\tilde{\alpha}(c)\circ \kappa_{(a, z)}$. Now, $\beta(c)\circ \kappa_{(a, z)}=\beta(c) \circ \tilde{F}(z)\circ \chi^{F}_{a}= \tilde{G}(z) \circ \beta(a) \circ \chi^{F}_{a}=\tilde{G}(z) \circ \chi^{G}_{a} \circ \alpha(a)=\kappa^{G}_{(a, z)}\circ \alpha(a)=\tilde{\alpha}(c)\circ \kappa^{F}_{(a, z)}$ (where the third equality follows from the naturality of $\beta$ and the last from Proposition \ref{alphaext}(a)), as required. 
 
\item Let $d$ be an object of $\cal D$ and $y_{\cal D}d=Hom_{\cal D}(-, d):{\cal D}^{\textrm{op}}\to \Set$ the representable on $\cal D$ associated to $d$. Then $\xi_{{\cal E}}(y_{\cal D}d)=y_{\cal C}d$, where $y_{\cal C}d=Hom_{\cal C}(-, d)$ is the representable on $\cal C$ associated to $d$ (considered here as an object of $\cal C$). Indeed, the flat functor $y_{\cal D}d$ corresponds to the geometric morphism whose inverse image is the evaluation functor $ev^{{\cal D}}_{d}:[{\cal D}, \Set]\to \Set$ at the object $d$, and the composite $ev^{{\cal D}}_{d} \circ E(i)^{\ast}$ coincides with the evaluation functor $ev^{{\cal C}}_{d}:[{\cal C}, \Set]\to \Set$ at the object $d$, which corresponds to the flat functor $y_{\cal C}d$.

\item Let $\cal E$ be a Grothendieck topos and $\gamma_{\cal E}:{\cal E}\to \Set$ the unique geometric morphism from $\cal E$ to $\Set$. Then by Proposition \ref{alphaext}(ii) and Remark \ref{cov}(v), for any object $d$ of $\cal D$, the functor $\xi_{{\cal E}}:\textbf{Flat}({\cal D}^{\textrm{op}}, {\cal E})\to \textbf{Flat}({\cal C}^{\textrm{op}}, {\cal E})$ sends the flat functor $\gamma_{\cal E}^{\ast}\circ y_{\cal D}d$ to the functor $\gamma_{\cal E}^{\ast}\circ y_{\cal C}d$.    
\end{enumerate}
\end{remarks}

\begin{proposition}
Under the natural equivalences $e_{\cal C}:\Ind{\cal C}\simeq \textbf{Flat}({\cal C}^{\textrm{op}}, \Set)$ and $e_{\cal D}:\Ind{\cal D}\simeq \textbf{Flat}({\cal D}^{\textrm{op}}, \Set)$, the functor $\xi_{\Set}$ corresponds to the functor $\Ind{i}:\Ind{\cal D}\to \Ind{\cal C}$.
\end{proposition}   

\begin{proofs}
Let us first recall the definition of the functor $\Ind{i}:\Ind{\cal D}\to \Ind{\cal C}$. For any flat functor $F:{\cal C}^{\textrm{op}}\to \Set$, $\Ind{i}(F)$ is the functor given by the colimit in $\Ind{\cal C}$ of the composite functor $i\circ \pi_{F}$, where $\pi_{F}:{\int F} \to {\cal D}$ is the canonical projection from the category of elements ${\int F}$ of $F$ to $\cal D$. For any $c\in {\cal C}$, we thus have $\Ind{i}(F)(c)=colim(ev_{c}\circ i \circ \pi_{F})$, where $ev_{c}:\Ind{\cal C} \to \Set$ is the evaluation functor at $c$.

Now, the functor $ev_{c}\circ i:{\cal D}\to \Set$ is equal to the functor $Hom_{{\cal C}}(c, -):{\cal D}\to \Set$ considered above, and $\tilde{F}(c)=colim(F\circ \pi_{c}):{\cal A}_{c}\to \Set$, where ${\cal A}_{c}$ is the category of elements of $ev_{c}\circ i$ and $\pi_{c}:{\cal A}_{c} \to {\cal D}^{\textrm{op}}$ the associated canonical projection. The commutativity of the tensor product between a presheaf and a covariant set-valued functor (cf. chapter VII of \cite{MM} or Theorem \ref{tensorpr} above) thus yields a natural isomorphism between the two sets, as required.
\end{proofs}

\begin{corollary}\label{logical}
Let ${\cal D}\hookrightarrow {\cal C}$ be an embedding of small categories and $F:{\cal D}^{\textrm{op}}\to {\cal E}$ a flat functor. With the above notation, for any object $c\in {\cal C}$, we have 
\[
\tilde{F}(c)\cong \mathbin{\mathop{\textrm{ $\coprod$}}\limits_{(a, z)\in {\cal A}_{c}}}F(c)\slash R_{c},
\]
where $R_{c}$ is the equivalence relation in $\cal E$ defined by saying that, for any objects $(a, z), (a',z')$ of the category ${\cal A}_{c}$, the geometric bi-sequent
\[
\begin{array}{lll}
( \name{R}(\name{\xi_{(a, z)}}(x), \name{\xi_{(a', z')}}(x')) \dashv\vdash_{x^{F(a)}, x'^{F(a')}} & & \\ 
\mathbin{\mathop{\textrm{\huge $\vee$}}\limits_{\stackrel{a\stackrel{f}{\rightarrow} b \stackrel{g}{\leftarrow} a' \textrm{|}}{{\scriptscriptstyle f\circ z=g\circ z'}}}} (\exists y^{F(b)})(\name{F(f)}(y)=x \wedge \name{F(g)}(y)=x')) & &    
\end{array}
\]
is valid in the tautological $\Sigma_{\cal E}$-structure ${\cal S}_{\cal E}$ (defined in section \ref{expcolim}), where for any object $(a, z)$ of the category ${\cal A}_{c}$, $\xi_{(a, z)}:F(a)\to \tilde{F}(c)$ is the canonical colimit arrow.

In particular, for any objects $(a, z), (a,z')$ of the category ${\cal A}_{c}$, the geometric bi-sequent
\[
( \name{R}(\name{\xi_{(a, z)}}(x), \name{\xi_{(a, z')}}(x)) \dashv\vdash_{x^{F(a)}} \mathbin{\mathop{\textrm{\huge $\vee$}}\limits_{\stackrel{a\stackrel{f}{\rightarrow} b \textrm{|}}{{\scriptscriptstyle f\circ z=f\circ z'}}}} (\exists y^{F(b)})(\name{F(f)}(y)=x)) 
\]
is valid in the $\Sigma_{\cal E}$-structure ${\cal S}_{\cal E}$.

Semantically, the relation $R_{c}$ can be characterized by saying that, for any objects $(a, z), (a',z')$ of the category ${\cal A}_{c}$ and any generalized elements $x:E\to F(a)$, $x':E\to F(a')$, $(\kappa_{(a, z)} \circ x, \kappa_{(a', z')} \circ x')\in R_{c}$ if and only if there exists an epimorphic family $\{e_{i}:E_{i}\to E \textrm{ | } i\in I\}$ and for each index $i\in I$ an object $b_{i}\in {\cal D}$, a generalized element $h_{i}:E_{i}\to F(b_{i})$ and two arrows $f_{i}:a\to b_{i}$ and $f_{i}':a'\to b_{i}$ in $\cal D$ such that $f_{i}'\circ z'=f_{i}\circ z$ and $<F(f_{i}), F(f'_{i})>\circ h_{i}=<x, x'>\circ e_{i}$.  

In particular, for any objects $(a, z)$ and $(a, z')$ of ${\cal A}_{c}$ and any generalized element $x:E\to F(a)$, $(\kappa_{(a, z)} \circ x, \kappa_{(a, z')} \circ x)\in R_{c}$ if and only if there exists an epimorphic family $\{e_{i}:E_{i}\to E \textrm{ | } i\in I\}$ and for each index $i\in I$ an object $b_{i}\in {\cal D}$, a generalized element $h_{i}:E_{i}\to F(b_{i})$ and an arrow $f_{i}:a\to b_{i}$ in $\cal D$ such that $f_{i}\circ z'=f_{i}\circ z$ and $F(f_{i}) \circ h_{i}=x \circ e_{i}$.  
\end{corollary}

\begin{proofs}
The corollary can be obtained by applying Proposition \ref{colim} (and its categorical reformulation provided by Proposition \ref{colimcat}) to the flat functor $F:{\cal D}^{\textrm{op}}\to {\cal E}$ and the functor $P:{\cal D}\to \Set$ given by $ev_{c}\circ i=Hom_{{\cal C}}(c, -):{\cal D}\to \Set$, whose category of elements coincides with ${\cal A}_{c}$ and whose associated fibration $\pi_{P}$ coincides with $\pi_{c}:{\cal A}_{c} \to {\cal D}^{\textrm{op}}$.
\end{proofs}

\subsection{Extensions from categories of set-based models to syntactic categories}\label{extsym}

Let $\mathbb T$ be a geometric theory and $\cal K$ a small category of set-based $\mathbb T$-models. Then the family of geometric morphisms $\Set \to \Sh({\cal C}_{\mathbb T}, J_{\mathbb T})$ corresponding to the $\mathbb T$-models in $\cal K$ induces a geometric morphism 
\[
p_{\cal K}:[{\cal K}, \Set] \to \Sh({\cal C}_{\mathbb T}, J_{\mathbb T})
\]
whose associated $\mathbb T$-model in $[{\cal K}, \Set]$ is given, at each sort, by the corresponding forgetful functor. 

We thus have, for each Grothendieck topos $\cal E$, an induced functor 
\[
u^{\mathbb T}_{({\cal K}, {\cal E})}:\mathbf{Flat}({\cal K}^{\textrm{op}}, {\cal E}) \to \mathbf{Flat}_{J_{\mathbb T}}({\cal C}_{\mathbb T}, {\cal E}),
\] 
which the following theorem describes explicitly. 

Before stating it, we need to introduce some notation. For any any object $\{\vec{x}. \phi\}$ of ${\cal C}_{\mathbb T}$, we write ${\cal A}^{\cal K}_{\{\vec{x}. \phi\}}$ or simply ${\cal A}_{\{\vec{x}. \phi\}}$ when the category $\cal K$ can be obviously inferred from the context, for the category whose objects are the pairs $(M, w)$, where $M\in {\cal K}$ and $w\in [[\vec{x}. \phi]]_{M}$ and whose arrows $(M, w)\to (N, z)$ are the $\mathbb T$-model homomorphisms $g:N\to M$ in $\cal K$ such that $g(w)=z$. We denote by  $\pi^{\cal K}_{\{\vec{x}. \phi\}}:{\cal A}_{\{\vec{x}. \phi\}}\to {\cal K}^{\textrm{op}}$ the canonical projection functor.  

\begin{theorem}\label{extsymthm}
Let $\mathbb T$ be a geometric theory. Then for any flat functor $F:{\cal K}^{\textrm{op}} \to {\cal E}$, the functor $\tilde{F}:=u^{\mathbb T}_{({\cal K}, {\cal E})}(F):{\cal C}_{\mathbb T}\to {\cal E}$ sends any object $\{\vec{x}. \phi\}$ of ${\cal C}_{\mathbb T}$ to the colimit $colim(F\circ \pi^{\cal K}_{\{\vec{x}. \phi\}})$ and acts on the arrows in the obvious way. In particular, for any formula $\{\vec{x}. \phi\}$ presenting a $\mathbb T$-model $M_{\{\vec{x}. \phi\}}$ in $\cal K$, $u^{\mathbb T}_{({\cal K}, {\cal E})}(F)(\{\vec{x}. \phi\})\cong F(M_{\{\vec{x}. \phi\}})$.
\end{theorem}

\begin{proofs}
Let $g_{F}:{\cal E}\to [{\cal K}, \Set]$ be the geometric morphism corresponding, via Diaconenscu's equivalence, to the flat functor $F$. Then the functor $u^{\mathbb T}_{({\cal K}, {\cal E})}(F)$ is equal to the composite $g_{F}^{\ast}\circ p_{\cal K}^{\ast} \circ y$, where $y:{\cal C}_{\mathbb T}\to \Sh({\cal C}_{\mathbb T}, J_{\mathbb T})$ is the Yoneda embedding.

Now, for any geometric formula $\{\vec{x}. \phi\}$ over $\Sigma$, $p_{\cal K}^{\ast}(y(\{\vec{x}. \phi\}))$ is the functor $F_{\{\vec{x}. \phi\}}$ sending to any model $M\in {\cal K}$ the set $[[\vec{x}. \phi]]_{M}$. This functor can be expressed as the colimit $colim(y'\circ \pi_{\{\vec{x}. \phi\}})$, where $y':{{\cal K}}^{\textrm{op}}\to [{\cal K}, \Set]$ is the Yoneda embedding, since the functor $\pi_{\{\vec{x}. \phi\}}$ coincides with the canonical projection from the category of elements of the functor $F_{\{\vec{x}. \phi\}}$ to ${{\cal K}}^{\textrm{op}}$. 

Therefore $u^{\mathbb T}_{({\cal K}, {\cal E})}(F)(\{\vec{x}. \phi\})=g_{F}^{\ast}(colim (y'\circ \pi_{\{\vec{x}. \phi\}}))\cong colim(g_{F}^{\ast} \circ y' \circ \pi_{\{\vec{x}. \phi\}})\cong colim(F\circ \pi_{\{\vec{x}. \phi\}})$, as required.    
\end{proofs}

\begin{remark}\label{remf}
From the proof of the theorem, it is clear that the isomorphism $u^{\mathbb T}_{({\cal K}, {\cal E})}(F)(\{\vec{x}. \phi\})\cong F(M_{\{\vec{x}. \phi\}})$ is natural in $\{\vec{x}. \phi\}$; that is, for any geometric formulae $\{\vec{x}. \phi\}$ and $\{\vec{y}. \psi\}$ respectively presenting $\mathbb T$-models $M_{\{\vec{x}. \phi\}}$ and $M_{\{\vec{y}. \psi\}}$ and any $\mathbb T$-provably functional formula $\theta:\{\vec{x}. \phi\} \to \{\vec{y}. \psi\}$ inducing a $\mathbb T$-model homomorphism $M_{\theta}:M_{\{\vec{y}. \psi\}} \to M_{\{\vec{x}. \phi\}}$, the arrow $u^{\mathbb T}_{({\cal K}, {\cal E})}(F)(\theta)$ is given by the image of the arrow $F(M_{\theta})$ across the isomorphisms $u^{\mathbb T}_{({\cal K}, {\cal E})}(F)(\{\vec{x}. \phi\})\cong F(M_{\{\vec{x}. \phi\}})$ and $u^{\mathbb T}_{({\cal K}, {\cal E})}(F)(\{\vec{y}. \psi\})\cong F(M_{\{\vec{y}. \psi\}})$.  
\end{remark}

\begin{corollary}\label{preservationinequalities}
Let $\mathbb T$ be a geometric theory, $\sigma = (\phi \vdash_{\vec{x}} \psi)$ a geometric sequent over the signature of $\mathbb T$ and $F:{{\cal K}}^{\textrm{op}} \to {\cal E}$ a flat functor. If $\sigma$ is valid in every $\mathbb T$-model in $\cal K$ then $\tilde{F}(\{\vec{x}. \phi\})\leq \tilde{F}(\{\vec{x}. \psi\})$ as subobjects of $\tilde{F}(\{\vec{x}. \top\})$ in $\cal E$. 
\end{corollary}

\begin{proofs}
By Theorem \ref{extsymthm}, we have that $\tilde{F}(\{\vec{x}. \phi\})=colim(F\circ \pi_{\{\vec{x}. \phi\}})$, $\tilde{F}(\{\vec{x}. \psi\})=colim(F\circ \pi_{\{\vec{x}. \psi\}})$ and $\tilde{F}(\{\vec{x}. \top\})=colim(F\circ \pi_{\{\vec{x}. \top\}})$. Now, ${\cal A}_{\{\vec{x}. \phi\}}$ and ${\cal A}_{\{\vec{x}. \psi\}}$ canonically embed as subcategories of ${\cal A}_{\{\vec{x}. \top\}}$, and if $\sigma$ is valid in every $\mathbb T$-model in $\cal K$ then we have a canonical functor $i: {\cal A}_{\{\vec{x}. \phi\}} \to {\cal A}_{\{\vec{x}. \psi\}}$ which commutes with these embeddings. It thus follows from the functoriality of colimits that $\tilde{F}(\{\vec{x}. \phi\})\leq \tilde{F}(\{\vec{x}. \psi\})$ as subobjects of $\tilde{F}(\{\vec{x}. \top\})$ in $\cal E$, as required.   
\end{proofs}

Let us now apply Proposition \ref{colimcat} in the context of extensions $F \to \tilde{F}$ of flat functors induced by the geometric morphism
\[
p_{\cal K}:[{\cal K}, \Set] \to \Sh({\cal C}_{\mathbb T}, J_{\mathbb T}).
\] 

The following characterization is obtained by applying it in conjunction with Proposition \ref{explicit}.

\begin{proposition}\label{descrcolim}
Let $\mathbb T$ be a geometric theory over a signature $\Sigma$, $\cal K$ a small category of set-based $\mathbb T$-models, $\cal E$ a Grothendieck topos with a separating set $S$ and $F:{\cal K}^{\textrm{op}}\to {\cal E}$ a flat functor. With the above notation, for any geometric formula-in-context $\phi(\vec{x})$ over $\Sigma$, we have 
\[
\tilde{F}(\{\vec{x}. \phi\})\cong (\mathbin{\mathop{\textrm{ $\coprod$}}\limits_{(M, z)\in {\cal A}_{\{\vec{x}. \phi\}}}}F(M))\slash R_{\{\vec{x}. \phi\}},
\]
where $R_{\{\vec{x}. \phi\}}$ is the equivalence relation in $\cal E$ defined by saying that for any objects $(M, z), (N, w)$ of the category ${\cal A}_{\{\vec{x}. \phi\}}$ and any generalized elements $x:E\to F(M)$, $x':E\to F(N)$ (where $E\in S$), we have $(\xi_{(M, z)} \circ x, \xi_{(N, w)} \circ x')\in R_{\{\vec{x}. \phi\}}$ if and only if there exists an epimorphic family $\{e_{i}:E_{i}\to E \textrm{ | } i\in I, E_{i}\in S\}$ and for each index $i\in I$ a $\mathbb T$-model $a_{i}$ in $\cal K$, a generalized element $h_{i}:E_{i}\to F(a_{i})$ and two $\mathbb T$-model homomorphisms $f_{i}:M\to a_{i}$ and $f_{i}':N\to a_{i}$ such that $f_{i}(z)=f_{i}'(w)$ and $<F(f_{i}), F(f'_{i})>\circ h_{i}=<x, x'>\circ e_{i}$ (where for any object $(M, z)$ of the category ${\cal A}_{\{\vec{x}. \phi\}}$, $\kappa_{(M, z)}:F(M)\to \tilde{F}(\{\vec{x}. \phi\})$ is the canonical colimit arrow).  

In particular, for any objects $(M, z)$ and $(M, w)$ of ${\cal A}_{\{\vec{x}. \phi\}}$ and any generalized element $x:E\to F(M)$ (where $E\in S$), $(\kappa_{(M, z)} \circ x, \kappa_{(M, w)} \circ x)\in R_{\{\vec{x}. \phi\}}$ if and only if there exists an epimorphic family $\{e_{i}:E_{i}\to E \textrm{ | } i\in I, E_{i}\in S\}$ and for each index $i\in I$ a $\mathbb T$-model $a_{i}$ in $\cal K$, a generalized element $h_{i}:E_{i}\to F(a_{i})$ and a $\mathbb T$-model homomorphism $f_{i}:M\to a_{i}$ such that $f_{i}(z)=f_{i}(w)$ and $F(f_{i}) \circ h_{i}=x \circ e_{i}$. 
\end{proposition}\qed

Let $M_{\mathbb T}$ the universal model of $\mathbb T$ in the syntactic category ${\cal C}_{\mathbb T}$. We can represent the model $\tilde{F}(M_{\mathbb T})$ as a $\cal E$-indexed filtered colimit of set-based models of $\mathbb T$. For simplicity, we shall first establish this representation in the case of the topos of sets, and then generalize it to an arbitrary topos.

Let $F:{\cal K}^{\textrm{op}}\to \Set$ be a flat functor. By Theorem \ref{extsymthm}, for any sort $A$ over $\Sigma$, $\tilde{F}(M_{\mathbb T})A=\tilde{F}(\{x^{A}. \top\})=colim(F\circ \pi_{\{x^{A}. \top\}})$, where $\pi_{\{x^{A}. \top\}}$ is the canonical projection functor ${\cal A}_{\{x^{A}. \top\}}\to {\cal K}^{\textrm{op}}$ to ${\cal K}^{\textrm{op}}$ from the category of elements ${\cal A}_{\{x^{A}. \top\}}$ of the functor $P_{A}:{\cal K}\to \Set$ which assigns any model $M$ in $\cal K$ to the set $MA$ and acts accordingly on the arrows. Now, it follows from the commutativity of the tensor product between a presheaf and a covariant set-valued functor (cf. chapter VII of \cite{MM} or section \ref{tensorpr} above) that the colimit $colim(F\circ \pi_{\{x^{A}. \top\}})$ is isomorphic to the colimit $colim(P_{A}\circ \pi_{F})$, where $\pi_{F}:{\int F} \to {\cal K}$ is the canonical projection functor from the category of elements of the functor $F$ to ${\cal K}$. Since the functor $F$ is flat, the category ${\int F} $ is filtered. Therefore, as filtered colimits in ${\mathbb T}\textrm{-mod}(\Set)$ are computed sortwise as in $\Set$, $\tilde{F}(M_{\mathbb T})\cong colim(i\circ \pi_{F})$, where $i:{\cal K}\hookrightarrow {\mathbb T}\textrm{-mod}(\Set)$ is the canonical inclusion. So for any object $(c, x)$ of the category ${\int F}$ we have a $\mathbb T$-model homomorphism $\xi_{(c, x)}:c\to \tilde{F}(M_{\mathbb T})$, which can be expressed in terms of the colimit arrows $\kappa_{(a,y)}:F(a)\to \tilde{F}(\{x^{A}. \top\})=\tilde{F}(M_{\mathbb T})A$ (for $y\in cA$ and $A$ sort over $\Sigma$) as follows: for any sort $A$ over $\Sigma$, $\xi_{(c,x)}A(y)=\kappa_{(c,y)}(x)$. The explicit description of filtered colimits in the category $\Set$ thus yields, for each sort $A$ over $\Sigma$, the following characterizing properties of the colimit $colim(P_{A}\circ \pi_{F})$ in terms of the arrows $\xi_{(c, x)}$:

\begin{enumerate}[(i)]
\item For any element $x$ of $\tilde{F}(M_{\mathbb T})A$ there exists an object $(c, x)$ of the category ${\int F}$ and an element $y$ of $cA$ such that $\xi_{(c, x)}A(y)=x$; 

\item For any objects $(c, x)$ and $(c', x')$ of the category ${\int F}$ and elements $y\in cA$ and $y'\in c'A$, we have $\xi_{(c, x)}A(y)=\xi_{(c', x')}A(y')$ if and only if there exists an object $(c'', x'')$ of ${\int F}$ and arrows $f:c\to c''$ and $g:c'\to c''$ in $\cal K$ such that $F(f)(x'')=x$, $F(g)(x'')=x'$ and $fA(y)=gA(y')$.    
\end{enumerate}

Moreover, the filteredness of the category ${\int F}$ implies the following `joint embedding property': for any objects $(c, x)$ and $(c', x')$ of the category ${\int F}$ there exists an object $(c'', x'')$ of ${\int F}$ and arrows $f:c\to c''$ and $g:c'\to c''$ in $\cal K$ such that $F(f)(x'')=x$, $F(g)(x'')=x'$ and $\xi_{(c, x)}\circ f=\xi_{(c, x)}\circ g$.

Let us now proceed to establish the $\cal E$-indexed generalization of this result.

Suppose that $F:{\cal K}^{\textrm{op}}\to {\cal E}$ is a flat functor, where ${\cal K}$ is a small subcategory of the category ${\mathbb T}\textrm{-mod}(\Set)$. By Theorem \ref{extsymthm}, for any sort $A$ over $\Sigma$, we have that $\tilde{F}(M_{\mathbb T})A=\tilde{F}(\{x^{A}. \top\})=colim(F\circ \pi_{\{x^{A}. \top\}})$; in particular, for any object $(a, y)$ of the category ${\cal A}_{\{x^{A}. \top\}}$ we have a colimit arrow $\kappa_{(a,y)}:F(a)\to \tilde{F}(\{x^{A}. \top\})=\tilde{F}(M_{\mathbb T})A$ in $\cal E$. More generally, for any formula-in-context $\phi(\vec{x})$ over $\Sigma$, any model $c$ in $\cal K$ and any element $\vec{y}$ of $[[\vec{x}. \phi]]_{c}$, we have a colimit arrow $\kappa_{(c, \vec{y})}:F(c)\to \tilde{F}(\{\vec{x}.\phi\})$.

\begin{proposition}
Let $\mathbb T$ be a geometric theory over a signature $\Sigma$, ${\cal K}$ a small subcategory of the category ${\mathbb T}\textrm{-mod}(\Set)$, $\cal E$ a Grothendieck topos and $F:{\cal K}^{\textrm{op}}\to {\cal E}$ a flat functor. With the above notation, for any pair $(c, x)$ consisting of an object $c$ of $\cal K$ and of a generalized element $x:E\to F(c)$, there is a $\Sigma$-structure homomorphism $\xi_{(c, x)}:c \to Hom_{\cal E}(E, \tilde{F}(M_{\mathbb T}))$ defined as follows: for any sort $A$ over $\Sigma$, the function $\xi_{(c, x)}:cA \to Hom_{\cal E}(E, \tilde{F}(\{x^{A}. \top\}))$ sends any element $y\in cA$ to the generalized element $\kappa_{(c, y)} \circ x:E \to \tilde{F}(\{x^{A}. \top\})$.  
\end{proposition}

\begin{proofs}
We have to verify that:

\begin{enumerate}[(1)]

\item For any function symbol $f:A_{1}, \ldots, A_{n}\to B$ over $\Sigma$, the diagram
\[  
\xymatrix {
cA_{1}\times \cdots \times cA_{n} \ar[d]^{\xi_{(c, x)}A_{1}\times \cdots \times \xi_{(c, x)}A_{n}} \ar[rrr]^{cf} & & & cB \ar[dd]^{\xi_{(c, x)}B}\\
Hom_{\cal E}(E, \tilde{F}(\{x^{A_{1}}. \top\}) \times \cdots \times \tilde{F}(\{x^{A_{n}}. \top\})) \ar[d]^{i} & & & \\
Hom_{\cal E}(E, \tilde{F}(\{x^{A_{1}}, \ldots, x^{A_{n}}. \top\})) \ar[rrr]_{Hom_{\cal E}(E, \tilde{F}([f]))} & & & Hom_{\cal E}(E, \tilde{F}(\{x^{B}. \top\})),} 
\] 
where $i$ is the canonical isomorphism $\tilde{F}(\{x^{A_{1}}. \top\}) \times \cdots \times \tilde{F}(\{x^{A_{n}}. \top\}) \cong \tilde{F}(\{x^{A_{1}}, \ldots, x^{A_{n}}. \top\})$ induced by the preservation of finite products by $\tilde{F}$, commutes;

\item For any relation symbol $R\mono A_{1}, \ldots, A_{n}$ over $\Sigma$, we have a commutative diagram
\[  
\xymatrix {
cR \ar[d] \ar[rrr] & & & Hom_{\cal E}(E, \tilde{F}(\{x^{A_{1}}, \ldots, x^{A_{n}}. R)) \ar[d]\\
cA_{1}\times \cdots \times cA_{n}  \ar[rrr] & & & Hom_{\cal E}(E, \tilde{F}(\{x^{A_{1}}, \ldots, x^{A_{n}}. \top\})),} 
\] 
where the unnamed vertical arrows are the canonical ones and the lower horizontal arrow is ${Hom_{\cal E}(E, i)\circ \xi_{(c, x)}A_{1}\times \cdots \times \xi_{(c, x)}A_{n}}$.
 
\end{enumerate} 

To prove $(1)$, we first observe that for any $n$-tuple $\vec{y}=(y_{1}, \ldots, y_{n})\in cA_{1}\times \cdots \times cA_{n}$, $i\circ <\kappa_{(c, y_{1})}, \ldots, \kappa_{(c, y_{n})}>=\kappa_{(c, \vec{y})}$. Indeed, for any $i\in \{1, \ldots, n\}$, the arrow $\tilde{F}(\pi_{i}) \circ i$, where $\pi_{i}:\{x^{A_{1}}, \ldots, x^{A_{n}}. \top\} \to \{x^{A_{i}}. \top\}$ is the canonical projection arrow in the syntactic category ${\cal C}_{\mathbb T}$, is equal to the $i$-th product projection $\tilde{F}(\{x^{A_{1}}. \top\}) \times \cdots \times \tilde{F}(\{x^{A_{n}}. \top\}) \to \tilde{F}(\{x^{A_{i}}. \top\})$, by Proposition \ref{alphaext}(iii). Therefore, to prove the required condition it is equivalent to verify that for any $n$-tuple $\vec{y}=(y_{1}, \ldots, y_{n})\in cA_{1}\times \cdots \times cA_{n}$, $\xi_{(c, x)}B(cf(\vec{y}))=\tilde{F}([f]) \circ \kappa_{(c,\vec{y})} \circ x$, where $[f]:\{x^{A_{1}}, \ldots, x^{A_{n}}. \top\}\to \{x^{B}. \top\}$ is the arrow in the syntactic category ${\cal C}_{\mathbb T}$ corresponding to the function symbol $f$. But $\xi_{(c, x)}B(cf(\vec{y}))=\kappa_{(c, cf(\vec{y}))} \circ x$, and we have that $\tilde{F}([f]) \circ \kappa_{(c,\vec{y})}=\kappa_{(c, cf(\vec{y}))}$, again by Proposition \ref{alphaext}(iii).  

To prove $(2)$, it suffices to observe that, by Proposition \ref{alphaext}(iii), for any $n$-tuple $\vec{y}=(y_{1}, \ldots, y_{n})\in cA_{1}\times \cdots \times cA_{n}$ in $cR$, $\kappa_{(c, \vec{y})}$ factors through the canonical subobject $\tilde{F}(R)\mono \tilde{F}(\{x^{A_{1}}, \ldots, x^{A_{n}}. \top\})$.  

\end{proofs}

The following lemma describes some basic properties of the homomorphisms $\xi_{(c, x)}$.

\begin{lemma}\label{lemmaxi}
Let $\mathbb T$ be a geometric theory over a signature $\Sigma$, ${\cal K}$ a small subcategory of the category ${\mathbb T}\textrm{-mod}(\Set)$ and $F:{\cal K}^{\textrm{op}}\to {\cal E}$ a flat functor with values in a Grothendieck topos $\cal E$. Then
\begin{enumerate}[(i)]
\item For any generalized element $x:E\to F(c)$ and any arrow $f:d\to c$ in $\cal K$, $\xi_{(c, x)} \circ f=\xi_{(c, F(f)\circ x)}$;

\item For any generalized element $x:E\to F(c)$ and any arrow $e:E' \to E$, $\xi_{(c, x\circ e)}=Hom(e, \tilde{F}(M_{\mathbb T})) \circ \xi_{(c, x)}$
\end{enumerate}
\end{lemma}

\begin{proofs}
These properties can be easily proved by using the definition of the arrows $\xi_{(c, x)}$ in terms of the arrows $\kappa_{(a, y)}$ and Proposition \ref{alphaext}(iii).
\end{proofs}

\begin{proposition}\label{res}
Let $\mathbb T$ be a geometric theory over a signature $\Sigma$, ${\cal K}$ a small subcategory of the category ${\mathbb T}\textrm{-mod}(\Set)$ and $F:{\cal K}^{\textrm{op}}\to {\cal E}$ a flat functor with values in a Grothendieck topos $\cal E$. Let $M$ be the $\mathbb T$-model $\tilde{F}(M_{\mathbb T})$ in $\cal E$. Then, for any sort $A$ over $\Sigma$, we have that

\begin{enumerate}[(i)] 
\item For any generalized element $x:E\to MA$ there exists an epimorphic family $\{e_{i}:E_{i}\to E \textrm{ | } i\in I\}$ in $\cal E$, for each index $i\in I$ a $\mathbb T$-model $a_{i}$ in $\cal K$, a generalized element $x_{i}:E_{i}\to F(a_{i})$ and an element $y_{i}\in a_{i}A$ such that $\xi_{(a_{i}, x_{i})}A(y_{i})=x\circ e_{i}$. 

\item For any pairs $(a, x)$ and $(b, x')$, where $a$ and $b$ are $\mathbb T$-models in $\cal K$ and $x:E\to F(a)$, $x':E\to F(b)$ are generalized elements, and any elements $y\in aA$ and $y'\in bA$, we have that $\xi_{(a, x)}A(y)=\xi_{(b, x')}A(y')$ if and only if there exists an epimorphic family $\{e_{i}:E_{i}\to E \textrm{ | } i\in I\}$ in $\cal E$ in $\cal E$ and for each index $i\in I$ a $\mathbb T$-model $c_{i}$ in $\cal K$, arrows $f_{i}:a\to c_{i}$, $g_{i}:b\to c_{i}$ and a generalized element $x_{i}:E_{i}\to F(c_{i})$ such that $<x, x'>\circ e_{i}=<F(f_{i}), F(g_{i})>\circ x_{i}$ for all $i\in I$ and $f_{i}A(y)=g_{i}A(y')$.
\end{enumerate}

Moreover, the following `joint embedding property' holds: for any pairs $(a, x)$ and $(b, x')$, where $a$ and $b$ are $\mathbb T$-models in $\cal K$ and $x:E\to F(a)$, $x':E\to F(b)$ are generalized elements, there exists an epimorphic family $\{e_{i}:E_{i}\to E \textrm{ | } i\in I\}$ in $\cal E$ in $\cal E$ and for each index $i\in I$ a $\mathbb T$-model $c_{i}$ in $\cal K$, arrows $f_{i}:a\to c_{i}$, $g_{i}:b\to c_{i}$ in $\cal K$ and a generalized element $x_{i}:E_{i}\to F(c_{i})$ such that $<x, x'>\circ e_{i}=<F(f_{i}), F(g_{i})>\circ x_{i}$ for all $i\in I$ and (by Lemma \ref{lemmaxi}) the following diagram commutes:

\[  
\xymatrix {
a \ar[dr]^{f} \ar[rr]^{\xi_{(a, x)}} & & Hom_{\cal E}(E, M) \ar[dr]^{Hom_{\cal E}(e_{i}, M)} & \\
& c_{i} \ar[rr]^{\xi_{(c_{i}, x_{i})}} & &  Hom_{\cal E}(E_{i}, M) \\
b \ar[ur]_{g} \ar[rr]^{\xi_{(b, x')}} & & Hom_{\cal E}(E, M) \ar[ur]_{Hom_{\cal E}(e_{i}, M)} &} 
\]

In fact, such family can be taken to be the pullback of the family of arrows $<F(f), F(g)>:F(c)\to F(a)\times F(b)$ (for all spans $(f:a\to c,\textrm{ }g:b\to c)$ in the category $\cal K$) along the arrow $<x,x'>:E\to F(a)\times F(b)$.  
\end{proposition} 
 
\begin{proofs}
Recall from Theorem \ref{extsymthm} that for any sort $A$ over $\Sigma$, $\tilde{F}(M_{\mathbb T})A=colim(F\circ \pi^{\cal K}_{\{x^{A}. \top\}})$; now, $\pi^{\cal K}_{\{x^{A}. \top\}}$ is precisely the functor $P_{A}:{\cal K}\to \Set$ sending a model $N$ in $\cal C$ to the set $NA$.
 
The proposition then straightforwardly follows from Theorem \ref{thmres}, applied to the functors $P_{A}:{\cal K}\to \Set$ (with $A$ varying among the sorts of $\Sigma$), in view of Theorems \ref{tensorpr} and \ref{altdes}. 
\end{proofs}

\subsection{A general adjunction}\label{generaladjunction}

Let $\cal C$ be a small category and $\cal E$ a Grothendieck topos. Recall from \cite{OC} that the indexed category $\underline{[{\cal C}, {\cal E}]}_{\cal E}$ is locally small, that is for any two functors $F, G:{\cal C} \to {\cal E}$ there exists an object $Hom^{\cal E}_{\underline{[{\cal C}, {\cal E}]}_{\cal E}}(F,G)$ of $\cal E$ satisfying the universal property that for any object $E$ of $\cal E$ the generalized elements $E \to Hom^{\cal E}_{\underline{[{\cal C}, {\cal E}]}_{\cal E}}(F,G)$ correspond bijectively, naturally in $E\in {\cal C}$, to the natural transformations $!_{E}^{\ast}\circ F \to !_{E}^{\ast}\circ G$.

The following theorem establishes a general adjunction between categories of $\cal E$-valued functors induced by a functor $P:{\cal C}\to [{\cal D}^{\textrm{op}}, \Set]$. 

Before stating it, we need to introduce some notation. For any object $c$ of $\cal C$, we denote by $\int P(c)$ the category of elements of the functor $P(c):{\cal D}^{\textrm{op}}\to \Set$ and by $\pi_{c}:\int P(c) \to {\cal D}$ the canonical projection functor. We denote by $y_{{\cal D}}:{\cal D}^{\textrm{op}}\to [{\cal D}, \Set]$ the Yoneda embedding. For any functor $H:{\cal D}\to {\cal E}$ and any object $(d, z)$ of the category $\int P(c)$, we denote by $\kappa_{(d, z)}:H(d)\to colim(H\circ \pi_{c})$ the canonical colimit arrow. 

\begin{theorem}\label{genadj}
Let $\cal C$ and $\cal D$ small categories, $\cal E$ a Grothendieck topos and $P:{\cal C}\to [{\cal D}^{\textrm{op}}, \Set]$ a functor. Let $\tilde{(-)}:[{\cal D}, {\cal E}]\to [{\cal C}, {\cal E}]$ be the functor sending to any functor $F:{\cal D}\to {\cal E}$ the functor $\tilde{F}$ defined by: 

for any $c\in {\cal C}$, $\tilde{F}(c)=colim(F\circ \pi_{c})$, and 

for any arrow $f:c\to c'$, $\tilde{F}(f):colim(F\circ \pi_{c}) \to colim(F\circ \pi_{c'})$ is defined by the conditions $\tilde{F}(f)\circ \kappa^{F}_{(d, z)}=\kappa^{F}_{(d, P(f)(z))}$ (for any object $(d, z)$ of the category $\int P(c)$),

and which acts on arrows in the obvious way.

Let $(-)_{r}:[{\cal C}, {\cal E}]\to [{\cal D}, {\cal E}]$ be the functor assigning to any functor $G:{\cal C}\to {\cal E}$ the functor $Hom^{\cal E}_{\underline{[{\cal C}, {\cal E}]}_{\cal E}}(\gamma_{\cal E}^{\ast} \circ \tilde{y_{{\cal D}}-}, G)$, and acting on the arrows in the obvious way. 

Then the functors 
\[
\tilde{(-)}:[{\cal D}, {\cal E}]\to [{\cal C}, {\cal E}]
\] 
and 
\[
(-)_{r}:[{\cal C}, {\cal E}]\to [{\cal D}, {\cal E}]
\]
are adjoint to each other ($\tilde{(-)}$ on the left and $(-)_{r}$ on the right). The unit $\eta:F \to (\tilde{F})_{r}$ and the counit $\epsilon^{G}:\tilde{(G_{r})}\to G$ are defined as follows:  

For any $d\in {\cal D}$, $\eta^{F}(d):F(d) \to (\tilde{F})_{r}(d)$ is the arrow in $\cal E$ defined by means of generalized elements by saying that it sends any generalized element $E\to F(d)$, regarded as a natural transformation $\gamma_{{\cal E}\slash E}^{\ast} \circ y_{\cal D}d \to !_{E}^{\ast} \circ F$, to the image of this arrow $\gamma_{{\cal E}\slash E}^{\ast} \circ \tilde{y_{\cal D}d} \to !_{E}^{\ast}\circ \tilde{F}$ under the functor $\tilde{(-)}$, regarded as a generalized element $E\to Hom^{\cal E}_{\underline{[{\cal C}, {\cal E}]}_{\cal E}}(\gamma_{\cal E}^{\ast} \circ \tilde{y_{{\cal D}}'d}, \tilde{F})=(\tilde{F})_{r}(d)$. 

For any $c\in {\cal C}$, $\epsilon^{G}(c): \tilde{(G_{r})}(c)=colim(G_{r}\circ \pi_{c}) \to G(c)$ is defined by setting, for each object $(d, z)$ of $\int P(c)$, $\epsilon^{G}(c) \circ \kappa^{G_{r}}_{(d,z)}$ equal to the arrow $G_{r}(d) \to G(c)$ defined by means of generalized elements by saying that a generalized element $x:E \to G_{r}(d)$, corresponding to a natural transformation $\overline{x}:!_{E}^{\ast} \circ \gamma^{\ast}_{\cal E} \circ \tilde{y_{\cal D}d} \cong \gamma_{{\cal E}\slash E}^{\ast} \circ \tilde{y_{\cal D}d} \to !_{E}^{\ast} \circ \tilde{G}$, is sent to the arrow $E\to G(c)$ obtained by composing $\overline{x}(c)$ with the component of $\gamma_{{\cal E}\slash E}^{\ast}(\tilde{y_{\cal D}d}(c))$ corresponding to the element of $\tilde{y_{\cal D}d}(c)$ given by the image of the identity on $d$ via the colimit arrow $\kappa^{ y_{\cal D}d }_{(d, z)}: (y_{\cal D}d)(d) \to \tilde{y_{\cal D}d}(c)$.      
\end{theorem}

\begin{proofs}
For simplicity we shall prove the result only in the case ${\cal E}=\Set$, the proof of the general case being entirely analogous (the only care that one has to take is to use generalized elements in place of standard set-theoretic elements).

We shall define a bijective correspondence between the natural transformations $\tilde{F}\to G$ and the natural transformations $F \to G_{r}$, natural in $F\in [{\cal D}, \Set]$ and $G\in [{\cal C}, \Set]$.

Given a natural transformation $\alpha: F \to G_{r}$, we define $\tau(\alpha):\tilde{F}\to G$ by setting, for each $c\in {\cal C}$, $\tau_{\alpha}(c):\tilde{F}(c)\to G(c)$ equal to the arrow defined as follows. As $\tilde{F}(c)$ is the colimit of the cone $\{\kappa^{F}_{(d, z)}:F(d) \to \tilde{F}(c) \textrm{ | } (d, z)\in \int P(c)\}$, it suffices to define an arrow $u_{(d, z)}:F(d)\to G(c)$ for each pair $(d, z)$ where $d\in {\cal D}$ and $z\in P(c)(d)$, checking that whenever $(d, z)$ and $(d', z')$ are pairs such that for an arrow $f:d\to d'$ in $\cal D$, $P(c)(f)(z')=z$ then $u_{(d', z')} \circ F(f)= u_{(d, z)}$; indeed, by the universal property of the colimit, such a family of arrows will induce a unique arrow $\tau_{\alpha}(c):\tilde{F}(c)\to G(c)$ such that $\tau_{\alpha}(c) \circ \kappa^{F}_{(d, z)}=u_{(d, z)}$ (for each $(d, z)$ in $\int P(c)$). We set $u_{(d, z)}:F(d)\to G(c)$ equal to the function sending every element $x\in F(d)$ to the element $\alpha(d)(x)(c)( \kappa^{y_{\cal D}'d}_{(d, z)}(1_{d}))$. Let us check that the compatibility condition is satisfied. Given an arrow $f:d\to d'$ such that $P(c)(f)(z')=z$, we have to verify that for every $x\in F(d)$, $\alpha(d')(Ff(x))(c)( \kappa^{y_{\cal D}'d'}_{(d', z')}(1_{d'}))=\alpha(d)(x)(c)( \kappa^{y_{\cal D}'d}_{(d, z)}(1_{d}))$. This identity immediately follows from the commutativity of the naturality diagram
\[  
\xymatrix {
F(d) \ar[d]^{F(f)} \ar[rrrr]^{\alpha(d)} & & & & G_{r}(d)=Hom_{[{\cal C}, \Set]}(\tilde{y_{\cal D}d}, G) \ar[d]^{G_{r}(f)=-\circ \tilde{y_{\cal D}f}}\\
F(d') \ar[rrrr]_{\alpha(d')} & & & & G_{r}(d')=Hom_{[{\cal C}, \Set]}(\tilde{y_{\cal D}d'}, G)} 
\] 
for $\alpha$ with respect to the arrow $f$ and that of the diagram
\[  
\xymatrix {
{\tilde{y_{\cal D}'d'}}(c)  \ar[r]^{f} & {\tilde{y_{\cal D}'d}}(c)   \\ \ar[u]^{\kappa^{y_{\cal D}d'}_{(d', z')}} y_{\cal D}d' \ar[r]_{(y_{\cal D}f)(d')} & y_{\cal D}d \ar[u]^{ \kappa^{y_{\cal D}d}_{(d', z')}}, } 
\] 
which is an instance of Proposition \ref{alphaext}(i). 

Indeed, $\alpha(d')(Ff(x))(c)( \kappa^{y_{\cal D}d'}_{(d', z')}(1_{d'}))=\alpha(d)(x)(\tilde{y_{\cal D}f}(c)(\kappa^{y_{\cal D}d'}_{(d', z')}(1_{d'})))=\alpha(d)(x)(\kappa^{y_{\cal D}'d}_{(d', z')}(f))=\alpha(d)(x)(c)( \kappa^{y_{\cal D}'d}_{(d, z)}(1_{d}))$, where the first equality follows from the commutativity of the first diagram, the second follows from the commutativity of the second diagram, and the last follows from the fact that the family $\{\kappa^{y_{\cal D}d}_{(d, z)} \textrm{ | } (d, z)\in \int P(c)\}$ is a cocone.  

To complete the definition of $\tau(\alpha)$, it remains to check that the assignment $c\to \tau(\alpha)(c)$ defines a natural transformation $\tilde{F}\to G$. We have to verify that for any arrow $f:c\to c'$ in $\cal C$, the diagram 
\[  
\xymatrix {
\tilde{F}(c) \ar[d]^{\tilde{F}(f)}  \ar[r]^{\tau(\alpha)(c)} & G(c) \ar[d]^{G(f)} \\ 
\tilde{F}(c') \ar[r]_{\tau(\alpha)(c')} & G(c') } 
\] 
commutes.

As, by definition of $\tilde{F}(f)$, the diagram
\[  
\xymatrix {
\tilde{F}(c)  \ar[r]^{\tilde{F}(f)} & \tilde{F}(c') \\
F(d) \ar[u]^{\kappa^{F}_{(d,z)}} \ar[ur]_{\kappa^{F}_{(d,P(f)(z))}} & } 
\]
commutes and the arrows $\kappa^{F}_{(d,z)}$ are jointly epimorphic (as they are colimit arrows), the square above commutes if and only if the arrows $G(f)\circ \tau(\alpha)(c)\circ \kappa^{F}_{(d,z)}$ and $\tau(\alpha)(c')\circ \kappa^{F}_{(d,P(f)(z))}$ are equal, that is if and only if they take the same values at any element $x\in F(d)$. Let us set $z'=P(f)(d)(z)$. By definition of $\tau(\alpha)$, $\tau(\alpha)(c)( \kappa^{F}_{(d,z)}(x))= \alpha(d)(x)(c)(\kappa^{y_{\cal D}d}_{(d, z)}(1_{d}))$, while $\tau(\alpha)(c')(\kappa^{F}_{(d, z')}(x))=\alpha(d)(x)(c')(\kappa^{y_{\cal D}d}_{(d, z')}(1_{d}))$. Now, for any $d\in {\cal D}$ and any $x\in F(d)$, $\alpha(d)(x)$ is a natural transformation $\tilde{yd}\to G$; in particular, the diagram
\[  
\xymatrix {
\tilde{yd}(c)  \ar[d]^{\tilde{yd}(f)} \ar[rr]^{\alpha(d)(x)(c)} & & G(c) \ar[d]^{G(f)} \\
\tilde{yd}(c')  \ar[rr]_{\kappa^{F}_{(d,P(f)(z))}} & & G(c')} 
\]
commutes. 

The commutativity of this diagram, together with that of the diagram
\[  
\xymatrix {
\tilde{yd}(c)  \ar[r]^{\tilde{yd}(f)} & \tilde{yd}(c') \\
yd(d) \ar[u]^{\kappa^{F}_{(d,z)}} \ar[ur]_{\kappa^{yd}_{(d,z')}} & } 
\]
(which follows by the definition of $\tilde{yd}(f)$) now immediately implies our thesis.  
    
Let us now define a function $\chi$ assigning to any natural transformation $\beta:\tilde{F}\to G$ a natural transformation $\chi(\beta):F\to G_{r}$. For any $d\in {\cal D}$, we define $\chi(\beta)(d):F(d)\to G_{r}(d)$ by setting, for each $x\in F(d)$, $\chi(\beta)(x)$ equal to the natural transformation $\beta \circ \tilde{a_{x}}:\tilde{yd}\to G$, where $a_{x}:yd\to F$ is the natural transformation corresponding, via the Yoneda lemma, to the element $x\in F(c)$.

To verify that $\chi(\beta)$ is well-defined, we have to check that for every arrow $g:d\to d'$ in $\cal D$, the diagram
\[  
\xymatrix {
F(d)  \ar[d]^{F(g)} \ar[r]^{\chi(\beta)(d)} & G_{r}(d) \ar[d]^{G_{r}(g)=-\circ \tilde{yd}} \\
F(d')  \ar[r]_{\chi(\beta)(d')} & G_{r}(d')} 
\]
commutes, i.e. for every $x\in F(d)$, $G_{r}(g)(\chi(\beta)(d)(x))=\chi(\beta)(d')(F(g)(x))$. Now, $G_{r}(g)(\chi(\beta)(d)(x))=\beta \circ \tilde{\alpha_{x}} \circ \tilde{yg}$, while $\chi(\beta)(d')(F(g)(x))=\beta \circ \tilde{\alpha_{F(g)(x)}}$; but $\alpha_{F(g)(x)}=\alpha_{x}\circ y_{\cal D}g$, from which our thesis follows. 

The proof of the fact that the correspondences $\tau$ and $\beta$ are natural in $F$ and $G$ is straightforward and left to the reader. 

To conclude the proof of the theorem, it thus remains to show that $\tau$ and $\chi$ are inverse to each other. The verification of the fact that the unit and counit of the adjunction coincide with those given in the statement of the theorem is straightforward and left to the reader.

Let us show that for any natural transformation $\alpha:F \to G_{r}$, $\chi(\tau(\alpha))=\alpha$. Let us set $\beta=\tau(\alpha)$. We have to prove that for any $d\in {\cal D}$, $x\in F(d)$ and $z\in P(c)(d')$, $\alpha(d)(x)(c)\circ \kappa^{y_{\cal D}d}_{(d',z)}= \chi(\beta)(d)(x)(c) \circ \kappa^{y_{\cal D}d}_{(d',z)}$ as functions $y_{\cal D}d(d')\to G(c)$, i.e. that for any element $g\in y_{\cal D}d(d')$, $\alpha(d)(x)(c)(\kappa^{y_{\cal D}d}_{(d',z)}(g))= (\chi(\beta)(d)(x)(c) \circ \kappa^{y_{\cal D}d}_{(d',z)})(g)$. Now, by definition of the functor $\tilde{(-)}$ and of the correspondence $\tau$, the diagram
\[  
\xymatrix {
\tilde{y_{\cal D}d}(c)  \ar[r]^{\tilde{a_{x}}(c)} & \tilde{F}(c) \ar[r]^{\beta(c)} & G(c) \\
y_{\cal D}d(d') \ar[u]^{\kappa^{y_{\cal D}d}_{(d',z)}}  \ar[r]_{a_{x}(d')} & F(d') \ar[u]^{\kappa^{F}_{(d',z)}} \ar[u]^{\kappa^{F}_{(d',z)}} \ar[ur]_{u_{(d',z)}} & } 
\]
commutes, and since $\chi(\beta)(x)=\beta \circ \tilde{a_{x}}:\tilde{y_{\cal D}d}\to G$, we have that $(\chi(\beta)(d)(x)(c) \circ \kappa^{y_{\cal D}d}_{(d',z)})(g)$ is equal to $\alpha(d')(F(g)(x))(c)(\kappa^{y_{\cal D}d'}_{(d',z)}(1_{d'}))$.

Now, the naturality diagram for $\alpha$ with respect to the arrow $g:d\to d'$ yields the equality $\alpha(d)(x)\circ \tilde{y_{\cal D}g}=\alpha(d')(Fg(x))$ and hence the equality $\alpha(d)(x)(c)\circ \tilde{y_{\cal D}g}(c) \circ \kappa^{y_{\cal D}d'}_{(d', z)}=\alpha(d')(Fg(x))(c) \circ \kappa^{y_{\cal D}d'}_{(d', z)}$.

But the commutativity of the diagram
\[  
\xymatrix {
\tilde{y_{\cal D}d'}(c)  \ar[r]^{\tilde{y_{\cal D}d}(g)(c)} & \tilde{y_{\cal D}d}(c') \\
y_{\cal D}d'(d') \ar[u]^{\kappa^{y_{\cal D}d'}_{(d',z)}} \ar[r]^{y_{\cal D}g(d')} & y_{\cal D}d(d') \ar[u]_{\kappa^{y_{\cal D}d}_{(d',z)}} } 
\]
(which follows by definition of the functor $\tilde{(-)}$) implies that $\alpha(d)(x)(c)\circ \tilde{y_{\cal D}g}(c) \circ \kappa^{y_{\cal D}d'}_{(d', z)}= \alpha(d)(x)(c)\circ \kappa^{y_{\cal D}d}_{(d', z)} \circ y_{\cal D}g$. Therefore $\alpha(d)(x)(c)\circ \kappa^{y_{\cal D}d}_{(d', z)} \circ y_{\cal D}g = \alpha(d')(Fg(x))(c) \circ \kappa^{y_{\cal D}d'}_{(d', z)}$ and evaluating at $1_{d'}$ yields the desired equality.

Finally, let us show that the composite $\tau \circ \chi$ is equal to the identity. Let $\beta:\tilde{F}\to G$ be a natural transformation. We have to show that $\beta=\tau(\chi(\beta))$. Let us set $\alpha=\chi(\beta)$. To prove that $\beta=\tau(\alpha)$ it is equivalent to verify that for every object $c\in {\cal C}$, any pair $(d, z)$ with $d\in {\cal D}$ and $z\in P(c)(d)$ and any element $x\in F(d)$, $\beta(c)( \kappa^{F}_{(d, z)}(x))=\tau(\alpha)(c)(\kappa^{F}_{(d, z)}(x))$. But $\tau(\alpha)(c)(\kappa^{F}_{(d, z)}(x))=\alpha(d)(x)(c)(\kappa^{y_{\cal D}d}_{(d, z)}(1_{d}))= \chi(\beta)(x)(c)(\kappa^{y_{\cal D}d}_{(d, z)}(1_{d}))=(\beta \circ \tilde{a_{x}})(c)(\kappa^{y_{\cal D}d}_{(d, z)}(1_{d}))=\beta(c)((\tilde{a_{x}})(c)(\kappa^{y_{\cal D}d}_{(d, z)}(1_{d})))$. Thus our thesis follows from the commutativity of the diagram 
\[  
\xymatrix {
\tilde{y_{\cal D}d}(c) \ar[r]^{\tilde{a_{x}}(c)} & \tilde{F}(c) \\
yd(d) \ar[u]^{\kappa^{y_{\cal D}d}_{(d,z)}} \ar[r]^{a_{x}(d)} & F(d) \ar[u]_{\kappa^{F}_{(d,z)}}, } 
\]
which is an immediate consequence of the definition of the functor $\tilde{(-)}$.   

\end{proofs}

\begin{remarks}\label{remadj}
\begin{enumerate}[(a)]
\item Let $f:{\cal D}\to {\cal C}$ be a functor between two small categories, and let $P:{\cal C}\to [{\cal D}^{\textrm{op}}, \Set]$ be the functor $y_{\cal C}(-)\circ f^{\textrm{op}}$. Notice that $P$ is the flat functor corresponding, via Diaconescu's equivalence, to the essential geometric morphism $[{\cal D}^{\textrm{op}}, \Set]\to [{\cal C}^{\textrm{op}}, \Set]$ induced by the functor $f^{\textrm{op}}:{\cal D}^{\textrm{op}}\to {\cal C}^{\textrm{op}}$. The functor $\tilde{(-)}:[{\cal D}, {\cal E}]\to [{\cal C}, {\cal E}]$ coincides with the left Kan extension functor along $f$, while the right adjoint functor $(-)_{r}$ coincides with the functor $-\circ f:[{\cal C}, {\cal E}]\to [{\cal D}, {\cal E}]$.     

\item Let $u:\Sh({\cal D}, K) \to \Sh({\cal C}, J)$ be a geometric morphism. Then for any Grothendieck topos $\cal E$, $u$ induces as in section \ref{ext} a functor
\[
\xi_{{\cal E}}:\textbf{Flat}_{K}({\cal D}, {\cal E})\to \textbf{Flat}_{J}({\cal C}, {\cal E}).
\]  
By Diaconescu's equivalence, the morphism $u$ corresponds to a flat functor ${\cal C}\to \Sh({\cal D}, K)$, which composed with the canonical geometric inclusion $\Sh({\cal D}, K)\hookrightarrow [{\cal D}^{\textrm{op}}, \Set]$ yields a functor $P:{\cal C}\to [{\cal D}^{\textrm{op}}, \Set]$. Then $\xi_{{\cal E}}$ coincides with the restriction of the functor $\tilde{(-)}:[{\cal D}, {\cal E}]\to [{\cal C}, {\cal E}]$ induced by $P$ as in the theorem to the full subcategories $\textbf{Flat}_{K}({\cal D}, {\cal E})\hookrightarrow  [{\cal D}, {\cal E}]$ and $\textbf{Flat}_{J}({\cal C}, {\cal E})\hookrightarrow  [{\cal C}, {\cal E}]$. In general, the right adjoint functor $(-)_{r}:[{\cal C}, {\cal E}]\to [{\cal D}, {\cal E}]$ does not restrict to these subcategories, but if it does, it becomes a right adjoint to the functor $\xi_{{\cal E}}:\textbf{Flat}_{K}({\cal D}, {\cal E})\to \textbf{Flat}_{J}({\cal C}, {\cal E})$.

This notably applies in the case of the geometric morphism
\[
p_{\cal K}:[{\cal K}, \Set] \to \Sh({\cal C}_{\mathbb T}, J_{\mathbb T})
\]
considered in section \ref{extsym}. In this case, ${\cal C}={\cal C}_{\mathbb T}$, ${\cal D}={{\cal K}}^{\textrm{op}}$ and $P$ is the functor ${\cal C}_{\mathbb T}\to [{\cal K}, \Set]$ sending any geometric formula $\phi(\vec{x})$ to the functor $M \to [[\vec{x}. \phi]]_{M}$.

\item By a basic property of adjoint functors, the (left adjoint) functor 
\[
\tilde{(-)}:[{\cal D}, {\cal E}]\to [{\cal C}, {\cal E}]
\]
is full and faithful if and only if the unit $\eta^{F}:F \to (\tilde{F})_{r}$ is an isomorphism for any $F$ in $[{\cal D}, {\cal E}]$. It easily follows (by purely formal considerations) that for every full subcategory $\cal H$ of $[{\cal D}, {\cal E}]$ with canonical embedding $i:{\cal H}\hookrightarrow [{\cal D}, {\cal E}]$, the composite functor $\tilde{(-)}\circ i$ is full and faithful if and only if the unit $\eta^{F}:F \to (\tilde{F})_{r}$ is an isomorphism for any $F$ in $\cal H$.    

\end{enumerate}
\end{remarks}

The following proposition will be useful in the sequel; see section \ref{filt} for the notation employed in it.

\begin{proposition}\label{finalcolimit}
Let $f:{\cal D}\to {\cal C}$ be a functor between small categories, $\cal E$ a Grothendieck topos and $F$ a functor ${\cal C}^{\textrm{op}}\to {\cal E}$. Then the $\cal E$-indexed functor $\underline{\int (F\circ f^{\textrm{op}})}_{\cal E} \to \underline{\int F}_{\cal E}$ sending any object $(c, x)$ of $\int (F\circ f^{\textrm{op}})_{E}$ to the object $(f(c), x)$ of $({\int F})_{E}$ is $\cal E$-final if and only if the $\cal E$-indexed functor $\underline{{\int \epsilon^{F}}}_{\cal E}:\underline{\int (\tilde{F\circ f^{\textrm{op}}})}_{\cal E} \to \underline{\int F}_{\cal E}$ induced by the natural transformation $\epsilon^{F}:\tilde{F\circ f} \to F$ of Theorem \ref{genadj} is $\cal E$-final. 
\end{proposition}

\begin{proofs}
From the general analysis of section \ref{gene} we know that for any $c\in {\cal C}$, $\tilde{F\circ f^{\textrm{op}}}(c)=colim(F\circ f^{\textrm{op}}\circ \pi_{c})$, where $\pi_{c}$ is the canonical projection to ${\cal D}^{\textrm{op}}$ from the category ${\cal A}_{c}$ whose objects are the pairs $(d, h)$, where $d$ is an object of the category $\cal D$ and $h$ is an arrow $c \to f(d)$ in $\cal C$, and whose arrows are the obvious ones. The arrow $\epsilon^{F}:\tilde{F\circ f} \to F$ is defined by the property that for any object $(d, h)$ of ${\cal A}_{c}$, $\epsilon^{F}(c)\circ \kappa^{F}_{(d, h)}=F(h)$, where $\kappa^{F}_{(d, h)}:F(f(d))\to colim(F\circ f^{\textrm{op}}\circ \pi_{c})$ is the canonical colimit arrow.
  
Now, the thesis follows immediately from the fact that for any object of the category ${\int (\tilde{F\circ f^{\textrm{op}}})}_{E}$ of the form $(c, \kappa^{F}_{(d, h)}\circ y )$, where $y$ is a generalized element $E\to F(f(d))$, $({\int \epsilon^{F}})_{E}((c, \kappa^{F}_{(d, h)}\circ y ))=(c, F(h)\circ y)$, invoking the fact that the colimit arrows are jointly epimorphic.     
\end{proofs}

\section{Preliminary results on theories of presheaf type}

In this section we establish some results on theories of presheaf type which will be important in the sequel.

\subsection{A syntactic description of finitely presentable models}\label{syntacticdescription}

For a theory of presheaf type $\mathbb T$, it is possible to give an explicit syntactic description of the finitely presentable $\mathbb T$-models; specifically, we have the following result.

\begin{theorem}\label{syn}
Let $\mathbb T$ be a theory of presheaf type over a signature $\Sigma$ and $\{\vec{x}. \phi\}$ be a formula over $\Sigma$ presenting a $\mathbb T$-model $U_{\{\vec{x}. \phi\}}$. Then $U_{\{\vec{x}. \phi\}}$ is isomorphic to the $\Sigma$-structure $M_{\{\vec{x}. \phi\}}$ defined as follows:
\begin{enumerate}[(i)]
\item for any sort $A$ over $\Sigma$, $M_{\{\vec{x}. \phi\}}A$ is equal to the set $Hom_{{\cal C}_{\mathbb T}}(\{\vec{x}. \phi\}, \{x^{A}. \top\})$ of $\mathbb T$-provably functional geometric formulae from $\{\vec{x}. \phi\}$ to $\{x^{A}. \top\}$;

\item for any function symbol $f:A_{1}, \ldots, A_{n}\to B$ over $\Sigma$, the function $M_{\{\vec{x}. \phi\}}f:Hom_{{\cal C}_{\mathbb T}}(\{\vec{x}. \phi\}, \{x^{A_{1}}. \top\})\times \cdots \times Hom_{{\cal C}_{\mathbb T}}(\{\vec{x}. \phi\}, \{x^{A_{n}}. \top\})\cong  Hom_{{\cal C}_{\mathbb T}}(\{\vec{x}. \phi\}, \{x^{A_{1}}, \ldots, x^{A_{n}}. \top\})  \to Hom_{{\cal C}_{\mathbb T}}(\{\vec{x}. \phi\}, \{x^{B}. \top\})$ is equal to $[f]\circ -$ (where $[f]:\{x^{A_{1}}, \ldots, x^{A_{n}}. \top\} \to \{x^{B}. \top\}$ is the morphism in ${\cal C}_{\mathbb T}$ corresponding to $f$);

\item for any relation symbol $R\mono A_{1}, \ldots, A_{n}$, $M_{\{\vec{x}. \phi\}}R$ is the subobject of $Hom_{{\cal C}_{\mathbb T}}(\{\vec{x}. \phi\}, \{x^{A_{1}}. \top\})\times \cdots \times Hom_{{\cal C}_{\mathbb T}}(\{\vec{x}. \phi\}, \{x^{A_{n}}. \top\})\cong  Hom_{{\cal C}_{\mathbb T}}(\{\vec{x}. \phi\}, \{x^{A_{1}}, \ldots, x^{A_{n}}. \top\})$ given by $[R]\circ -$, (where $[R]:\{x^{A_{1}}, \ldots, x^{A_{n}}. R\}\mono \{x^{A_{1}}, \ldots, x^{A_{n}}. \top\}$ is the subobject in ${\cal C}_{\mathbb T}$ corresponding to $R$).  
\end{enumerate}

\end{theorem}

\begin{proofs}
Recall that we have a canonical equivalence of categories
\[
\mathbf{Flat}_{J_{\mathbb T}}({\cal C}_{\mathbb T}, \Set)\simeq {\mathbb T}\textrm{-mod}(\Set),
\]
sending any flat $J_{\mathbb T}$-continuous functor $F:{\cal C}_{\mathbb T}\to \Set$ to the $\mathbb T$-model $F(M_{\mathbb T})$, where $M_{\mathbb T}$ is the universal model of $\mathbb T$ in ${\cal C}_{\mathbb T}$.

We know from \cite{OCS} that for any theory of presheaf type $\mathbb T$ over a signature $\Sigma$, any formula-in-context $\{\vec{x}. \phi\}$ over $\Sigma$ which presents a $\mathbb T$-model is $\mathbb T$-irreducible, in the sense that every $J_{\mathbb T}$-covering sieve on $\{\vec{x}. \phi\}$ in ${\cal C}_{\mathbb T}$ is maximal. From this it easily follows that the (flat) representable functor $Hom_{{\cal C}_{\mathbb T}}(\{\vec{x}. \phi\}, -):{\cal C}_{\mathbb T}\to \Set$ is $J_{\mathbb T}$-continuous; indeed, for any $J_{\mathbb T}$-covering sieve $S$ on an object $\{\vec{y}. \psi\}$ of ${\cal C}_{\mathbb T}$, any arrow $\gamma:\{\vec{x}. \phi\}\to \{\vec{y}. \psi\}$ in ${\cal C}_{\mathbb T}$ factors through one of the arrows belonging to $S$ as the pullback of $S$ along $\gamma$ is $J_{\mathbb T}$-covering and hence maximal. The image $Hom_{{\cal C}_{\mathbb T}}(\{\vec{x}. \phi\}, M_{\mathbb T})$ under this functor of the model $M_{\mathbb T}$, which clearly coincides with the $\Sigma$-structure $M_{\{\vec{x}. \phi\}}$ in the statement of the theorem, is therefore a $\mathbb T$-model. In order to deduce our thesis, it thus remains to verify that the model $M_{\{\vec{x}. \phi\}}$ satisfies the universal property of the $\mathbb T$-model presented by the formula $\{\vec{x}. \phi\}$, i.e. that for any $\mathbb T$-model $N$ in $\Set$, the $\mathbb T$-model homomorphisms $M_{\{\vec{x}. \phi\}}=Hom_{{\cal C}_{\mathbb T}}(\{\vec{x}. \phi\}, M_{\mathbb T}) \to N$ are in natural bijection with the elements of the set $[[\vec{x}. \phi]]_{N}$. But the $\mathbb T$-model homomorphisms $Hom_{{\cal C}_{\mathbb T}}(\{\vec{x}. \phi\}, M_{\mathbb T}) \to N$, are in natural bijection, by the equivalence $\mathbf{Flat}_{J_{\mathbb T}}({\cal C}_{\mathbb T}, \Set)\simeq {\mathbb T}\textrm{-mod}(\Set)$, with the natural transformations $Hom_{{\cal C}_{\mathbb T}}(\{\vec{x}. \phi\}, -)\to F_{N}$, that is, by the Yoneda lemma, with the elements of the set $F_{N}(\{\vec{x}. \phi\})=[[\vec{x}. \phi]]_{N}$, as required.   
\end{proofs}

\begin{remarks}\label{remblasce}
\begin{enumerate}[(a)]
\item If $\mathbb T$ is a universal Horn theory (in the sense of \cite{blasce}) and $\phi(\vec{x})$ is a finite conjunction of atomic formulas in a context $\vec{x}=(x^{A_{1}}, \ldots, x^{A_{n}})$ then the set $M_{\{\vec{x}. \phi\}}A$ can be identified with the set of equivalence classes of terms over $\Sigma$ of type $A_{1}, \ldots, A_{n} \to A$ modulo the equivalence relation which identifies two terms $t_{1}$ and $t_{2}$ precisely when the sequent $(\phi \vdash_{\vec{x}} t_{1}=t_{2})$ is provable in $\mathbb T$. Indeed, it is shown in \cite{blasce} (cf. p. 120 therein) that any $\mathbb T$-provably functional geometric formula $\theta(\vec{x}, \vec{y})$ between Horn formulae over $\Sigma$ is $\mathbb T$-provably equivalent to a formula of the form $\vec{y}=\vec{t}(\vec{x})$, where $\vec{t}$ is a sequence of terms of the appropriate sorts in the context $\vec{x}$. 

\item Let $\mathbb T$ be a geometric theory over a signature $\Sigma$, $\cal K$ a small category of set-based $\mathbb T$-models and $\phi(\vec{x})$ a geometric formula over $\Sigma$ presenting a $\mathbb T$-model in $\cal K$. If the geometric morphism 
\[
p:[{\cal K}, \Set] \to \Sh({\cal C}_{\mathbb T}, J_{\mathbb T})
\]
has the property that its inverse image $p^{\ast}$ is full and faithful (for instance, if $p$ is hyperconnected) then $\phi(\vec{x})$ is $\mathbb T$-irreducible and the argument in the proof of Theorem \ref{syn} applies yielding a syntactic description of the model presented by $\phi(\vec{x})$ as specified in the statement of the theorem. 
Indeed, denoted by $y$ and $y'$ the Yoneda embeddings respectively of ${\cal C}_{\mathbb T}$ into $\Sh({\cal C}_{\mathbb T}, J_{\mathbb T})$ and of ${{\cal K}}^{\textrm{op}}$ into $[{\cal K}, \Set]$, we have that $p^{\ast}(y\{\vec{x}. \phi\})\cong y'(M_{\{\vec{x}. \phi\}})$. Now, as $y'(M_{\{\vec{x}. \phi\}})$ is an irreducible object in the topos $[{\cal K}, \Set]$ and the property of an object of a topos to be irreducible is reflected by full and faithful inverse images of geometric morphisms, $y\{\vec{x}. \phi\}$ is an irreducible object of the topos $\Sh({\cal C}_{\mathbb T}, J_{\mathbb T})$, equivalently $\phi(\vec{x})$ is $\mathbb T$-irreducible, as required. 
\end{enumerate}
\end{remarks}

\subsection{Objects of homomorphisms}\label{objhom}

For any first-order signature $\Sigma$ and any Grothendieck topos $\cal E$, we have a $\cal E$-indexed category $\underline{\Sigma\textrm{-}\mathbf{str}({\cal E})}$ whose fibre at $E\in {\cal E}$ is the category $\Sigma\textrm{-}\mathbf{str}({\cal E})$ and whose `change of base' functors are the obvious pullback functors. For any sort $A$ over $\Sigma$, we have a $\cal E$-indexed forgetful functor $U_{A}:\underline{\Sigma\textrm{-}\mathbf{str}({\cal E})} \to \underline{{\cal E}}_{\cal E}$ assigning to any $\Sigma$-structure $M$ in the topos ${\cal E}\slash E$ the object $MA$. It is easy to see, by adapting the classical proof in the case ${\cal E}=\Set$ and exploiting Theorem \ref{fincol}, that the $\cal E$-indexed functors $U_{A}$ jointly create colimits of $\cal E$-indexed diagrams defined on $\cal E$-final and $\cal E$-filtered subcategories $\underline{{\cal A}}_{\cal E}$ of a small $\cal E$-indexed category. Indeed, the notion of $\Sigma$-structure only involves finite set-indexed limits, and, as we have remarked above, the colimit functor $\underline{colim}_{\cal E}:[\underline{{\cal A}}_{\cal E}, \underline{\cal E}_{{\cal E}}] \to \underline{\cal E}_{{\cal E}}$ preserves them. Moreover, the structure used in interpreting geometric formulae over $\Sigma$ is all derived from set-indexed finite limits and arbitrary colimits, which are both preserved by $\underline{colim}_{\cal E}$ (the fact that colimits commute with colimits is obvious, while the commutation with finite limits has been observed above). This implies that for any geometric theory $\mathbb T$ over $\Sigma$, the $\cal E$-indexed full subcategory $\underline{{\mathbb T}\textrm{-mod}(\cal E)}$ of $\underline{\Sigma\textrm{-}\mathbf{str}({\cal E})}$ is closed in $\underline{\Sigma\textrm{-}\mathbf{str}({\cal E})}$ under $\cal E$-indexed colimits of diagrams defined on $\cal E$-final and $\cal E$-filtered subcategories $\underline{{\cal A}}_{\cal E}$ of a small $\cal E$-indexed category.     

The following result asserts that the indexed category of models of a geometric theory in a Grothendieck topos is locally small.

\begin{theorem}\label{objectmodelhomo}
Let $\mathbb T$ be a geometric theory. Then for any $\mathbb T$-models $M$ and $N$ in a Grothendieck topos $\cal E$ there exists an object $Hom_{\underline{{\mathbb T}\textrm{-mod}(\cal E)}}^{\cal E}(M,N)$ of $\cal E$, called the `object of $\mathbb T$-model homomorphisms from $M$ to $N$', satisfying the universal property that for any object $E$ of $\cal E$ the generalized elements $E\to Hom_{\underline{{\mathbb T}\textrm{-mod}(\cal E)}}^{\cal E}(M,N)$ are in bijective correspondence, naturally in $E\in {\cal E}$, with the $\mathbb T$-model homomorphisms $!_{E}^{\ast}(M)\to !_{E}^{\ast}(N)$ in ${\mathbb T}\textrm{-mod}({\cal E}\slash E)$.
\end{theorem}

\begin{proofs}
We recall from \cite{OC} that for any small category $\cal C$ and any Grothendieck topos $\cal E$, for any two functors $F,G:{\cal C}\to {\cal E}$ there exists an object $Hom^{\cal E}(F,G)$ satisfying the property that for any object $E$ of $\cal E$ the generalized elements $E\to Hom^{\cal E}(F,G)$ are in bijective correspondence, naturally in $E\in {\cal E}$, with the arrows $!_{E}^{\ast}\circ F\to !_{E}^{\ast}\circ G$ in $[{\cal C}, {\cal E}\slash E]$, that is with the natural transformations $!_{E}^{\ast}\circ F\to !_{E}^{\ast}\circ G$. This implies that for any Grothendieck topos $\cal E$ the indexed category $\underline{[{\cal C}, {\cal E}]}_{\cal E}$ of functors on $\cal C$ with values in $\cal E$ and natural transformations between them is locally small. Clearly, it follows at once that any $\cal E$-indexed full subcategory of $\underline{[{\cal C}, {\cal E}]}_{\cal E}$, such as for example the indexed category of $J$-continuous flat functors on $\cal C$ with values in $\cal E$, is also locally small.  

Now, every geometric theory $\mathbb T$ is Morita-equivalent to the theory of flat $J_{\mathbb T}$-continuous functors ${\cal C}_{\mathbb T}$ (cf. \cite{El}); in other words, the $\cal E$-indexed category $\underline{{\mathbb T}\textrm{-mod}(\cal E)}$ is equivalent to the $\cal E$-indexed category of flat $J_{\mathbb T}$-continuous functors on ${\cal C}_{\mathbb T}$ with values in $\cal E$. Hence $\underline{{\mathbb T}\textrm{-mod}(\cal E)}$ is locally small, as required.       
\end{proofs}

\begin{remarks}\label{functoriality}
\begin{enumerate}[(a)]
\item The assignment $(M, N)\to Hom_{\underline{{\mathbb T}\textrm{-mod}(\cal E)}}^{\cal E}(M,N)$ is functorial both in $M$ and $N$; that is, any homomorphism $f: M \to M'$ in ${\mathbb T}\textrm{-mod}(\cal E)$ (resp. any homomorphism $g:N' \to N$ in ${\mathbb T}\textrm{-mod}(\cal E)$) induces an arrow 
\[
Hom_{\underline{{\mathbb T}\textrm{-mod}(\cal E)}}^{\cal E}(f,N):Hom_{\underline{{\mathbb T}\textrm{-mod}(\cal E)}}^{\cal E}(M',N)\to Hom_{\underline{{\mathbb T}\textrm{-mod}(\cal E)}}^{\cal E}(M,N)   
\]
(resp. an arrow
\[
Hom_{\underline{{\mathbb T}\textrm{-mod}(\cal E)}}^{\cal E}(M,g):Hom_{\underline{{\mathbb T}\textrm{-mod}(\cal E)}}^{\cal E}(M,N')\to Hom_{\underline{{\mathbb T}\textrm{-mod}(\cal E)}}^{\cal E}(M,N))   
\]
functorially in $f$ (resp. functorially in $g$).

\item For any Grothendieck topos $\cal E$ and $\mathbb T$-models $M$ and $N$ in $\cal E$, we have a canonical embedding 
\[
Hom_{\underline{{\mathbb T}\textrm{-mod}(\cal E)}}^{\cal F}(M,N) \to \mathbin{\mathop{\textrm{ $\prod$}}\limits_{A \textrm{ sort of } \Sigma}} NA^{MA}
\]
induced by arrows 
\[
\pi_{A}:Hom_{\underline{{\mathbb T}\textrm{-mod}(\cal E)}}^{\cal E}(M,N) \to NA^{MA}
\]
(for any sort $A$ of $\Sigma$) defined in terms of generalized elements as follows: $\pi A$ sends any arrow $E\to Hom_{\underline{{\mathbb T}\textrm{-mod}(\cal E)}}^{\cal E}(M,N)$ in $\cal E$, corresponding to a $\mathbb T$-model homomorphism $r:!_{E}^{\ast}(M) \to !_{E}^{\ast}M$ in ${\cal E}\slash E$, to the arrow $E\to NA^{MA}$ whose transpose is the `evaluation' $!_{E}^{\ast}(MA)\to !_{E}^{\ast}(NA)$ at $A$ of $r$.

\item For any $\mathbb T$-models $M$ and $N$ in a Grothendieck topos $\cal E$ and any geometric morphism $f:{\cal F}\to {\cal E}$, there is a canonical arrow 
\[
f^{\ast}(Hom_{\underline{{\mathbb T}\textrm{-mod}(\cal E)}}^{\cal E}(M,N)) \to Hom_{\underline{{\mathbb T}\textrm{-mod}(\cal F)}}^{\cal F}(f^{\ast}(M),f^{\ast}(N)).
\]
Indeed, this arrow corresponds, by the universal property of the object $Hom_{\underline{{\mathbb T}\textrm{-mod}(\cal F)}}^{\cal F}(f^{\ast}(M),f^{\ast}(N))$ to the $\mathbb T$-model homomorphism 
\[
f^{\ast}(Hom_{\underline{{\mathbb T}\textrm{-mod}(\cal E)}}^{\cal E}(M,N))\times f^{\ast}(M) \to f^{\ast}(Hom_{\underline{{\mathbb T}\textrm{-mod}(\cal E)}}^{\cal E}(M,N))\times f^{\ast}(N).
\]
in the topos ${\cal F}\slash (f^{\ast}(Hom_{\underline{{\mathbb T}\textrm{-mod}(\cal E)}}^{\cal F}(M,N)))$ whose first component, at any sort $A$, is the canonical projection and whose second component at $A$ is the arrow 
\[
f^{\ast}(Hom_{\underline{{\mathbb T}\textrm{-mod}(\cal E)}}^{\cal E}(M,N))\times f^{\ast}(MA)\cong f^{\ast}(Hom_{\underline{{\mathbb T}\textrm{-mod}(\cal E)}}^{\cal E}(M,N) \times MA) \to f^{\ast}(NA).
\]
obtained by taking the image under $f^{\ast}$ of the arrow 
\[
Hom_{\underline{{\mathbb T}\textrm{-mod}(\cal E)}}^{\cal E}(M,N) \times MA \to NA
\]
given by the transpose of the arrow $\pi_{A}: Hom_{\underline{{\mathbb T}\textrm{-mod}(\cal E)}}^{\cal E}(M,N)\to NA^{MA}$ defined above.       

\item Since the $\cal E$-indexed category $\underline{{\mathbb T}\textrm{-mod}(\cal E)}$ is locally small (cf. the proof of Theorem \ref{objectmodelhomo}), we have a $\cal E$-indexed hom functor
\[
Hom_{\underline{{\mathbb T}\textrm{-mod}(\cal E)}}^{\cal E}(-,-): \underline{{\mathbb T}\textrm{-mod}(\cal E)} \times \underline{{\mathbb T}\textrm{-mod}(\cal E)} \to \underline{{\cal E}}_{\cal E}
\]
\end{enumerate}
\end{remarks}

The following proposition provides an explicit description of the generalized elements of the objects of homomorphism of Theorem \ref{objectmodelhomo} in a particular case of interest. In the statement of the proposition, for a $\mathbb T$-model $M$ in a Grothendieck topos $\cal E$, we denote by $Hom_{\cal E}(E, M)$ the $\Sigma$-structure in $\Set$ obtained as the image of $M$ under the product-preserving functor $Hom_{\cal E}(E,-):{\cal E}\to \Set$.

\begin{proposition}\label{explicit}
Let $\mathbb T$ be a geometric theory over a signature $\Sigma$, $M$ a model of $\mathbb T$ in a Grothendieck topos $\cal E$, $c$ a set-based $\mathbb T$-model and $E$ an object of $\cal E$. Then the generalized elements $x:E\to Hom_{\underline{{\mathbb T}\textrm{-mod}(\cal E)}}^{\cal E}(\gamma_{\cal E}^{\ast}(c),M)$ correspond bijectively to the $\Sigma$-structure homomorphisms $\xi_{x}:c\to Hom_{\cal E}(E, M)$.
\end{proposition} 

\begin{proofs}
By definition of $Hom_{\underline{{\mathbb T}\textrm{-mod}_{m}(\cal E)}}^{\cal E}(\gamma_{\cal E}^{\ast}(c),M)$, a generalized element $E\to Hom_{\underline{{\mathbb T}\textrm{-mod}_{m}(\cal E)}}^{\cal E}(\gamma_{\cal E}^{\ast}(c),M)$ corresponds precisely to a $\mathbb T$-model homomorphism ${\gamma}_{{\cal E}\slash E}^{\ast}(c)\to !_{E}^{\ast}(M)$ in the topos ${\cal E}\slash E$. Concretely, such a $\mathbb T$-model homomorphism consists of a family of arros $\tau_{A}:{\gamma}_{{\cal E}\slash E}^{\ast}(cA)\to !_{E}^{\ast}(MA)$ in ${\cal E}\slash E$ indexed by the sorts $A$ over $\Sigma$ which satisfies the preservation conditions defining the notion of $\Sigma$-structure homomorphism. Now, each of the arrows $\tau_{A}:{\gamma}_{{\cal E}\slash E}^{\ast}(cA)\to !_{E}^{\ast}(MA)$ corresponds, via the adjunction between ${\gamma}_{{\cal E}\slash E}^{\ast}$ and the global section functor on the topos ${\cal E}\slash E$, to a function $\xi_{A}:cA\to Hom_{\cal E}(E, MA)$, and it is immediate to see that the above-mentioned preservation conditions translate precisely into the requirement that the arrows $\xi_{A}:cA\to Hom_{\cal E}(E, M)$ should yield a $\Sigma$-structure homomorphism $c\to Hom_{\cal E}(c, M)$. 
\end{proofs}

\begin{remark}\label{remmono}
If the model $c$ is finitely presentable also as a ${\mathbb T}_{c}$-model, where ${\mathbb T}_{c}$ is the cartesianization of $\mathbb T$ as defined in section \ref{secquotients} (notice that this is always the case if $\Sigma$ is finite and $\mathbb T$ has only a finite number of axioms, cf. Theorem 6.4 \cite{OCG}) then, if $\phi(\vec{x})$ is a formula over $\Sigma$ which presents it, the $\Sigma$-structure homomorphisms $c\to Hom_{\cal E}(E, M)$ can be identified with the elements of the interpretation of the formula $\phi(\vec{x})$ in the $\Sigma$-structure $Hom_{\cal E}(E, M)$. 
\end{remark}

\subsection{Strong finite presentability}\label{strong}

We shall show in this section that the finitely presentable models of a theory $\mathbb T$ of presheaf type enjoy a strong form of finite presentability by a geometric formula with respect to the models of $\mathbb T$ in arbitrary Grothendieck toposes.

Let $\mathbb T$ be a geometric theory over a signature $\Sigma$, $c$ a set-based model of $\mathbb T$, $\phi(\vec{x})$ a geometric formula over $\Sigma$ and $\vec{u}$ an element of $[[\vec{x}. \phi]]_{c}$. Then for any model $M$ of $\mathbb T$ in a Grothendieck topos $\cal E$ there is a canonical arrow
\[
\tau_{\phi(\vec{x}), \vec{u}}^{M}:Hom_{\underline{{\mathbb T}\textrm{-mod}(\cal E)}}^{\cal E}(\gamma_{\cal E}^{\ast}(c),M) \to [[\vec{x}. \phi]]_{M}
\] 
in $\cal E$, defined on generalized elements as follows: $\tau_{\phi(\vec{x}), \vec{u}}^{M}$ sends a generalized element 
\[
E \to Hom_{\underline{{\mathbb T}\textrm{-mod}(\cal E)}}^{\cal E}(\gamma_{\cal E}^{\ast}(c),M),
\]
corresponding under the bijection of Proposition \ref{explicit} to a $\Sigma$-structure homomorphism $f:c \to Hom_{\cal E}(E, M)$, to the image of $\vec{u}$ under $f$; notice that such element indeed belongs to $Hom_{\cal E}(E, [[\vec{x}. \phi]]_{M})$ since, as $f$ is a $\Sigma$-structure homomorphism, the image of $[[\vec{x}. \phi]]_{c}$ under $f$ is contained in $[[\vec{x}. \phi]]_{Hom_{\cal E}(E, M)}$, which is contained in $Hom_{\cal E}(E, [[\vec{x}. \phi]]_{M})$ since the functor $Hom_{\cal E}(E, -)$ is cartesian.

\begin{definition}\label{strongfp}
Let $\mathbb T$ be a geometric theory over a signature $\Sigma$, and $c$ a model of $\mathbb T$ in $\Set$. The $\mathbb T$-model $c$ is said to be \emph{strongly finitely presented} if there exists a geometric formula $\phi(\vec{x})$ over $\Sigma$ and a finite string of elements $\vec{u}$ of $[[\vec{x}. \phi]]_{c}$, called the \emph{strong generators} of $c$, such that for any $\mathbb T$-model $M$ in a Grothendieck topos $\cal E$ the arrow
\[
\tau_{\phi(\vec{x}), \vec{u}}^{M}:Hom_{\underline{{\mathbb T}\textrm{-mod}(\cal E)}}^{\cal E}(\gamma_{\cal E}^{\ast}(M_{\{\vec{x}. \phi\}}),M) \to [[\vec{x}. \phi]]_{M}
\]
is an isomorphism (equivalently, the $\Sigma$-structure homomorphisms $\xi:c\to Hom_{\cal E}(E, M)$ are in natural bijection with the generalized elements $E \to [[\vec{x}. \phi]]_{M}$ via the assignment $\xi \to \xi(\vec{u})$).
\end{definition}

\begin{remark}
If the latter condition in the definition is satisfied by all models $M$ of $\mathbb T$ inside Grothendieck toposes for $E=1_{\cal E}$ then it is true in general, by the localization principle.
\end{remark}

\begin{theorem}\label{correspondence}
Let $\mathbb T$ be a theory of presheaf type classified by the topos $[{\cal K}, \Set]$, where $\cal K$ is a small subcategory of ${\mathbb T}\textrm{-mod}(\Set)$, $M$ a $\mathbb T$-model in a Grothendieck topos $\cal E$ and $c$ a $\mathbb T$-model in $\cal K$. Then there is a natural bijective correspondence between the $\Sigma$-structure homomorphisms $c\to Hom_{\cal E}(E, M)$ and the elements of the set $Hom_{\cal E}(E, F_{M}c)$, where $F_{M}$ is the flat functor ${{\cal K}}^{\textrm{op}}\to {\cal E}$
corresponding to the model $M$ via the canonical Morita-equivalence for $\mathbb T$.  
\end{theorem}

\begin{proofs}
We can clearly suppose without loss of generality $E=1_{\cal E}$. The adjunction between $\gamma_{\cal E}^{\ast}$ and the global sections functor $\Gamma_{\cal E}:{\cal E}\to \Set$ provides a natural bijective correspondence between the elements of the set $Hom_{\cal E}(E, F_{M}c)$ and the natural transformations $\gamma_{\cal E}^{\ast}\circ yc \to F_{M}$. By the canonical Morita-equivalence 
\[
\tau_{\cal E}:\mathbf{Flat}({{\cal K}}^{\textrm{op}}, {\cal E})\simeq {\mathbb T}\textrm{-mod}(\cal E),
\] 
for $\mathbb T$, such natural transformations are in natural bijective correspondence with the $\mathbb T$-model homomorphisms $\gamma_{\cal E}^{\ast}(c)\cong \tau(\gamma_{\cal E}^{\ast}\circ yc)\to \tau(F_{M})=M$. But these homomorphisms are, by Proposition \ref{explicit}, in natural bijective correspondence with the $\Sigma$-structure homomorphisms $c\to Hom_{\cal E}(1_{\cal E}, M)$, as required. 
\end{proofs}

\begin{corollary}\label{reprtopos}
Let $\mathbb T$ be a theory of presheaf type and $\phi(\vec{x})$ a formula presenting a $\mathbb T$-model $M_{\{\vec{x}. \phi\}}$. Then $M_{\{\vec{x}. \phi\}}$ is strongly finitely presented by $\phi(\vec{x})$, i.e. for any $\mathbb T$-model $M$ in a Grothendieck topos $\cal E$ and any object $E$ of $\cal E$, the $\Sigma$-structure homomorphisms $M_{\{\vec{x}. \phi\}} \to Hom_{\cal E}(E, M)$ correspond bijectively to the elements of the set $Hom_{\cal E}(E, [[\vec{x}. \phi]]_{M})$, via the assignment sending a $\Sigma$-structure homomorphism $M_{\{\vec{x}. \phi\}} \to Hom_{\cal E}(E, M)$ to the image under it of the generators of $M_{\{\vec{x}. \phi\}}$. 
\end{corollary}

\begin{proofs}
Clearly, we can suppose without loss of generality $E=1_{\cal E}$.

By the results of \cite{OC}, the canonical Morita-equivalence 
\[
\mathbf{Flat}({{\textrm{f.p.} {\mathbb T}\textrm{-mod}(\Set)}}^{\textrm{op}}, {\cal E})\simeq {\mathbb T}\textrm{-mod}(\cal E)
\] 
for $\mathbb T$ can be described as the correspondence sending, on one hand, to any $\mathbb T$-model $M$ the flat functor $F_{M}:=Hom_{\underline{{\mathbb T}\textrm{-mod}(\cal E)}}(\gamma_{\cal E}^{\ast}(c), M)$ and conversely to any flat functor $F:{{\textrm{f.p.} {\mathbb T}\textrm{-mod}(\Set)}}^{\textrm{op}}\to {\cal E}$ the model $\tilde{F}(M_{\mathbb T})$, where $\tilde{F}:{\cal C}_{\mathbb T}\to {\cal E}$ is the extension of $F$ to the syntactic category ${\cal C}_{\mathbb T}$ (in the sense of section \ref{extsym}). By Theorem \ref{extsymthm}, for any $\mathbb T$-model $M$ in a Grothendieck topos $\cal E$, there is an isomorphism $z_{(M, \{\vec{x}. \phi\})}:F_{M}(M_{\{\vec{x}. \phi\}}) \cong \tilde{F_{M}}({\{\vec{x}. \phi\}})=[[\vec{x}. \phi]]_{M}$. Thus, applying Theorem \ref{correspondence} to the model $c=M_{\{\vec{x}. \phi\}}$, we obtain a bijective correspondence between the $\Sigma$-structure homomorphisms $M_{\{\vec{x}. \phi\}} \to Hom_{\cal E}(1_{\cal E}, M)$ and the elements of the set $Hom_{\cal E}(1_{\cal E}, [[\vec{x}. \phi]]_{M})$. It remains to show that this correspondence can be identified with the assignment sending a $\Sigma$-structure homomorphism $M_{\{\vec{x}. \phi\}} \to Hom_{\cal E}(1_{\cal E}, M)$ to the image under it of the generators of $M_{\{\vec{x}. \phi\}}$. Recall that the bijection of Theorem \ref{correspondence} can be described as follows: any $\Sigma$-structure homomorphism $f:M_{\{\vec{x}. \phi\}} \to Hom_{\cal E}(1_{\cal E}, M)$, corresponding to a $\mathbb T$-model homomorphism $\gamma_{\cal E}^{\ast}(M_{\{\vec{x}. \phi\}})\to M$, and hence, via the canonical Morita-equivalence for $\mathbb T$, to a natural transformation $\gamma_{\cal E}^{\ast} \circ y M_{\{\vec{x}. \phi\}}\cong F_{\gamma_{\cal E}^{\ast}(M_{\{\vec{x}. \phi\}} )} \to F_{M}$ (where $y$ is the Yoneda embedding ${{\textrm{f.p.} {\mathbb T}\textrm{-mod}(\Set)}}^{\textrm{op}}\hookrightarrow [{{\textrm{f.p.} {\mathbb T}\textrm{-mod}(\Set)}}, \Set]$), is sent to the global element of $F_{M}(M_{\{\vec{x}. \phi\}})$ corresponding via the Yoneda lemma to this transformation. To deduce our thesis, it thus remains to verify that the canonical isomorphism of functors $\gamma_{\cal E}^{\ast} \circ y M_{\{\vec{x}. \phi\}} = \gamma_{\cal E}^{\ast} \circ F_{M_{\{\vec{x}. \phi\}}}\cong F_{\gamma_{\cal E}^{\ast}(M_{\{\vec{x}. \phi\}})}$, when evaluated in $M_{\{\vec{x}. \phi\}}$ and composed with the isomorphism $z_{(\gamma_{\cal E}^{\ast}(M_{\{\vec{x}. \phi\}}), \{\vec{x}. \phi\})}:F_{\gamma_{\cal E}^{\ast}(M_{\{\vec{x}. \phi\}})    }(M_{\{\vec{x}. \phi\}})\cong [[\vec{x}. \phi]]_{\gamma_{\cal E}^{\ast}(M_{\{\vec{x}. \phi\}})}\cong \gamma_{\cal E}^{\ast}(  [[\vec{x}. \phi]]_{M_{\{\vec{x}. \phi\}}})$ sends the coproduct component of $\gamma_{\cal E}^{\ast}( yM_{\{\vec{x}. \phi\}} (M_{\{\vec{x}. \phi\}}) )$ corresponding to the identity on $M_{\{\vec{x}. \phi\}}$ to the coproduct component of $\gamma_{\cal E}^{\ast}(  [[\vec{x}. \phi]]_{M_{\{\vec{x}. \phi\}}})$ corresponding to the generators of $M_{\{\vec{x}. \phi\}}$. By the naturality in $\cal E$ of the Morita-equivalence for $\mathbb T$ we can clearly suppose, without loss of generality, $\cal E$ equal to $\Set$. But from the proof of Theorem \ref{syn}, it is clear that the generators of $M_{\{\vec{x}. \phi\}}$ correspond to the identity on $\{\vec{x}. \phi\})$ via the Yoneda lemma; hence they correspond to the identity on $M_{\{\vec{x}. \phi\}}$ via the above-mentioned bijection, as required.

\begin{remark}\label{indexedrepresentability}
We can express the bijective correspondence of Corollary \ref{reprtopos} by saying that for any Grothendieck topos $\cal E$ the $\cal E$-indexed functor $[[\vec{x}. \phi]]_{-}:\underline{{\mathbb T}\textrm{-mod}(\cal E)} \to \underline{{\cal E}}_{\cal E}$ assigning to any $\mathbb T$-model $M$ the interpretation of $\phi(\vec{x})$ in $M$ is represented as a $\cal E$-indexed functor by the object $\gamma_{\cal E}^{\ast}(M_{\{\vec{x}. \phi\}})$, in the sense of being naturally isomorphic to the $\cal E$-representable functor $Hom_{\underline{{\mathbb T}\textrm{-mod}(\cal E)}}^{\cal E}(\gamma_{\cal E}^{\ast}(M_{\{\vec{x}. \phi\}}), -)$.
\end{remark}

\end{proofs}

\subsection{Semantic $\cal E$-finite presentability}\label{secfp}

In this section, we introduce a semantic notion of $\cal E$-finite presentability of a model $M$ of a geometric theory $\mathbb T$ in a Grothendieck topos $\cal E$, which generalizes the classical notion in the context of finitely accessible categories, and show that all the `constant' finitely presentable models of a theory of presheaf type in a Grothendieck topos $\cal E$ are $\cal E$-finitely presentable.

For any Grothendieck topos $\cal E$, we can make any small category $\cal C$ into a $\cal E$-indexed category $\underline{{\cal C}}_{\cal E}$ defined by: ${\cal C}_{E}={\cal C}$ for all $E\in {\cal E}$ and ${\cal C}_{\alpha}=1_{{\cal C}}$ for all arrows $\alpha$ in $\cal E$. To any (set-indexed) diagram $D:{\cal C}\to {\cal E}$ corresponds a $\cal E$-indexed functor $\underline{D}_{\cal E}:\underline{{\cal C}}_{\cal E} \to \underline{\cal E}_{{\cal E}}$ defined by: for any object $E$ of $\cal E$, $\underline{D}_{\cal E}E=!_{E}^{\ast}\circ D$, where $!_{E}^{\ast}:{\cal E}\to {\cal E}\slash E$ is the pullback functor along the unique arrow $E\to 1_{\cal E}$ in $\cal E$. Since the pullback functors preserve all small limits and colimits, giving a colimiting cocone (resp. a limiting cone) on the diagram $D$ in the classical sense is equivalent to giving a $\cal E$-indexed colimiting cocone (resp. limiting cone) over the $\cal E$-indexed diagram $\underline{D}_{\cal E}$. Since the $\cal E$-indexed category $\underline{\cal E}_{{\cal E}}$ is complete, for any $\cal E$-indexed category $\underline{{\cal A}}_{\cal E}$, the $\cal E$-indexed functor category $[\underline{{\cal A}}_{\cal E}, \underline{\cal E}_{{\cal E}}]$ is also complete (since limits are computed pointwise); in particular, it has limits of diagrams defined on $\cal E$-indexed categories of the form $\underline{{\cal C}}_{\cal E}$. On the other hand, if $\underline{{\cal A}}_{\cal E}$ is a $\cal E$-final $\cal E$-filtered subcategory of a small $\cal E$-indexed category, the cocompleteness of $\underline{\cal E}_{{\cal E}}$ as a $\cal E$-indexed category ensures that there exists a well-defined $\cal E$-indexed colimit functor $\underline{colim}_{\cal E}:[\underline{{\cal A}}_{\cal E}, \underline{\cal E}_{{\cal E}}] \to \underline{\cal E}_{{\cal E}}$. It is easy to see, by mimicking the classical proof of the fact that finite limits commute with filtered colimits, that this functor preserves limits of diagrams defined on $\cal E$-indexed categories of the form $\underline{{\cal C}}_{\cal E}$ for a finite category $\cal C$.    

\begin{definition}\label{sfp}
Let $\mathbb T$ be a geometric theory and $M$ be a model of $\mathbb T$ in a Grothendieck topos $\cal E$. Then $M$ is said to be \emph{$\cal E$-finitely presentable} if the $\cal E$-indexed functor $Hom_{\underline{{\mathbb T}\textrm{-mod}(\cal E)}}^{\cal E}(M,-):\underline{{\mathbb T}\textrm{-mod}(\cal E)} \to \underline{{\cal E}}_{\cal E}$ of section \ref{objhom} preserves $\cal E$-filtered colimits (of $\cal E$-final and $\cal E$-filtered subcategories of a small $\cal E$-indexed category). 
\end{definition} 

\begin{theorem}\label{semfinpres}
Let $\mathbb T$ be a theory of presheaf type and $c$ a finitely presentable $\mathbb T$-model. Then for any Grothendieck topos $\cal E$, the $\mathbb T$-model $\gamma_{\cal E}^{\ast}(c)$ is $\cal E$-finitely presentable.
\end{theorem}

\begin{proofs}
If $c$ is presented by a geometric formula $\phi(\vec{x})$ over the signature $\Sigma$ of $\mathbb T$, by Corollary \ref{reprtopos}, the functor $Hom_{\underline{{\mathbb T}\textrm{-mod}(\cal E)}}^{\cal E}(\gamma_{\cal E}^{\ast}(c),M)$ is naturally isomorphic to the functor $[[\vec{x}. \phi]]_{-}:\underline{{\mathbb T}\textrm{-mod}(\cal E)} \to \underline{{\cal E}}_{\cal E}$. 

This latter functor preserves $\cal E$-filtered colimits (of $\cal E$-final and $\cal E$-filtered subcategories of a small $\cal E$-indexed category) since the structure used in interpreting geometric formulae over $\Sigma$ is all derived from set-indexed finite limits and arbitrary colimits, which, as remarked above, are all preserved by $\cal E$-indexed functors of the form $\underline{colim}_{\cal E}:[\underline{{\cal A}}_{\cal E}, \underline{\cal E}_{{\cal E}}] \to \underline{\cal E}_{{\cal E}}$ where $\underline{{\cal A}}_{\cal E}$ is a $\cal E$-final $\cal E$-filtered subcategory of a small $\cal E$-indexed category.  
\end{proofs}

\begin{remarks}
\begin{enumerate}[(a)]
\item The proof of the theorem shows that, more generally, for any strongly finitely presentable model $c$ of a geometric theory $\mathbb T$ in the sense of Definition \ref{strongfp} (for instance, a finite model of $\mathbb T$ if the signature of $\mathbb T$ is finite - note that such a model is finitely presented with respect to the empty theory over the signature of $\mathbb T$ by Theorem 6.4 \cite{OCG}) and any Grothendieck topos $\cal E$, the $\mathbb T$-model $\gamma_{\cal E}^{\ast}(c)$ is $\cal E$-finitely presentable. 

\item We could have more strongly required in Definition \ref{sfp} the preservation of all existing colimits of diagrams defined on $\cal E$-filtered $\cal E$-indexed categories, in the sense of Definition \ref{Efiltered}. The theorem remains true with respect to this stronger notion, but a smallness condition for the domain category, such as the requirement for it to be a $\cal E$-final and $\cal E$-filtered $\cal E$-indexed subcategory of a small $\cal E$-indexed category, is necessary to dispose of the explicit characterization of filtered colimits provided by Corollary \ref{corchar}.
\end{enumerate} 
\end{remarks}

The following proposition provides an explicit characterization of the $\cal E$-finitely presentable models of a geometric theory $\mathbb T$.

\begin{proposition}\label{propsemfp}
Let $\mathbb T$ be a geometric theory, $\cal E$ a Grothendieck topos and $c$ a $\mathbb T$-model in $\cal E$. Then $c$ is $\cal E$-finitely presentable if and only if for every $\cal E$-indexed diagram $D:\underline{{\cal A}}_{\cal E} \to \underline{{\mathbb T}\textrm{-mod}(\cal E)}$ defined on a $\cal E$-filtered $\cal E$-final subcategory $\underline{{\cal A}}_{\cal E}$ of a small $\cal E$-indexed category with $\cal E$-indexed colimiting cocone $(M, \mu)$ (we denote by $\mu_{(E, x)}:D_{E}(x)\to !_{E}^{\ast}(M)$ the colimit arrows), the following conditions are verified:
\begin{enumerate}[(a)]
\item For any object $E$ of $\cal E$ and $\mathbb T$-model homomorphism $h:!_{E}^{\ast}(c)\to !_{E}^{\ast}(M)$ in the topos ${\cal E}\slash E$ there exists an epimorphic family $\{e_{i}:E_{i}\to E \textrm{ | } i\in I\}$ in $\cal E$ and for each $i\in I$ an object $x_{i}$ of ${\cal A}_{E_{i}}$ and a $\mathbb T$-model homomorphism $\alpha_{i}:!_{E_{i}}^{\ast}(c)\to D_{E_{i}}(x_{i})$ in the topos ${\cal E}\slash E_{i}$ such that for all $i\in I$, $\mu_{(E_{i}, x_{i})}\circ \alpha_{i}=e_{i}^{\ast}(h)$; 

\item For any pairs $(x, y)$ and $(x', y')$, where $x$ and $x'$ are objects of ${\cal A}_{E}$, $f$ is an arrow $F\to E$ in $\cal E$, $y$ is a $\mathbb T$-model homomorphism $!_{F}^{\ast}(c)\to f^{\ast}(D_{E}(x))$ in ${\cal E}\slash F$ and $y'$ is a $\mathbb T$-model homomorphism $!_{F}^{\ast}(c)\to f^{\ast}(D_{E}(x'))$ in ${\cal E}\slash F$, we have $f^{\ast}(\mu_{E}(x))\circ y=f^{\ast}(\mu_{E}(x'))\circ y'$ if and only if there exists an epimorphic family $\{f_{i}:F_{i}\to F \textrm{ | } i\in I\}$ in $\cal E$ and for each $i\in I$ arrows $g_{i}:{\cal A}_{f\circ f_{i}}(x) \to z_{i}$ and $h_{i}:{\cal A}_{f\circ f_{i}}(x') \to z_{i}$ in the category ${\cal A}_{F_{i}}$ such that $D_{F_{i}}(g_{i})\circ f_{i}^{\ast}(y)=D_{F_{i}}(h_{i})\circ f_{i}^{\ast}(y')$.  
\end{enumerate}    
\end{proposition}

\begin{proofs}
The proposition follows as an an immediate consequence of Corollary \ref{corchar}, applied to the composite $\cal E$-indexed functor $Hom_{\underline{{\mathbb T}\textrm{-mod}(\cal E)}}^{\cal E}(c,-) \circ D$.  
\end{proofs}

\begin{remark}
By the construction of (small) $\cal E$-indexed colimits in the $\cal E$-indexed category $\underline{{\mathbb T}\textrm{-mod}(\cal E)}$, a $\cal E$-indexed cocone $(M, \mu)$ over a diagram $D:\underline{{\cal A}}_{\cal E} \to \underline{{\mathbb T}\textrm{-mod}(\cal E)}$ as in the statement of the proposition is colimiting if and only if for every sort $A$ over the signature of $\mathbb T$, $UA((M, \mu))$ is a colimiting cocone over the diagram $U_{A}\circ D$ in $\underline{\cal E}_{{\cal E}}$.
\end{remark}

\section{Semantic criteria for a theory to be of presheaf type}

In this section we establish our main characterization theorem providing necessary and sufficient conditions for a geometric theory to be of presheaf type. These conditions are entirely expressed in terms of the models of the theory in arbitrary Grothendieck toposes. 

We shall first prove the theorem and then proceed to reformulate its conditions in more concrete terms so for them to be directly applicable in practice.

\subsection{The characterization theorem}
 
Recall from \cite{OC} that the classifying topos of a theory of presheaf type $\mathbb T$ can be canonically represented as the topos $[\textrm{f.p.} {\mathbb T}\textrm{-mod}(\Set), \Set]$ of set-valued functors on the category ${\textrm{f.p.} {\mathbb T}\textrm{-mod}(\Set)}$ of finitely presentable $\mathbb T$-models in $\Set$. This is not the only possible representation of the classifying topos of $\mathbb T$ as a presheaf topos, but for any small category $\cal K$, $\mathbb T$ is classified by the topos $[{\cal K}, \Set]$ if and only if the Cauchy-completion of $\cal K$ is equivalent to $\textrm{f.p.} {\mathbb T}\textrm{-mod}(\Set)$. 

\begin{theorem}\label{main}

Let $\mathbb T$ be a geometric theory over a signature $\Sigma$ and let $\cal K$ be a full subcategory of ${\textrm{f.p.} {\mathbb T}\textrm{-mod}(\Set)}$. Then $\mathbb T$ is a theory of presheaf type classified by the topos $[{\cal K}, \Set]$ if and only if all of the following conditions are satisfied:
\begin{enumerate}[(i)]

\item For any $\mathbb T$-model $M$ in a Grothendieck topos $\cal E$, the functor $H_{M}:=Hom_{\underline{{\mathbb T}\textrm{-mod}(\cal E)}}^{\cal E}(\gamma_{\cal E}^{\ast}(-),M):{\cal K}^{\textrm{op}} \to {\cal E}$ is flat;  

\item The extension $\tilde{H_{M}}:{\cal C}_{\mathbb T}\to {\cal E}$ of the functor $H_{M}:{\cal K}^{\textrm{op}} \to {\cal E}$ to the syntactic category ${\cal C}_{\mathbb T}$ (in the sense of section \ref{extsym}) satisfies the property that the canonical morphism $\tilde{H_{M}}(M_{\mathbb T})\to M$ is an isomorphism;

\item Any of the following conditions (equivalent, under assumptions $(i)$ and $(ii)$) is satisfied:

\begin{enumerate}[(a)]

\item For any model $c$ in $\cal K$, $\mathbb T$-model $M$ in a Grothendieck topos $\cal E$ and geometric morphism $f:{\cal F}\to {\cal E}$, the canonical morphism 
\[
f^{\ast}(Hom_{\underline{{\mathbb T}\textrm{-mod}(\cal E)}}^{\cal E}(\gamma_{\cal E}^{\ast}(c),M))\to Hom_{\underline{{\mathbb T}\textrm{-mod}(\cal F)}}^{\cal F}(\gamma_{\cal F}^{\ast}(c),f^{\ast}(M))
\]
provided by Remark \ref{functoriality}(c) using the identification $\gamma_{\cal F}^{\ast}(c)\cong f^{\ast}(\gamma_{\cal E}^{\ast}(c))$, is an isomorphism;  

\item For any flat functor $F:{\cal K}^{\textrm{op}} \to {\cal E}$, the canonical natural transformation 
\[
\eta^{F}:F\to Hom_{\underline{{\mathbb T}\textrm{-mod}(\cal E)}}^{\cal E}(\gamma_{\cal E}^{\ast}(-),\tilde{F}(M_{\mathbb T}))\cong Hom_{\underline{\mathbf{Flat}_{J_{\mathbb T}}({\cal C}_{\mathbb T}, {\cal E})}}^{\cal E}(\gamma_{\cal E}^{\ast} \circ \tilde{y_{\cal K}(-)},\tilde{F})
\]
of Theorem \ref{genadj} is an isomorphism, where $y_{\cal K}:{\cal K}^{\textrm{op}}\hookrightarrow \mathbf{Flat}({\cal K}^{\textrm{op}}, {\Set})$ is the Yoneda embedding;

\item The functor 
\[
u^{\mathbb T}_{({\cal K}, {\cal E})}:\mathbf{Flat}({\cal K}^{\textrm{op}}, {\cal E}) \to \mathbf{Flat}_{J_{\mathbb T}}({\cal C}_{\mathbb T}, {\cal E}) \simeq {\mathbb T}\textrm{-mod}({\cal E})
\]
of section \ref{extsym} is full and faithful.
\end{enumerate}

\end{enumerate}

\end{theorem}

\begin{proofs}
Let us begin by proving that each of the three listed conditions are necessary.

By the results in \cite{OC}, if $\mathbb T$ is of presheaf type classified by the topos $[{\cal K}, \Set]$ then we have a Morita-equivalence 
\[
\tau_{\cal E}:\mathbf{Flat}({\cal K}^{\textrm{op}}, {\cal E}) \simeq {\mathbb T}\textrm{-mod}({\cal E})
\] 
which can be supposed canonical without loss of generality, i.e. which sends any finitely presentable $\mathbb T$-model $c$ in $\cal K$ to the functor $\gamma_{\cal E}^{\ast}\circ y_{\cal K}c$. It follows that for any $\mathbb T$-model $M$ in a Grothendieck topos $\cal E$, with corresponding flat functor $F_{M}$ under this Morita-equivalence, the object $Hom_{\underline{{\mathbb T}\textrm{-mod}(\cal E)}}^{\cal E}(\gamma_{\cal E}^{\ast}(c),M)$ is isomorphic to the object $Hom_{\underline{\mathbf{Flat}({\cal K}^{\textrm{op}}, {\cal E})}}^{\cal E}(\gamma_{\cal E}^{\ast} \circ y_{\cal K}c,F_{M})\cong F_{M}(c)$, naturally in $M$ and in $c$. Therefore the functor
\[
Hom_{\underline{{\mathbb T}\textrm{-mod}(\cal E)}}^{\cal E}(\gamma_{\cal E}^{\ast}(-),M):{\cal K}^{\textrm{op}}\to {\cal E}
\]
is flat, it being isomorphic to $F_{M}$. This proves that condition $(i)$ of the theorem is satisfied. 

Next, we notice that the left-to-right functor forming the canonical Morita-equivalence $\tau_{\cal E}$
for $\mathbb T$ can be described as follows: for any flat functor $F:{\cal K}^{\textrm{op}} \to {\cal E}$, the $\mathbb T$-model corresponding to it is naturally isomorphic to the model $\tilde{F}(M_{\mathbb T})$. Indeed, as the Morita-equivalence for $\mathbb T$ is canonical, it is induced by the canonical geometric morphism (in fact, equivalence) $p_{\cal K}:[{\cal K}, \Set]\to \Sh({\cal C}_{\mathbb T}, J_{\mathbb T})$, i.e. it is given by the composite of the induced equivalence $\mathbf{Flat}({\cal K}^{\textrm{op}}, {\cal E}) \simeq \mathbf{Flat}_{J_{\mathbb T}}({\cal C}_{\mathbb T}, {\cal E})$ with the canonical equivalence $\mathbf{Flat}_{J_{\mathbb T}}({\cal C}_{\mathbb T}, {\cal E})\simeq {\mathbb T}\textrm{-mod}({\cal E})$ sending any flat $J_{\mathbb T}$-continuous functor $G$ on ${\cal C}_{\mathbb T}$ to the $\mathbb T$-model $G(M_{\mathbb T})$.

The fact that $\tau_{\cal E}$ is an equivalence thus implies that the canonical arrow $\tilde{H_{M}}\to M$ is an isomorphism. This shows that condition $(ii)$ of the theorem is satisfied.

The fact that condition $(iii)(c)$ is satisfied also immediately follows from the fact that $\tau_{\cal E}$ is an equivalence.   

We have thus proved that conditions $(i)$, $(ii)$ and $(iii)(c)$ of the theorem are necessary. 

Let us now show that conditions $(i)$, $(ii)$ and $(iii)(b)$ are, all together, sufficient for the theory $\mathbb T$ to be classified by the topos $[{\cal K}, \Set]$.

First, we notice that condition $(iii)(b)$ implies condition $(iii)(a)$ under the assumption that condition $(i)$ holds. Indeed, for any $\mathbb T$-model $M$ in a Grothendieck topos $\cal E$, condition $(i)$ ensures that the functor $H_{M}:{\cal K}^{\textrm{op}}\to {\cal E}$ is flat. On the other hand, for any geometric morphism $f:{\cal F}\to {\cal E}$ condition $(iii)(b)$ yields, in view of the naturality in $\cal E$ of the operation $\tilde{(-)}$, a natural isomorphism between the flat functor $f^{\ast}\circ H_{M}$ and the functor $H_{f^{\ast}(M)}$. This ensures that the requirement of condition $(iii)(a)$ is satisfied.

Under conditions $(i)$, $(ii)$ and $(iii)(b)$, we shall construct, for any Grothen-dieck topos $\cal E$, a categorical equivalence
\[
\mathbf{Flat}({\cal K}^{\textrm{op}}, {\cal E}) \simeq {\mathbb T}\textrm{-mod}({\cal E})
\]
natural in $\cal E$.
 
We shall define two functors $G_{\cal E}:\mathbf{Flat}({\cal K}^{\textrm{op}}, {\cal E}) \to {\mathbb T}\textrm{-mod}({\cal E})$ and $H_{\cal E}:{\mathbb T}\textrm{-mod}({\cal E}) \to \mathbf{Flat}({\cal K}^{\textrm{op}}, {\cal E})$ which are natural in $\cal E$ and categorical inverses to each other (up to isomorphism).

We set $H_{\cal E}(M)$ equal to the functor 
\[
Hom_{\underline{{\mathbb T}\textrm{-mod}(\cal E)}}^{\cal E}(\gamma_{\cal E}^{\ast}(-),M):{\cal K}^{\textrm{op}} \to {\cal E}.
\]
This assignment is natural in $M$ (cf. Remark \ref{functoriality}) and hence defines a functor $H_{\cal E}:{\mathbb T}\textrm{-mod}({\cal E}) \to \mathbf{Flat}({\cal K}^{\textrm{op}}, {\cal E})$, which is natural in $\cal E$ by condition $(iii)(a)$.  

In the converse direction, for any flat functor $F:{\cal K}^{\textrm{op}}\to {\cal E}$ we set $G_{\cal E}(F)=\tilde{F}(M_{\mathbb T})$. Clearly, this assignment is natural in $F$ and defines a functor $G_{\cal E}:\mathbf{Flat}({\cal K}^{\textrm{op}}, {\cal E}) \to {\mathbb T}\textrm{-mod}({\cal E})$.

For each Grothendieck topos $\cal E$, the functors $G_{\cal E}$ and $H_{\cal E}$ are categorical inverses to each other (up to isomorphism). Indeed, condition $(iii)(b)$ ensures that $H_{\cal E} \circ G_{\cal E}$ is naturally isomorphic to the identity, while condition $(ii)$ ensures that $G_{\cal E} \circ H_{\cal E}$ is naturally isomorphic to the identity.  

Now, the functors 
\[
G_{\cal E}:\mathbf{Flat}({\cal K}^{\textrm{op}}, {\cal E}) \to {\mathbb T}\textrm{-mod}({\cal E})
\]
and 
\[
H_{\cal E}:{\mathbb T}\textrm{-mod}({\cal E}) \to \mathbf{Flat}({\cal K}^{\textrm{op}}, {\cal E})
\]
are natural in $\cal E$ and therefore induce geometric morphisms
\[
G: [{\textrm{f.p.} {\mathbb T}\textrm{-mod}(\Set)}, \Set] \to \Sh({\cal C}_{\mathbb T}, J_{\mathbb T})
\]
and
\[
H: \Sh({\cal C}_{\mathbb T}, J_{\mathbb T}) \to [{\cal K}, \Set]. 
\]
The fact that $G$ and $H$ are two halves of a categorical equivalence between $[{\cal K}, \Set]$ and $\Sh({\cal C}_{\mathbb T}, J_{\mathbb T})$ follows immediately from the fact that for every Grothendieck topos $\cal E$ there are natural isomorphisms between $H_{\cal E}\circ G_{\cal E}$ and the identity functor on $\mathbf{Flat}({\cal K}^{\textrm{op}}, {\cal E})$ and between $G_{\cal E}\circ H_{\cal E}$ and the identity functor on ${\mathbb T}\textrm{-mod}({\cal E})$, provided respectively by condition $(iii)(b)$ and condition $(ii)$.
 
Finally, let us show that, under the assumption of conditions $(i)$ and $(ii)$ of the theorem, conditions $(iii)(a)$, $(iii)(b)$, $(iii)(c)$ are all equivalent.

Conditions $(iii)(b)$ and $(iii)(c)$ are equivalent by Remark \ref{remadj}(c). The necessity of condition $(iii)(a)$ follows from Remark \ref{indexedrepresentability} and the fact that any finitely presentable model of a theory of presheaf type is finitely presented (cf. \cite{OCS}). 

Having already verified that conditions $(i)$, $(ii)$ and $(iii)(b)$ imply all together that $\mathbb T$ is classified by the topos $[{\cal K}, \Set]$, and that conditions $(iii)(a)$ and $(iii)(d)$ are both necessary conditions for $\mathbb T$ to be of presheaf type, it remains to check that, under conditions $(i)$ and $(ii)$ of the theorem, condition $(iii)(a)$ implies condition $(iii)(b)$. We shall do so by verifying that conditions $(i)$, $(ii)$ and $(iii)(a)$ imply all together that $\mathbb T$ is classified by the presheaf topos $[{\cal K}, \Set]$.

Under conditions $(i)$, $(ii)$ and $(iii)(a)$, we clearly have functors 
\[
G_{\cal E}:\mathbf{Flat}({\cal K}^{\textrm{op}}, {\cal E}) \to {\mathbb T}\textrm{-mod}({\cal E})
\]
and 
\[
H_{\cal E}:{\mathbb T}\textrm{-mod}({\cal E}) \to \mathbf{Flat}({\cal K}^{\textrm{op}}, {\cal E})
\]
which are natural in $\cal E$ and therefore induce geometric morphisms
\[
G: [{\cal K}, \Set] \to \Sh({\cal C}_{\mathbb T}, J_{\mathbb T})
\]
and
\[
H: \Sh({\cal C}_{\mathbb T}, J_{\mathbb T}) \to [{\cal K}, \Set]. 
\]

Notice that the morphism $G$ coincides with the morphism $p_{\cal K}$ canonically induced by the universal property of the classifying topos for $\mathbb T$ by the $\mathbb T$-models in $\cal K$.

Notice that for any $\mathbb T$-model $P$ in $\cal E$, $H_{\cal E}(P)=f_{P}^{\ast} \circ H^{\ast} \circ y_{{\cal K}}$, where $f_{P}:{\cal E}\to \Sh({\cal C}_{\mathbb T}, J_{\mathbb T})$ is the geometric morphism corresponding to $P$ via the universal property of the classifying topos and $y_{\cal K}:{\cal K}^{\textrm{op}}\to [{\cal K}, \Set]$ is the Yoneda embedding.   

Under our assumptions, we have to prove that
\[
H \circ G \cong 1_{[{\cal K}, \Set]}
\]
or equivalently, that $G^{\ast} \circ H^{\ast} \circ y_{\cal K} \cong y_{\cal K}$. 

Let us consider the geometric morphisms 
\[
e_{N}:\Set \to [{\cal K}, \Set]
\]
corresponding to the models $N$ of $\mathbb T$ in $\cal K$. As the $e_{N}$ are jointly surjective, it is equivalent to prove that $e_{N}^{\ast} \circ G^{\ast} \circ H^{\ast} \circ y_{\cal K} \cong e_{N}^{\ast} \circ y_{\cal K}=Hom_{\cal K}(-, N)$ naturally in $N\in {\cal K}$.

Let us denote by $M_{G}$ the $\mathbb T$-model in $[{\cal K}, \Set]$ corresponding to the geometric morphism $G$. Clearly, for any $\mathbb T$-model $N$ in $\cal K$, $e_{N}^{\ast}(M_{G})\cong N$.

We have that $G^{\ast} \circ H^{\ast} \circ y_{\cal K}=H_{[{\cal K}, \Set]}(M_{G})$. But $H_{[{\cal K}, \Set]}(M_{G})= Hom(\gamma_{[{\cal K}, \Set]}^{\ast}(-),M_{G})$, and, by condition $(iii)(a)$, 
\[
e_{N}^{\ast}(Hom(\gamma_{[{\cal K}, \Set]}^{\ast}(-),M_{G}))\cong Hom_{{\mathbb T}\textrm{-mod}(\Set)}(-, e_{N}^{\ast}(M_{G}))\cong Hom_{{\cal K}}(-, N), 
\]
(we have omitted the subscripts and superscripts in the $Hom$s above to lighten the notation), as required.

On the other hand, if condition $(ii)$ holds then $G\circ H$ is isomorphic to the identity, so we can conclude that $\mathbb T$ is classified by the topos $[{\cal K}, \Set]$. In particular, $\mathbb T$ satisfies condition $(iii)(b)$. This completes the proof of the theorem.
\end{proofs}

\begin{remarks}\label{conditionstar}
\begin{enumerate}[(a)]
\item The following condition is sufficient, together with conditions $(i)$ and $(ii)$ of the theorem (or equivalently, together with condition $(i)$ and the requirement that the models in $\cal K$ should be jointly conservative for $\mathbb T$), to ensure that $\mathbb T$ is classified by the topos $[{\cal K}, \Set]$ but necessary only if one assumes the axiom of choice:

($\ast$) There is an assignment $M \to \phi_{M}(\vec{x_{M}})$ to a $\mathbb T$-model in $\cal K$ of a formula $\phi_{M}(\vec{x_{M}})$ presenting it such that every $\mathbb T$-model homomorphism $M \to N$ between models in $\cal K$ is induced by a $\mathbb T$-provably functional formula from $\phi_{N}(\vec{x_{N}})$ to $\phi_{M}(\vec{x_{M}})$.   

This condition can be easier to verify in practice than the original condition since it allows to work with distinguished presentations (rather than with \emph{all} of them) of finitely presentable models of the theory.

The necessity of condition $(\ast)$, under the axiom of choice, was established in \cite{OCS}. In the converse direction, it suffices by Theorem \ref{main} to prove that, under conditions $(i)$ and $(ii)$ (or equivalently, under condition $(i)$ and the assertion that the models in $\cal K$ are jointly conservative for $\mathbb T$), condition ($\ast$) implies condition $(iii)(c)$. First, we remark that under either of these assumptions, the canonical geometric morphism $p_{\cal K}:[{\cal K}, \Set]\to \Sh({\cal C}_{\mathbb T}, J_{\mathbb T})$ is a surjection. From this it easily follows that for any geometric formulae $\{\vec{x}. \phi\}$ and $\{\vec{y}. \psi\}$ over $\Sigma$ respectively presenting models $c$ and $d$ in $\cal K$, there can be at most one provably functional formula $\{\vec{x}. \phi\}\to \{\vec{y}. \psi\}$ over $\Sigma$, up to $\mathbb T$-provable equivalence, inducing a given homomorphism of $\mathbb T$-models $d\to c$. Condition ($\ast$) thus implies that we have a well-defined full and faithful functor ${\cal K}\to {\cal C}_{\mathbb T}$. From this it is immediate to see, by invoking Theorem \ref{extsymthm} and Remark \ref{remf}, that condition $(iii)(c)$ is satisfied.

\item If conditions $(i)$ and $(ii)$ in the theorem are satisfied then the canonical geometric morphism
\[
p_{\cal K}:[{\cal K}, \Set] \to \Sh({\cal C}_{\mathbb T}, J_{\mathbb T})
\]
is a surjection; in other words, the models in $\cal K$ are jointly conservative for $\mathbb T$. Indeed, by condition $(i)$ the functor $H_{M}$, where $M$ is the universal model of $\mathbb T$ lying in its classifying topos, is flat. By applying Corollary \ref{preservationinequalities} to it we thus obtain that for any geometric sequent $\sigma=(\phi \vdash_{\vec{x}} \psi)$ over the signature of $\mathbb T$, $\tilde{H_{M}}(\{\vec{x}. \phi\})\leq \tilde{H_{M}}(\{\vec{x}. \psi\})$; condition $(ii)$, combined with the conservativity of $M$, thus allows to conclude that $\sigma$ is provable in $\mathbb T$, as required.  

\item Under condition $(i)$, if condition $(iii)(a)$ is satisfied and the models in $\cal K$ are jointly conservative for $\mathbb T$ (that is, every geometric sequent over $\Sigma$ which is valid in every model in $\cal K$ is provable in $\mathbb T$) then $\mathbb T$ satisfies condition $(ii)$. Indeed, the assertion that the models in $\cal K$ should be enough for the theory $\mathbb T$ is precisely equivalent to the requirement that the geometric morphism $p_{\cal K}$ should be a surjection. By the naturality in $\cal E$ of the operation $\tilde{(-)}$ and that of the functor $H_{\cal E}$ (notice that the latter follows from condition $(iii)(a)$), it suffices to verify, since $p_{\cal K}$ is surjective, that the canonical morphism $\tilde{H_{M}}(M_{\mathbb T})\to M$ is an isomorphism for $M$ equal to a model in $\cal K$. But the fact that this condition holds is obvious, as required.

\item Conditions $(ii)$ and $(iii)(c)$ in Theorem \ref{main} admit invariant formulations, which can be profitably exploited in presence of different representations for the classifying topos of $\mathbb T$. Indeed, they can both be entirely reformulated in terms of the extension of flat functors operation (in the sense of section \ref{gene}) along the canonical geometric morphism  
\[
p_{\cal K}:[{\cal K}, \Set] \to \Set[{\mathbb T}]
\]
to the classifying topos $\Set[{\mathbb T}]$ for $\mathbb T$ induced by the $\mathbb T$-models in $\cal K$. Specifically, any site of definition $({\cal C}, J)$ for $\Set[{\mathbb T}]$ gives rise to a functor
\[
G^{({\cal C}, J)}_{\cal E}: \mathbf{Flat}({\cal K}^{\textrm{op}}, {\cal E}) \to \mathbf{Flat}_{J}({\cal C}, {\cal E}).
\]
Condition $(iii)(c)$ asserts that this functor, whose explicit description is given section \ref{gene}, is full and faithful, while condition $(ii)$ asserts that, denoting by
\[
u_{{\cal E}}:\mathbf{Flat}_{J}({\cal C}, {\cal E}) \simeq {\mathbb T}\textrm{-mod}({\cal E})
\]
the equivalence canonically induced by the universal property of the classifying topos of $\mathbb T$, for any $\mathbb T$-model $M$ in a Grothendieck topos $\cal E$, the canonical morphism $u_{\cal E}(G^{({\cal C}, J)}_{\cal E}(H_{M})) \to M$ is an isomorphism.  

\item Condition $(iii)(b)$ can be split into two separate conditions: $\eta^{F}$ is pointwise monic and $\eta^{F}$ is pointwise epic. We shall refer to the first condition as to condition $(iii)(b)\textrm{-}(1)$ and to the second as to condition $(iii)(b)\textrm{-}(2)$. By Theorem \ref{genadj} and Lemma \ref{ad} below, condition $(iii)(b)\textrm{-}(1)$ is equivalent to the requirement that the functor $u^{\mathbb T}_{({\cal K}, {\cal E})}$ should be faithful.
\end{enumerate}
\end{remarks}

\subsection{Concrete reformulations}

In this section we shall give `concrete' reformulations of the conditions of Theorem \ref{main}, in full generality as well as in some particular cases in which they admit relevant simplifications.

\subsubsection{Condition $(i)$}

In this section we shall give an explicit reformulation of condition $(i)$ in Theorem \ref{main}.

\begin{theorem}\label{thcond1}
Let $\mathbb T$ be a geometric theory, $\cal K$ a small category of set-based models of $\mathbb T$, $\cal E$ a Grothendieck topos with a separating set $S$ and $M$ a $\mathbb T$-model in $\cal E$. Then the following conditions are equivalent:
\begin{enumerate}[(i)]

\item The functor 
\[
H_{M}:=Hom_{\underline{{\mathbb T}\textrm{-mod}(\cal E)}}^{\cal E}(\gamma_{\cal E}^{\ast}(-),M):{\cal K}^{\textrm{op}} \to {\cal E}
\]
is flat.

\item 

\begin{enumerate}[(a)]
\item There exists an epimorphic family $\{E_{i} \to 1_{\cal E} \textrm{ | }i\in I, E_{i}\in S\}$ and for each $i\in I$ a $\mathbb T$-model $c_{i}$ in $\cal K$ and a $\Sigma$-structure homomorphism $c_{i} \to Hom_{\cal E}(E_{i}, M)$;
\item For any $\mathbb T$-models $c$ and $d$ in $\cal K$ and $\Sigma$-structure homomorphisms $x:c\to Hom_{\cal E}(E, M)$ (where $E\in S$) and $y:d\to Hom_{\cal E}(E, M)$, there exists an epimorphic family $\{e_{i}:E_{i} \to E \textrm{ | } i\in I, E_{i}\in S\}$ and for each $i\in I$ a $\mathbb T$-model $b_{i}$ in $\cal K$, $\mathbb T$-model homomorphisms $u_{i}:c\to b_{i}$, $v_{i}:d\to b_{i}$ and a $\Sigma$-structure homomorphism $z_{i}:b_{i}\to Hom_{\cal E}(E_{i}, M)$ such that $Hom_{\cal E}(e_{i}, M)\circ x=z_{i}\circ u_{i}$ and $Hom_{\cal E}(e_{i}, M)\circ y=z_{i}\circ v_{i}$. 
\item For any two parallel homomorphisms $u, v: d \to c$ of $\mathbb T$-models in $\cal K$ and any $\Sigma$-structure homomorphism $x:c \to Hom_{\cal E}(E, M)$ in $\cal E$ (where $E\in S$) for which $x\circ u=x\circ v$, there is an epimorphic family $\{e_{i}: E_{i} \to E \textrm{ | } i\in I, E_{i}\in S \}$ in $\cal E$ and for each index $i$ a homomorphism of $\mathbb T$-models in $\cal K$ $w_{i}:c \to b_{i}$ and a $\Sigma$-structure homomorphism $y_{i}:b_{i}\to Hom_{\cal E}(E_{i}, M)$ such that $w_{i}\circ u=w_{i}\circ v$ and $y_{i} \circ w_{i}=Hom_{\cal E}(e_{i}, M) \circ x$.
\end{enumerate} 
\end{enumerate} 

\end{theorem}

\begin{proofs}
This follows immediately from Proposition \ref{explicit} in view of the characterization of flat functors as filtering functors established in chapter VII of \cite{MM} and reported in section \ref{expcolim}.
\end{proofs}

\begin{remarks}\label{remcond1}
\begin{enumerate}[(a)]
\item If for any $\mathbb T$-model $M$ in a Grothendieck topos $\cal E$ and any object $E$ of $\cal E$ the $\Sigma$-structure $Hom_{\cal E}(E, M)$ is a $\mathbb T$-model then the three conditions $(a)$, $(b)$ and $(c)$ are satisfied if they are satisfied in the case ${\cal E}=\Set$.

\item Condition $(ii)(b)$ implies condition $(ii)(c)$ if all the $\mathbb T$-model homomorphisms in any Grothendieck topos are monic.

\item Condition $(ii)(a)$ follows from condition $(ii)(a)$ of Theorem \ref{thmcond2} below if the signature of $\mathbb T$ contains at least one constant. Indeed, for any $\mathbb T$-model $M$ in a Grothendieck topos, the interpretation of such constant in $M$ will be an arrow $1\to MA$ in the topos, where $A$ is the sort of the constant; applying part $(a)$ of condition $(ii)$ of Theorem \ref{thmcond2} thus yields an epimorphic family satisfying the required property.

\item If all the $\mathbb T$-models in $\cal K$ are finitely generated as $\Sigma$-structures and all the $\Sigma$-structure homomorphisms of the form $c\to Hom_{\cal E}(E, M)$ (where $c$ is a $\mathbb T$-model in $\cal K$ and $M$ is a $\mathbb T$-model in $\cal E$) are injective if $E\ncong 0$, a sufficient condition for property $(ii)(b)$ to hold is that, for any context $\vec{x}$, the disjunction $(\top \vdash_{\vec{x}} \mathbin{\mathop{\textrm{ $\bigvee$}}\limits_{\phi(x)\in {\cal I}^{\vec{x}}_{\cal K}}}\phi(\vec{x}))$ be provable in $\mathbb T$, where ${\cal I}^{\vec{x}}_{\cal K}$ is the set of geometric formulae in the context $\vec{x}$ which strongly finitely present a $\mathbb T$-model in $\cal K$ (in the sense of Definition \ref{strongfp}). Indeed, for any two $\Sigma$-structure homomorphisms $x:c\to Hom_{\cal E}(E, M)$ and $y:d\to Hom_{\cal E}(E, M)$, where $c$ and $d$ are finitely generated $\Sigma$-structures and $E\ncong 0$, the substructure $r_{e}:e\hookrightarrow Hom_{\cal E}(E, M)$ of $Hom_{\cal E}(E, M)$ generated by the images of $c$ under $x$ and of $y$ under $y$ is finitely generated, say by elements $\xi_{1}, \ldots, \xi_{n}$; by choosing a context $\vec{x}$ of the same length as the number of generators of $e$, we obtain an epimorphic family $\{e_{i}:E_{i}\to E \textrm{ | } i\in I, E_{i}\in S, E_{i}\ncong 0\}$ with the property that for each $i\in I$ there exists a geometric formula $\phi_{i}(\vec{x})$ strongly presenting a $\mathbb T$-model $u_{i}$ in $\cal K$ such that $(\xi_{1}\circ e_{i}, \ldots, \xi_{n}\circ e_{i})$ factors through $[[\vec{x}. \phi_{i}]]_{M}$. Therefore, by the universal property of $u_{i}$ as $\mathbb T$-model strongly presented by the formula $\phi_{i}(\vec{x})$, for each $i\in I$ we have a $\Sigma$-structure homomorphism $z_{i}:u_{i}\to Hom_{\cal E}(E_{i}, M)$ which sends the generators of $u_{i}$ to the element $(\xi_{1}\circ e_{i}, \ldots, \xi_{n}\circ e_{i})$. Since $z_{i}$ is injective by our hypothesis and its image contains a set of generators for $e$, we have a factorization of $Hom_{\cal E}(e_{i}, M)\circ r_{e}$ through $r_{i}$. Therefore both $Hom_{\cal E}(e_{i}, M)\circ x$ and $Hom_{\cal E}(e_{i}, M) \circ y$ factor through $z_{i}$. This shows that condition $(ii)(b)$ is satisfied. 
\end{enumerate}
\end{remarks}
 
Let us suppose that every $\mathbb T$-model in $\cal K$ is strongly finitely presented by a formula over $\Sigma$ (in the sense of Definition \ref{strongfp}) and that every $\mathbb T$-model homomorphism between two models in $\cal K$ is induced by a $\mathbb T$-provably functional formula between formulas which present them. Then, in light of the discussion preceding Definition \ref{strongfp}, we can express the conditions for the functor $H_{M}$ to be flat (equivalently, filtering) in terms of the satisfaction by $M$ of certain geometric sequents involving these formulas. Specifically, we have the following result.

\begin{theorem}\label{flatcond}
Let $\mathbb T$ be a geometric theory over a signature $\Sigma$, $\cal K$ a small full category of the category of set-based $\mathbb T$-models and $\cal P$ a family of geometric formulae over $\Sigma$ such that every $\mathbb T$-model in $\cal K$ is strongly presented by a formula in $\cal P$ and for any two formulae $\phi(\vec{x})$ and $\psi(\vec{y})$ in $\cal P$ presenting respectively models $c$ and $d$ in $\cal K$, any $\mathbb T$-model homomorphism $d\to c$ is induced by a $\mathbb T$-provably functional formula from $\phi(\vec{x})$ to $\psi(\vec{y})$. Then $\mathbb T$ satisfies condition $(i)$ of Theorem \ref{main} with respect to the category $\cal K$ if and only if the following conditions are satisfied:
\begin{enumerate}[(i)]

\item The sequent 
\[
(\top \vdash_{[]} \mathbin{\mathop{\textrm{ $\bigvee$}}\limits_{\phi(\vec{x})\in {\cal P}}} (\exists \vec{x})\phi(\vec{x}))
\]
is valid in $M$;

\item For any formulae $\phi(\vec{x})$ and $\psi(\vec{y})$ in $\cal P$, the sequent  
\[
\phi(\vec{x}) \wedge \psi(\vec{y}) \vdash_{\vec{x}, \vec{y}} \mathbin{\mathop{\textrm{ $\bigvee$}}\limits_{\substack{ \chi(\vec{z})\in {\cal P}, \\ \{\vec{x}. \phi\} \stackrel{\theta_{1}}{\leftarrow}  \{\vec{z}. \chi\}     \stackrel{\theta_{2}}{\rightarrow} \{\vec{y}.\psi\} \textrm{ in } {\cal C}_{\mathbb T}}} (\exists \vec{z})(\theta_{1}(\vec{z},\vec{x}) \wedge \theta_{2}(\vec{z},\vec{y}))},
\]
is valid in $M$;

\item For any $\mathbb T$-provably functional formulae $\theta_{1}, \theta_{2}:\phi(\vec{x})\to \psi(\vec{y})$ between two formulae $\phi(\vec{x})$ and $\psi(\vec{y})$ in $\cal P$, the sequent
\[
\theta_{1}(\vec{x}, \vec{y}) \wedge \theta_{2}(\vec{x},\vec{y}) \vdash_{\vec{x}, \vec{y}} \mathbin{\mathop{\textrm{ $\bigvee$}}\limits_{\substack{\chi(\vec{z})\in {\cal P}, \\ \{\vec{z}. \chi\}\stackrel{\tau}{\rightarrow} \{\vec{x}. \phi\} \textrm{ in } {\cal C}_{\mathbb T}, \\ \tau \wedge \theta_{1} \dashv\vdash_{\mathbb T} \tau \wedge \theta_{2}}} (\exists \vec{z})\tau(\vec{z}, \vec{x})}
\]
is valid in $M$.

\end{enumerate}
\end{theorem}

\begin{proofs}
The fact that every model $c$ in $\cal K$ is strongly finitely presented by a formula $\phi(\vec{x})$ in $\cal P$ ensures that for any object $E$ of $\cal E$ and any $\mathbb T$-model $M$ in $\cal E$ the $\Sigma$-homomorphism $c\to Hom_{\cal E}(E, M)$ are in bijective correspondence with the generalized elements $E\to [[\vec{x}. \phi]]_{M}$ in $\cal E$. The thesis then follows from the Kripke-Joyal semantics for the topos $\cal E$, noticing that by our hypotheses for any two formulae $\phi(\vec{x})$ and $\psi(\vec{y})$ in $\cal P$ presenting respectively models $c$ and $d$ in $\cal K$, any $\mathbb T$-model homomorphism $d\to c$ is induced by a $\mathbb T$-provably functional formula from $\phi(\vec{x})$ to $\psi(\vec{y})$.
\end{proofs}

\begin{remarks}\label{remfg1}
\begin{enumerate}[(a)]
\item For any geometric theory $\mathbb T$ and category $\cal K$ satisfying the hypotheses of the theorem, all the sequents in the statement of the theorem are satisfied by every model in $\cal K$. Therefore, adding them to $\mathbb T$ yields a quotient of $\mathbb T$ satisfying condition $(i)$ of Theorem \ref{main} with respect to the category $\cal K$. 

\item Under the alternative hypothesis that every model of $\cal K$ is both strongly finitely presentable and finitely generated (with respect to the same generators), for any two formulae $\phi(\vec{x})$ and $\psi(\vec{y})$ in $\cal P$ presenting respectively models $c$ and $d$ in $\cal K$, the $\mathbb T$-model homomorphisms $c\to d$ are in bijection (via the evaluation of such homomorphisms at the generators of $c$) with the tuples of elements of $d$ which satisfy $\phi$, each of which has the form $(t_{1}(\vec{z}), \ldots, t_{n}(\vec{z}))$ (where $n$ is the length of the context $\vec{x}$) for some terms $t_{1}, \ldots, t_{n}$ in the context $\vec{z}$ over the signature $\Sigma$. Conditions $(ii)$ and $(iii)$ of Theorem \ref{flatcond} can thus be reformulated more explicitly as follows:

\begin{enumerate}[]
\item $(ii')$ For any formulae $\phi(\vec{x})$ and $\psi(\vec{y})$ in $\cal P$, where $\vec{x}=(x_{1}^{A_{1}}, \ldots, x_{n}^{A_{n}})$ and $\vec{y}=(y_{1}^{B_{1}}, \ldots, y_{m}^{B_{m}})$, the sequent  
\[
(\phi(\vec{x}) \wedge \psi(\vec{y}) \vdash_{\vec{x}, \vec{y}} \mathbin{\mathop{\textrm{ $\bigvee$}}\limits_{\substack{\chi(\vec{z})\in {\cal P}, t_{1}^{A_{1}}(\vec{z}), \ldots, t_{n}^{A_{n}}(\vec{z}) \\ s_{1}^{B_{1}}(\vec{z}), \ldots, s_{m}^{B_{m}}(\vec{z})} } (\exists \vec{z})(\chi(\vec{z}) \wedge \mathbin{\mathop{\textrm{ $\bigwedge$}}\limits_{\substack{i\in \{1, \ldots, n\}, \\ {j\in \{1, \ldots, m\}}}} (x_{i}=t_{i}(\vec{z}) \wedge y_{j}=s_{j}(\vec{z})))}}),
\]
where the disjunction is taken over all the formulae $\chi(\vec{z})$ in $\cal P$ and all the sequences of terms $t_{1}^{A_{1}}(\vec{z}), \ldots, t_{n}^{A_{n}}(\vec{z})$ and $s_{1}^{B_{1}}(\vec{z}), \ldots, s_{m}^{B_{m}}(\vec{z})$ whose output sorts are respectively $A_{1}, \ldots, A_{n}, B_{1}, \ldots, B_{m}$ and such that, denoting by $\vec{\xi}$ the set of generators of the model $M_{\{\vec{z}. \chi\}}$ (strongly) finitely presented by the formula $\chi(\vec{z})$, $(t_{1}^{A_{1}}(\vec{\xi}), \ldots, t_{n}^{A_{n}}(\vec{\xi}))\in [[\vec{x}. \phi]]_{M_{\{\vec{z}. \chi\}}}$ and $(s_{1}^{B_{1}}(\vec{\xi}), \ldots, s_{m}^{B_{m}}(\vec{\xi}))\in [[\vec{y}. \psi]]_{M_{\{\vec{z}. \chi\}}}$, is valid in $M$;

\item $(iii')$ For any formulae $\phi(\vec{x})$ and $\psi(\vec{y})$ in $\cal P$, where $\vec{x}=(x_{1}^{A_{1}}, \ldots, x_{n}^{A_{n}})$ and $\vec{y}=(y_{1}^{B_{1}}, \ldots, y_{m}^{B_{m}})$, and any terms $t_{1}^{A_{1}}(\vec{y}), s_{1}^{A_{1}}(\vec{y}), \ldots, t_{n}^{A_{n}}(\vec{y}), s_{n}^{A_{n}}(\vec{y})$ whose output sorts are respectively $A_{1}, \ldots, A_{n}$, the sequent 

\begin{equation*}
\begin{split}
(\mathbin{\mathop{\textrm{ $\bigwedge$}}\limits_{i\in \{1, \ldots, n\}} (t_{i}(\vec{y})=s_{i}(\vec{y}))} \wedge \phi(t_{1}\slash x_{1}, \ldots, t_{n}\slash x_{n}) \wedge \phi(s_{1}\slash x_{1}, \ldots, s_{n}\slash x_{n}) \wedge \psi(\vec{y}) \\
\vdash_{\vec{y}} \mathbin{\mathop{\textrm{ $\bigvee$}}\limits_{\chi(\vec{z})\in {\cal P}, u_{1}^{B_{1}}(\vec{z}), \ldots, u_{m}^{B_{m}}(\vec{z})}} ((\exists \vec{z})(\chi(\vec{z}) \wedge \mathbin{\mathop{\textrm{ $\bigwedge$}}\limits_{j\in \{1, \ldots, m\}}} (y_{j}=u_{j}(\vec{z}))),
\end{split}
\end{equation*}
where the disjunction is taken over all the formulae $\chi(\vec{z})$ in $\cal P$ and all the sequences of terms $u_{1}^{B_{1}}(\vec{z}), \ldots, u_{m}^{B_{m}}(\vec{z})$ whose output sorts are respectively $B_{1}, \ldots, B_{m}$ and such that, denoting by $\vec{\xi}$ the set of generators of the model $M_{\{\vec{z}. \chi\}}$ (strongly) finitely presented by the formula $\chi(\vec{z})$, $(u_{1}^{B_{1}}(\vec{\xi}), \ldots, u_{m}^{B_{m}}(\vec{\xi}))\in [[\vec{y}. \psi]]_{M_{\{\vec{z}. \chi\}}}$ and $t_{i}(u_{1}(\vec{\xi}), \ldots, u_{m}(\vec{\xi}))=s_{i}(u_{1}(\vec{\xi}), \ldots, u_{m}(\vec{\xi}))$ in $M_{\{\vec{z}. \chi\}}$ for all $i\in \{1, \ldots, n\}$, is valid in $M$.

\end{enumerate}
\end{enumerate}
\end{remarks}

\subsubsection{Condition $(ii)$}

In this section we shall give concrete reformulations of condition $(ii)$ of Theorem \ref{main}, under the assumption that condition $(i)$ of Theorem \ref{main} is satisfied.

First, notice that the canonical morphism $\tilde{H_{M}}(M_{\mathbb T})\to M$ considered in condition $(ii)$ is an isomorphism if and only if for every sort $A$ over $\Sigma$, the induced arrow $\tilde{H_{M}}(\{x^{A}. \top\}) \to MA$ is an isomorphism in $\cal E$. To understand this condition more concretely, we apply Proposition \ref{descrcolim} to the geometric formula $\{x^{A}. \top\}$ and the flat functor $F=H_{M}:{\cal K}^{\textrm{op}}\to {\cal E}$. Recall that the category ${\cal A}^{\cal K}_{\{x^{A}. \top \}}$ defined in that context has as objects the pairs $(c, z)$ where $c$ is a $\mathbb T$-model in $\cal K$ and $z$ is an element of the set $cA$ and as arrows $(c, z)\to (d, w)$ the $\mathbb T$-model homomorphisms $g:d\to c$ in $\cal K$ such that $gA(w)=z$, and that the equivalence relation $R_{\{x^{A}. \top\}}$ is defined by the condition that for any generalized elements $x:E \to F(c)$ and $x':E \to F(d)$, $(x, x')\in R_{\{x^{A}. \top\}}$ if and only if there exists an epimorphic family $\{e_{i}:E_{i}\to E \textrm{ | } i\in I\}$ and for each index $i\in I$ a $\mathbb T$-model $a_{i}$ in $\cal K$, a generalized element $h_{i}:E_{i}\to F(b_{i})$ and two $\mathbb T$-model homomorphisms $f_{i}:c\to b_{i}$ and $f_{i}':d\to b_{i}$ in $\cal K$ such that $f_{i}A(z)=f_{i}'A(w)$ and $<F(f_{i}), F(f'_{i})>\circ h_{i}=<x, x'>\circ e_{i}$.

This yields that the canonical arrow $\tilde{H_{M}}(\{x^{A}. \top\}) \to MA$ is an isomorphism if and only if (using the notation of Proposition \ref{descrcolim}):

\begin{enumerate}[(1)]
\item the canonical arrows $\kappa_{(c, z)}:H_{M}(c) \to MA$ for $(c, z)\in {\cal A}_{\{x^{A}. \top \}}$ are jointly epimorphic and

\item for any two objects $(c, z)$ and $(d, w)$ of the category ${\cal A}_{\{x^{A}. \top \}}$ and any generalized elements $x:E \to H_{M}(c)$ and $x':E \to H_{M}(d)$ (where $E\in S$), $\kappa_{(c, z)} \circ x=\kappa_{(d, w)} \circ x'$ if and only if $(x, x')\in R_{\{x^{A}. \top\}}$. 
\end{enumerate} 

Thanks to the identification between the generalized elements $x:E \to H_{M}(c)$ and the $\Sigma$-structure homomorphisms $f_{x}:c\to Hom_{\cal E}(E, M)$ provided by Proposition \ref{explicit}, we can rewrite conditions $(1)$ and $(2)$ more explicitly. 

To this end, we preliminarily notice that for any object $(c, z)$ of the category ${\cal A}_{\{x^{A}. \top \}}$, the canonical arrow $\kappa_{(c, z)}:H_{M}(c)=Hom_{\underline{{\mathbb T}\textrm{-mod}(\cal E)}}^{\cal E}(\gamma_{\cal E}^{\ast}(c),M) \to MA$ can be described in terms of generalized elements as the arrow sending any generalized element $x:E \to H_{M}(c)$, corresponding via the identification of Proposition \ref{explicit} to a $\Sigma$-structure homomorphism $f_{x}:c \to Hom_{\cal E}(E, M)$, to the generalized element $E \to MA$ of $MA$ given by $f_{x}A(z)$. Therefore, for any generalized elements $x:E \to H_{M}(c)$ and $x':E \to H_{M}(d)$, corresponding respectively to $\Sigma$-structure homomorphisms $f_{x}:c \to Hom_{\cal E}(E, M)$ and $f_{x'}:d \to Hom_{\cal E}(E, M)$, we have that $\kappa_{(c, z)} \circ x=\kappa_{(d, w)} \circ x'$ if and only if $f_{x}A(z)=f_{x'}A(w)$.   

Now, condition $(1)$ can be formulated by saying that for any generalized element $x:E\to MA$ (where $E\in S$) there exists an epimorphic family $\{e_{i}:E_{i}\to E \textrm{ | } i\in I, E_{i}\in S\}$, a family $\{(c_{i}, z_{i}) \textrm{ | } i\in I\}$ of objects of the category ${\cal A}_{\{x^{A}. \top \}}$ and generalized elements $\{y_{i}:E_{i}\to H_{M}(c_{i}) \textrm{ | } i\in I\}$ such that $x\circ e_{i}=\xi_{(c_{i}, z_{i})} \circ y_{i}$ for every $i\in I$. 

Under the identification of Proposition \ref{explicit}, condition $(1)$ thus rewrites as follows: for any generalized element $x:E\to MA$ (where $E\in S$) there exists an epimorphic family $\{e_{i}:E_{i}\to E \textrm{ | } i\in I, E_{i}\in S\}$, a family $\{(c_{i}, z_{i}) \textrm{ | } i\in I\}$ of objects of the category ${\cal A}_{\{x^{A}. \top \}}$ and $\Sigma$-structure homomorphisms $f_{i}:c_{i}\to Hom_{\cal E}(E_{i}, M)$ such that $f_{i}A(z_{i})=x \circ e_{i}$.     
 
Concerning condition $(2)$, we observe that for any objects $(c, z)$ and $(d, w)$ of the category ${\cal A}_{\{x^{A}. \top \}}$ and any $\Sigma$-structure homomorphisms $f_{x}:c \to Hom_{\cal E}(E, M)$ and $f_{x'}:d \to Hom_{\cal E}(E, M)$ (where $E\in S$) corresponding respectively to generalized elements $x:E \to H_{M}(c)$ and $x':E \to H_{M}(d)$ via the identification of Proposition \ref{explicit}, we have that $(x, x')\in R_{\{x^{A}. \top\}}$ if and only if there exists an epimorphic family $\{e_{i}:E_{i}\to E \textrm{ | } i\in I, E_{i}\in S\}$ and for each index $i\in I$ a $\mathbb T$-model $a_{i}$ in $\cal K$, a $\Sigma$-structure homomorphism $h_{i}:b_{i}\to Hom_{\cal E}(E_{i}, M)$ and two $\mathbb T$-model homomorphisms $f_{i}:c\to b_{i}$ and $f_{i}':d\to b_{i}$ in $\cal K$ such that $f_{i}A(z)=f_{i}'A(w)$ and $h_{i}\circ f_{i}=Hom_{\cal E}(e_{i}, M)\circ f_{x}$ and $h_{i}\circ f_{i}'=Hom_{\cal E}(e_{i}, M)\circ f_{x'}$.  
  
Summarizing, we have the following result.    
  
\begin{theorem}\label{thmcond2}
Let $\mathbb T$ be a geometric theory, $\cal K$ a small full subcategory of the category of set-based models of $\mathbb T$, $\cal E$ a Grothendieck topos with a separating set $S$ and $M$ a $\mathbb T$-model in $\cal E$. Then the following conditions are equivalent:

\begin{enumerate}[(i)]
\item The extension $\tilde{H_{M}}:{\cal C}_{\mathbb T}\to {\cal E}$ of the functor $H_{M}:{\cal K}^{\textrm{op}} \to {\cal E}$ to the syntactic category ${\cal C}_{\mathbb T}$ (in the sense of section \ref{extsym}) satisfies the property that the canonical morphism $\tilde{H_{M}}(M_{\mathbb T})\to M$ is an isomorphism;

\item For any sort $A$ over $\Sigma$, the following conditions are satisfied:
\begin{enumerate}[(a)]

\item For any generalized element $x:E\to MA$ (where $E\in S$) there exists an epimorphic family $\{e_{i}:E_{i}\to E \textrm{ | } i\in I, E_{i}\in S\}$ and for each index $i\in I$ a $\mathbb T$-model $c_{i}$ in $\cal K$, an element $z_{i}$ of $c_{i}A$ and a $\Sigma$-structure homomorphism $f_{i}:c_{i}\to Hom_{\cal E}(E_{i}, M)$ such that $(f_{i}A)(z_{i})=x \circ e_{i}$;     

\item For any two pairs $(c, z)$ and $(d, w)$ consisting of $\mathbb T$-models $c$ and $d$ in $\cal K$ and elements $z\in cA, w\in dA$, and any $\Sigma$-structure homomorphisms $f:c \to Hom_{\cal E}(E, M)$ and $f':d \to Hom_{\cal E}(E, M)$ (where $E$ is an object of $S$), we have that $fA(z)=f'A(w)$ if and only if there exists an epimorphic family $\{e_{j}:E_{j}\to E \textrm{ | } j\in J, E_{j}\in S\}$ and for each index $j\in J$ a $\mathbb T$-model $b_{j}$ in $\cal K$, a $\Sigma$-structure homomorphism $h_{j}:b_{j}\to Hom_{\cal E}(E_{j}, M)$ and two $\mathbb T$-model homomorphisms $f_{j}:c\to b_{j}$ and $f_{j}':d\to b_{j}$ in $\cal K$ such that $f_{j}A(z)=f_{j}A'(w)$, $h_{j}\circ f_{j}=Hom_{\cal E}(e_{j}, M)\circ f$ and $h_{j}\circ f_{j}'=Hom_{\cal E}(e_{j}, M)\circ f'$.  
\end{enumerate}

\end{enumerate}

\end{theorem}

\begin{remarks}\label{mon}
\begin{enumerate}[(a)]
\item If all the $\mathbb T$-model homomorphisms in any Grothendieck topos are monic and ${\mathbb T}$ satisfies condition $(ii)(b)$ of Theorem \ref{thcond1} then condition $(ii)(b)$ of Theorem \ref{thmcond2} is automatically satisfied.

\item A sufficient condition for condition $(ii)(a)$ to hold is that the disjunction $(\top \vdash_{x} \mathbin{\mathop{\textrm{ $\bigvee$}}\limits_{\phi(x)\in {\cal I}^{x}_{\cal K}}}\phi(x))$ be provable in $\mathbb T$, where ${\cal I}^{x}_{\cal K}$ is the set of geometric formulae in one variable which strongly finitely present a $\mathbb T$-model in $\cal K$ (in the sense of Definition \ref{strongfp}).

\item If for every sort $A$ over $\Sigma$ the formula $\{x^{A}.\top\}$ strongly presents a $\mathbb T$-model $F_{A}$ in $\cal K$ then condition $(i)$ of the theorem is automatically satisfied; indeed, under this hypothesis for any sort $A$ over $\Sigma$ the canonical arrow $\tilde{H_{M}}(M_{\mathbb T})A=Hom_{\underline{{\mathbb T}\textrm{-mod}(\cal E)}}^{\cal E}(\gamma_{\cal E}^{\ast}(F_{A}),M) \to MA$ is an isomorphism (cf. section \ref{strong}). 
\end{enumerate}
\end{remarks}

The following result provides an explicit formulation of condition $(ii)$ of Theorem \ref{main} holding for theories $\mathbb T$ with respect to small categories $\cal K$ such that every model of $\cal K$ is both strongly finitely presentable and finitely generated (with respect to the same generators). 

\begin{theorem}\label{cond2fg}
Let $\mathbb T$ be a geometric theory over a signature $\Sigma$ and $\cal K$ a small full subcategory of the category of set-based $\mathbb T$-models such that every model in $\cal K$ is both strongly finitely presentable and finitely generated (with respect to the same generators). Then $\mathbb T$ satisfies condition $(ii)$ of Theorem \ref{main} with respect to the category $\cal K$ if and only if for every model $M$ of $\mathbb T$ in a Grothendieck topos, the following conditions are satisfied (where $\cal P$ denotes the set of formulae over $\Sigma$ which present a model in $\cal K$): 
\begin{enumerate}[(i)]

\item For any sort $A$ over $\Sigma$, the sequent 
\[
(\top \vdash_{x_{A}} \mathbin{\mathop{\textrm{ $\bigvee$}}\limits_{\chi(\vec{z})\in {\cal P}, t^{A}(\vec{z})}} (\exists \vec{z})(\chi(\vec{z}) \wedge x=t(\vec{z}))),
\]
where the the disjunction is taken over all the formulae $\chi(\vec{z})$ in $\cal P$ and all the terms $t^{A}(\vec{z})$  whose output sort is $A$;

\item For any sort $A$ over $\Sigma$, any formulae $\phi(\vec{x})$ and $\psi(\vec{y})$ in $\cal P$, where $\vec{x}=(x_{1}^{A_{1}}, \ldots, x_{n}^{A_{n}})$ and $\vec{y}=(y_{1}^{B_{1}}, \ldots, y_{m}^{B_{m}})$, and any terms $t^{A}(\vec{x})$ and $s^{A}(\vec{y})$, the sequent 

\begin{equation*}
\begin{split}
(\phi(\vec{x}) \wedge \psi(\vec{y}) \wedge t(\vec{x})=s(\vec{y})
\vdash_{\vec{x}, \vec{y}}
\mathbin{\mathop{\textrm{ $\bigvee$}}\limits_{\substack{\chi(\vec{z})\in {\cal P}, p_{1}^{A_{1}}(\vec{z}), \ldots, p_{n}^{A_{n}}(\vec{z}) \\ q_{1}^{B_{1}}(\vec{z}), \ldots, q_{m}^{B_{m}}(\vec{z})} } (\exists \vec{z})(\chi(\vec{z})} \wedge \\ \wedge \mathbin{\mathop{\textrm{ $\bigwedge$}}\limits_{\substack{i\in \{1, \ldots, n\},\\ {j\in \{1, \ldots, m\}}}} (x_{i}=p_{i}(\vec{z}) \wedge y_{j}=q_{j}(\vec{z}))) },
\end{split}
\end{equation*}
where the disjunction is taken over all the formulae $\chi(\vec{z})$ in $\cal P$ and all the sequences of terms $p_{1}^{A_{1}}(\vec{z}), \ldots, p_{n}^{A_{n}}$ and $q_{1}^{B_{1}}(\vec{z}), \ldots, q_{m}^{B_{m}}(\vec{z})$ whose output sorts are respectively $A_{1}, \ldots, A_{n}, B_{1}, \ldots, B_{m}$ and such that, denoting by $\vec{\xi}$ the set of generators of the model $M_{\{\vec{z}. \chi\}}$ (strongly) finitely presented by the formula $\chi(\vec{z})$, $(p_{1}^{A_{1}}(\vec{\xi}), \ldots, p_{n}^{A_{n}}(\vec{\xi}))\in [[\vec{x}. \phi]]_{M_{\{\vec{z}. \chi\}}}$ and $(q_{1}^{B_{1}}(\vec{\xi}), \ldots, q_{m}^{B_{m}}(\vec{\xi}))\in [[\vec{y}. \psi]]_{M_{\{\vec{z}. \chi\}}}$ and $t(p_{1}(\vec{\xi}), \ldots, p_{n}(\vec{\xi}))=s(q_{1}(\vec{\xi}), \ldots, q_{m}(\vec{\xi}))$ in $M_{\{\vec{z}. \chi\}}$.
\end{enumerate}    
\end{theorem}

\begin{proofs}
The proof is analogous to that of Theorem \ref{thcond1} and left to the reader.
\end{proofs}

\begin{remark}\label{remfg2}
All the sequents in the statement of the theorem are satisfied by every model in $\cal K$. Therefore, adding them to any theory $\mathbb T$ satisfying the hypotheses of the theorem yields a quotient of $\mathbb T$ satisfying condition $(ii)$ of Theorem \ref{main} with respect to the category $\cal K$. 
\end{remark}

\subsubsection{Condition $(iii)$}

In this section, we shall give concrete reformulations of conditions $(iii)(a)\textrm{-}(b)\textrm{-}(c)$ of Theorem \ref{main}.

Let us begin our analysis with a proposition which shows that, under some natural assumptions which are often verified in practice, condition $(iii)(a)$ of Theorem \ref{main} is satisfied.

\begin{proposition}\label{criteriongeometricity}
Let $\mathbb T$ be a geometric theory over a signature $\Sigma$ and let $\cal K$ be a small category of ${\mathbb T}\textrm{-mod}(\Set)$. Then

\begin{enumerate}[(i)]
\item If $\mathbb T$ is a quotient of a theory $\mathbb S$ satisfying condition $(iii)(a)$ of Theorem \ref{main} with respect to a category $\cal H$ of set-based $\mathbb S$-models and $\cal K$ is a subcategory of $\cal H$ then $\mathbb T$ satisfies property $(iii)(a)$ of Theorem \ref{main} with respect to $\cal K$;

\item If every $\mathbb T$-model in $\cal K$ is strongly finitely presented (in the sense of section \ref{strong}) then $\mathbb T$ satisfies condition $(iii)(a)$ of Theorem \ref{main} with respect to the category $\cal K$; 
  
\item If for every $\mathbb T$-model $c$ in $\cal K$ the object $Hom_{\underline{{\mathbb T}\textrm{-mod}(\cal E)}}^{\cal E}(\gamma_{\cal E}^{\ast}(c),M)$ can be built from $c$ and $M$ by only using geometric constructions (i.e. constructions only involving finite limits and arbitrary small colimits) then $\mathbb T$ satisfies condition $(iii)(a)$ of Theorem \ref{main} with respect to the category $\cal K$.
\end{enumerate}

\end{proposition}

\begin{proofs}
If $\mathbb T$ is a quotient of $\mathbb S$ then for every Grothendieck topos $\cal E$, the category ${\mathbb T}\textrm{-mod}({\cal E})$ is a full subcategory of the category ${\mathbb S}\textrm{-mod}({\cal E})$. This clearly implies that for any $\mathbb T$-models $M$ and $N$ in a Grothendieck topos $\cal E$, $Hom_{\underline{{\mathbb T}\textrm{-mod}(\cal E)}}^{\cal E}(M,N)\cong Hom_{\underline{{\mathbb S}\textrm{-mod}(\cal E)}}^{\cal E}(M,N)$; in particular, for any $\mathbb T$-model $c$ in $\cal K$ and any $\mathbb T$-model $M$ in a Grothendieck topos $\cal E$, $Hom_{\underline{{\mathbb T}\textrm{-mod}(\cal E)}}^{\cal E}(\gamma_{\cal E}^{\ast}(c),M)\cong Hom_{\underline{{\mathbb S}\textrm{-mod}(\cal E)}}^{\cal E}(\gamma_{\cal E}^{\ast}(c),M)$. The fact that $\mathbb S$ satisfies condition $(iii)(a)$ of Theorem \ref{main} with respect to the category $\cal H$ thus implies that $\mathbb T$ does with respect to the category $\cal K$, as required. 

If every $\mathbb T$-model $c$ in $\cal K$ is strongly finitely presented (in the sense of section \ref{strong}) by a formula $\phi(\vec{x})$ over the signature of $\mathbb T$ then for any $\mathbb T$-model $M$ in a Grothendieck topos, $Hom_{\underline{{\mathbb T}\textrm{-mod}(\cal E)}}^{\cal E}({\gamma}_{\cal E}^{\ast}(c),M)\cong [[\vec{x}. \phi]]_{M}$ (cf. section \ref{strong}). Therefore, as the interpretation of geometric formulae is always preserved by inverse image functors of geometric morphisms, condition $(iii)(a)$ of Theorem \ref{main} is satisfied by the theory $\mathbb T$ with respect to the category $\cal K$.  

The fact that condition $(iii)$ implies condition $(iii)(a)$ of Theorem \ref{main} follows immediately from the fact that geometric constructions are preserved by inverse image functors of geometric morphisms. 
\end{proofs}

We shall now proceed to giving concrete reformulations of the conditions $(iii)(b)\textrm{-}(1)$ and $(iii)(b)\textrm{-}(2)$ of Theorem \ref{main} introduced in Remark \ref{conditionstar}(e), in order to make them more easily verifiable in practice.

First, let us explicitly describe, for any object $c$ of $\cal K$, the arrow 
\[
\eta^{F}(c):F(c)\to Hom_{\underline{{\mathbb T}\textrm{-mod}(\cal E)}}^{\cal E}(\gamma_{\cal E}^{\ast}(c),\tilde{F}(M_{\mathbb T}))
\] 
of condition $(iii)(b)$ of Theorem \ref{main} in terms of generalized elements. 

For any generalized element $x:E\to F(c)$, $\eta^{F}(c)(x)$ corresponds under the identification of Proposition \ref{explicit} to the $\Sigma$-structure homomorphism $z_{x}:c\to Hom_{\cal E}(E, \tilde{F}(M_{\mathbb T}))$ defined at each sort $A$ over $\Sigma$ as the function $cA \to Hom_{\cal E}(E, \tilde{F}(\{x^{A}. \top\}))$ sending any element $y\in cA$ to the generalized element $E\to \tilde{F}(\{x^{A}. \top\})$ obtained by composing the canonical colimit arrow $\kappa^{F}_{(c,y)}:F(c)\to \tilde{F}(\{x^{A}. \top\})$ with the generalized element $x:E\to F(c)$. 

It follows that $\eta^{F}(c)(x)$ is a monomorphism if and only if for any generalized elements $x,x':E\to F(c)$, $\kappa^{F}_{(c,y)}\circ x=\kappa^{F}_{(c,y)}\circ x'$ for every sort $A$ over $\Sigma$ and element $y\in cA$ implies $x=x'$. By Proposition \ref{descrcolim}, the condition `$\kappa^{F}_{(c,y)}\circ x=\kappa^{F}_{(c,y)}\circ x'$' is satisfied if and only if there exists an epimorphic family $\{e_{i}:E_{i}\to E \textrm{ | } i\in I\}$ in $\cal E$ and for each index $i\in I$ a $\mathbb T$-model $a_{i}$ in $\cal K$, a generalized element $h_{i}:E_{i}\to F(a_{i})$ and two $\mathbb T$-model homomorphisms $f_{i}, f_{i}':c\to a_{i}$ in $\cal K$ such that $f_{i}A(y)=f_{i}'A(y)$ and $<F(f_{i}), F(f'_{i})>\circ h_{i}=<x, x'>\circ e_{i}$.   

To obtain an explicit characterization of the condition for $\eta^{F}(c)(x)$ to be an epimorphism, we notice that an arrow $f:A\to B$ in a Grothendieck topos $\cal E$ is an epimorphism if and only if for every generalized element $x:E\to B$ of $B$ there exists an epimorphic family $\{e_{i}:E_{i}\to E \textrm{ | } i\in I\}$ in $\cal E$ and for each index $i\in I$ a generalized element $y_{i}:E_{i}\to A$ such that $f\circ y_{i}=x\circ e_{i}$. Applying this characterization to the arrow $\eta^{F}(c)(x)$ modulo the identification of Proposition \ref{explicit}, we obtain the following criterion: $\eta^{F}(c)(x)$ is an epimorphism if and only if for every object $E$ of $\cal E$ and any $\Sigma$-structure homomorphism $v:c\to Hom_{\cal E}(E, \tilde{F}(M_{\mathbb T}))$, there exists an epimorphic family $\{e_{i}:E_{i}\to E \textrm{ | } i\in I\}$ in $\cal E$ and for each index $i\in I$ a generalized element $x_{i}:E_{i}\to F(c)$ such that $z_{x_{i}}=Hom_{\cal E}(e_{i}, \tilde{F}(M_{\mathbb T}))\circ v$ for all $i$. The condition `$z_{x_{i}}=Hom_{\cal E}(e_{i}, \tilde{F}(M_{\mathbb T}))\circ v$' can be explicitly reformulated as the requirement that for every sort $A$ over $\Sigma$ and every element $y\in cA$, $\kappa^{F}_{(c,y)} \circ x_{i}=vA(y)\circ e_{i}$.  

Summarizing, we have the following 

\begin{theorem}\label{cond2}
Let $\mathbb T$ be a geometric theory over a signature $\Sigma$, $\cal K$ be a small subcategory of ${\mathbb T}\textrm{-mod}(\Set)$, $\cal E$ a Grothendieck topos with a separating set $S$ and $F:{\cal K}^{\textrm{op}} \to {\cal E}$ a flat functor. Then

\begin{enumerate}[(i)]

\item $F$ satisfies condition $(iii)(b)\textrm{-}(1)$ of Theorem \ref{main} if and only if for any $\mathbb T$-model $c$ in $\cal K$ and any generalized elements $x,x':E\to F(c)$ (where $E\in S$), if for every sort $A$ over $\Sigma$ and any element $y\in cA$ there exists an epimorphic family $\{e_{i}:E_{i}\to E \textrm{ | } i\in I, E_{i}\in S\}$ and for each index $i\in I$ a $\mathbb T$-model $a_{i}$ in $\cal K$, a generalized element $h_{i}:E_{i}\to F(a_{i})$ and two $\mathbb T$-model homomorphisms $f_{i}, f_{i}':c\to a_{i}$ in $\cal K$ such that $f_{i}A(y)=f_{i}'A(y)$ and $<F(f_{i}), F(f'_{i})>\circ h_{i}=<x, x'>\circ e_{i}$ then $x=x'$.  

\item $F$ satisfies condition $(iii)(b)\textrm{-}(2)$ of Theorem \ref{main} if and only if for any $\mathbb T$-model $c$ in $\cal K$, any object $E$ of $S$ and any $\Sigma$-structure homomorphism $v:c\to Hom_{\cal E}(E, \tilde{F}(M_{\mathbb T}))$, there exists an epimorphic family $\{e_{i}:E_{i}\to E \textrm{ | } i\in I, E_{i}\in S\}$ and for each index $i\in I$ a generalized element $x_{i}:E_{i}\to F(c)$ such that for every sort $A$ over $\Sigma$ and any element $y\in cA$, $\kappa^{F}_{(c,y)} \circ x_{i}=vA(y)\circ e_{i}$.  
\end{enumerate}

\end{theorem}\qed

Under the hypothesis that for any sort $A$ over $\Sigma$ the formula $\{x^{A}. \top\}$ presents a $\mathbb T$-model $F_{A}$, the model $\tilde{F}(M_{\mathbb T})$ is isomorphic to the model interpreting each sort $A$ with the object $F(F_{A})$, and the homomorphism $z_{x}:c\to Hom_{\cal E}(E, \tilde{F}(M_{\mathbb T}))$ corresponding to a generalized element $x:E\to F(c)$ can be described as follows: for any sort $A$ over $\Sigma$, $z_{x}A:cA \to Hom_{\cal E}(E, F(P_{A}))$ assigns to any element $y\in c_{A}$, corresponding to a $\mathbb T$-model homomorphism $s_{y}:F_{A}\to c$ via the universal property of $F_{A}$, the generalized element $F(s_{y})\circ x$. Therefore, given two generalized elements $x, x':E\to F(c)$ and an element $y\in cA$, $z_{x}A(y)=z_{x'}A(y)$ if and only if $F(s_{y})\circ x=F(s_{y})\circ x'$. 

Summarizing, we have the following

\begin{theorem}
Let $\mathbb T$ be a geometric theory over a signature $\Sigma$, such that for any sort $A$ over $\Sigma$ the formula $\{x^{A}. \top\}$ presents a $\mathbb T$-model $F_{A}$, $\cal K$ be a small subcategory of ${{\mathbb T}\textrm{-mod}(\Set)}$ containing the models $F_{A}$ and $F:{\cal K}^{\textrm{op}} \to {\cal E}$ a flat functor with values in a Grothendieck topos $\cal E$. Then $F$ satisfies condition $(iii)(b)\textrm{-}(1)$ of Theorem \ref{main} with respect to the category $\cal K$ if and only if for every $\mathbb T$-model $c$ in $\cal K$ the family of arrows $\{F(s_{y}):F(c)\to F(F_{A}) \textrm{ | } A\in \Sigma, \textrm{ }, y\in cA\}$ is jointly monic.  
\end{theorem}\qed

The representation of $\tilde{F}(M_{\mathbb T})$ as a filtered $\cal E$-indexed colimit established in section \ref{extsym} allows us to obtain a different reformulation, in terms of the $\Sigma$-structure homomorphisms $\xi_{(c, x)}$ defined in that context, of Theorem \ref{cond2}:

\begin{theorem}\label{cond3}
Let $\mathbb T$ be a geometric theory over a signature $\Sigma$, $\cal K$ a small full subcategory of ${\mathbb T}\textrm{-mod}(\Set)$ and $F:{\cal K}^{\textrm{op}} \to {\cal E}$ be a flat functor with values in a Grothendieck topos $\cal E$. Then

\begin{enumerate}[(i)]

\item $F$ satisfies condition $(iii)(b)\textrm{-}(1)$ of Theorem \ref{main} if and only if for any $\mathbb T$-model $c$ in $\cal K$ and any generalized elements $x,x':E\to F(c)$, the $\Sigma$-structure homomorphisms $\xi_{(c, x)}$ and $\xi_{(c,x')}$ are equal if and only if $x=x'$.  

\item $F$ satisfies condition $(iii)(b)\textrm{-}(2)$ of Theorem \ref{main} if and only if for any $\mathbb T$-model $c$ in $\cal K$, any object $E$ of $\cal E$ and any $\Sigma$-structure homomorphism $z:c\to Hom_{\cal E}(E, \tilde{F}(M_{\mathbb T}))$ there exists an epimorphic family $\{e_{i}:E_{i}\to E \textrm{ | } i\in I\}$ in $\cal E$ and for each index $i\in I$ a generalized element $x_{i}:E_{i}\to F(c)$ such that $Hom_{\cal E}(e_{i}, \tilde{F}(M_{\mathbb T}))\circ z=\xi_{(c, x_{i})}$ for all $i\in I$.
\end{enumerate}

\end{theorem}

\begin{proofs}
$(i)$ The equality $\xi_{(c, x)}=\xi_{(c,x')}$ holds if and only if for every sort $A$ over $\Sigma$ and any $y\in cA$, $\xi_{(c, x)}A(y)=\xi_{(c,x')}A(y)$, i.e. $\kappa_{(c,y)}\circ x=\kappa_{(c, y)}\circ x'$. Now, by Proposition \ref{descrcolim}, this latter condition is satisfied if and only if there exists an epimorphic family $\{e_{i}:E_{i}\to E \textrm{ | } i\in I\}$ in $\cal E$ and for each index $i\in I$ a $\mathbb T$-model $a_{i}$ in $\cal K$, a generalized element $h_{i}:E_{i}\to F(a_{i})$ and two $\mathbb T$-model homomorphisms $f_{i}, f_{i}':c\to a_{i}$ such that $f_{i}A(y)=f_{i}'A(y)$ and $<F(f_{i}), F(f'_{i})>\circ h_{i}=<x, x'>\circ e_{i}$, as required.     

$(ii)$ It suffices to notice that the condition that for every sort $A$ over $\Sigma$ and every element $y\in cA$, $\kappa^{F}_{(c,y)} \circ x_{i}=zA(y)\circ e_{i}$ can be reformulated as the requirement that $Hom_{\cal E}(e_{i}, \tilde{F}(M_{\mathbb T}))\circ z=\xi_{(c, x_{i})}$ (for any $i\in I$).   
\end{proofs}

Thanks to this different reformulation of condition $(iii)(b)\textrm{-}(1)$ of Theorem \ref{main}, we can prove that if all the arrows in $\cal C$ are sortwise monic then ${\mathbb T}$ always satisfies the condition. 

First, we need a lemma.

\begin{lemma}\label{injcolim}
Under the hypotheses specified above, if all the homomorphisms in the category $\cal K$ are sortwise monic then for any pair $(c, x)$ consisting of an object $c$ of $\cal K$ and of a generalized element $x:E\to F(c)$ such that $E\ncong 0_{\cal E}$, the $\Sigma$-structure homomorphism $\xi_{(c, x)}:a \to Hom_{\cal E}(E, \tilde{F}(M_{\mathbb T}))$ is sortwise injective. 
\end{lemma}

\begin{proofs}
For any sort $A$ over $\Sigma$, the function $\xi_{(c, x)}A:cA \to Hom_{\cal E}(E, MA)$ sends any element $y\in cA$ to the generalized element $\kappa_{(c, y)}\circ x$. Now, for any $y_{1}, y_{2}\in cA$, we have, by Proposition \ref{descrcolim}, that $\kappa_{(c, y_{1})}\circ x=\kappa_{(c, y_{2})}\circ x$ if and only if there exists an epimorphic family $\{e_{i}:E_{i}\to E \textrm{ | } i\in I\}$ in $\cal E$ and for each index $i\in I$ a $\mathbb T$-model $a_{i}$ in $\cal K$, a generalized element $h_{i}:E_{i}\to F(a_{i})$ and a $\mathbb T$-model homomorphism $f_{i}:c\to a_{i}$ in $\cal K$ such that $f_{i}A(y_{1})=f_{i}A(y_{2})$ and $F(f_{i}) \circ h_{i}=x \circ e_{i}$. If $E\ncong 0_{\cal E}$ then the set $I$ is non-empty, that is there exists $i\in I$; thus we have that $f_{i}A(y_{1})=f_{i}A(y_{2})$, which entails $y_{1}=y_{2}$ as $f_{i}$ is sortwise monic by our hypothesis.  
\end{proofs}

\begin{corollary}\label{corcondition3}
Let $\mathbb T$ be a geometric theory over a signature $\Sigma$, $\cal K$ a small full subcategory of ${{\mathbb T}\textrm{-mod}(\Set)}$ whose arrows are all sortwise monic homomorphisms. Then $\mathbb T$ satisfies condition $(iii)(b)\textrm{-}(1)$ of Theorem \ref{main}.
\end{corollary}  

\begin{proofs}
By Proposition \ref{res}, for any $x:E\to F(c)$ and $x':E\to F(c)$, the following `joint embedding property' holds: there exists an epimorphic family $\{e_{i}:E_{i}\to E \textrm{ | } i\in I\}$ in $\cal E$ and for each index $i\in I$ a $\mathbb T$-model $c_{i}$ in $\cal K$, homomorphisms $f_{i}:c\to c_{i}$, $g_{i}:c\to c_{i}$ in $\cal K$ and a generalized element $x_{i}:E_{i}\to F(c_{i})$ such that $<x, x'>\circ e_{i}=<F(f_{i}), F(g_{i})>\circ x_{i}$, $Hom_{\cal E}(E_{i}, M)\circ \xi_{(c, x)}=\xi_{(c_{i}, x_{i})}\circ f_{i}$ and $Hom_{\cal E}(E_{i}, M)\circ \xi_{(c, x')}=\xi_{(c_{i}, x_{i})}\circ g_{i}$ (for all $i\in I$). Clearly, we can suppose without loss of generality all the objects $E_{i}$ to be non-zero.

Now, since all the arrows in $\cal K$ are sortwise monic homomorphisms, by Lemma \ref{injcolim} for each $i\in I$ the arrow $\xi_{(c_{i}, x_{i})}$ is sortwise monic and hence $\xi_{(c, x)}=\xi_{(c, x')}$ implies $f_{i}=g_{i}$. Therefore $x\circ e_{i}=x'\circ e_{i}$ for all $i\in I$ (since $<x, x'>\circ e_{i}=<F(f_{i}), F(g_{i})>\circ x_{i}$), and hence $x=x'$, as required. 
\end{proofs}

We shall now proceed to identifying some natural sufficient conditions for a theory $\mathbb T$ to satisfy condition $(iii)(b)\textrm{-}(2)$ of Theorem \ref{main}. Before stating the relevant theorem, we need a number of preliminary results.

\begin{lemma}\label{reflectivereltopos}
Let $\Sigma$ be a signature without relation symbols and $\mathbb T$ a geometric theory over a signature $\Sigma'$ obtained from $\Sigma$ by solely adding relation symbols whose interpretation in any $\mathbb T$-model in a Grothendieck topos is the complement of the interpretation of a geometric formula over $\Sigma$ (for instance, the injectivization of a geometric theory over $\Sigma$ in the sense of Definition \ref{injectiviz}). Let $f:M\to N$ and $g:P\to N$ be homomorphisms of $\mathbb T$-models in a Grothendieck topos $\cal E$, and let $k$ be an assignment to any sort $A$ over $\Sigma$ of an arrow $kA:MA \to PA$ such that $gA \circ kA=fA$.
Then, if $g$ is sortwise monic, $k$ defines a $\mathbb T$-model homomorphism $M\to P$ such that $g\circ k=f$.
\end{lemma}

\begin{proofs}
The fact that $k$ preserves the interpretation of function symbols over $\Sigma'$ follows from the fact that $f$ and $g$ do, as $g$ is sortwise monic. It remains to prove that for any relation symbol $R\mono A_{1}, \ldots, A_{n}$ over $\Sigma'$ and any generalized element $x:E\to MA_{1}\times \cdots \times MA_{n}$ in $\cal E$, if $x$ factors through $R_{M}\mono MA_{1}\times \cdots \times MA_{n}$ then $(kA_{1}\times \cdots \times kA_{n}) \circ x$ factors through $R_{P}\mono PA_{1}\times \cdots \times PA_{n}$. Now, if $R_{P}$ is the complement of the interpretation $i_{\phi}:[[\phi(\vec{x})]]_{P}\mono PA_{1}\times \cdots \times PA_{n}$ of a geometric formula $\phi(\vec{x})=\phi(x^{A_{1}}, \ldots, x^{A_{n}})$ in the model $P$, this latter condition is equivalent to the requirement that the equalizer $e$ of $(kA_{1}\times \cdots \times kA_{n}) \circ x$ and $i_{\phi}$ be zero; but, since $g$ is a $\Sigma'$-structure homomorphism, the subobject $e$ is contained in the equalizer of the arrows $(gA_{1}\times \cdots \times gA_{n})\circ (kA_{1}\times \cdots \times kA_{n}) \circ x=(fA_{1}\times \cdots \times fA_{n}) \circ x$ and $[[\phi(\vec{x})]]_{N}\mono NA_{1}\times \cdots \times NA_{n}$, which is zero since $f$ is a $\Sigma'$-structure homomorphism and $R_{N}$ is the complement of $[[\phi(\vec{x})]]_{N}$. Therefore $e\cong 0_{\cal E}$, as required.  
\end{proofs}

\begin{lemma}\label{faithfulness}
Let $\cal E$ be a Grothendieck topos. Then the inverse image functor $\gamma_{\cal E}^{\ast}:\Set \to {\cal E}$ of the unique geometric morphism $\gamma_{\cal E}:{\cal E}\to \Set$ is faithful if and only if $\cal E$ is non-trivial (i.e., $1_{\cal E}\ncong 0_{\cal E}$).
\end{lemma}

\begin{proofs}
It is clear that if $\cal E$ is trivial then $\gamma_{\cal E}^{\ast}:\Set \to {\cal E}$ is not faithful, it being the constant functor with value $0_{\cal E}$.

In the converse direction, suppose that $\cal E$ is non-trivial. Given two functions $f,g:A\to B$ in $\Set$, the arrows $\gamma_{\cal E}^{\ast}(f), \gamma_{\cal E}^{\ast}(g):  \mathbin{\mathop{\textrm{ $\coprod$}}\limits_{a\in A}}1_{\cal E} \to \mathbin{\mathop{\textrm{ $\coprod$}}\limits_{b\in B}}1_{\cal E}$ in $\cal E$ are characterized by the following identities: for any $a\in A$ $\gamma_{\cal E}^{\ast}(f) \circ s_{a}=t_{f(a)}$ and $\gamma_{\cal E}^{\ast}(g) \circ s_{a}=t_{g(a)}$, where $s_{a}:1\to \mathbin{\mathop{\textrm{ $\coprod$}}\limits_{a\in A}}1_{\cal E}$ and $t_{b}:1\to \mathbin{\mathop{\textrm{ $\coprod$}}\limits_{b\in B}}1_{\cal E}$ are respectively the $a$-th and $b$-th coproduct arrows (for any $a\in A$ and $b\in B$); so $\gamma_{\cal E}^{\ast}(f)= \gamma_{\cal E}^{\ast}(g)$ if and only if for every $a\in A$, $t_{f(a)}=t_{g(a)}$. Now, since a subobject with domain $1_{\cal E}$ in a non-trivial topos $\cal E$ cannot be disjoint from itself, it follows that $f(a)=g(a)$ for every $a\in A$. Therefore $f=g$, as required.     
\end{proofs}

\begin{corollary}\label{faithmod}
Let $\mathbb T$ be a geometric theory over a signature $\Sigma$, $\cal E$ a non-trivial Grothendieck topos (i.e., $1_{\cal E}\ncong 0_{\cal E}$), $M$ and $N$ two $\mathbb T$-models in $\Set$ and $f$ a function sending each sort $A$ over $\Sigma$ to a map $fA:MA\to NA$ in such a way that the assignment $A \to \gamma_{\cal E}^{\ast}(fA):\gamma_{\cal E}^{\ast}(MA) \to \gamma_{\cal E}^{\ast}(NA)$ is a $\mathbb T$-model homomorphism in $\cal E$. Then the assignment $A \to fA$ is a $\mathbb T$-model homomorphism $M \to N$ in $\Set$.
\end{corollary}

\begin{proofs}
We have to prove that $f$ preserves the interpretation of function symbols over $\Sigma$ and the satisfaction of atomic relations over $\Sigma$. Since the preservation of the interpretation of function symbols over $\Sigma$ can be expressed as the commutativity of a certain square involving finite products of sets of the form $MA$ and $NA$, its satisfaction follows from that of $\gamma_{\cal E}^{\ast}(f)$ by virtue of Lemma \ref{faithfulness}. It remains to prove that $f$ preserves the satisfaction of atomic relations over $\Sigma$. Let $R\mono A_{1}, \ldots, A_{n}$ be a relation symbol over $\Sigma$. From the fact that $\gamma_{\cal E}^{\ast}(f):\gamma_{\cal E}^{\ast}(M) \to \gamma_{\cal E}^{\ast}(N)$ preserves the satisfaction by $R$ it follows that for any $n$-tuple $\vec{a}=(a_{1}, \ldots, a_{n})\in R_{M}$, the coproduct arrow $:1_{\cal E}\to \mathbin{\mathop{\textrm{ $\coprod$}}\limits_{(b_{1}, \ldots, b_{n})\in NA_{1}\times \cdots \times NA_{n}}}1_{\cal E}$ corresponding to the $n$-tuple $(fA_{1}(a_{1}), \ldots, fA_{n}(a_{n}))$ factors through the subobject $\mathbin{\mathop{\textrm{ $\coprod$}}\limits_{(b_{1}, \ldots, b_{n})\in R_{N}}}1_{\cal E} \mono \mathbin{\mathop{\textrm{ $\coprod$}}\limits_{(b_{1}, \ldots, b_{n})\in NA_{1}\times \cdots \times NA_{n}}}1_{\cal E}$. But this immediately implies, since distinct coproduct arrows are disjoint from each other and the topos $\cal E$ is non-trivial, that $f(\vec{a})\in R_{N}$, as required.  
\end{proofs}

\begin{corollary}\label{corhom}
Let $\Sigma$ be a signature without relation symbols and $\mathbb T$ a geometric theory over a signature $\Sigma'$ obtained from $\Sigma$ by solely adding relation symbols whose interpretation in any $\mathbb T$-model in a Grothendieck topos is the complement of the interpretation of a geometric formula over $\Sigma$ (for instance, the injectivization of a geometric theory over $\Sigma$ in the sense of Definition \ref{injectiviz}). Let $\cal K$ be a small full subcategory of ${\mathbb T}\textrm{-mod}(\Set)$, $a$ and $b$ $\mathbb T$-models in $\cal K$, $M$ a $\mathbb T$-model in a Grothendieck topos $\cal E$ and $f:a\to Hom_{\cal E}(E, M)$, $g:b \to Hom_{\cal E}(E, M)$ $\Sigma'$-structure homomorphisms, where $E$ is an object of $\cal E$. Let $k$ be an assignment sending any sort $A$ over $\Sigma$ to a function $kA:aA \to bA$ such that $gA \circ kA=fA$. Suppose that either
\begin{enumerate}[(i)] 
\item there are no relation symbols in $\Sigma'$ except possibly for a binary relation symbol which is $\mathbb T$-provably complemented to the equality relation on a sort of $\Sigma$ and the $\Sigma'$-structure homomorphisms $f$ and $g$ are sortwise monic; or

\item $E\ncong 0_{\cal E}$ and the $\mathbb T$-model homomorphism $\tilde{g}:\gamma_{{\cal E}\slash E}^{\ast}(b)\to M$ corresponding to $g$ under the identification of Proposition \ref{explicit} is sortwise monic.
\end{enumerate}
Then the assignment $A \to kA$ defines a $\mathbb T$-model homomorphism $a\to b$ in $\Set$.
\end{corollary}

\begin{proofs}
We have to verify that $k$ preserves the interpretation of function symbols over $\Sigma'$ as well as the satisfaction by any binary relation symbol which is $\mathbb T$-provably complemented to the equality relation on some sort over $\Sigma$.

By Lemma \ref{mondec} below, $k$ preserves the interpretation of any binary relation symbol which is $\mathbb T$-provably complemented to the equality relation on some sort over $\Sigma$ if and only if it is sortwise monic. This holds under hypothesis $(i)$ since $f$ is sortwise monic and for every sort $A$ over $\Sigma$, $gA\circ kA=fA$. On the other hand, under hypothesis $(i)$, the fact that $k$ preserves the interpretation of function symbols over $\Sigma'$ follows from the fact that $g$ and $f$ do, since $g$ is sortwise monic and $gA\circ kA=fA$ for every sort $A$ over $\Sigma$. This shows that our thesis is satisfied under hypothesis $(i)$.

Let us now suppose that hypothesis $(ii)$ holds. Consider the $\mathbb T$-model homomorphisms $\tilde{f}:\gamma_{{\cal E}\slash E}^{\ast}(a)\to M$ and $\tilde{g}:\gamma_{{\cal E}\slash E}^{\ast}(b)\to M$ corresponding to $f$ and $g$ under the identification of Proposition \ref{explicit}. Clearly, the assignment $A \to \tilde{k}A=\gamma_{{\cal E}\slash E}^{\ast}(kA)$ satisfies $\tilde{g}A\circ \tilde{k}A=\tilde{f}A$ for all sorts $A$ over $\Sigma$. It thus follows from Lemma \ref{reflectivereltopos} that $A \to \tilde{k}A$ is a $\mathbb T$-model homomorphism $\gamma_{{\cal E}\slash E}^{\ast}(a)\to \gamma_{{\cal E}\slash E}^{\ast}(b)$; but this in turn implies, by Corollary \ref{faithmod} (notice that $E\ncong 0_{\cal E}$ if and only if the topos ${\cal E}\slash E$ is non-trivial), that the assignment $A \to kA$ is a $\mathbb T$-model homomorphism $a\to b$, as required.  
\end{proofs}

\begin{theorem}\label{condition3}
Let $\Sigma$ be a signature without relation symbols and $\mathbb T$ a geometric theory over a signature $\Sigma'$ obtained from $\Sigma$ by solely adding relation symbols whose interpretation in any $\mathbb T$-model in a Grothendieck topos is the complement of the interpretation of a geometric formula over $\Sigma$ (for instance, the injectivization of a geometric theory over $\Sigma$ in the sense of Definition \ref{injectiviz}). Let $\cal K$ be a small full subcategory of ${\mathbb T}\textrm{-mod}(\Set)$ whose objects are all finitely generated $\mathbb T$-models. Suppose that 
\begin{enumerate}[(i)]
\item either all the arrows in $\cal K$ are sortwise monic homomorphisms and there are no relation symbols except possibly for a binary relation symbol which is $\mathbb T$-provably complemented to the equality relation on a given sort of $\Sigma$ or
\item all the $\mathbb T$-model homomorphisms in any Grothendieck topos are sortwise monic.   
\end{enumerate}

Then ${\mathbb T}$ satisfies condition $(iii)(b)\textrm{-}(2)$ of Theorem \ref{main} with respect to the category $\cal K$.
\end{theorem}

\begin{proofs}
By Theorem \ref{cond3}, we have to verify that for any $\mathbb T$-model $c$ in $\cal K$, any object $E$ of $\cal E$ and any $\Sigma$-structure homomorphism $z:c\to Hom_{\cal E}(E, \tilde{F}(M_{\mathbb T}))$ there exists an epimorphic family $\{e_{i}:E_{i}\to E \textrm{ | } i\in I\}$ in $\cal E$ and for each index $i\in I$ a generalized element $x_{i}:E_{i}\to F(c)$ such that $Hom_{\cal E}(e_{i}, M)\circ z=\xi_{(c, x_{i})}$ for all $i\in I$.

If $E\cong 0_{\cal E}$ then the condition is trivially satisfied; indeed, one can take $I=\emptyset$. We shall therefore suppose $E\cong 0_{\cal E}$ or, equivalently, that the topos ${\cal E}\slash E$ is non-trivial. 

From now on, we shall suppose for simplicity that $\Sigma$ is one-sorted, but all our arguments can be straightforwardly extended to the general case.

Let $\{r_{1}, \ldots, r_{n}\}$ be a set of generators for the model $c$. Consider their images $z(r_{1}), \ldots, z(r_{n}):E\to M$ under the homomorphism $z$. By Proposition \ref{res}, for any $i\in \{1, \ldots, n\}$ there exists an epimorphic family $\{e^{i}_{j}:E^{i}_{j}\to E \textrm{ | } j\in I_{i}\}$ in $\cal E$ and for each index $j\in J_{i}$ a $\mathbb T$-model $a^{i}_{j}$ in $\cal K$, a generalized element $x^{i}_{j}:E^{i}_{j}\to F(a^{i}_{j})$ and an element $y^{i}_{j}\in a^{i}_{j}$ such that $\xi_{(a^{i}_{j}, x^{i}_{j})}(y^{i}_{j})=z(r_{i})\circ e^{i}_{j}$. 

For any tuple $\vec{k}=(k_{1}, \ldots, k_{n})\in J_{1}\times \cdots \times J_{n}$, consider the iterated pullback $e_{\vec{k}}:E_{\vec{k}}=:E^{1}_{k_{1}}\times_{E} \cdots \times_{E} E^{n}_{k_{n}}\to E$. The family of arrows $\{e_{\vec{k}}:E_{\vec{k}} \to E \textrm{ | } \vec{k}\in J_{1}\times \cdots \times J_{n}\}$ is clearly epimorphic. For any $i\in \{1, \ldots, n\}$ and $k_{i}\in J_{i}$, set $x^{i}_{\vec{k}}:E_{\vec{k}} \to F(a^{i}_{k_{i}})$ equal to the composite of the generalized element $x^{i}_{k_{i}}:E^{i}_{k_{i}}\to F(a^{i}_{k_{i}})$ with the canonical pullback arrow $p_{k_{i}}^{i}:E_{\vec{k}}\to E^{i}_{k_{i}}$. For any fixed $\vec{k}\in J_{1}\times \cdots \times J_{n}$, by inductively applying the joint embedding property of Proposition \ref{res}, we can find an epimorphic family $\{u^{\vec{k}}_{l}:U^{\vec{k}}_{l}\to E_{\vec{k}} \textrm{ | } l\in L_{\vec{k}}\}$ and for each index $l\in L_{\vec{k}}$ a $\mathbb T$-model $d^{\vec{k}}_{l}$ in $\cal K$, a generalized element $x^{\vec{k}}_{l}:U^{\vec{k}}_{l} \to F(d^{\vec{k}}_{l})$, and arrows $f^{\vec{k}}_{l}:a^{i}_{k_{i}} \to d^{\vec{k}}_{l}$ in $\cal K$ such that $x^{i}_{\vec{k}} \circ u^{\vec{k}}_{l}=F(f^{\vec{k}}_{l}) \circ x^{\vec{k}}_{l}$ (for all $i\in \{1, \ldots, n\}$ and $l\in L_{\vec{k}}$). 

Let us prove that $z(r_{i})\circ e_{\vec{k}}\circ u^{\vec{k}}_{l}=\xi_{(d^{\vec{k}}_{l}, x^{\vec{k}}_{l})}(f^{\vec{k}}_{l}(y^{i}_{k_{i}}))$ (for any $i, \vec{k}$ and $l$).

We have already observed that for any $i\in \{1, \ldots, n\}$, we have $z(r_{i})\circ e^{i}_{k_{i}}=\xi_{(a^{i}_{k_{i}}, x^{i}_{k_{i}})}(y^{i}_{k_{i}})$. Composing both sides of the equation with $p_{k_{i}}^{i}$ and applying Lemma \ref{lemmaxi}(ii) yields the equality $z(r_{i})\circ e_{\vec{k}}=\xi_{(a^{i}_{k_{i}}, x^{i}_{\vec{k}})}(y^{i}_{k_{i}})$. On the other hand, by Lemma \ref{lemmaxi}(i) we have that $\xi_{(a^{i}_{k_{i}}, x^{i}_{\vec{k}})}(y^{i}_{k_{i}}) \circ u^{\vec{k}}_{l}=\xi_{(d^{\vec{k}}_{l}, x^{\vec{k}}_{l})}(f^{\vec{k}}_{l}(y^{i}_{k_{i}}))$. Our thesis thus follows by combining the former identity with the one obtained by composing both sides of the latter identity with $u^{\vec{k}}_{l}$.

Let us now consider the $\Sigma'$-structure homomorphisms
\[
Hom_{\cal E}(e_{\vec{k}}\circ u^{\vec{k}}_{l}, \tilde{F}(M_{\mathbb T})) \circ z:c \to Hom_{\cal E}(U^{\vec{k}}_{l}, \tilde{F}(M_{\mathbb T}))
\]
and 
\[
\xi_{(d^{\vec{k}}_{l}, x^{\vec{k}}_{l})}:d^{\vec{k}}_{l} \to Hom_{\cal E}(U^{\vec{k}}_{l}, \tilde{F}(M_{\mathbb T})).
\]
Since, under either assumption $(i)$ or $(ii)$, the arrows of the category $\cal K$ are sortwise monic homomorphisms, by Proposition \ref{injcolim} the homomorphism $\xi_{(d^{\vec{k}}_{l}, x^{\vec{k}}_{l})}$ is injective. Therefore, since the image of all the generators $r_{1}, \ldots, r_{n}$ of $c$ under the homomorphism $Hom_{\cal E}(e_{\vec{k}}\circ u^{\vec{k}}_{l}, \tilde{F}(M_{\mathbb T})) \circ z$ is contained in the image of the homomorphism $\xi_{(d^{\vec{k}}_{l}, x^{\vec{k}}_{l})}$, an easy induction on the structure $t(r_{1},\ldots, r_{n})$ of the elements of $c$ shows that the image of all elements of $c$ belongs to the image of $\xi_{(d^{\vec{k}}_{l}, x^{\vec{k}}_{l})}$; in other words, there exists a factorization $k^{\vec{k}}_{l}:c\to d^{\vec{k}}_{l}$ of $Hom_{\cal E}(e_{\vec{k}}\circ u^{\vec{k}}_{l}, \tilde{F}(M_{\mathbb T})) \circ z$ across $\xi_{(d^{\vec{k}}_{l}, x^{\vec{k}}_{l})}$. By Corollary \ref{corhom}, such factorization is a $\mathbb T$-model homomorphism. We have thus found an epimorphic family on $E$, namely $\{e_{\vec{k}}\circ u^{\vec{k}}_{l} \textrm{ | } \vec{k}\in J_{1}\times \cdots \times J_{n}, \textrm{ }l\in L_{\vec{k}}\}$, and for every $\vec{k}$ and $l$ a $\mathbb T$-model homomorphism $k^{\vec{k}}_{l}:c\to d^{\vec{k}}_{l}$ in $\cal K$ such that
\[
Hom_{\cal E}(e_{\vec{k}}\circ u^{\vec{k}}_{l}, \tilde{F}(M_{\mathbb T})) \circ z=\xi_{(d^{\vec{k}}_{l}, x^{\vec{k}}_{l})} \circ  k^{\vec{k}}_{l}.
\]     

The composites ${x^{\vec{k}}_{l}}':U^{\vec{k}}_{l} \to F(c)$ of the generalized elements $x^{\vec{k}}_{l}:U^{\vec{k}}_{l} \to F(d^{\vec{k}}_{l})$ with the arrows $F(k^{\vec{k}}_{l})$ thus yield generalized elements such that $Hom_{\cal E}(e_{\vec{k}}\circ u^{\vec{k}}_{l}, \tilde{F}(M_{\mathbb T})) \circ z=\xi_{(d^{\vec{k}}_{l}, x^{\vec{k}}_{l})} \circ  k^{\vec{k}}_{l}=\xi_{(c, {x^{\vec{k}}_{l}}')}$, where the last equality follows by Lemma \ref{lemmaxi}(i). This completes our proof.
\end{proofs}

\begin{remark}
The assumption that all the arrows in $\cal K$ be sortwise monic in condition $(i)$ is weaker in general than the requirement of condition $(ii)$ that all the $\mathbb T$-model homomorphisms in any Grothendieck topos be sortwise monic. Anyway, condition $(ii)$ is a necessary condition for ${\mathbb T}$ to be classified by the topos $[{\cal K}, \Set]$ if all the arrows in $\cal K$ are sortwise monic (cf. Corollary \ref{monicarr}). 
\end{remark}

\begin{corollary}
Let $\mathbb T$ be a geometric theory satisfying the hypotheses of Theorem \ref{condition3} with respect to a small full subcategory $\cal K$ of ${\mathbb T}\textrm{-mod}(\Set)$. Then $\mathbb T$ is of presheaf type classified by the topos $[{\cal K}, \Set]$ if and only if it satisfies condition $(i)$ of Theorem \ref{main} and condition $(ii)(a)$ of Theorem \ref{thmcond2}.
\end{corollary}

\begin{proofs}
Condition $(iii)$ of Theorem \ref{main} is satisfied by Theorem \ref{condition3}, while the first part of condition $(ii)$ holds by Remark \ref{mon}(a). 
\end{proofs}

\subsection{Abstract reformulation}\label{reint}

Thanks to the general theory of indexed filtered colimits developed in section \ref{indexed}, we can reformulate the conditions of Theorem \ref{main} in more abstract, though less explicit, terms as follows.

Let $\mathbb T$ be a geometric theory over a signature $\Sigma$ and $\cal K$ a small full subcategory of ${\mathbb T}\textrm{-mod}(\Set)$.

Condition $(i)$ of Theorem \ref{main} for $\mathbb T$ with respect to $\cal K$ can be reformulated as the requirement that for every model $M$ of $\mathbb T$ in a Grothendieck topos $\cal E$, the $\cal E$-indexed category of elements $\int^{f}\overline{H_{M}}$ of the functor $H_{M}$ (or equivalently its $\cal E$-final subcategory $\underline{\int H_{M}}_{\cal E}$, cf. Theorem \ref{finalfunc}) be $\cal E$-filtered.

For any sort $A$ over $\Sigma$, we have a functor $P_{A}:{\cal K}\to \Set$ of evaluation of models in $\cal K$ at the sort $A$. Using the notation of Theorem \ref{tensorpr}, we have an associated internal diagram $(P_{A})_{\cal E}$ in $[{\mathbb K}, {\cal E}]$ and hence $\cal E$-indexed functor $\underline{{\mathbb K}}_{\cal E}\to \underline{{\cal E}}_{\cal E}$. These functors, with $A$ varying among the sorts over $\Sigma$, clearly lift to a $\cal E$-indexed functor $P:\underline{{\mathbb K}}_{\cal E}\to \underline{{\mathbb T}\textrm{-mod}(\cal E)}$. 

In these terms, condition $(ii)$ of Theorem \ref{main} for the model $M$ can be reformulated as the requirement that the canonical $\cal E$-indexed cone over the $\cal E$-indexed functor $P\circ \pi^{f}_{\overline{H_{M}}}$ be colimiting (cf. the proof of Proposition \ref{res}).

Conditions $(i)$ and $(ii)$ can thus be interpreted by saying that every $\mathbb T$-model $M$ in a Grothendieck topos $\cal E$ is canonically a $\cal E$-filtered colimit of constant finitely presentable $\mathbb T$-models in $\Set$ which belong to $\cal K$.

Under the assumption that conditions $(i)$ and $(ii)$ are satisfied, condition $(iii)$ is equivalent to the requirement that for any $\mathbb T$-model $c$ in $\cal K$ and any Grothendieck topos $\cal E$, the $\mathbb T$-model $\gamma_{\cal E}^{\ast}(c)$ be $\cal E$-finitely presentable (in the sense of section \ref{secfp}) or, more weakly, that the internal $\cal E$-indexed hom functor $Hom_{\underline{{\mathbb T}\textrm{-mod}(\cal E)}}^{\cal E}(\gamma_{\cal E}^{\ast}(c), -):\underline{{\mathbb T}\textrm{-mod}(\cal E)} \to \underline{\cal E}_{\cal E}$ preserve $\cal E$-filtered colimits of diagrams of the form $\gamma_{\cal E}\circ D$, where $D$ is a diagram defined on a ($\cal E$-final subcategory of a) small internal filtered category. Indeed, the fact that this condition is necessary for $\mathbb T$ to be of presheaf type follows from Theorem \ref{semfinpres}, while the fact that it is sufficient, together with conditions $(i)$ and $(ii)$ of Theorem \ref{main}, for $\mathbb T$ to be classified by the topos $[{\cal K}, \Set]$ can be proved by showing that it implies condition $(iii)(a)$ of Theorem \ref{main}. In fact, this immediately follows from the fact that inverse image functors of geometric morphisms preserve internal (filtered) colimits in view of the fact that, by condition $(ii)$, every $\mathbb T$-model $M$ in a Grothendieck topos $\cal E$ is canonically a $\cal E$-filtered colimit of `constant' finitely presentable $\mathbb T$-models in $\cal K$. 

Overall, we can conclude that a geometric theory $\mathbb T$ is of presheaf type if and only if every $\mathbb T$-model in any Grothendieck topos $\cal E$ is a $\cal E$-filtered colimit of its canonical diagram made of `constant' finitely presentable models which are $\cal E$-finitely presentable.

\section{Faithful interpretations of theories of presheaf type}

In this section we introduce the notion of faithful interpretation of geometric theories and establish sufficent criteria for faithful interpretations of theories of presheaf type to be again of presheaf type. We shall treat conditions (i), (ii) and (iii) of Theorem \ref{main} separately.

\subsection{General results}\label{faithint}

\begin{definition}\label{faith}
A geometric morphism $a:\Set[{{\mathbb T}}]\to \Set[{{\mathbb T}'}]$ between the classifying toposes of two geometric theories  is said to be a \emph{faithful interpretation} of ${\mathbb T}'$ into ${\mathbb T}$ if the induced morphism 
\[
a_{\cal E}:{{{\mathbb T}}\textrm{-mod}({\cal E})} \to {{{\mathbb T}'}\textrm{-mod}({\cal E})}  
\]
of categories of models is faithful, reflects isomorphisms and the equality relation between objects. 
\end{definition}

For any faithful interpretation of ${\mathbb T}'$ into ${\mathbb T}$ and any Grothendieck topos $\cal E$, we have a subcategory $Im(a_{\cal E})$ of ${{{\mathbb T}'}\textrm{-mod}({\cal E})}$ whose objects (resp. arrows) are exactly the models (resp. the model homomorphisms) in the image of the functor $a_{\cal E}$. For any  functor $F:{\cal A} \to {{{\mathbb T}}\textrm{-mod}({\cal E})}$, $F$ is full and faithful if and only if its composite $a_{\cal E} \circ F$ with $a_{\cal E}$ is full and faithful as a functor ${\cal A} \to Im(a_{\cal E})$.  

\begin{remarks}
\begin{enumerate}[(a)] 

\item Any \emph{quotient} $\mathbb S$ of a geometric theory $\mathbb T$ defines a faithful interpretation of $\mathbb T$ into $\mathbb S$;

\item Any \emph{expansion} (in the sense of section \ref{exp}) ${\mathbb T}$ of a geometric theory ${\mathbb T}'$ over a signature which does not contain new sorts with respect to the signature of $\mathbb T$ defines a faithful interpretation of ${\mathbb T}'$ into ${\mathbb T}$; in particular, the \emph{injectivization} ${\mathbb T}_{m}$ (in the sense of Definition \ref{injectiviz}) of a geometric theory $\mathbb T$ defines a faithful interpretation of $\mathbb T$ into ${\mathbb T}_{m}$.
\end{enumerate}
\end{remarks}

Given a faithful interpretation of ${\mathbb T}'$ into $\mathbb T$ as above, we can reformulate the conditions of Theorem \ref{main} for the theory $\mathbb T$ with respect to a small full subcategory $\cal K$ of ${\mathbb T}\textrm{-mod}(\Set)$ in alternative ways, as follows. 

Condition $(iii)(c)$ of Theorem \ref{main} for $\mathbb T$ with respect to $\cal K$ can be formulated as follows: the composite functor 
\[ 
q_{{\cal E}}=a_{{\cal E}} \circ u_{\cal K}: \mathbf{Flat}({\cal K}^{\textrm{op}}, {\cal E}) \to  {{\mathbb T}'}\textrm{-mod}({\cal E}) 
\]  
is full and faithful into the image $Im(a_{\cal E})$ of $a_{\cal E}$, where $u_{\cal K}$ is the functor $\mathbf{Flat}({\cal K}^{\textrm{op}}, {\cal E}) \to {{\mathbb T}}\textrm{-mod}({\cal E})$ induced by the canonical geometric morphism $p_{\cal K}:[{\cal K}, \Set]\to \Set[{\mathbb T}]$. 

Condition $(ii)$ of Theorem \ref{main} for $\mathbb T$ with respect to $\cal K$ can be reformulated by saying that for any $\mathbb T$-model $M$ in a Grothendieck topos $\cal E$, the canonical morphism
\[
q_{{\cal E}}(H_{M}) \to a_{\cal E}(M)
\]
is an isomorphism.  

Suppose moreover that the functor
\[
a_{\Set}:{{{\mathbb T}}\textrm{-mod}({\Set})} \to {{{\mathbb T}'}\textrm{-mod}({\Set})}  
\] 
restricts to a functor
\[
f:{\cal K}\to {\cal H}, 
\]
where $\cal H$ is a small full subcategory of ${{{\mathbb T}'}\textrm{-mod}({\Set})}$.

We have a commutative diagram
\[  
\xymatrix {
[{\cal K}, \Set] \ar[d]^{p_{\cal K}} \ar[r]^{[f, \Set]} & [{\cal H}, \Set] \ar[d]^{p_{{\cal H}}}\\
\Set[{\mathbb T}] \ar[r]_{a} & \Set[{\mathbb T}'],} 
\] 
where $p_{\cal K}$ (resp. $p_{{\cal H}}$) is the canonical geometric morphism determined by the universal property of the classifying topos of $\mathbb T$ (resp. of ${\mathbb T}'$). 
 
Also, for any ${\mathbb T}'$-model $M$ in a Grothendieck topos $\cal E$, the functors $a_{\cal E}$ induce a $\cal E$-indexed functor
\[
{\int a}:{\int H^{{\mathbb T}}_{M}} \to {\int H^{{\mathbb T}'}_{a_{{\cal E}}(M)}}
\]
between the $\cal E$-indexed categories of elements ${\int H^{{\mathbb T}}_{M}}$ and ${\int H^{{\mathbb T}'}_{a_{{\cal E}}(M)}}$ of the functors
\[
H^{{\mathbb T}}_{M}:=Hom_{\underline{{\mathbb T}\textrm{-mod}(\cal E)}}^{\cal E}(\gamma_{\cal E}^{\ast}(-),M):{\cal K}^{\textrm{op}} \to {\cal E}
\]
and
\[
H^{{\mathbb T}'}_{a_{{\cal E}}(M)}:=Hom_{\underline{{\mathbb T}'\textrm{-mod}(\cal E)}}^{\cal E}(\gamma_{\cal E}^{\ast}(-),a_{{\cal E}}(M)):{\cal H}^{\textrm{op}} \to {\cal E}.
\]

Since $a$ is by our hypothesis a faithful interpretation of theories, ${\int a}$ is faithful and reflects equalities and isomorphisms. Suppose now that ${\int a}$ is moreover $\cal E$-full (in the sense of Definition \ref{fulldefinition}). Since for any functor $F$ with values in a Grothendieck topos, $F$ is flat if and only if its $\cal E$-indexed category of elements is $\cal E$-filtered, we conclude by Proposition \ref{finalfiltered} that if ${\int a}$ is $\cal E$-final then the theory $\mathbb {\mathbb T}$ satisfies condition (i) of Theorem \ref{main} with respect to the category ${\cal K}$ if ${\mathbb T}'$ does with respect to the category $\cal H$. Moreover, Proposition \ref{finalcolimit} and Theorem \ref{fincol} ensure that if ${\mathbb T}'$ satisfies condition (ii) of Theorem \ref{main} with respect to the category $\cal H$ then ${\mathbb T}$ does with respect to the category $\cal K$.  

Summarizing, we have the following result.

\begin{theorem}\label{finalitythm}
Let $a:\Set[{{\mathbb T}}]\to \Set[{{\mathbb T}'}]$ be a faithful interpretation of geometric theories and let $\cal K$ and $\cal H$ be small subcategories respectively of ${{{\mathbb T}}\textrm{-mod}({\Set})}$ and of ${{{\mathbb T}'}\textrm{-mod}({\Set})}$ such that the functor 
\[
a_{\Set}:{{{\mathbb T}}\textrm{-mod}({\Set})} \to {{{\mathbb T}'}\textrm{-mod}({\Set})}  
\]
restricts to a functor ${\cal K}\to {\cal H}$ and the $\cal E$-indexed functor
\[
{\int a}:{\int H^{{\mathbb T}}_{M}} \to {\int H^{{\mathbb T}'}_{a_{{\cal E}}(M)}}  
\]
is $\cal E$-full and satisfies the equivalent conditions of Proposition \ref{finality}. 

Then the theory ${\mathbb T}$ satisfies condition $(i)$ (resp. condition $(ii)$) of Theorem \ref{main} with respect to the category ${\cal K}$ if ${\mathbb T}'$ satisfies condition $(i)$ (resp. condition $(ii)$) of Theorem \ref{main} with respect to the category ${\cal H}$.
\end{theorem}\qed

\begin{remark}\label{appl}
If ${\mathbb T}$ is the injectivization of ${\mathbb T}'$ (in the sense of Definition \ref{injectiviz}) or if ${\mathbb T}$ is a quotient of ${\mathbb T}'$ then the functor $\int a$ is ${\cal E}$-full. Indeed, in the first case this follows from the fact that if the composite of two arrows is monic then the first one is monic, while in the second it follows from the fact that the category of models of ${\mathbb T}$ is a full subcategory of the category of models of ${\mathbb T}'$.
\end{remark}

The following proposition provides an explicit rephrasing of the finality condition for the functor $\int a$. 

\begin{proposition}\label{finality}
Let $a:\Set[{{\mathbb T}}]\to \Set[{{\mathbb T}'}]$ be a faithful interpretation of geometric theories let $M$ be a ${\mathbb T}$-model in a Grothendieck topos $\cal E$ such that the functor $H^{{\mathbb T}'}_{a_{{\cal E}}(M)}$ is flat. Then the following conditions are equivalent:
\begin{enumerate}[(a)]
\item the functor ${\int a}$ is $\cal E$-final; 

\item For any object $E$ of $\cal E$, any ${\mathbb T}'$-model $c$ in $\cal H$, any ${\mathbb T}$-model $M$ in $\cal E$ and any ${\mathbb T}'$-model homomorphism $x:\gamma_{{\cal E}\slash E}^{\ast}(c)\to !_{E}^{\ast}(a_{\cal E}(M))$, there exists an epimorphic family $\{e_{i}:E_{i}\to E \textrm{ | } i\in I\}$ in $\cal E$ and for each $i\in I$ a ${\mathbb T}'$-model homomorphism $f_{i}:c\to a_{\Set}(c_{i})$, where $c_{i}$ is a ${\mathbb T}$-model in $\cal K$, and a ${\mathbb T}'$-model homomorphism $x_{i}:\gamma_{{\cal E}\slash E_{i}}^{\ast}(a_{\Set}(c_{i}))\to !_{E_{i}}(a_{\cal E}(M))$ such that $x_{i}\circ \gamma_{{\cal E}\slash E_{i}}^{\ast}(f_{i})=e_{i}^{\ast}(x)$ for all $i\in I$; 

\item For any object $E$ of $\cal E$ and any $\Sigma$-structure homomorphism $x:c\to Hom_{\cal E}(E, a_{\cal E}(M))$, where $c$ is a ${\mathbb T}'$-model in $\cal H$, there exists an epimorphic family $\{e_{i}:E_{i}\to E \textrm{ | } i\in I\}$ in $\cal E$ and for each $i\in I$ a ${\mathbb T}'$-model homomorphism $f_{i}:c\to a_{\Set}(c_{i})$, where $c_{i}$ is a ${\mathbb T}$-model in $\cal K$, and a $\Sigma$-structure homomorphism $x_{i}:a_{\Set}(c_{i})\to Hom_{{\cal E}}(E_{i}, a_{\cal E}(M))$ such that $x_{i} \circ f_{i}=Hom_{\cal E}(e_{i}, a_{\cal E}(M))\circ x$ for all $i\in I$.
\end{enumerate}
\end{proposition}

\begin{proofs}
The equivalence between the first two conditions immediately follows from the fact that, $H_{M}$ being flat, the $\cal E$-indexed category ${\int H_{M}}$ is $\cal E$-filtered and hence Remark \ref{fullfinal}(c) applies, while the equivalence between the second and the third condition follows from Proposition \ref{explicit}.    
\end{proofs}

\begin{remark}
If ${\cal E}=\Set$ and ${\mathbb T}$ is a quotient of ${\mathbb T}'$ the condition rewrites as follows: for any ${\mathbb T}$-model $M$ in $\Set$, any model $c$ in $\cal H$ and any ${\mathbb T}'$-model homomorphism $f:c\to M$, there exists a ${\mathbb T}$-model $d$ in $\cal K$ and homomorphisms $g:c\to d$ and $h:d\to M$ such that $h\circ g=f$.

In the particular case where ${\mathbb T}'$ is the empty theory over a signature $\Sigma$, ${\mathbb T}$ is a ${\cal F}$-finitary geometric theory over $\Sigma$ (in the sense of \cite{vickers}), $\cal H$ is the category of finite $\Sigma$-structures and $\cal K$ is the category of finite $\mathbb T$-models, the condition, required for every model $M$ of ${\mathbb T}$ in $\Set$, specializes precisely to the `finite structure condition' of \cite{vickers}. 
\end{remark}

The following result shows a relation between the action on points of a morphism of classifying toposes and a related induced action on flat functors.

\begin{theorem}\label{injt}
Let $a:\Set[{{\mathbb T}}]\to \Set[{{\mathbb T}'}]$ be an interpretation of geometric theories and $\cal K$ and $\cal H$ small subcategories respectively of ${{{\mathbb T}}\textrm{-mod}({\Set})}$ and of ${{{\mathbb T}'}\textrm{-mod}({\Set})}$ such that the functor 
\[
a_{\Set}:{{{\mathbb T}}\textrm{-mod}({\Set})} \to {{{\mathbb T}'}\textrm{-mod}({\Set})}  
\]
restricts to a functor $f:{\cal K}\to {\cal H}$. Then for any Grothendieck topos $\cal E$ the extension functor
\[
\textbf{Flat}({\cal H}^{\textrm{op}}, {\cal E})\to \textbf{Flat}({\cal K}^{\textrm{op}}, {\cal E}) \to {\mathbb T}\textrm{-mod}({\cal E})
\] 
along the geometric morphism $[f, \Set]:[{\cal K}, \Set]\to [{\cal H}, \Set]$ takes values in the full subcategory $Im(a_{\cal E})$ of ${\mathbb T}\textrm{-mod}({\cal E})$.   
\end{theorem}

\begin{proofs}
The diagram
\[  
\xymatrix {
[{\cal H}, \Set]  \ar[rrr]^{p_{\cal H}} & & & \Set[{\mathbb T}'] \\
[{\cal K}, \Set] \ar[u]^{[f, \Set]} \ar[rrr]_{p_{\cal K}} & & & \Set[{\mathbb T}] \ar[u]^{a}} 
\] 
clearly commutes by definition of the functor $f$. Therefore, in view of Diaconescu's theorem, for any Grothendieck topos $\cal E$ we have a commutative diagram
\[  
\xymatrix {
\textbf{Flat}({\cal H}^{\textrm{op}}, {\cal E})  \ar[rrr]^{u^{{\mathbb T}'}_{({\cal H}, {\cal E})}} & & & {\mathbb T}'\textrm{-mod}({\cal E}) \\
\textbf{Flat}({\cal K}^{\textrm{op}}, {\cal E}) \ar[u]^{\xi_{\cal E}} \ar[rrr]^{u^{{\mathbb T}}_{({\cal K}, {\cal E})}} & & & {{\mathbb T}}\textrm{-mod}({\cal E}) \ar[u]^{a_{\cal E}},} 
\] 
where $u^{\mathbb T}_{({\cal K}, {\cal E})}$ and $u^{{\mathbb T}'}_{({\cal H}, {\cal E})}$ are the functors of section \ref{extsym} and $\xi_{\cal E}$ is the extension functor induced by the geometric morphism $[f, \Set]$, from which our thesis immediately follows.
\end{proofs}

The following corollary can be obtained as a particular case of the theorem by taking ${\mathbb T}'$ to be the geometric theory of flat functors on ${\cal C}^{\textrm{op}}$, ${\mathbb T}$ to be its injectivization (in the sense of Definition \ref{injectiviz}), $\cal H$ equal to $\cal C$ and $\cal K$ equal to $\cal D$. 

\begin{corollary}\label{monic}
Let $\cal D$ be a small subcategory of a small category $\cal C$ whose arrows are all monic (in ${\cal C}$), and $\cal E$ a Grothendieck topos. Then

\begin{enumerate}[(a)]
\item For any object $d\in {\cal D}$, $\tilde{F}(d)$ is a decidable object of $\cal E$;

\item For any natural transformation $\alpha:F\to G$ between two flat functors $F, G:{\cal D}^{\textrm{op}}\to {\cal E}$, $\tilde{\alpha}(c):\tilde{F}(c)\to \tilde{G}(c)$ is monic in $\cal E$ for every $c\in {\cal C}$. In particular, $\alpha(d):F(d)\to G(d)$ is monic in $\cal E$ for every $d\in {\cal D}$. 
\end{enumerate}
\end{corollary}\qed

The following corollary is obtained by applying the theorem in the case ${\mathbb T}$ is equal to the injectivization of ${\mathbb T}'$ (in the sense of Definition \ref{injectiviz}) and ${\cal H}={\cal K}={\textrm{f.p.} {\mathbb T}'\textrm{-mod}(\Set)}$.

\begin{corollary}\label{monicarr}
Let $\mathbb T$ be a theory of presheaf type such that all the $\mathbb T$-model homomorphisms between its finitely presentable (set-based) models are sortwise injective. Then every $\mathbb T$-model homomorphism in every Grothendieck topos is sortwise monic. 
\end{corollary}\qed

\subsection{Finitely presentable and finitely generated models}\label{fpfgmodels}

In this section, we discuss, for the purpose of establishing our main criteria for the injectivization of a theory of presheaf type to be again of presheaf type, the relationship between finitely presentable and finitely generated models of a given geometric theory.

Throughout this section, we shall assume for simplicity all our first-order signatures to be one-sorted, if not otherwise stated. It is certainly possible to extend our definitions and results to the general situation, but we shall not embark in the straightforward task of making this explicit.

Let $\Sigma$ be a first-order signature. Recall from \cite{Hodges} that, given a $\Sigma$-structure $M$ and a subset $A\subseteq M$, the $\Sigma$-structure $<A>$ generated by $A$ is the smallest $\Sigma$-substructure of $M$ containing $A$; by the proof of Theorem 1.2.3 \cite{Hodges}, $<A>$ can be concretely represented as the union $\mathbin{\mathop{\textrm{ $\bigcup$}}\limits_{n\in {\mathbb N}}}A_{n}$, where the subsets $A_{n}$ are defined inductively as follows: $A_{0}=A \cup \{c^{M}\}$, where $c^{M}$ are the interpretations in $M$ of the constants over $\Sigma$, and $A_{n+1}=A_{n}\cup \{f^{M}(\vec{a}): \textrm{for some $n\gt 0$, $f$ is a $n$-ary function symbol of $\Sigma$ and}$ \textrm{$\vec{a}$ is a $n$-tuple of elements of $A_{n}$}\}. A $\Sigma$-structure $M$ is said to be \emph{finitely generated} if there exists a finite subset $A\subseteq M$ such that $M=<A>$. 

In the category ${\mathbb T}\textrm{-mod}(\Set)$ of set-based models of a geometric theory $\mathbb T$, we say that an object $B$ is a \emph{quotient} of an object $A$ if there exists a surjective $\mathbb T$-model homomorphism $A \to B$. 

\begin{proposition}\label{quotient}
Let $\mathbb T$ be a geometric theory over a signature $\Sigma$. Then any quotient in ${\mathbb T}\textrm{-mod}(\Set)$ of a finitely generated $\mathbb T$-model is finitely generated.
\end{proposition} 

\begin{proofs}
Let $q:M\to N$ be a surjective homomorphism of $\mathbb T$-models. Suppose that $M$ is finitely generated, by a finite subset $A\subseteq M$. Set $B=q(A)$. This set is clearly finite, and $N$ is equal to $<B>$. To prove this we show, by induction on $n$, that for any $n\in {\mathbb N}$, $q(A_{n})=B_{n}$. As $q$ is a $\Sigma$-structure homomorphism, we have that $q(A_{0})=B_{0}$ and if $q(A_{n})=B_{n}$ then clearly  $q(A_{n+1})=B_{n+1}$. Now, as $q$ is surjective, $N=q(A)=q(\mathbin{\mathop{\textrm{ $\bigcup$}}\limits_{n\in {\mathbb N}}}A_{n})=\mathbin{\mathop{\textrm{ $\bigcup$}}\limits_{n\in {\mathbb N}}}q(A_{n})=\mathbin{\mathop{\textrm{ $\bigcup$}}\limits_{n\in {\mathbb N}}}B_{n}$, as required.  
\end{proofs}

\begin{proposition}\label{coincid}
Let $\mathbb T$ be a geometric theory over a signature $\Sigma$ whose axioms are all of the form $(\phi \vdash_{\vec{x}} \psi)$, where $\psi$ is a quantifier-free geometric formula. Then any finitely presentable $\mathbb T$-model is finitely generated.
\end{proposition} 

\begin{proofs}
First, we notice that any geometric theory over a signature $\Sigma$ such that all its axioms are of the form $(\phi \vdash_{\vec{x}} \psi)$, where $\psi$ is a quantifier-free geometric formula, satisfies the property that any substructure of a model of $\mathbb T$ is a model of $\mathbb T$.

Let $M$ be a finitely presentable $\mathbb T$-model. We can clearly represent $M$ as a filtered colimit (actually, directed union) of its finitely generated substructures (equivalently, submodels, cf. the above remark). By definition of finite presentability, the identity on $M$ factors through one of the embeddings of a finitely generated $\mathbb T$-submodel of $M$ into $M$; such embedding is thus necessarily an isomorphism, and $M$ coincides with its domain. In particular, $M$ is finitely generated, as required.
\end{proofs}

\begin{proposition}\label{pr}
Let $\mathbb T$ be a geometric theory over a signature $\Sigma$ whose category of set-based models is finitely accessible (in particular, a theory of presheaf type). Then any finitely generated $\mathbb T$-model is a quotient of a finitely presentable models of $\mathbb T$ in the category ${\mathbb T}\textrm{-mod}(\Set)$.

If moreover all the axioms of $\mathbb T$ are of the form $(\phi \vdash_{\vec{x}} \psi)$, where $\psi$ is a quantifier-free geometric formula, then the finitely generated $\mathbb T$-models are precisely the quotients of the finitely presentable models of $\mathbb T$ in the category ${\mathbb T}\textrm{-mod}(\Set)$.
\end{proposition} 

\begin{proofs}
As the category ${\mathbb T}\textrm{-mod}(\Set)$ is finitely accessible, $M$ can be expressed as a filtered colimit $(M=colim(N_{i}), \{ J_{i}:N_{i}\to M \textrm{ | } i\in {\cal I}\}) $ of finitely presentable models $N_{i}$. Let $a_{1}, \ldots, a_{n}$ a finite set of generators for $M$. As the family of arrows $\{J_{i}:N_{i}\to M \textrm{ | } i\in {\cal I}\}$ is jointly surjective (since filtered colimits in ${\mathbb T}\textrm{-mod}(\Set)$ are computed pointwise as in $\Set$), for any generator $a_{j}$ there exists an index $k_{j}\in {\cal I}$ and an element $d_{j}$ of $N_{k_{j}}$ such that $J_{k_{j}}(d_{j})=a_{j}$. Now, the set of generators $a_{j}$ being finite and the index category $\cal I$ being filtered, we can suppose without loss of generality that all the $k_{j}$ are equal. We thus have an index $k\in {\cal I}$ and a string of elements of $N_{k}$ which are sent by $J_{k}$ to the generators of $M$; hence the arrow $J_{k}:N_{k}\to M$ is surjective. 

The second part of the proposition follows from Propositions \ref{coincid} and \ref{quotient}, which ensure that every quotient of a finitely presentable model is finitely generated.   
\end{proofs}

\begin{proposition}\label{reflectiverel}
Let $\Sigma$ be a signature without relation symbols and $\mathbb T$ a geometric theory over a larger signature $\Sigma'$ obtained from $\Sigma$ by solely adding relation symbols whose interpretation in any set-based $\mathbb T$-model coincides with the complement of a geometric formula over $\Sigma$ (for instance, the injectivization of a geometric theory over $\Sigma$). Suppose that all the $\mathbb T$-model homomorphisms between set-based $\mathbb T$-models are injective. Then any finitely generated $\mathbb T$-model is finitely presentable.  
\end{proposition}

\begin{proofs}
From Proposition \ref{coincid}, we know that every finitely presentable $\mathbb T$-model is finitely generated. It thus remains to prove the converse. 

We observed in \cite{OCG} (Lemma 6.2) that, given a category $\cal D$ with filtered colimits, an object $M$ of $\cal D$ is finitely presentable in $\cal D$ if and only if for any filtered diagram $D:{\cal I}\to {\cal D}$, any arrow $M\to colim(D)$ factors through one of the canonical colimit arrows $J_{i}:D(i)\to colim(D)$ and for any arrows $f:M\to D(i)$ and $g:M\to D(j)$ in $\cal D$ such that $J_{i}\circ f=J_{j}\circ g$, there exists $k\in {\cal I}$ and two arrows $s:i\to k$ and $t:j\to k$ such that $D(s)\circ f=D(t)\circ g$. Now, the latter condition is automatically satisfied if all the arrows of $\cal D$ are monic (since the category $\cal I$ is filtered); in particular, applying this criterion to the case of our theory $\mathbb T$ yields the following characterization of the finitely presentable $\mathbb T$-models: a set-based $\mathbb T$-model $M$ is finitely presentable if and only if for any filtered colimit $colim(N_{i})$ of set-based $\mathbb T$-models, with colimit arrows $J_{i}:N_{i}\to colim(N_{i})$, any $\mathbb T$-model homomorphism $f:M\to colim(N_{i})$ factors through at least one arrow $J_{i}$ in the category ${\mathbb T}\textrm{-mod}(\Set)$. 

Now, if $M$ is finitely generated by elements $a_{1}, \ldots, a_{n}$, for any of generator $a_{i}$ of $M$ there exists an index $k_{i}$ of $\cal I$ and an element $b_{i}\in N_{k_{i}}$ such that $f(a_{i})=J_{i}(b_{i})$. As the category $\cal I$ is filtered, we can suppose without loss of generality that all the $k_{i}$ are equal to each other. Therefore there exists an index $k\in {\cal I}$ and elements $b_{1}, \ldots, b_{n}$ of $N_{k}$ such that $f(a_{i})=J_{k}(b_{i})$ for all $i$. The image of $M$ under $f$ is thus entirely contained in $N_{k}$ and hence we have a function $g:M \to N_{k}$ such that $J_{k}\circ g=f$. Let us prove that $g$ is a $\Sigma'$-structure homomorphism. The fact that $g$ is a $\Sigma$-structure homomorphism, i.e. that it commutes with the interpretation of the function symbols over $\Sigma$, follows from the fact that $f$ is by the injectivity of $J_{k}$. Concerning the preservation of atomic relations over $\Sigma'$, this is guaranteed by the fact that, were they not preserved, $f$ would preserve the relations of satisfaction by elements of the geometric formulae over $\Sigma$ defining the complements of such relations, contradicting the fact that $f$ is a $\Sigma'$-structure homomorphism.
\end{proofs}

A useful sufficient condition for condition $(iii)(b)\textrm{-}(1)$ of Theorem \ref{main} to hold, which is applicable to categories $\cal K$ whose arrows are not necessarily monomorphisms is the following.

\begin{theorem}
Let $\mathbb T$ be a geometric theory over a signature $\Sigma$, and $\cal K$ a small full subcategory of ${{\mathbb T}\textrm{-mod}(\Set)}$ such that every model in $\cal K$ is a quotient of a $\mathbb T$-model which is finitely presented by a geometric formula over $\Sigma$. Then $\mathbb T$ satisfies condition $(iii)(b)\textrm{-}(1)$ of Theorem \ref{main} with respect to the category $\cal K$.  
\end{theorem}

\begin{proofs}
We have to show that the extension functor
\[
u^{\mathbb T}_{({\cal K}, {\cal E})}:\mathbf{Flat}({\cal K}^{\textrm{op}}, {\cal E}) \to \mathbf{Flat}_{J_{\mathbb T}}({\cal C}_{\mathbb T}, {\cal E})
\]
of section \ref{extsym} is faithful. 

By Theorem \ref{extsymthm}, for any formula $\{\vec{x}. \phi\}$ presenting a $\mathbb T$-model $M_{\{\vec{x}. \phi\}}$, $u^{\mathbb T}_{({\cal K}, {\cal E})}(F)(\{\vec{x}. \phi\})\cong F(M_{\{\vec{x}. \phi\}})$. The thesis thus follows immediately from Proposition \ref{monoepi}. 
\end{proofs}

\subsection{Reformulations of condition $(iii)$ of Theorem \ref{main}}

Let us work in the context of a faithful interpretation of theories as in the last section.

Let 
\[
u^{\mathbb T}_{({\cal K}, {\cal E})}:\textbf{Flat}({\cal K}^{\textrm{op}}, {\cal E}) \to {\mathbb T}\textrm{-mod}({\cal E}) 
\]
and
\[
u^{{\mathbb T}'}_{({\cal H}, {\cal E})}: \textbf{Flat}({\cal H}^{\textrm{op}}, {\cal E}) \to {\mathbb T}'\textrm{-mod}({\cal E}) 
\]
be the functors of section \ref{extsym}.

Suppose that ${\mathbb T}'$ is of presheaf type and that $\cal H$ contains a category $\cal P$ whose Cauchy-completion coincides with the category of finitely presentable ${\mathbb T}'$-models (so that ${\mathbb T}'$ is classified by the topos $[{\cal P}, \Set]$). Notice that, since the functor $f:{\cal K}\to {\cal H}$ is a restriction of the functor $a_{\Set}$, $\cal K$ gets identified under it with a subcategory of $\cal H$ and hence, by Remark \ref{cov}(c), the extension functor 
\[
\xi_{\cal E}:\textbf{Flat}({\cal K}^{\textrm{op}}, {\cal E})  \to \textbf{Flat}({\cal H}^{\textrm{op}}, {\cal E}) 
\]   
along the geometric morphism $[f, \Set]$ is faithful.

From the fact that ${\mathbb T}'$ faithfully interprets in $\mathbb T$ (in the sense of Definition \ref{faith}) it follows that the functor $u^{\mathbb T}_{({\cal K}, {\cal E})}$ is faithful (resp. full and faithful) if and only if the composite functor $u^{{\mathbb T}'}_{({\cal H}, {\cal E})}\circ
\xi_{\cal E}$ is, when regarded as a functor with values in the subcategory $Im(a_{\cal E})$ of ${\mathbb T}'\textrm{-mod}({\cal E})$. The interest of this reformulation lies in the alternative description of the functor 
\[
u^{{\mathbb T}'}_{({\cal H}, {\cal E})}:\textbf{Flat}({\cal H}^{\textrm{op}}, {\cal E}) \to {\mathbb T}'\textrm{-mod}({\cal E})
\]
which is available under our hypotheses. Indeed, the following proposition holds.

\begin{proposition}\label{restriction}
Under the hypotheses specified above, the functor 
\[
u^{{\mathbb T}'}_{({\cal H}, {\cal E})}:\textbf{Flat}({\cal H}^{\textrm{op}}, {\cal E}) \to {\mathbb T}'\textrm{-mod}({\cal E})\simeq \textbf{Flat}({\cal P}^{\textrm{op}}, {\cal E})
\]
sends any flat functor $F:{\cal H}^{\textrm{op}} \to {\cal E}$ to its restriction $F|{\cal P}:{\cal P}^{\textrm{op}} \to {\cal E}$, and acts accordingly on the natural transformations. 
\end{proposition}

\begin{proofs}
This immediately follows from Theorem \ref{extsymthm} in view of the fact that $\cal H$ contains $\cal P$ and all the objects of $\cal P$ are finitely presented ${\mathbb T}'$-models (since ${\mathbb T}'$ is by our hypothesis of presheaf type classified by the topos $[{\cal P}, \Set]$).  
\end{proofs}

As a corollary of the proposition, we immediately obtain the following result.

\begin{theorem}\label{thmsimpl}
Under the hypotheses specified above, we have that
\begin{enumerate}[(i)]
\item ${\mathbb T}$ satisfies condition $(iii)(b)\textrm{-}(1)$ of Theorem \ref{main} with respect to the category $\cal K$ if and only if for any Grothendieck topos $\cal E$, every natural transformation between flat functors in the image of the functor $\xi_{\cal E}$ is determined by its restriction to the category ${\cal P}$;

\item ${\mathbb T}$ satisfies condition $(iii)(b)\textrm{-}(2)$ of Theorem \ref{main} with respect to the category $\cal K$ if and only if for any Grothendieck topos $\cal E$, every natural transformation between restrictions to ${\cal P}$ of flat functors in the image of the functor $\xi_{\cal E}$ can be extended to a natural transformation between the functors themselves. 
\end{enumerate}
\end{theorem}
  
Let us now consider the satisfaction of condition $(iii)(b)$ of Theorem \ref{main} by the injectivization of a geometric theory $\mathbb T$ satisfying the hypotheses of Propositions \ref{coincid} and \ref{reflectiverel}. Since by the proposition every finitely generated $\mathbb T$-model is finitely presentable in ${\mathbb T}\textrm{-mod}(\Set)$ and every model of $\mathbb T$ is a directed union of its finitely generated submodels, the category ${\mathbb T}\textrm{-mod}(\Set)$ is equivalent to $\mathbf{Flat}({{\textrm{f.g.} {\mathbb T}\textrm{-mod}(\Set)}}^{\textrm{op}}, \Set)$, where ${\textrm{f.g.} {\mathbb T}\textrm{-mod}(\Set)}$ is the full subcategory of ${\mathbb T}\textrm{-mod}(\Set)$ on the finitely generated $\mathbb T$-models. The following theorem shows that condition $(iii)(b)\textrm{-}(1)$ is satisfied. Before stating it, we need a lemma.

\begin{lemma}\label{monoepi}
Let $\cal C$ be a small category, $\cal E$ a topos and $F:{\cal C}^{\textrm{op}}\to {\cal E}$ a flat functor. Then for any epimorphism $f:a\to b$ in ${\cal C}$, the arrow $F(f):F(b)\to F(a)$ is monic in $\cal E$.
\end{lemma}

\begin{proofs}
The thesis follows at once from the fact that flat functors preserve all the finite limits which exist in the domain category using the well-known characterization of monomorphisms in terms of pullbacks. 
\end{proofs}
 
\begin{theorem}\label{thminj1}
Let $\mathbb T$ be the injectivization of a theory of presheaf type ${\mathbb T}'$ (in the sense of Definition \ref{injectiviz}) such that the every finitely presentable ${\mathbb T}'$-model is finitely generated (cf. for instance Proposition \ref{coincid}) and the finitely generated ${\mathbb T}'$-models are all quotients (in the sense of section \ref{fpfgmodels}) of finitely presentable ${\mathbb T}'$-models (cf. for instance Proposition \ref{pr}). Then $\mathbb T$ satisfies condition $(iii)(b)\textrm{-}(1)$ of Theorem \ref{main} with respect to the category of finitely generated ${\mathbb T}'$-models. 
\end{theorem}

\begin{proofs}
It suffices to apply Proposition \ref{restriction} and Lemma \ref{monoepi} in conjunction with Theorem \ref{thmsimpl} by taking ${\cal K}$ equal to the category of finitely generated ${\mathbb T}'$-models and sortwise monic homomorphisms between them, $\cal H$ equal to the category of finitely generated ${\mathbb T}'$-models and homomorphisms between them and $\cal P$ equal to the category of finitely presentable ${\mathbb T}'$-models.
\end{proofs}

The following proposition can be useful in verifying that condition $(ii)$ of Theorem \ref{thmsimpl} holds.

\begin{proposition}\label{longprop}
Let ${\cal D}$ be a full subcategory of a category $\cal C$ with filtered colimits such that 
\begin{itemize}
\item for any object $c$ of $\cal C$ there exists an epimorphism $a\to c$ in $\cal C$ from an object $a$ of $\cal D$ to $c$,
\item every object of $\cal D$ is finitely presentable in $\cal C$ and
\item every object $c$ of $\cal C$ can be canonically expressed as a filtered colimit $colim(d_{i})$ of objects of $\cal D$, with colimit arrows $J_{i}:d_{i} \to c$ (for $i\in {\cal I}$).
\end{itemize}
Let $F$ and $G$ be flat functors ${\cal C}^{\textrm{op}}\to {\cal E}$. Then a natural transformation $\beta:F|_{{\cal D}}\to G|_{{\cal D}}$ lifts (uniquely) to a natural transformation $F \to G$ if and only if for every epimorphism $q:a \to b$ of the form $J_{i}$, the arrow $\beta(a):F(a)\to G(a)$ restricts to an arrow $F(b)\to G(b)$ along the monomorphisms $F(q):F(b)\to F(a)$ and $G(q):G(b)\to G(a)$.  
\end{proposition}

\begin{proofs}
Clearly, the `only if' direction is trivially satisfied, so we only have to care about the `if' direction. The uniqueness of the extended natural transformation follows from Proposition \ref{monoepi}, so we just have to prove the existence.

Notice that for every object $c$ of $\cal C$ there exists an $i\in {\cal I}$ such that the arrow $J_{i}:d_{i}\to c$ is an epimorphism. Indeed, by our hypotheses there exists an epimorphism $q:a\to c$ to $c$ from an object $a$ of $\cal D$ to $c$ and, since $a$ is finitely presentable in $\cal C$, $q$ necessarily factors through an arrow of the form $J_{i}$, which is necessarily epic as $q$ is.

Let $\beta:F|_{{\cal D}}\to G|_{{\cal D}}$ be a natural transformation. Choose an $i\in {\cal I}$ such that the arrow $J_{i}:d_{i}\to c$ is an epimorphism (such arrow exists by the above remark), and set $\alpha(c)$ equal to the restriction of $\beta(d_{i})$ along the arrows $F(J_{i})$ and $G(J_{i})$.  

First, let us verify that this definition is well-posed, i.e. that it does not depend on the choice of $i\in {\cal I}$. For any other index $j$ such that $J_{j}$ is an epimorphism with codomain $c$, the filteredness of $\cal I$ ensures the existence of an index $k$ and two arrows $\xi:i\to k$ and $\chi:j\to k$ in $\cal I$, whence $J_{k}\circ \xi=J_{i}$ and $J_{k}\circ \chi=J_{j}$. Now, if we denote by $u$ the restriction of the arrow $\beta(d_{k})$ along the arrows $F(J_{k})$ and $G(J_{k})$ we obtain, by invoking the naturality of $\beta$ with respect to the arrows $\xi$ and $\chi$ in $\cal I$, that $u$ is equal to the unique restriction of $\beta(d_{i})$ along the arrows $F(J_{i})$ and $G(J_{i})$ and to the unique restriction of $\beta(d_{j})$ along the arrows $F(J_{j})$ and $G(J_{j})$; in particular, these two restrictions are equal, as required.   

It remains to show that the assignment $c\to \alpha(c)$ is natural in $c$, i.e. that for any arrow $f:c\to c'$ in $\cal C$, the following naturality square commutes:
\[  
\xymatrix {
F(c)  \ar[r]^{\alpha(c)} & G(c) \\
F(c') \ar[u]^{F(f)} \ar[r]_{\alpha(c')} & G(c') \ar[u]_{G(f)}.} 
\] 
To prove this we observe that, using the fact that there exists an epimorphism $J_{i}:d_{i}\to c$ from an object of $\cal D$ to $c$, we can suppose without loss of generality $c$ to lie in $\cal D$; indeed, since $G$ sends epimorphisms to monomorphisms, the commutativity of the naturality square above is equivalent, by definition of $\alpha(c)$, to the commutativity of the naturality square relative to the arrow $f\circ J_{i}$. Now, consider the canonical representation of $c'$ as a filtered colimit $colim(e_{j})$ (for $j\in {\cal J}$) of objects in $\cal D$, with colimit arrows $K_{j}:e_{j}\to c'$ ($j\in {\cal J}$). Since $c$ lies in $\cal D$, $c$ is finitely presentable in $\cal C$ by our hypotheses and hence there exists an index $j\in {\cal J}$ and an arrow $r:c\to e_{j}$ such that $K_{j}\circ r=f$. The commutativity of our diagram thus follows from the commutativity of the naturality square of $\beta$ relative to the arrow $r$ and the definition of $\alpha(c')$ in terms of $\beta(e_{j})$.   
\end{proofs}

The following result, obtained by combining Theorem \ref{monoflat} and Proposition \ref{longprop}, gives a criterion for condition $(iii)(b)\textrm{-}(2)$ to hold. 

\begin{theorem}
Let $\mathbb T$ be the injectivization of a theory of presheaf type ${\mathbb T}'$ such that every finitely presentable ${\mathbb T}'$-model is finitely generated (cf. for instance Proposition \ref{coincid}), the finitely generated ${\mathbb T}'$-models are all quotients (in the sense of section \ref{fpfgmodels}) of finitely presentable ${\mathbb T}'$-models (cf. for instance Proposition \ref{pr}) and every monic ${\mathbb T}'$-model homomorphism between finitely generated models of $\mathbb T$ is sortwise monic. Then $\mathbb T$ satisfies condition $(iii)(b)\textrm{-}(2)$ of Theorem \ref{main} with respect to the category $\cal C$ of finitely generated models of ${\mathbb T}'$ and sortwise monic homomorphisms between them if and only if for any flat functors $F,G:{\cal C}^{\textrm{op}}\to {\cal E}$, denoting by $\tilde{F}$ and $\tilde{G}$ their extension to the category of finitely generated ${\mathbb T}'$-models and homomorphisms between them and by $\cal P$ the full subcategory of $\cal C$ on the finitely presentable ${\mathbb T}'$-models, any natural transformation $\alpha:F|_{\cal P}\to G|_{\cal P}$ satisfies the following property: for any quotient $\mathbb T$-model homomorphism $q:a\to b$, where $a$ is finitely presentable and $b$ is finitely generated, the arrow $\alpha(b)$ restricts along the monic arrows $F(q)$ and $G(q)$.   
\end{theorem}\qed

\begin{remark}
Notice that if $F$ and $G$ are representable functors then the condition of the theorem is trivially satisfied, the category ${\mathbb T}'\textrm{-mod}(\Set)$ being finitely accessible.
\end{remark}

\subsection{Quotient theories}\label{secquotients}

Recall from \cite{OC6} that a \emph{quotient} of a geometric theory $\mathbb T$ over a signature $\Sigma$ is a geometric theory ${\mathbb T}'$ over $\Sigma$ such that every axiom of $\mathbb T$ is provable in ${\mathbb T}'$. 

It is often useful, in investigating whether a given geometric theory $\mathbb S$ is of presheaf type, to consider it in relation to a theory $\mathbb T$ of presheaf type of which $\mathbb S$ is a quotient.

The following result is a corollary of Theorem \ref{finalitythm} and Proposition \ref{finality}.

\begin{corollary}\label{corquotients}
Let $\mathbb T$ be a theory of presheaf type over a signature $\Sigma$ and $\mathbb S$ be a quotient of $\mathbb T$ such that all the finitely presentable $\mathbb S$-models are finitely presentable as $\mathbb T$-models (for instance, $\mathbb T$ can be the empty theory over a finite signature and $\mathbb S$ can be any geometric theory over $\Sigma$ whose finitely presentable models are all finite, cf. Theorem 6.4 \cite{OCG}). Suppose moreover that for any object $E$ of a Grothendieck topos $\cal E$ and any $\Sigma$-structure homomorphism $x:c \to Hom_{\cal E}(E, M)$, where $c$ is a finitely presentable $\mathbb T$-model and $M$ is a $\mathbb S$-model in $\cal E$, there exists an epimorphic family $\{e_{i}:E_{i}\to E \textrm{ | } i\in I\}$ in $\cal E$ and for each $i\in I$ a $\mathbb T$-model homomorphism $f_{i}:c\to c_{i}$, where $c_{i}$ is a finitely presentable $\mathbb S$-model, and a $\Sigma$-structure homomorphism $x_{i}:c_{i}\to Hom_{\cal E}(E_{i}, M)$ such that $x_{i}\circ f_{i}=Hom_{\cal E}(e_{i}, M)\circ x$ for all $i\in I$. Then $\mathbb S$ is of presheaf type.
\end{corollary}

\begin{proofs}
Condition $(iii)$ of Theorem \ref{main} satisfied by Proposition \ref{criteriongeometricity}, while the fact that conditions $(i)$ and $(ii)$ of Theorem \ref{main} are satisfied follows immediately from Theorem \ref{finalitythm} in view of Proposition \ref{finality} (take ${\mathbb T}$ equal to $\mathbb S$, ${\mathbb T}'$ equal to $\mathbb T$, ${\cal K}$ equal to the category of finitely presentable ${\mathbb S}$-models and ${\cal H}$ equal to the category of finitely presentable ${\mathbb T}$-models - the fact that ${\int a}$ is $\cal E$-full is clear). 
\end{proofs}

We shall say that the family of homomorphisms $x_{i}$ in the statement of the corollary defines a `localization' of the homomorphism $x$. 

\begin{remarks}\label{cartesianiz}
\begin{enumerate}[(a)]
\item Corollary \ref{corquotients} generalizes Theorem 3.6 \cite{vickers}, whose hypotheses, which are expressed syntactically in terms of geometric logic, once interpreted in topos-theoretic semantics, are stronger than those of Corollary \ref{corquotients} in the particular case where $\mathbb T$ is the empty theory over a finite signature and $\mathbb S$ is any geometric theory over $\Sigma$ whose finitely presentable models are all finite.

\item Theorem \ref{closed} follows as an immediate consequence of the corollary, observing that if $\mathbb T$ and ${\mathbb T}'$ are theories as in the hypotheses of the theorem then any $\mathbb T$-model homomorphism having as codomain a ${\mathbb T}'$-model also has as domain a ${\mathbb T}'$-model whence the hypotheses of the corollary are trivially satisfied.   
 
\item Assuming the first hypothesis of the corollary, the second is actually necessary for $\mathbb S$ to be of presheaf type. Indeed, if $\mathbb S$ is of presheaf type then every model $M$ of $\mathbb S$ is a $\cal E$-filtered colimit of the associated canonical diagram of `constant' finitely presentable $\mathbb S$-models (cf. sections \ref{secfp} and \ref{reint}), which by the first hypothesis are also finitely presentable $\mathbb T$-models; the $\cal E$-finite presentability of $c$ thus implies the second hypothesis of the corollary (cf. Proposition \ref{propsemfp}). 

\item Corollary \ref{corquotients} is often applied to pairs of the form $({\mathbb S}, {\mathbb T}={\mathbb S}_{c})$, where $\mathbb S$ is a geometric theory over a signature $\Sigma$ and $\mathbb T$ equal to the \emph{cartesianization} (or finite-limit part) of ${\mathbb S}_{c}$ of $\mathbb S$, namely the cartesian theory over $\Sigma$ consisting of all the cartesian sequents over $\Sigma$ which are provable in $\mathbb S$ (recall that a $\mathbb S$-cartesian sequent is a sequent over $\Sigma$ involving $\mathbb S$-cartesian formulae over $\Sigma$, that is formulae built from atomic formulae by only using finitary conjunctions and $\mathbb S$-provably unique existential quantifications). Notice that the cartesianization of $\mathbb S$ is the biggest cartesian theory over the signature of $\mathbb S$ of which $\mathbb S$ is a quotient.

The question of whether for any theory of presheaf type $\mathbb T$, every finitely presentable $\mathbb T$-model is finitely presentable as a ${\mathbb T}_{c}$-model is still open and will be addressed in full generality in another paper. For the moment, we limit ourselves to remarking that this property is satisfied by all the examples of theories of presheaf type considered in this paper (cf. section \ref{examples}). 
\end{enumerate}
\end{remarks}

Corollary \ref{corquotients} can also be applied to pairs consisting of a geometric theory over a signature $\Sigma$ and the empty theory over the same signature, provided that the former satisfies appropriate hypotheses. In order to apply the corollary in this context, we need the following lemma. Below, when we say that a $\Sigma$-structure $c$ is finite we mean that for every sort $A$ over $\Sigma$, $cA$ is finite.  

\begin{lemma}\label{lemmafinite}
Let $\mathbb T$ be a geometric theory over a finite signature $\Sigma$ with a finite number of axioms each of which is of the form $(\top \vdash_{\vec{x}} \mathbin{\mathop{\textrm{ $\bigvee$}}\limits_{i\in I}} \phi_{i})$, where the $\phi_{i}$ are finite conjunctions of atomic formulae. Then for any $\mathbb T$-model $M$ in a Grothendieck topos $\cal E$, any object $E$ of $\cal E$ and any $\Sigma$-structure homomorphism $f:c\to Hom_{\cal E}(E, M)$ from a finite $\Sigma$-structure $c$ there exists an epimorphic family $\{e_{i}:E_{i}\to E \textrm{ | } i\in I\}$ in $\cal E$ and for each $i\in I$ a finite $\Sigma$-substructure $c_{i}$ of $Hom_{\cal E}(E_{i}, M)$ which is a model of $\mathbb T$ and a $\Sigma$-structure homomorphism $f_{i}:c\to c_{i}$ such that $Hom_{\cal E}(e_{i}, M)\circ f=j_{i}\circ f_{i}$ (where $j_{i}$ is the canonical inclusion of $c_{i}$ into $Hom_{\cal E}(E_{i}, M)$).
\end{lemma}

\begin{proofs}
For each axiom $\sigma$ of the form $(\top \vdash_{\vec{x_{\sigma}}} \mathbin{\mathop{\textrm{ $\bigvee$}}\limits_{i\in I_{\sigma}}} \phi^{\sigma}_{i})$ and any finite tuple $\vec{\xi}$ of elements of $c$ of the same type as $\vec{x_{\sigma}}$, since $M$ is by our hypotheses a $\mathbb T$-model, there exists an epimorphic family $\{e^{\vec{\xi}}_{i}: E_{i}^{\vec{\xi}}\to E \textrm{ | } i\in I_{\sigma}\}$ in $\cal E$ such that for any $i\in I_{\sigma}$, $f(\vec{\xi})\circ \vec{e_{i}^{\xi}}$ factors through $[[\vec{x_{\sigma}}. \phi_{i}]]_{M}$. Since there is only a finite number of axioms of $\mathbb T$ and of elements of $c$, and the fibered product of a finite number of epimorphic families is again an epimorphic family, there exists an epimorphic family $\{e_{k}:E_{k}\to E \textrm{ | } k\in K\}$ such that for any $k\in K$, any axiom $\sigma$ of $\mathbb T$ and any tuple $\vec{\xi}$ of elements of $c$ of the appropriate sorts, the element $f(\vec{\xi})\circ e_{k}\in Hom_{\cal E}(E_{k}, M)$ satisfies the formula on the right-hand side of the sequent $\sigma$. Hence if we consider, for each $k\in K$, the surjection-inclusion factorization $i_{k}\circ f_{k}$ of the homomorphism $Hom_{\cal E}(e_{k}, M)\circ f:c\to Hom_{\cal E}(E_{k}, M)$ (in the sense of Lemma \ref{lemmafact}), we obtain that $c_{k}$ is finite, since $f_{k}$ is surjective, and that all the axioms of $\mathbb T$ are satisfied in $c_{k}$ (since the formulae $\phi_{i}$ are finite conjunctions of atomic formulae and $c_{k}$ is a substructure of $Hom_{\cal E}(E_{k}, M)$), that is all the $c_{k}$ are models of $\mathbb T$. This proves our thesis.    
\end{proofs}

\begin{corollary}\label{fincor}
Let $\mathbb T$ be a geometric theory over a finite signature $\Sigma$ with a finite number of axioms each of which is of the form $(\top \vdash_{\vec{x}} \mathbin{\mathop{\textrm{ $\bigvee$}}\limits_{i\in I}} \phi_{i})$, where the $\phi_{i}$ are finite conjunctions of atomic formulae. Suppose that for any $\mathbb T$-model $M$ in a Grothendieck topos $\cal E$ and object $E$ of $\cal E$, every finitely generated $\Sigma$-substructure of $Hom_{\cal E}(E, M)$ has sortwise only a finite number of elements besides the constants (for instance, when the signature $\Sigma$ does not contain function symbols except for a finite number of constants). Then $\mathbb T$ is of presheaf type, classified by the category of covariant set-valued functors on the category of finite models of $\mathbb T$.  
\end{corollary}

\begin{proofs}
As the signature of $\mathbb T$ is finite, every model which  sortwise contains only finitely many elements besides the constants is finitely presentable as a model of the empty theory over the signature of $\mathbb T$ (cf. Theorem 6.4 \cite{OCG}). This remark ensures, by Lemma \ref{lemmafinite}, that the theory $\mathbb T$ satisfies the hypotheses of Corollary \ref{corquotients} with respect to the empty theory over its signature. Since the latter theory is clearly of presheaf type, it follows that $\mathbb T$ is of presheaf type as well, as required. 
\end{proofs}

\subsubsection{Presheaf-type quotients and rigid topologies}\label{rigidtopologies}

In this section we shall analyze the presheaf-type quotients of a presheaf-type theory $\mathbb T$ in terms of the associated subtoposes of the classifying topos for $\mathbb T$.

We know from \cite{OC6} that every quotient ${\mathbb T}'$ a theory of presheaf type $\mathbb T$ corresponds to a unique Grothendieck topology $J$ on the category ${\textrm{f.p.} {\mathbb T}\textrm{-mod}(\Set)}^{\textrm{op}}$ (under the duality theorem of \cite{OC6}).

Recall from \cite{El} (Definition C2.2.18) that a Grothendieck topology $J$ on a category $\cal C$ is said to be \emph{rigid} if for every object $c$ of $\cal C$, the family of all the arrows to $c$ in $\cal C$ from $J$-irreducible objects of $\cal C$ (i.e., objects $d$ with the property that the only $J$-covering sieve on $d$ is the maximal one), generates a $J$-covering sieve.

We can characterize the topologies $J$ such that the corresponding subtopos $\Sh({\textrm{f.p.} {\mathbb T}\textrm{-mod}(\Set)}^{\textrm{op}}, J)\hookrightarrow [{\textrm{f.p.} {\mathbb T}\textrm{-mod}(\Set)}, \Set]$ is an essential geometric inclusion (that is, a geometric inclusion whose inverse image which admits a left adjoint) by using the following site characterizations:

$(1)$ A canonical geometric inclusion $\Sh({\cal C}, J)\hookrightarrow [{\cal C}^{\textrm{op}}, \Set]$ from a topos $\Sh({\cal C}, J)$ which is equivalent to a presheaf topos and such that $\cal C$ is a Cauchy-complete category, is essential if and only if the topology $J$ is rigid;

$(2)$ A geometric morphism (resp. geometric inclusion) $[{\cal C}, \Set]\to [{\cal D}, \Set]$, where $\cal D$ is a Cauchy-complete category, is essential if and only if it is induced by a functor (resp. by a full and faithful functor) ${\cal C}\to {\cal D}$.  

The latter characterization is well-known (cf. Lemma A4.1.5 \cite{El} and Example A4.2.12(b)), while the former can be proved as follows. If $J$ is rigid then, denoting by $\cal D$ the full subcategory of $\cal C$ on the $J$-irreducible objects, the Comparison Lemma yields an equivalence $\Sh({\cal C}, J)\simeq [{\cal D}^{\textrm{op}}, \Set]$ which makes the canonical inclusion $\Sh({\cal C}, J)\hookrightarrow [{\cal C}^{\textrm{op}}, \Set]$ isomorphic to the canonical geometric inclusion $[{\cal D}^{\textrm{op}}, \Set] \hookrightarrow [{\cal C}^{\textrm{op}}, \Set]$ induced by the full inclusion of categories ${\cal D}^{\textrm{op}} \hookrightarrow {\cal C}^{\textrm{op}}$; in particular, the morphism  $\Sh({\cal C}, J)\hookrightarrow [{\cal C}^{\textrm{op}}, \Set]$ is essential. Conversely, if the canonical geometric inclusion $i_{J}:\Sh({\cal C}, J)\hookrightarrow [{\cal C}^{\textrm{op}}, \Set]$ is essential then if $\Sh({\cal C}, J)$ is equivalent to $[{\cal D}^{\textrm{op}}, \Set]$, by property $(2)$ $i_{J}$ is isomorphic to a geometric inclusion $[{\cal D}^{\textrm{op}}, \Set]\to [{\cal C}^{\textrm{op}}, \Set]$ induced by a full embedding ${\cal D}^{\textrm{op}}\hookrightarrow {\cal C}^{\textrm{op}}$. Now, the topology $J_{\cal D}$ on $\cal C$ defined by saying that a sieve $R$ on $c$ is $J$-covering if and only if it contains all the morphisms from objects of $\cal D$ to $c$ is clearly rigid, and the Comparison Lemma yields an equivalence $\Sh({\cal C}, J_{\cal D})\simeq [{\cal D}^{\textrm{op}}, \Set]$ which makes the geometric morphism $[{\cal D}^{\textrm{op}}, \Set]\to [{\cal C}^{\textrm{op}}, \Set]$ isomorphic to the canonical inclusion $\Sh({\cal C}, J_{\cal D}) \hookrightarrow [{\cal C}^{\textrm{op}}, \Set]$ (cf. Example A2.2.4(d) \cite{El}). It thus follows that $J=J_{\cal D}$; in particular, $J$ is rigid. 

These site characterizations define the arches of a `bridge' (in the sense of \cite{OCU}) leading to the following result (cf. also Theorem 6.8 \cite{OCG} for a related result).

\begin{theorem}\label{rigidity}
Let ${\mathbb T}'$ be a quotient of a theory of presheaf type $\mathbb T$, corresponding to a Grothendieck topology $J$ on the category ${\textrm{f.p.} {\mathbb T}\textrm{-mod}(\Set)}^{\textrm{op}}$ under the duality theorem of \cite{OC6}. Suppose that ${\mathbb T}'$ is itself of presheaf type. Then every finitely presentable ${\mathbb T}'$-model is finitely presentable also as a $\mathbb T$-model if and only if the topology $J$ is rigid.
\end{theorem}

\begin{proofs}
If every finitely presentable ${\mathbb T}'$-model is finitely presentable also as a $\mathbb T$-model then the geometric inclusion corresponding via the duality theorem to the quotient ${\mathbb T}'$ of $\mathbb T$ is induced by the inclusion of the respective full subcategory of ${\textrm{f.p.} {\mathbb T}\textrm{-mod}(\Set)}$ on the ${\mathbb T}'$-models. In view of the above characterizations, the topology $J$ is rigid. Conversely, if $J$ is rigid then $\Sh({\textrm{f.p.} {\mathbb T}\textrm{-mod}(\Set)}^{\textrm{op}}, J)$ is equivalent to the presheaf topos $[{\cal D}, \Set]$, where ${\cal D}^{\textrm{op}}$ is the full subcategory of ${\textrm{f.p.} {\mathbb T}\textrm{-mod}(\Set)^{\textrm{op}}}$ on the $J$-irreducible objects. Clearly, this subcategory is Cauchy-complete (as ${\textrm{f.p.} {\mathbb T}\textrm{-mod}(\Set)}$ is Cauchy-complete and the condition of $J$-irreducibility is stable under retracts) and hence, since ${\mathbb T}'$ is of presheaf type, it is equivalent to the category of finitely presentable ${\mathbb T}'$-models; in particular, any such model is finitely presentable as a $\mathbb T$-model.  
\end{proofs}

\begin{remark}\label{qpf}

\begin{enumerate}[(i)]	
\item Under the hypothesis that every finitely presentable ${\mathbb T}'$-model is finitely presentable as a $\mathbb T$-model, the full subcategory ${\textrm{f.p.} {{\mathbb T}'}\textrm{-mod}(\Set)}$ of ${\textrm{f.p.} {\mathbb T}\textrm{-mod}(\Set)}$ and the topology $J$ can be defined in terms of each other as follows (cf. the proof of the theorem): the objects of ${\textrm{f.p.} {{\mathbb T}'}\textrm{-mod}(\Set)}$ are precisely the $J$-irreducible objects of ${\textrm{f.p.} {{\mathbb T}}\textrm{-mod}(\Set)}$, and a sieve $S$ in ${\textrm{f.p.} {{\mathbb T}'}\textrm{-mod}(\Set)}$ on an object $c$ is $J$-covering if and only if it contains all the arrows in ${\textrm{f.p.} {\mathbb T}\textrm{-mod}(\Set)}$ from finitely presentable ${\mathbb T}'$-models to $c$.

\item If ${\mathbb T}'$ is a presheaf-type quotient of a theory of presheaf type ${\mathbb T}$ then for any finitely presentable $\mathbb T$-model $c$, any ${\mathbb T}'$-model $M$ and any $\mathbb T$-model homomorphism $f:c\to M$, there exists a $\mathbb T$-model homomorphism $g:c\to c'$ to a finitely presentable ${\mathbb T}'$-model $c'$ and a $\mathbb T$-model homomorphism $h:c' \to M$ such that $h\circ g=f$. Indeed, $M$ can be expressed as a filtered colimit of finitely presentable ${\mathbb T}'$-models $c'$; therefore, as $c$ is finitely presentable as a $\mathbb T$-model, $f$ necessarily factors through a colimit arrow $c'\to M$.

\item Under the hypothesis that every finitely presentable ${\mathbb T}'$-model is finitely presentable as a $\mathbb T$-model, by the duality theorem of \cite{OC6}, the syntactic description of ${\mathbb T}'$ in terms of $J$ given therein and the description of $J$ in terms of ${\textrm{f.p.} {{\mathbb T}'}\textrm{-mod}(\Set)}$ given in Remark \ref{qpf}(a), ${\mathbb T}'$ can be characterized as the quotient of $\mathbb T$ obtained by adding all the sequents of the form $( \phi(\vec{x}) \vdash_{\vec{x}} \mathbin{\mathop{\textrm{ $\bigvee$}}\limits_{i\in I}}\exists \vec{x_{i}}\theta_{i}(\vec{x_{i}}, \vec{x}))$, where $\{\theta_{i}(\vec{x_{i}}, \vec{x}):\phi_{i}(\vec{x_{i}}) \to \phi(\vec{x})\}$ is the family of $\mathbb T$-provably functional formulae from the formulae $\phi_{i}(\vec{x_{i}})$ presenting a ${\mathbb T}'$-model to a $\mathbb T$-irreducible formula (equivalently, a formula presenting a $\mathbb T$-model) $\phi(\vec{x})$.  
\end{enumerate} 
\end{remark}

The following result provides a natural class of presheaf-type quotients of presheaf-type theories whose associated topologies are rigid.

\begin{theorem}\label{closed}
Let $\mathbb T$ be a theory of presheaf type over a signature $\Sigma$. Then any quotient ${\mathbb T}'$ of $\mathbb T$ obtained from $\mathbb T$ by adding sequents of the form $\phi \vdash_{\vec{x}} \bot$, where $\phi(\vec{x})$ is a geometric formula over $\Sigma$, is classified by the topos $[{\cal T}, \Set]$, where $\cal T$ is the full subcategory of ${\textrm{f.p.} {\mathbb T}\textrm{-mod}(\Set)}$ on the ${\mathbb T}'$-models. 
\end{theorem}

\begin{proofs}
By covering each $\phi(\vec{x})$ by $\mathbb T$-irreducible formulae in the syntactic category ${\cal C}_{\mathbb T}$ of $\mathbb T$, we can suppose without loss of generality all of the $\phi(\vec{x})$ to be $\mathbb T$-irreducible, that is to present a $\mathbb T$-model $M_{\phi}$. Then, by the results of \cite{OC6}, the quotient ${\mathbb T}'$ is classified by the topos $\Sh({\textrm{f.p.} {\mathbb T}\textrm{-mod}(\Set)}^{\textrm{op}}, J)$, where $J$ is the smallest Grothendieck topology on ${\textrm{f.p.} {\mathbb T}\textrm{-mod}(\Set)}^{\textrm{op}}$ which contains all the empty sieves on the models $M_{\phi}$ presented by the formulae $\phi$ involved in the axioms $\phi \vdash_{\vec{x}} \bot$ added to $\mathbb T$ to form ${\mathbb T}'$. Let $\cal T$ be the full subcategory of ${\textrm{f.p.} {\mathbb T}\textrm{-mod}(\Set)}$ on the ${\mathbb T}'$-models. Then $\cal T$ is $J$-dense; indeed, any object not in $\cal T$ admits an arrow from a model of the form $M_{\phi}$ and hence is covered by the empty sieve. Further, for any object of $\cal T$, the $J$-covering sieves on it are exactly the maximal ones; therefore, the Comparison Lemma yields an equivalence $\Sh({\textrm{f.p.} {\mathbb T}\textrm{-mod}(\Set)}^{\textrm{op}}, J)\simeq [{\cal T}, \Set]$, as required. 
\end{proofs}

\subsubsection{Finding theories classified by a given presheaf topos}\label{findingpresheaf}

The following theorem provides a method for constructing theories of presheaf type whose categories of finitely presentable models are equivalent, up to Cauchy-completion, to a given small category of structures. 

\begin{theorem}\label{crit}
Let $\mathbb T$ be a theory of presheaf type and $\cal A$ a full subcategory of ${\textrm{f.p.} {\mathbb T}\textrm{-mod}(\Set)}$. Then the $\cal A$-completion ${\mathbb T}'$ of ${\mathbb T}$ (i.e., the set of all geometric sequents over the signature of $\mathbb T$ which are valid in all models in $\cal A$) is of presheaf type classified by the topos $[{\cal A}, \Set]$; in particular, every finitely presentable ${\mathbb T}'$-model is a retract of a model in $\cal A$.
\end{theorem}

\begin{proofs}
Since every model in $\cal A$ is finitely presentable as a $\mathbb T$-model, we have a geometric inclusion $i:[{\cal A}, \Set] \hookrightarrow \Set[{\mathbb T}]\simeq [\textrm{f.p.} {{\mathbb T}}\textrm{-mod}(\Set), \Set]$ induced by the canonical inclusion ${\cal A} \hookrightarrow \textrm{f.p.} {{\mathbb T}}\textrm{-mod}(\Set)$. This subtopos corresponds, by the duality theorem of \cite{OC6}, to a quotient ${\mathbb T}'$ of $\mathbb T$ classified by the topos $[{\cal A}, \Set]$, which can be characterized as the collection of all geometric sequents which hold in every model in $\cal A$, that is, as the $\cal A$-completion ${\mathbb T}'$ of ${\mathbb T}$ (recall that theories of presheaf type have enough set-based models). Therefore ${\mathbb T}'$ is of presheaf type classified by the topos $[{\cal A}, \Set]$; but the finitely presentable ${\mathbb T}'$-models are the finitely presentable objects of $\Ind{\cal A}$, that is the retracts of objects of $\cal A$ in $\Ind{\cal A}\simeq {{\mathbb T}'}\textrm{-mod}({\Set})$. This completes the proof of the theorem. 
\end{proofs}

\begin{remarks}
\begin{enumerate}[(a)]
\item Theorem \ref{crit} is a generalization of Joyal and Wraith's recognition theorem (Theorem 1.1 \cite{beke}), as well as of Proposition 5.3 \cite{diers}. Indeed, the former theorem can be obtained as the particular case of Theorem \ref{crit} when ${\mathbb T}'$ has enough set-based models and every set-based model of ${\mathbb T}'$ is a filtered colimit of finitely presentable models of ${\mathbb T}'$ which are also finitely presentable as $\mathbb T$-models, while Proposition 5.3 \cite{diers} can be obtained as the specialization of this latter theorem to the case when ${\mathbb T}'$ is a disjunctive theory and $\mathbb T$ is a cartesian theory (in view of the fact, proved in \cite{diers}, that the category of set-based models of a disjunctive theory is multiply finitely presentable and hence that every set-based model of a disjunctive theory can be expressed as a filtered colimit of finitely presentable models of it).

\item It is worth to compare Theorem \ref{crit} with Theorem \ref{main}. If ${\mathbb T}'$ is a quotient of a theory of presheaf type $\mathbb T$ and $\cal K$ is a small category of set-based ${\mathbb T}'$-models such that every model in $\cal K$ is finitely presentable as a $\mathbb T$-model then, in order to conclude that ${\mathbb T}'$ is of presheaf type classified by the topos $[{\cal K}, \Set]$, we can either verify conditions $(i)$ and $(ii)$ of Theorem \ref{main} or to verify that the models in $\cal K$ are jointly conservative for ${\mathbb T}'$; in fact, as observed in Remark \ref{conditionstar}(c), conditions $(i)$ and $(ii)$ of Theorem \ref{main} together imply that the models in $\cal K$ are jointly conservative for ${\mathbb T}'$.
\end{enumerate}
\end{remarks}

\begin{corollary}
Assuming the axiom of choice, every coherent theory (or, more generally, any theory in a countable fragment of geometric logic as in the hypotheses of Theorem 5.1.7 \cite{MR}) $\mathbb T$ over a finite relational signature whose axioms do not contain existential quantifications is of presheaf type.
\end{corollary}

\begin{proofs}
Since the axioms of $\mathbb T$ do not contain quantifications, every substructure of a model of $\mathbb T$ is a model of $\mathbb T$. Moreover, since the signature of $\mathbb T$ is relational, every finitely generated substructure over the signature of $\mathbb T$ contains only a finite number of elements besides the constants. Therefore every finitely presentable model of $\mathbb T$ contains only a finite number of elements besides the constants (as any model of $\mathbb T$ is a filtered union of its finite submodels); on the other hand, since the signature of $\mathbb T$ is finite, every such model is finitely presentable as a model of the empty theory over the  signature of $\mathbb T$ (cf. Theorem 6.4 \cite{OCG}). Therefore condition $(iii)$ of Theorem \ref{main} is satisfied (cf. Proposition \ref{criteriongeometricity}(i)). Now, under the axiom of choice, every theory satisfying the hypotheses of Theorem 5.1.7 \cite{MR}, in particular any coherent theory, has enough models, whence our thesis follows from Theorem \ref{crit}.   
\end{proofs}

It is nonetheless important, from a constructive viewpoint, to be able to prove that a theory is of presheaf type without invoking the axiom of choice. Theorem \ref{main} allows us to do so in a great variety of cases.

Theorem \ref{crit} shows that a good first step in constructing a geometric theory classified by a presheaf topos $[{\cal K}, \Set]$ consists in finding a theory of presheaf type $\mathbb T$ such that the category $\cal K$ can be identified as a full subcategory of the category of finitely presentable models of $\mathbb T$. Indeed, under these hypotheses, Theorem \ref{crit} ensures the existence of a quotient ${\mathbb T}_{\cal K}$ of $\mathbb T$ classified by the topos $[{\cal K}, \Set]$, which can be characterized as the theory consisting of all the geometric sequents over the signature of $\mathbb T$ which are valid in every model in $\cal K$. In most cases, if one has a natural candidate $\mathbb S$ for a theory classified by the topos $[{\cal K}, \Set]$, the theory $\mathbb T$ can be chosen to be the Horn part of $\mathbb S$ or the cartesianization ${\mathbb S}_{c}$ of $\mathbb S$. 
 
Of course, the abstract characterization of ${\mathbb T}_{\cal K}$ as the theory consisting of all the geometric sequents over the signature of $\mathbb T$ which are valid in every model in $\cal K$ is not very useful in specific contexts, where one looks for an axiomatization of ${\mathbb T}_{\cal K}$ as simple and `economical' as possible. To this end, we observe that if every $M$ in $\cal K$ is strongly finitely presented as a model of $\mathbb T$ as well as finitely generated (in the sense that for any sort $A$ over the signature of $\mathbb T$ the elements of the set $MA$ are precisely given by the interpretations in $M$ of terms $t^{A}(\vec{x})$ (or more generally of $\mathbb T$-provably functional predicates) over the signature of $\mathbb T$, where $\vec{x}$ are the generators of $M$ as strongly finitely presented model of $\mathbb T$) then we dispose of an explicit axiomatization of the theory ${\mathbb T}_{\cal K}$, as given by the following 

\begin{theorem}\label{extending}
Let $\mathbb T$ be a geometric theory over a signature $\Sigma$ and $\cal K$ a full subcategory of the category set-based $\mathbb T$-models such that every $\mathbb T$-model in $\cal K$ is both strongly finitely presentable and finitely generated (with respect to the same generators). Then the following sequents (where we denote by $\cal P$ the set of geometric formulae over $\Sigma$ which strongly present a $\mathbb T$-model in $\cal K$), added to the axioms of $\mathbb T$, yield an axiomatization of a quotient of $\mathbb T$ classified by the topos $[{\cal K}, \Set]$ via a Morita-equivalence induced by the canonical geometric morphism $[{\cal K}, \Set]\to \Sh({\cal C}_{\mathbb T}, J_{\mathbb T})$.  

In particular, if the theory $\mathbb T$ is of presheaf type (whence every finitely presentable $\mathbb T$-model is strongly finitely presentable) and all the models in the full subcategory $\cal K$ of ${\mathbb T}\textrm{-mod}(\Set)$ are finitely presented and finitely generated (with respect to the same generators), the following sequents yield an axiomatization of the theory ${\mathbb T}_{\cal K}$ defined above: 

\begin{enumerate}[(i)]
\item The sequent 
\[
(\top \vdash_{[]} \mathbin{\mathop{\textrm{ $\bigvee$}}\limits_{\phi(\vec{x})\in {\cal P}}} (\exists \vec{x})\phi(\vec{x}));
\]

\item For any formulae $\phi(\vec{x})$ and $\psi(\vec{y})$ in $\cal P$, where $\vec{x}=(x_{1}^{A_{1}}, \ldots, x_{n}^{A_{n}})$ and $\vec{y}=(y_{1}^{B_{1}}, \ldots, y_{m}^{B_{m}})$, the sequent  
\[
(\phi(\vec{x}) \wedge \psi(\vec{y}) \vdash_{\vec{x}, \vec{y}} \mathbin{\mathop{\textrm{ $\bigvee$}}\limits_{\substack{\chi(\vec{z})\in {\cal P}, t_{1}^{A_{1}}(\vec{z}), \ldots, t_{n}^{A_{n}}(\vec{z}) \\ s_{1}^{B_{1}}(\vec{z}), \ldots, s_{m}^{B_{m}}(\vec{z})} } (\exists \vec{z})(\chi(\vec{z}) \wedge \mathbin{\mathop{\textrm{ $\bigwedge$}}\limits_{\substack{i\in \{1, \ldots, n\}, \\ {j\in \{1, \ldots, m\}}}} (x_{i}=t_{i}(\vec{z}) \wedge y_{j}=s_{j}(\vec{z})))}}),
\]
where the disjunction is taken over all the formulae $\chi(\vec{z})$ in $\cal P$ and all the sequences of terms $t_{1}^{A_{1}}(\vec{z}), \ldots, t_{n}^{A_{n}}(\vec{z})$ and $s_{1}^{B_{1}}(\vec{z}), \ldots, s_{m}^{B_{m}}(\vec{z})$ whose output sorts are respectively $A_{1}, \ldots, A_{n}, B_{1}, \ldots, B_{m}$ and such that, denoting by $\vec{\xi}$ the set of generators of the model $M_{\{\vec{z}. \chi\}}$ (strongly) finitely presented by the formula $\chi(\vec{z})$, $(t_{1}^{A_{1}}(\vec{\xi}), \ldots, t_{n}^{A_{n}}(\vec{\xi}))\in [[\vec{x}. \phi]]_{M_{\{\vec{z}. \chi\}}}$ and $(s_{1}^{B_{1}}(\vec{\xi}), \ldots, s_{m}^{B_{m}}(\vec{\xi}))\in [[\vec{y}. \psi]]_{M_{\{\vec{z}. \chi\}}}$;

\item For any formulae $\phi(\vec{x})$ and $\psi(\vec{y})$ in $\cal P$, where $\vec{x}=(x_{1}^{A_{1}}, \ldots, x_{n}^{A_{n}})$ and $\vec{y}=(y_{1}^{B_{1}}, \ldots, y_{m}^{B_{m}})$, and any terms $t_{1}^{A_{1}}(\vec{y}), s_{1}^{A_{1}}(\vec{y}), \ldots, t_{n}^{A_{n}}(\vec{y}), s_{n}^{A_{n}}(\vec{y})$ whose output sorts are respectively $A_{1}, \ldots, A_{n}$, the sequent 

\begin{equation*}
\begin{split}
(\mathbin{\mathop{\textrm{ $\bigwedge$}}\limits_{i\in \{1, \ldots, n\}} (t_{i}(\vec{y})=s_{i}(\vec{y}))} \wedge \phi(t_{1}\slash x_{1}, \ldots, t_{n}\slash x_{n}) \wedge \phi(s_{1}\slash x_{1}, \ldots, s_{n}\slash x_{n}) \wedge \psi(\vec{y}) \\
\vdash_{\vec{y}} \mathbin{\mathop{\textrm{ $\bigvee$}}\limits_{\chi(\vec{z})\in {\cal P}, u_{1}^{B_{1}}(\vec{z}), \ldots, u_{m}^{B_{m}}(\vec{z})}} ((\exists \vec{z})(\chi(\vec{z}) \wedge \mathbin{\mathop{\textrm{ $\bigwedge$}}\limits_{j\in \{1, \ldots, m\}}} (y_{j}=u_{j}(\vec{z}))),
\end{split}
\end{equation*}
where the disjunction is taken over all the formulae $\chi(\vec{z})$ in $\cal P$ and all the sequences of terms $u_{1}^{B_{1}}(\vec{z}), \ldots, u_{m}^{B_{m}}(\vec{z})$ whose output sorts are respectively $B_{1}, \ldots, B_{m}$ and such that, denoting by $\vec{\xi}$ the set of generators of the model $M_{\{\vec{z}. \chi\}}$ (strongly) finitely presented by the formula $\chi(\vec{z})$, $(u_{1}^{B_{1}}(\vec{\xi}), \ldots, u_{m}^{B_{m}}(\vec{\xi}))\in [[\vec{y}. \psi]]_{M_{\{\vec{z}. \chi\}}}$ and $t_{i}(u_{1}(\vec{\xi}), \ldots, u_{m}(\vec{\xi}))=s_{i}(u_{1}(\vec{\xi}), \ldots, u_{m}(\vec{\xi}))$ in $M_{\{\vec{z}. \chi\}}$ for all $i\in \{1, \ldots, n\}$;

\item For any sort $A$ over $\Sigma$, the sequent 
\[
(\top \vdash_{x_{A}} \mathbin{\mathop{\textrm{ $\bigvee$}}\limits_{\chi(\vec{z})\in {\cal P}, t^{A}(\vec{z})}} (\exists \vec{z})(\chi(\vec{z}) \wedge x=t(\vec{z}))),
\]
where the the disjunction is taken over all the formulae $\chi(\vec{z})$ in $\cal P$ and all the terms $t^{A}(\vec{z})$  whose output sort is $A$;

\item For any sort $A$ over $\Sigma$, any formulae $\phi(\vec{x})$ and $\psi(\vec{y})$ in $\cal P$, where $\vec{x}=(x_{1}^{A_{1}}, \ldots, x_{n}^{A_{n}})$ and $\vec{y}=(y_{1}^{B_{1}}, \ldots, y_{m}^{B_{m}})$, and any terms $t^{A}(\vec{x})$ and $s^{A}(\vec{y})$, the sequent 

\begin{equation*}
\begin{split}
(\phi(\vec{x}) \wedge \psi(\vec{y}) \wedge t(\vec{x})=s(\vec{y})
\vdash_{\vec{x}, \vec{y}}
\mathbin{\mathop{\textrm{ $\bigvee$}}\limits_{\substack{\chi(\vec{z})\in {\cal P}, p_{1}^{A_{1}}(\vec{z}), \ldots, p_{n}^{A_{n}}(\vec{z}) \\ q_{1}^{B_{1}}(\vec{z}), \ldots, q_{m}^{B_{m}}(\vec{z})} } (\exists \vec{z})(\chi(\vec{z})} \wedge \\ \wedge \mathbin{\mathop{\textrm{ $\bigwedge$}}\limits_{\substack{i\in \{1, \ldots, n\},\\ {j\in \{1, \ldots, m\}}}} (x_{i}=p_{i}(\vec{z}) \wedge y_{j}=q_{j}(\vec{z}))) },
\end{split}
\end{equation*}
where the disjunction is taken over all the formulae $\chi(\vec{z})$ in $\cal P$ and all the sequences of terms $p_{1}^{A_{1}}(\vec{z}), \ldots, p_{n}^{A_{n}}$ and $q_{1}^{B_{1}}(\vec{z}), \ldots, q_{m}^{B_{m}}(\vec{z})$ whose output sorts are respectively $A_{1}, \ldots, A_{n}, B_{1}, \ldots, B_{m}$ and such that, denoting by $\vec{\xi}$ the set of generators of the model $M_{\{\vec{z}. \chi\}}$ (strongly) finitely presented by the formula $\chi(\vec{z})$, $(p_{1}^{A_{1}}(\vec{\xi}), \ldots, p_{n}^{A_{n}}(\vec{\xi}))\in [[\vec{x}. \phi]]_{M_{\{\vec{z}. \chi\}}}$ and $(q_{1}^{B_{1}}(\vec{\xi}), \ldots, q_{m}^{B_{m}}(\vec{\xi}))\in [[\vec{y}. \psi]]_{M_{\{\vec{z}. \chi\}}}$ and $t(p_{1}(\vec{\xi}), \ldots, p_{n}(\vec{\xi}))=s(q_{1}(\vec{\xi}), \ldots, q_{m}(\vec{\xi}))$ in $M_{\{\vec{z}. \chi\}}$.
\end{enumerate} 
\end{theorem}

\begin{proofs}
Let $\mathbb R$ be the geometric theory obtained from $\mathbb T$ by adding all the sequents specified above. The objects in $\cal K$ are clearly models of $\mathbb R$ (cf. Remarks \ref{remfg1}(a) and \ref{remfg2}). From the fact that they are strongly finitely presentable it follows by Proposition \ref{criteriongeometricity}(ii) that condition $(iii)$ of Theorem \ref{main} is satisfied by the theory $\mathbb R$ with respect to the category $\cal K$. By Theorems \ref{thcond1} and \ref{cond2fg} and Remarks \ref{remfg1}(a) and \ref{remfg2}, $\mathbb R$ also satisfies conditions $(i)$ and $(ii)$ of Theorem \ref{main} with respect to $\cal K$. Theorem \ref{main} thus implies that $\mathbb R$ is of presheaf type classified by the topos $[{\cal K}, \Set]$. Moreover, it can be readily seen that the resulting Morita-equivalence is induced by the canonical geometric morphism from $[{\cal K}, \Set]$ to the classifying topos for $\mathbb T$; this ensures in particular that if $\mathbb T$ is of presheaf type then $\mathbb R$ is equal to the theory ${\mathbb T}_{\cal K}$.     
\end{proofs}

\subsection{Injectivizations of theories}\label{inj}

For any geometric theory $\mathbb T$, we can slightly modify its syntax so to obtain a geometric theory whose models in $\Set$ are the same as those of $\mathbb T$ and whose homomorphisms between them are precisely the sortwise injective $\mathbb T$-model homomorphisms.

This construction is useful in many contexts. For instance, in \cite{OC6} we showed that if the category $\textrm{f.p.} {\mathbb T}\textrm{-mod}(\Set)$ of finitely presentable models of a theory of presheaf type $\mathbb T$ satisfies the amalgamation property then the quotient of $\mathbb T$ corresponding to the subtopos $\Sh({\textrm{f.p.} {\mathbb T}\textrm{-mod}(\Set)}^{\textrm{op}}, J_{at})$ of $[\textrm{f.p.} {\mathbb T}\textrm{-mod}(\Set), \Set]$ (where $J_{at}$ is the atomic topology) via the duality theorem of \cite{OC6} axiomatizes the homogeneous $\mathbb T$-models (in the sense of \cite{OC}) in any Grothendieck topos. Now, the notion of homogeneous $\mathbb T$-model, which is strictly related to the notion of weakly homogeneous model considered in classical Model Theory (cf. \cite{Hodges}), is mostly interesting when the arrows of the category ${\mathbb T}\textrm{-mod}(\Set)$ are all monic; indeed, as shown in \cite{OCG}, a necessary condition for $\mathbb T$ to admit an associated `concrete' Galois theory, is that all the arrows in ${\textrm{f.p.} {\mathbb T}\textrm{-mod}(\Set)}$ should be strict monomorphisms.  

This motivates the following formal definition.

\begin{definition}\label{injectiviz}
Let $\mathbb T$ be a geometric theory over a signature $\Sigma$. The \emph{injectivization} ${\mathbb T}_{m}$ of $\mathbb T$ is the geometric theory obtained from $\mathbb T$ by adding a binary predicate $D_{A}\mono A, A$ for each sort $A$ over $\Sigma$ and the coherent sequents
\[
(D_{A}(x^{A}, y^{A}) \wedge x^{A}=y^{A}) \vdash_{x^{A}, y^{A}} \bot)
\]
and 
\[
(\top \vdash_{x^{A}, y^{A}} D_{A}(x^{A}, y^{A}) \vee x^{A}=y^{A}).
\]
\end{definition}

The models of ${\mathbb T}_{m}$ in an arbitrary topos $\cal E$ coincide with the models $M$ of $\mathbb T$ in $\cal E$ which are sortwise decidable, in the sense that for every sort $A$ over the signature of $\mathbb T$ the object $MA$ of $\cal E$ is \emph{decidable}, i.e. the diagonal subobject of $MA$ is complemented (by the interpretation of $D_{A}$ in $M$).

As shown by the following lemma, the arrows $M\to N$ in the category ${\mathbb T}_{m}\textrm{-mod}(\Set)$ are precisely the $\mathbb T$-model homomorphisms $f:M\to N$ such that for every sort $A$ over the signature of $\mathbb T$, $fA:MA\to NA$ is a monomorphism in $\cal E$.

\begin{lemma}\label{mondec}
Let $A$ and $B$ be decidable objects in a topos $\cal E$ and $f:A\to B$ an arrow in $\cal E$. Let $D_{A}\mono A\times A$ and $D_{B}\mono B\times B$ denote respectively the complements of the diagonal subobjects $\delta_{A}:A\mono A\times A$ and $\delta_{B}:B\mono B\times B$. Then $f$ is a monomorphism if and only if $f\times f:A\times A \to B\times B$ restricts to an arrow $D_{A}\to D_{B}$. 
\end{lemma}    

\begin{proofs}
It is immediate to see that $f:A\to B$ is a monomorphism if and only if the diagram
\[  
\xymatrix {
A \ar[d]^{\delta_{A}} \ar[r]^{f} & B \ar[d]^{\delta_{B}}\\
A\times A \ar[r]_{f\times f} & B\times B, } 
\] 
is a pullback. 

Since pullback functors preserve arbitrary unions and intersections of subobjects in a topos (they having both a left and a right adjoint), we have that $(f\times f)^{\ast}(D_{B})\cong \neg (f\times f)^{\ast}(\delta_{B})$. Now, $f\times f:A\times A \to B\times B$ restricts to an arrow $D_{A}\to D_{B}$ if and only if $D_{A}\leq (f\times f)^{\ast}(D_{B})$. But this condition holds if and only if $D_{A}\cap (f\times f)^{\ast}(\delta_{B})\cong 0$, i.e. if and only if $(f\times f)^{\ast}(\delta_{B})\leq \delta_{A}$, which is equivalent to the condition $(f\times f)^{\ast}(\delta_{B})\cong \delta_{A}$ (as $\delta_{A} \leq (f\times f)^{\ast}(\delta_{B})$).
\end{proofs}

Several injectivizations of theories of presheaf type have been considered in the literature, e.g. in \cite{MM} and in \cite{El}; see also \cite{OC2} and \cite{OCG} for some applications of this type of theories in the context of topos-theoretic Galois-type equivalences.    

As we shall see below, under certain conditions, the injectivization of a theory of presheaf type is again of presheaf type.

The following proposition is a corollary of Theorem \ref{objectmodelhomo}.
 
\begin{proposition}\label{objectmonomorphisms}
Let $\mathbb T$ be a geometric theory. Then for any $\mathbb T$-models $M$ and $N$ in a Grothendieck topos $\cal E$ which are sortwise decidable there exists an object $Hom_{\underline{{\mathbb T}_{m}\textrm{-mod}(\cal E)}}^{\cal E}(M,N)$ of $\cal E$ satisfying the following universal property: for any object $E$ of $\cal E$ we have an equivalence
\[
Hom_{\cal E}(E, Hom_{\underline{{\mathbb T}_{m}\textrm{-mod}(\cal E)}}^{\cal E}(M,N))\cong Hom_{{\mathbb T}_{m}\textrm{-mod}({\cal E}\slash E)}(!_{E}^{\ast}(M), !_{E}^{\ast}(N))
\]
natural in $E\in {\cal E}$. 
\end{proposition}
 
\begin{remarks}
\begin{enumerate}[(a)]
\item For any $\cal E$, $M$ and $N$, $Hom_{\underline{{\mathbb T}_{m}\textrm{-mod}(\cal E)}}^{\cal E}(M,N)$ embeds canonically as a subobject of $Hom_{\underline{{\mathbb T}\textrm{-mod}(\cal E)}}^{\cal E}(M,N)$.

\item Let $\mathbb T$ be a geometric theory over a signature $\Sigma$, $c$ a finitely presentable $\mathbb T$-model and $M$ a sortwise decidable model of $\mathbb T$ in a Grothendieck topos $\cal E$. Then the subobject 
\[
Hom_{\underline{{\mathbb T}_{m}\textrm{-mod}(\cal E)}}^{\cal E}(\gamma_{\cal E}^{\ast}(c),M)\mono Hom_{\underline{{\mathbb T}\textrm{-mod}(\cal E)}}^{\cal E}(\gamma_{\cal E}^{\ast}(c),M)
\]
can be identified with the interpretation of the formula 
\[
\chi_{c}:=\mathbin{\mathop{\textrm{$\bigwedge$}}\limits_{\substack{{A \textrm{ sort over } \Sigma,}\\ {x,y\in cA, \textrm{ } x\neq y }}}} (\pi_{A}(f(\gamma_{\cal E}^{\ast}(\overline{x})), \pi_{A}(f(\gamma_{\cal E}^{\ast}(\overline{y})))\in D_{MA},
\]  
written in the internal language of the topos $\cal E$, where 
\[
\pi_{A}:Hom_{\underline{{\mathbb T}\textrm{-mod}(\cal E)}}^{\cal E}(M,N) \to NA^{MA}
\]
is the arrow defined in Remark \ref{functoriality}(b) and $\overline{x}, \overline{y}:1 \to c$ are the arrows in $\Set$ corresponding respectively to the elements $x$ and $y$ of $cA$. Indeed, a $\mathbb T$-model homomorphism $f:\gamma_{\cal E}^{\ast}(c) \to M$ is sortwise monic if and only if for every sort $A$ over $\Sigma$ and any distinct elements $x, y\in cA$, $fA(\gamma_{\cal E}^{\ast}(\overline{x}))$ and $ fA(\gamma_{\cal E}^{\ast}(\overline{y}))$ are disjoint, in other words they satisfy the relation $D_{MA}$.

It follows that if for every finitely presentable $\mathbb T$-model $c$, the formula $\chi_{c}$ is equivalent to a geometric formula over the signature of ${\mathbb T}_{m}$, for instance when $\Sigma$ only contains a finite number of sorts and for any sort $A$ and any finitely presentable $\mathbb T$-model the set $c_{A}$ is finite, then the theory ${\mathbb T}_{m}$ satisfies condition $(iii)(a)$ of Theorem \ref{main} with respect to the category of finitely presentable $\mathbb T$-models if $\mathbb T$ does. 

\item Let $\mathbb T$ be a geometric theory over a signature $\Sigma$, $c$ a finitely presentable $\mathbb T$-model and $M$ a sortwise decidable model of $\mathbb T$ in a Grothendieck topos $\cal E$. If $c$ is strongly finitely presented by a geometric formula $\phi(\vec{x})$ over $\Sigma$ then the subobject $Hom_{\underline{{\mathbb T}\textrm{-mod}_{m}(\cal E)}}^{\cal E}(\gamma_{\cal E}^{\ast}(c),M)$ of $Hom_{\underline{{\mathbb T}\textrm{-mod}(\cal E)}}^{\cal E}(\gamma_{\cal E}^{\ast}(c),M)\cong [[\vec{x}. \phi]]_{M}$ can be identified with the intersection of the interpretations in $M$ of the formulae of the form $\neg (\exists x^{A})(\theta_{1}(\vec{x}, x^{A}) \wedge \theta_{2}(\vec{x}, x^{A}))$, where $A$ is a sort over $\Sigma$ and $\theta_{1}, \theta_{2}$ are ${\mathbb T}$-provably functional formulae from $\{\vec{x}.\phi\}$ to $\{x^{A}. \top\}$. 
\end{enumerate}
\end{remarks}

\begin{lemma}\label{explicit2}
Let $\mathbb T$ be a geometric theory over a signature $\Sigma$, $c$ a set-based $\mathbb T$-model and $M$ a $\mathbb T$-model in a Grothendieck topos $\cal E$. A generalized element $x:E\to G_{M}(c)=Hom_{\underline{{\mathbb T}_{m}\textrm{-mod}(\cal E)}}^{\cal E}(\gamma_{\cal E}^{\ast}(c),M)$ can be identified with a $\Sigma$-homomorphism $\xi_{x}:c\to Hom_{\cal E}(E, M)$ which is sortwise disjunctive in the sense that for every sort $A$ over $\Sigma$, the function $\xi_{x}A:cA\to Hom_{\cal E}(E, MA)$ has the property that for any distinct elements $z,w\in cA$, the arrows $\xi_{x}A(z), \xi_{x}A(w):E\to MA$ have equalizer zero in $\cal E$.  
\end{lemma}

\begin{proofs}
The thesis follows from Proposition \ref{explicit}, observing that, by Lemmas \ref{remarklemmamono} and \ref{mono}, an arrow $\tau_{A}:{\gamma}_{{\cal E}\slash E}^{\ast}(cA)\to !_{E}^{\ast}(MA)$ in ${\cal E}\slash E$ is monic in ${\cal E}\slash E$ if and only if the corresponding arrow $\xi_{A}:cA\to Hom_{\cal E}(E, MA)$ satisfies the property that for any distinct elements $z,w\in cA$, the arrows $\xi_{A}(z):E\to MA$ and $\xi_{A}(w):E\to MA$ have equalizer zero (notice that two arrows $s,t:(a:A\to E)\to (b:B\to E)$ in ${\cal E}\slash E$ have equalizer zero in ${\cal E}\slash E$ if and only if the arrows $s:A\to B$ and $t:A\to B$ have equalizer zero in $\cal E$).
\end{proofs}

\begin{remark}
Suppose that the formula $\{x^{A}. \top\}$ strongly presents a $\mathbb T$-model $F_{A}$. Then a $\Sigma$-structure homomorphism $s:F_{A}\to Hom_{\cal E}(E, M)$, corresponding to a generalized element $z:E\to MA$, is sortwise disjunctive if and only if for any sort $B$ over $\Sigma$ and any two $\mathbb T$-provably inequivalent $\mathbb T$-provably functional geometric formulae $\theta_{1}$ and $\theta_{2}$ from $\{x^{A}. \top\}$ to $\{x^{B}. \top\}$, the generalized elements $[[\theta_{1}]]_{M}\circ z$ and $[[\theta_{2}]]_{M}\circ z$ are disjoint. Indeed, from the proof of Theorem \ref{syn} we know that that for any sort $B$ and $\mathbb T$-provably functional geometric formula $\theta$ from $\{x^{A}. \top\}$ to $\{x^{B}. \top\}$, the function $[[\theta]]_{F_{A}}:F_{A}A=Hom_{{\cal C}_{\mathbb T}}(\{x^{A}. \top\}, \{x^{A}. \top\}) \to F_{A}B= Hom_{{\cal C}_{\mathbb T}}(\{x^{A}. \top\}, \{x^{B}. \top\})$ coincides with $[\theta]\circ -:Hom_{{\cal C}_{\mathbb T}}(\{x^{A}. \top\}, \{x^{A}. \top\}) \to Hom_{{\cal C}_{\mathbb T}}(\{x^{A}. \top\}, \{x^{B}. \top\})$, while the generator $u_{A}$ of $FA$ is precisely the identity arrow on $\{x^{A}. \top\}$. Now, it is immediate to see that for any $\Sigma$-structure homomorphism $s:F_{A}\to Hom_{\cal E}(E, M)$ and any sort $B$ over $\Sigma$, the diagram 
\[  
\xymatrix {
F_{A}A \ar[d]^{[\theta]\circ -} \ar[rr]^{sA} & & Hom_{{\cal E}}(E, MA) \ar[d]^{[[\theta]]_{M} \circ -}\\
F_{A}B \ar[rr]^{sB} & & Hom_{{\cal E}}(E, MB)} 
\]
commutes, i.e. $sA(u_{A})=sB(\theta)$. From these remarks, our claim immediately follows. 
 
Note that if $\mathbb T$ is a universal Horn theory (in the sense of \cite{blasce}) then we can suppose without loss of generality $\theta_{1}$ and $\theta_{2}$ to be functional formulae of the form $x^{B}=t(x^{A})$, where $t$ is a term of type $A\to B$ over $\Sigma$ (cf. Remark \ref{remblasce}(a)). 
\end{remark}

\subsubsection{Condition (iii) of Theorem \ref{main} for injectivizations}

In this section we shall establish a general result about injectivizations of theories of presheaf type, namely that for any theory of presheaf type $\mathbb T$ such that all the monic arrows in the category ${{\textrm{f.p.} {\mathbb T}\textrm{-mod}(\Set)}}$ are sortwise monic, the injectivization of $\mathbb T$ satisfies condition $(iii)$ of Theorem \ref{main} with respect to the category ${{\textrm{f.p.} {\mathbb T}\textrm{-mod}(\Set)}}$.

Before proving this theorem, we need a series of preliminary results.

\begin{lemma}\label{factorization}
Let $r:R\mono A\times A$ a subobject in a Grothendieck topos $\cal E$, $e:E\to A\times A$ an arrow in $\cal E$ and $\{f_{i}:E_{i}\to E \textrm{ | } i\in I\}$ an epimorphic family in $\cal E$. If for every $i\in I$ the arrow $e\circ f_{i}$ factors through $r$ then $e$ factors through $r$.
\end{lemma}

\begin{proofs}
The arrow $\mathbin{\mathop{\textrm{ $\coprod$}}\limits_{i\in I}}f_{i}:\mathbin{\mathop{\textrm{ $\coprod$}}\limits_{i\in I}}E_{i}\to E$ induced by the universal property of the coproduct is an epimorphism, since by our hypothesis the family $\{f_{i}:E_{i}\to E \textrm{ | } i\in I\}$ is epimorphic. The factorizations $b_{i}:E_{i}\to R$ of the arrows $e\circ f_{i}$ through $r$ induce an arrow $b:=\mathbin{\mathop{\textrm{ $\coprod$}}\limits_{i\in I}}b_{i}:\mathbin{\mathop{\textrm{$\coprod$}}\limits_{i\in I}}E_{i} \to R$, such that $r\circ b= e\circ \mathbin{\mathop{\textrm{ $\coprod$}}\limits_{i\in I}}f_{i}$. Consider the epi-mono factorization $b=h\circ k$ in $\cal E$ of the arrow $b$, where $k:\mathbin{\mathop{\textrm{ $\coprod$}}\limits_{i\in I}}E_{i} \epi U$ and $h:U\mono R$, and the epi-mono factorization $e=m\circ n$ of the arrow $e$ in $\cal E$, where $n:E\epi T$ and $m:T\mono A$. Clearly, the arrow $e\circ \mathbin{\mathop{\textrm{ $\coprod$}}\limits_{i\in I}}f_{i}$ factorizes both as $m\circ (n \circ \mathbin{\mathop{\textrm{ $\coprod$}}\limits_{i\in I}}f_{i})$ and as $(r\circ h) \circ k$. Now, as $m$, $r\circ h$ are monic and $n \circ \mathbin{\mathop{\textrm{ $\coprod$}}\limits_{i\in I}}f_{i}$, $k$ are epic, the uniqueness of the epi-mono factorization of a given arrow in a topos implies that there exists an isomorphism $i:T\cong U$ such that $r\circ h\circ i=m$ and $i\circ n \circ \mathbin{\mathop{\textrm{ $\coprod$}}\limits_{i\in I}}f_{i}=k$. The arrow $h\circ i \circ n$ thus provides a factorization of $e$ through $r$, as required.
\end{proofs}

\begin{lemma}\label{disjointness}
Let $\{f_{i}:A_{i}\to B \textrm{ | } i\in I\}$ be a family of arrows in a Grothendieck topos $\cal E$. Then the arrow $\mathbin{\mathop{\textrm{ $\coprod$}}\limits_{i\in I}}f_{i} \to \mathbin{\mathop{\textrm{ $\coprod$}}\limits_{i\in I}} A_{i} \to B$ is monic if and only if for every $i\in I$, $f_{i}$ is monic and for every $i, i'\in I$, either $i=i'$ or the subobjects $f_{i}:A_{i}\mono B$ and $f_{i}':A_{i'}\mono B$ are disjoint.
\end{lemma}

\begin{proofs}
The `if' direction follows from Proposition IV 7.6 \cite{MM}, while the `only if' one follows from the fact that coproducts in a topos are always disjoint (cf. Corollary IV 10.5 \cite{MM}) and the composite of a given monomorphism with two disjoint subobjects yields two disjoint subobjects.
\end{proofs}

\begin{remark}\label{remarklemmamono}
If in the statement of Lemma \ref{disjointness} the objects $A_{i}$ are all equal to the terminal object $1_{\cal E}$ of $\cal E$ then the arrows $f_{i}$ are automatically monic and any two of them are disjoint if and only if their equalizer is zero.
\end{remark}

The notation employed in the statements and proofs of the following results is borrowed from section \ref{seccov}.

\begin{lemma}\label{factel}
Let $F:{\cal D}^{\textrm{op}}\to {\cal E}$ be a flat functor from a subcategory $\cal D$ of a small category $\cal C$ to a Grothendieck topos $\cal E$ and $x:E\to \tilde{F}(d)$ a generalized element which factors through $\chi^{F}_{d}:F(d)\to \tilde{F}(d)$ as $\chi^{F}_{d} \circ x'$. Then the natural transformation $\alpha_{x}:\gamma_{{\cal E}\slash E}^{\ast} \circ y_{\cal C}d \to !_{E}^{\ast} \circ \tilde{F}$ corresponding to $x$ is equal to $\tilde{\alpha_{x'}}$, where $\alpha_{x'}$ is the natural transformation $\gamma_{{\cal E}\slash E}^{\ast} \circ y_{\cal D}d \to !_{E}^{\ast} \circ F$ corresponding to $x'$.
\end{lemma}

\begin{proofs}
Straightforward from the results of section \ref{seccov}.
\end{proofs}

Given a small category $\cal C$, we shall write ${\cal C}_{m}$ for the category whose objects are the objects of $\cal C$ and whose arrows are the monic arrows in $\cal C$ between them. The following results concern extensions $\tilde{F}$ of flat functors $F$ along the embedding ${\cal C}_{m}\hookrightarrow {\cal C}$.

Below, for a given flat functor $H:{\cal C}^{\textrm{op}}\to {\cal E}$ with values in a Grothendieck topos $\cal E$ and any object $E$ of $\cal E$, we write $H_{E}$ for the flat functor $!_{E}^{\ast}\circ H:{\cal C}^{\textrm{op}}\to {\cal E}\slash E$.

\begin{proposition}\label{mono}
Let $\cal C$ be a small category and and $F:{{\cal C}_{m}}^{\textrm{op}}\to {\cal E}$ a flat functor. Then for any object $c\in {\cal C}$ and any generalized element $x:E\to \tilde{F}(c)$, $x$ factors through $\chi^{F}_{c}:F(c)\to \tilde{F}(c)$ if and only if the corresponding natural transformation $\gamma_{{\cal E}\slash E}^{\ast} \circ y_{\cal C}c \to \tilde{!_{E}^{\ast} \circ F}$ is pointwise monic.  
\end{proposition}

\begin{proofs}
The `only if' direction follows at once from Lemma \ref{factel} and Theorem \ref{monic}. To prove the `if' one, thanks to the localization technique, we can suppose without loss of generality $E=1_{\cal E}$.  

Suppose that the natural transformation $\alpha_{x}:\gamma_{{\cal E}}^{\ast} \circ y_{\cal C}c \to \tilde{F}$ corresponding to the generalized element $x:1\to \tilde{F}(c)$ is pointwise monic. We want to prove that $x$ factors through $\chi^{F}_{c}:F(c)\to \tilde{F}(c)$. Recall that $\alpha_{x}(d):   \mathbin{\mathop{\textrm{ $\coprod$}}\limits_{f\in Hom_{\cal C}(d, c)}}1_{\cal E}\to \tilde{F}(c)$ is defined as the arrow which sends the component of the coproduct indexed by $f$ to the generalized element $\tilde{F}(f)\circ x:1\to \tilde{F}(d)$. 

By Lemma \ref{disjointness}, if $\alpha_{x}$ is pointwise monic then for any object $d\in {\cal C}$ and arrows $f,g:d\to c$ in $\cal C$, either $f=g$ or the equalizer of $\tilde{F}(f)\circ x$ and $\tilde{F}(g)\circ x$ is zero. Consider the pullbacks of the jointly epimorphic arrows $\kappa_{(a, z)}:F(a)\to \tilde{F}(c)$ along the arrow $x:1_{\cal E}\to \tilde{F}(c)$:
\[  
\xymatrix {
E_{(a, z)} \ar[rr]^{h_{(a, z)}} \ar[d]^{e_{(a, z)}} & &  F(a) \ar[d]^{\kappa_{(a, z)}} \\
1  \ar[rr]_{x} & & \tilde{F}(c)} 
\]      
We shall prove that the epimorphic family $\{e_{(a, z)}:E_{(a, z)}\to 1_{\cal E} \textrm{ | } (a,z)\in {\cal A}_{c}\}$ satisfies the condition that for any $(a, z)\in {\cal A}_{c}$, the composite arrow $x\circ e_{(a, z)}$ factors through $\chi^{F}_{c}:F(c)\to \tilde{F}(c)$; this will imply our thesis by Lemma \ref{factorization}. As if $E_{(a, z)}\cong 0$ then $x\circ e_{(a, z)}$ factors through $\chi^{F}_{c}:F(c)\to \tilde{F}(c)$, we can suppose $E_{(a, z)}\ncong 0$. Under this hypothesis, the arrow $z:c\to a$ is monic in $\cal C$. Indeed, for any two arrows $f,g:b\to c$ such that $z\circ f=z\circ g$, $\tilde{F}(f)\circ x \circ e_{(a, z)}=\tilde{F}(f)\circ \kappa_{(a, z)}\circ h_{(a, z)}=\tilde{F}(f) \circ \tilde{F}(z) \circ \chi^{F}_{a}\circ h_{(a, z)}=\tilde{F}(g) \circ \tilde{F}(z) \circ \chi^{F}_{a}\circ h_{(a, z)}=\tilde{F}(g)\circ \kappa_{(a, z)}\circ h_{(a, z)}=\tilde{F}(g)\circ x \circ e_{(a, z)}$, whence either $f=g$ or the the equalizer of $\tilde{F}(f)\circ x$ and $\tilde{F}(g)\circ x$ is zero. But if the equalizer of $\tilde{F}(f)\circ x$ and $\tilde{F}(g)\circ x$ were zero then $E_{(a, z)}$ would also be isomorphic to zero (since it would admit an arrow to zero), contrary to our hypothesis; so $f=g$, as required.    
\end{proofs}

\begin{remark}
By Proposition \ref{objectmonomorphisms}, Proposition \ref{mono} can be reformulated as follows: for any object $c\in {\cal C}$, $F(c)\cong Hom_{\underline{{{\mathbb T}_{\cal C}}_{m}\textrm{-mod}(\cal E)}}^{\cal E}(\gamma_{{\cal E}}^{\ast} \circ y_{\cal C}c, \tilde{F})$, where ${\mathbb T}_{{\cal C}_{m}}$ is the geometric theory of flat functors on ${\cal C}_{m}^{\textrm{op}}$.  
\end{remark}

\begin{theorem}\label{monoflat}
Let $\cal C$ be a small category and $\textbf{Flat}_{m}({\cal C}^{\textrm{op}}, {\cal E})$ the subcategory of $\textbf{Flat}({\cal C}^{\textrm{op}}, {\cal E})$ whose objects are the same as $\textbf{Flat}({\cal C}^{\textrm{op}}, {\cal E})$ and whose arrows are the natural transformations between them which are pointwise monic. Then the extension functor 
\[
\xi_{{\cal E}}:\textbf{Flat}({{\cal C}_{m}}^{\textrm{op}}, {\cal E})\to \textbf{Flat}({\cal C}^{\textrm{op}}, {\cal E}),
\]
along the embedding ${\cal C}_{m}\hookrightarrow {\cal C}$, which by Theorem \ref{monic} takes values into $\textbf{Flat}_{m}({\cal C}^{\textrm{op}}, {\cal E})$, is full on this latter category.    
\end{theorem}

\begin{proofs}
Let $F,G$ be flat functors ${{\cal C}_{m}}^{\textrm{op}}\to {\cal E}$ with values in a Grothendieck topos $\cal E$, and let $\beta:\tilde{F}\to \tilde{G}$ be a pointwise monic natural transformation between them. We want to prove that there exists a natural transformation $\alpha:F\to G$ such that $\beta=\tilde{\alpha}$. It suffices to show that for any $c\in {\cal C}$, $\beta(c):\tilde{F}(c)\to \tilde{G}(c)$ restricts (along the arrows $\chi^{F}_{c}:F(c)\to \tilde{F}(c)$ and $\chi^{G}_{c}:G(c)\to \tilde{G}(c)$) to an arrow $F(c)\to G(c)$. To this end, we define a function $\gamma_{E}:Hom_{\cal E}(E, F(c))\to Hom_{\cal E}(E, G(c))$ natural in $E\in {\cal E}$. By Remark \ref{cov}(iii), the set $Hom_{\cal E}(E, F(c))$ (resp. the set $Hom_{\cal E}(E, G(c))$) can be identified with the set of arrows $E\to \tilde{F}(c)$ (resp. with the set of arrows $E\to \tilde{G}(c)$) which factor through $\chi^{F}_{c}$ (resp. through $\chi^{G}_{c}$). By Proposition \ref{mono}, for any flat functor $H:{\cal C}^{\textrm{op}}\to {\cal E}$, the arrows $E \to \tilde{H}(c)$ which factor through $\chi_{c}^{H}$ correspond precisely to the natural transformations $\gamma_{{\cal E}\slash E}^{\ast} \circ y_{\cal C}c \to \tilde{H_{E}}$ which are pointwise monic. Now, since $\beta$ is pointwise monic then $!_{E}^{\ast} \circ \beta:\tilde{F_{E}}\to \tilde{G_{E}}$ is also pointwise monic (since the functor $!_{E}^{\ast}$ preserves monomorphisms, it being the inverse image of a geometric morphism); hence any generalized element $E\to \tilde{F}(c)$ which factors through $\chi^{F}_{c}$ gives rise, by composition with $!_{E}^{\ast}\circ \beta$ of the corresponding, pointwise monic, natural transformation $\gamma_{{\cal E}\slash E}^{\ast} \circ y_{\cal C}c \to \tilde{F_{E}}$, to a poitwise monic natural transformation $\gamma_{{\cal E}\slash E}^{\ast} \circ y_{\cal C}c \to \tilde{F_{E}}$, that is to a generalized element $E\to \tilde{G}(c)$ which factors through $\chi^{G}_{c}$. But this generalized element is precisely $\beta_{c}\circ x$.
So $\beta(c):\tilde{F}(c)\to \tilde{G}(c)$ restricts to an arrow $F(c)\to G(c)$, as required.  
\end{proofs}

We can now prove the following

\begin{theorem}\label{unif}
Let $\mathbb T$ be a theory of presheaf type over a signature $\Sigma$. Then the injectivization of $\mathbb T$ satisfies condition $(iii)(b)\textrm{-}(1)$ of Theorem \ref{main} with respect to the category ${{{\textrm{f.p.} {\mathbb T}\textrm{-mod}(\Set)}}_{m}}$, and condition $(iii)(b)\textrm{-}(2)$ with respect to the same category if every monic arrow in ${\textrm{f.p.} {\mathbb T}\textrm{-mod}(\Set)}$ is sortwise monic (for instance, by Proposition \ref{freemono}, if for every sort $A$ over $\Sigma$ there exists the free $\mathbb T$-model on $A$).
\end{theorem}

\begin{proofs}
By Theorem \ref{injt}, the composite functor
\[
\xi_{{\cal E}}:\textbf{Flat}({{{\textrm{f.p.} {\mathbb T}\textrm{-mod}(\Set)}}_{m}}^{\textrm{op}}, {\cal E})\to \textbf{Flat}({{\textrm{f.p.} {\mathbb T}\textrm{-mod}(\Set)}}^{\textrm{op}}, {\cal E})\simeq {\mathbb T}\textrm{-mod}({\cal E}),
\] 
takes values in the subcategory ${\mathbb T}_{m}\textrm{-mod}({\cal E})$ of ${\mathbb T}\textrm{-mod}({\cal E})$. The functor $\xi_{{\cal E}}$ is faithful by Proposition \ref{cov}(iii) and full on ${\mathbb T}_{m}\textrm{-mod}({\cal E})$ by Theorem \ref{monoflat} in view of the fact that for any sortwise monic $\mathbb T$-model homomorphism $f:M\to N$ in $\cal E$, the natural transformation $\alpha_{f}:F_{M}\to F_{N}$ corresponding to it under the Morita-equivalence
\[
\mathbf{Flat}({{\textrm{f.p.} {\mathbb T}\textrm{-mod}(\Set)}}^{\textrm{op}}, {\cal E})\simeq {\mathbb T}\textrm{-mod}({\cal E})
\]
for $\mathbb T$ is pointwise monic. This latter fact can be proved as follows. For any object $D$ of the category ${\textrm{f.p.} {\mathbb T}\textrm{-mod}(\Set)}$, the value of $\alpha_{f}$ at $D$ can be identified with the arrow $[[\vec{x}. \phi]]_{M}\to [[\vec{x}. \phi]]_{N}$ canonically induced by $f$, where $\phi(\vec{x})$ is `the' formula which presents the model $D$; therefore if $f$ is sortwise monic then $\alpha_{f}(D)$ is monic in $\cal E$, it being the restriction of a monic arrow (namely $fA_{1}\times \cdots \times fA_{n}$, where $\vec{x}=(x^{A_{1}}, \ldots, x^{A_{n}})$) along two subobjects.
\end{proofs} 

The following proposition identifies a class of theories of presheaf type $\mathbb T$ with the property that the monic arrows of the category ${\textrm{f.p.} {\mathbb T}\textrm{-mod}(\Set)}$ are sortwise monic.

\begin{proposition}\label{freemono}
Let $\mathbb T$ be a theory of presheaf type over a signature $\Sigma$ in which for every sort $A$ over $\Sigma$ the formula $\{x^{A}.\top\}$ presents a $\mathbb T$-model, and let $f:M\to N$ be a homomorphism of finitely presentable $\mathbb T$-models $M$ and $N$. Then $f$ is monic as an arrow of ${\textrm{f.p.} {\mathbb T}\textrm{-mod}(\Set)}$ (equivalently, as an arrow of ${\mathbb T}\textrm{-mod}(\Set)$) if and only if it is sortwise monic.   
\end{proposition}

\begin{proofs}
Let $A$ be a sort over $\Sigma$. As $\{x^{A}.\top\}$ presents a $\mathbb T$-model $F_{A}$ then for any model $P$ of $\mathbb T$ in $\Set$ we have an equivalence $Hom_{{{\mathbb T}\textrm{-mod}(\Set)}}(F_{A}, P)\cong PA$ natural in $P$. In particular we have equivalences $Hom_{{{\mathbb T}\textrm{-mod}(\Set)}}(F_{A}, M)\cong MA$ and $Hom_{{{\mathbb T}\textrm{-mod}(\Set)}}(F_{A}, N)\cong NA$ under which the function $f\circ - :Hom_{{{\mathbb T}\textrm{-mod}(\Set)}}(F_{A}, M) \to Hom_{{{\mathbb T}\textrm{-mod}(\Set)}}(F_{A}, N)$ corresponds to the function $fA:MA\to NA$. Now, if $f$ is monic then the function $f\circ -:Hom_{{{\mathbb T}\textrm{-mod}(\Set)}}(F_{A}, M) \to Hom_{{{\mathbb T}\textrm{-mod}(\Set)}}(P_{A}, N)$ is injective, equivalently $fA:MA\to NA$ is injective, as required. 
\end{proofs}

\subsubsection{A criterion for injectivizations}

In this section we shall establish a result providing a sufficient condition for the injectivization of a theory of presheaf type of a certain form to be again of presheaf type. Before stating it, we need some preliminaries.

The following definition gives a natural topos-theoretic generalization of the standard notion of congruence on a set-based structure.

\begin{definition}
Let $\Sigma$ be a one-sorted first-order signature and $M$ a $\Sigma$-structure in a Grothendieck topos $\cal E$. An equivalence relation $R\mono M\times M$ on $M$ in $\cal E$ is said to be a \emph{congruence} if for any function symbol $f$ over $\Sigma$ of arity $n$, we have a commutative diagram
\[  
\xymatrix {
R^{n} \ar[d] \ar[r] & (M\times M)^{n}\cong M^{n}\times M^{n} \ar[d]^{f_{M}\times f_{M}} \\
R \ar[r] & M\times M.} 
\]
\end{definition}

\begin{proposition}
Let $\Sigma$ be a one-sorted first-order signature and $M$ a $\Sigma$-structure in a Grothendieck topos $\cal E$. For any congruence $R$ on $M$ there exists a $\Sigma$-structure $M\slash R$ on $\cal E$ whose underlying object is the quotient in $\cal E$ of $M$ by the relation $R$, and a $\Sigma$-structure epimorphism $p_{R}:M\to M\slash R$ given by the canonical projection. Conversely, for any $\Sigma$-structure epimorphism $q:M\to N$, the kernel pair $R_{q}$ of $q$ is a congruence on $M$ such that $q$ is isomorphic to $M\slash R_{q}$.   
\end{proposition}

\begin{proofs}
The proof of the proposition is immediate by using the exactness properties of Grothendieck toposes relating epimorphisms and equivalence relations.
\end{proofs}

\begin{proposition}\label{congruence}
Let $\mathbb T$ be a geometric theory over a one-sorted signature $\Sigma$ and $M$ a $\mathbb T$-model in a Grothendieck topos $\cal E$. If the axioms of $\mathbb T$ are all of the form $(\phi \vdash_{\vec{x}} \psi)$, where $\phi$ does not contain any conjunctions, then for any congruence $R$ on $M$ the $\Sigma$-structure $M\slash R$ is a $\mathbb T$-model.
\end{proposition}

\begin{proofs}
The thesis can be easily proved by induction on the structure of geometric formulae over $\Sigma$, using the fact that the action of the canonical projection homomorphism $p_{R}:M\to M\slash R$ on subobjects (of powers of $M$) preserves the top subobject, the natural order on subobjects, image factorizations and arbitrary unions.  
\end{proofs}

The following lemma shows that one can always perform image factorization of homomorphisms of structures in regular categories.

\begin{lemma}\label{lemmafact}
Let $\Sigma$ be a first-order signature and $\cal C$ a regular category. Then any $\Sigma$-structure homomorphism $f:M\to N$ in $\cal C$ can be factored as $h\circ g$, where $h:N'\mono N$ is a $\Sigma$-substructure of $N$ and $g:M\to N'$ is sortwise a cover.
\end{lemma}

\begin{proofs}
First, we notice that in any regular category finite products of covers are covers; indeed, composite of covers are covers (this can be easily proved by using the definition of cover as an arrow orthogonal to the class of monomorphisms), and the product of two covers $f\times g$, where $f:A\to B$ and $g:C\to D$, is equal to the composite $(1_{B}\times g) \circ (f\times 1_{C})$ of two arrows which are pullbacks of covers.

For every sort $A$ over $\Sigma$ we set $N'A$ equal to $Im(fA)$, $gA$ equal to the canonical cover $MA\to Im(fA)$ and $hA$ equal to the canonical subobject $N'A\mono NA$. For any function symbol $\xi:A_{1}, \ldots A_{n}\to B$ over $\Sigma$, we set $N'\xi$ equal to the restriction $N'A_{1}\times \cdots \times N'A_{n}\to N'B$ of $N\xi:NA_{1}\times \cdots \times NA_{n}\to NB$. This restriction actually exists (and is unique) since $f$ is a $\Sigma$-structure homomorphism and $gA_{1}\times \cdots \times gA_{n}$ is a cover. For any relation symbol $R$ over $\Sigma$ of type $A_{1}, \ldots, A_{n}$, we set $N'R$ equal to the intersection of $NR$ with the canonical subobject $N'A_{1}\times \cdots \times N'A_{n}\mono NA_{1}\times \cdots \times NA_{n}$. It is clear that $f=h\circ g$, that $g$ is sortwise a cover and that $h$ is a substructure homomorphism, as required.       
\end{proofs}

The following proposition, giving an explicit characterization of decidable objects in terms of their generalized elements, was stated in \cite{MM} as Exercise VIII.8(a).  

\begin{proposition}\label{dec}
Let $\cal E$ be a cocomplete (in particular, a Grothendieck) topos and $A$ an object of $\cal E$. Then $A$ is decidable if and only if for any generalized elements $x,y:E\to A$ there exists an epimorphic family (possibly consisting of just two elements) $\{e_{i}:E_{i}\to E \textrm{ | } i\in I\}$ such that for any $i\in I$, either $x\circ e_{i}=y\circ e_{i}$ or the equalizer of $x\circ e_{i}$ and $y\circ e_{i}$ is zero.   
\end{proposition}

\begin{proofs}
Let us suppose that $A$ is decidable. Let $p:P\mono A\times A$ be the complement of the diagonal subobject $\Delta:A\mono A\times A$. Consider the pullback of $<x, y>:E\to A\times A$ along $p$:
\[  
\xymatrix {
E' \ar[d]^{s} \ar[r]^{u} & P \ar[d]^{p} \\
E \ar[r]^{<x,y>}  & A\times A .} 
\]
The equalizer $i:R\mono E'$ of $x\circ s$ and $y\circ s$ is zero; indeed, by definition of $P$, the diagram
\[  
\xymatrix {
0 \ar[r] \ar[d] & A \ar[d]^{\Delta} \\
P  \ar[r]_{p} & A\times A} 
\]
is a pullback, and the arrows $z:=x\circ s \circ i=y\circ s\circ i:R\to A$ and $u\circ i:R\to P$ satisfy the condition $\Delta \circ z=p\circ u\circ i$.

Let us denote by $t:E''\mono E$ the pullback of the subobject $\Delta:A\mono A\times A$ along $<x,y>$. The arrows $s:E'\to E$ and $t:E''\to E$ are jointly epimorphic, they being the pullbacks of arrows which are jointly epimorphic, namely $\Delta$ and $p$. 

Therefore, by setting $I=\{0,1\}$, $E_{0}=E'$, $E_{1}=E''$, $e_{0}:E_{0}\to E$ equal to $s$ and $e_{1}:E_{1}\to E$ equal to $t$, we have that the family $\{e_{i}:E_{i}\to E \textrm{ | } i\in I\}$ satisfies the condition in the statement of the proposition. 

Conversely, let us suppose that the condition in the statement of the proposition is satisfied. To prove that $A$ is decidable, we have to show that the subobject $\Delta \vee \neg \Delta \mono A\times A$ is an isomorphism, in other words that any subobject $<x, y>:E\mono A\times A$ factors through it. 

Notice that for any generalized elements $x', y':E'\to A$, the arrow $<x', y'>:E'\to A\times A$ factors through $\neg \Delta \mono A\times A$ if and only if the equalizer of $x'$ and $y'$ is zero. Indeed, $<x',y'>$ factors through $\neg \Delta$ if and only if its image does, and by definition of $\neg \Delta$ this holds if and only if the pullback of it (or equivalently, of $<x', y'>$) along $\Delta$ is zero, i.e. if and only if the equalizer of $x'$ and $y'$ is zero.

Now, by our hypotheses, there exists an epimorphic family $\{e_{i}:E_{i}\to E \textrm{ | } i\in I\}$ such that for any $i\in I$, either $x\circ e_{i}=y\circ e_{i}$ or the equalizer of $x\circ e_{i}$ and $y\circ e_{i}$ is zero. By Lemma \ref{factorization}, to prove that $<x, y>$ factors through $\Delta \vee \neg \Delta \mono A\times A$ it suffices to prove that for any $i\in I$, $<x\circ e_{i}, y\circ e_{i}>$ factors through $\Delta \vee \neg \Delta \mono A\times A$. But by our hypothesis, for any given $i\in I$, either $x\circ e_{i}=y\circ e_{i}$ (which implies that $<x\circ e_{i}, y\circ e_{i}>$ factors through $\Delta:A\mono A\times A$) or the equalizer of $x\circ e_{i}$ and $y\circ e_{i}$ is zero (which implies, as we have just seen, that $<x\circ e_{i}, y\circ e_{i}>$ factors through $\neg \Delta\mono A\times A\mono A$); therefore for any $i\in I$, $<x\circ e_{i}, y\circ e_{i}>$ factors through $\Delta \vee \neg \Delta \mono A\times A$, as required.        
\end{proofs}

\begin{proposition}\label{findecid}
Let $\Sigma$ be a one-sorted signature, $c$ a finite $\Sigma$-structure in $\Set$, $M$ a $\Sigma$-structure in a Grothendieck topos $\cal E$ whose underlying object is decidable and $E$ an object of $\cal E$. Then for any $\Sigma$-structure homomorphism $f:c\to Hom_{\cal E}(E, M)$ there exists an epimorphic family $\{e_{i}:E_{i}\to E \textrm{ | } i\in I\}$ in $\cal E$ and for each $i\in I$ a quotient map $q_{i}:c\to c_{i}$, where $c_{i}$ is a finite $\Sigma$-structure, and a disjunctive $\Sigma$-structure homomorphism (in the sense of Lemma \ref{explicit2}) $J_{i}:c_{i}\mono Hom_{\cal E}(E_{i}, M)$ such that $J_{i}\circ q_{i}=Hom_{\cal E}(e_{i}, M)\circ f$ for all $i\in I$:
\[  
\xymatrix {
c  \ar[rr]^{f} \ar[d]^{q_{i}} & & Hom_{\cal E}(E, M) \ar[d]^{f} \\
c_{i}  \ar[rr]_{J_{i}} & & Hom_{\cal E}(E_{i}, M).} 
\]
\end{proposition}

\begin{proofs}
Let us suppose that $c$ has $n$ elements $x_{1}, \ldots, x_{n}$. We know from Proposition \ref{dec} that for any pair $(i, j)$ where $i,j \in \{1, \ldots, n\}$, there exist arrows $e_{(i, j)}:E_{(i, j)}\to E$ and $e'_{(i, j)}:E'_{(i, j)}\to E$ such that $e_{(i, j)}$ and $e'_{(i, j)}$ are jointly epimorphic, $f(x_{i})\circ e_{(i, j)}=f(x_{j})\circ e_{(i, j)}$ and $f(x_{i})\circ e'_{(i, j)}, f(x_{j})\circ e'_{(i, j)}$ are disjoint. The iterated fibered product of all these epimorphic families (corresponding to the pairs $(i, j)$ such that $i,j \in \{1, \ldots, n\}$) thus yields an epimorphic family $\{u_{k}:U_{k}\to E \textrm{ | } k\in K\}$ such that for every $k\in K$ there exists a subset $S_{k}\subseteq \{1, \ldots, n\}\times \{1, \ldots, n\}$ such that for every $(i, j)\in S_{k}$,   $f(x_{i})\circ u_{k}=f(x_{j})\circ u_{k}$ and for every $(i, j)\notin S_{k}$, $f(x_{i})\circ u_{k}$ and $f(x_{j})\circ u_{k}$ are disjoint. For each $k\in K$, consider the quotient $q_{k}:c\to c_{q}$ of $c$ by the congruence generated by the pairs of the form $(x_{i}, x_{j})$ for $(i, j)\in S_{k}$. By definition of this congruence relation, the $\Sigma$-structure homomorphism $Hom_{\cal E}(e_{k}, M)\circ f$ factors through $q_{k}$, and the resulting factorization is a disjunctive $\Sigma$-structure homomorphism. We have thus found a set of data satisfying the condition in the statement of the proposition, as required. 
\end{proofs}

\begin{proposition}\label{propinjfin}
Let $\mathbb T$ be a theory of presheaf type over a signature $\Sigma$ such that the finitely presentable $\mathbb T$-models coincide with the finitely presentable ${\mathbb T}_{m}$-models. Suppose that the following condition is satisfied: for any Grothendieck topos $\cal E$, object $E$ of $\cal E$ and $\Sigma$-structure homomorphism $x:c\to Hom_{\cal E}(E, M)$, where $c$ is a finitely presentable ${\mathbb T}$-model and $M$ is a sortwise decidable $\mathbb T$-model, there exists an epimorphic family $\{e_{i}:E_{i}\to E \textrm{ | } i\in I\}$ in $\cal E$ and for each $i\in I$ a $\mathbb T$-model homomorphism $f_{i}:c\to c_{i}$ of finitely presentable $\mathbb T$-models and a sortwise disjunctive $\Sigma$-structure homomorphism $x_{i}:c_{i} \to Hom_{{\cal E}}(E_{i}, M)$ such that $x_{i} \circ f_{i}=Hom_{\cal E}(e_{i}, M)\circ x$ for all $i\in I$.
Then the injectification ${\mathbb T}_{m}$ of $\mathbb T$ satisfies conditions $(i)$ and $(ii)$ of Theorem \ref{main} (with respect to its category of finitely presentable models). 
\end{proposition}

\begin{proofs}
The proposition represents the particular case of Theorem \ref{finalitythm} for the faithful interpretation of ${\mathbb T}$ into its injectivization, with ${\cal K}={\cal H}$ equal to the category of finitely presentable $\mathbb T$-models. The hypotheses of the theorem are satisfied since the functor $\int a$ is full (cf. Remark \ref{appl}).
\end{proofs}

\begin{remark}
Any geometric theory $\mathbb T$ such that the finitely presentable $\mathbb T$-models are precisely the finitely generated ones and all its axioms are of the form $(\phi \vdash_{\vec{x}} \psi)$, where $\psi$ is a quantifier-free geometric formula, satisfies the first of the hypotheses of the Corollary. Indeed, every finitely presentable $\mathbb T$-model is clearly finitely presentable as a ${\mathbb T}_{m}$-model, and since by Proposition \ref{coincid} the substructures of models of $\mathbb T$ are again models of $\mathbb T$, every model of $\mathbb T$ can be expressed as the directed union of its finitely generated submodels; so the finitely presentable ${\mathbb T}_{m}$-models are precisely the finitely generated $\mathbb T$-models.  
\end{remark}

\begin{corollary}\label{injcorollary}
Let $\mathbb T$ be a theory of presheaf type over a signature $\Sigma$ such that the finitely presentable $\mathbb T$-models coincide with the finitely presentable ${\mathbb T}_{m}$-models and the monic $\mathbb T$-model homomorphisms in $\Set$ are precisely the homomorphisms which are sortwise monic. Suppose that the following condition is satisfied: for any Grothendieck topos $\cal E$, object $E$ of $\cal E$ and $\Sigma$-structure homomorphism $x:c\to Hom_{\cal E}(E, M)$, where $c$ is a finitely presentable ${\mathbb T}$-model and $M$ is a sortwise decidable $\mathbb T$-model, there exists an epimorphic family $\{e_{i}:E_{i}\to E \textrm{ | } i\in I\}$ in $\cal E$ and for each $i\in I$ a $\mathbb T$-model homomorphism $f_{i}:c\to c_{i}$ of finitely presentable $\mathbb T$-models and a sortwise disjunctive $\Sigma$-structure homomorphism $x_{i}:c_{i} \to Hom_{{\cal E}}(E_{i}, M)$ such that $x_{i} \circ f_{i}=Hom_{\cal E}(e_{i}, M)\circ x$ for all $i\in I$.
Then the injectification ${\mathbb T}_{m}$ of $\mathbb T$ is of presheaf type.
\end{corollary}

\begin{proofs}
By Theorem \ref{unif}, ${\mathbb T}_{m}$ satisfies condition $(iii)$ of Theorem \ref{main} (since all the the monic arrows in the category ${\textrm{f.p.} {\mathbb T}\textrm{-mod}(\Set)}$ are injective functions), while Proposition \ref{propinjfin} ensures that conditions $(i)$ and $(ii)$ of Theorem \ref{main} are satisfied. We can thus conclude that the theory ${\mathbb T}_{m}$ is of presheaf type.
\end{proofs}

As a consequence of Corollary \ref{injcorollary}, we obtain the following result. 
 
\begin{corollary}\label{corinjectivizations}
Let $\mathbb T$ be a geometric theory over a one-sorted signature $\Sigma$ such that any quotient of a finitely presentable $\mathbb T$-model is a $\mathbb T$-model (for instance, a theory whose axioms are all of the form $(\phi \vdash_{\vec{x}} \psi)$, where $\phi$ does not contain any conjunctions - cf. Proposition \ref{congruence}). Suppose that the finitely presentable $\mathbb T$-models are exactly the finite $\mathbb T$-models and that all the the monic arrows in the category ${\textrm{f.p.} {\mathbb T}\textrm{-mod}(\Set)}$ are injective functions. Then if $\mathbb T$ is of presheaf type, ${\mathbb T}_{m}$ is of presheaf type as well.   
\end{corollary}

\begin{proofs}
Propositions \ref{congruence} and \ref{findecid} ensure that all the conditions of Corollary \ref{injcorollary} are satisfied. We can thus conclude that ${\mathbb T}_{m}$ is of presheaf type, as required.   
\end{proofs}

\section{Expansions of theories}

\subsection{General theory}\label{genth}

In this section we introduce the syntactic notion of expansion of a geometric theory, and show that it corresponds in a natural way to having a geometric morphisms between the respective classifying toposes. Further, we characterize the expansions of theories whose corresponding geometric morphisms are localic and hyperconnected. This section has been inspired by section 6.2 of \cite{Coumans}, which contains a brief informal discussion of this topic.

\begin{definition}
Let $\mathbb T$ be a geometric theory over a signature $\Sigma$.
\begin{enumerate}[(i)]
\item A geometric \emph{expansion} of $\mathbb T$ is a geometric theory obtained from $\mathbb T$ by adding sorts, relation or function symbols to the signature $\Sigma$ and geometric axioms over the resulting extended signature; equivalently, a geometric expansion of $\mathbb T$ is a geometric theory ${\mathbb T}'$ over a signature containing $\Sigma$ such that every axiom of $\mathbb T$, regarded as a geometric sequent in the signature of ${\mathbb T}'$, is provable in ${\mathbb T}'$.

\item A geometric expansion ${\mathbb T}'$ of $\mathbb T$ is said to be \emph{localic} if no new sorts are added to $\Sigma$ to obtain the signature of ${\mathbb T}'$.

\item A geometric expansion ${\mathbb T}'$ of $\mathbb T$ is said to be \emph{hyperconnected} if no new function or relation symbols which only involve sorts over $\Sigma$ are added to $\Sigma$ to form the signature of ${\mathbb T}'$, and for any geometric sequent $\sigma$ over $\Sigma$, $\sigma$ is provable in ${\mathbb T}'$ if and only if it is provable in $\mathbb T$.  
\end{enumerate}
\end{definition}

Notice that if ${\mathbb T}'$ is an expansion of a geometric theory $\mathbb T$ then there is a canonical morphism of sites $({\cal C}_{\mathbb T}, J_{\mathbb T})\to ({\cal C}_{{\mathbb T}'}, J_{{\mathbb T}'})$ inducing a geometric morphism $p_{\mathbb T}^{{\mathbb T}'}:\Set[{\mathbb T}']\to \Set[{\mathbb T}]$.

We say that a geometric expansion ${\mathbb T}'$ of $\mathbb T$ is \emph{faithful} (resp. \emph{fully faithful}) if for every Grothendieck topos $\cal E$, the induced functor 
\[
(p_{\mathbb T}^{{\mathbb T}'})_{\cal E}:{\mathbb T}'\textrm{-mod}({\cal E})\to {\mathbb T}\textrm{-mod}({\cal E})
\]
is faithful (resp. full and faithful). 

We know from Theorem 9.1 \cite{OC6} that the inclusion part of the surjection-inclusion factorization of the geometric morphism $p_{\mathbb T}^{{\mathbb T}'}:\Set[{\mathbb T}']\to \Set[{\mathbb T}]$ can be identified with the classifying topos of the quotient of $\mathbb T$ consisting of all the sequents over $\Sigma$ which are provable in the theory ${\mathbb T}'$.

Recall (cf. for instance section A4.6 of \cite{El}) that a geometric morphism $f:{\cal F}\to {\cal E}$ is \emph{localic} if every object of $\cal F$ is a subquotient (i.e. a quotient of a subobject) of an object of the form $f^{\ast}(a)$, where $a$ is an object of $\cal E$; the morphism $f$ is \emph{hyperconnected} if  
$f^{\ast}$ is full and faithful, and its image is closed under subobjects in $\cal F$. 

The following technical lemma will be useful in the sequel.

\begin{lemma}\label{hypercon}
Let $f:{\cal F}\to {\cal E}$ a geometric morphism between Grothendieck toposes and ${\cal C}$ a full subcategory of $\cal E$ which is separating for $\cal E$. Suppose that the following conditions are satisfied:
\begin{enumerate}[(i)]
\item $\cal C$ is closed in $\cal E$ under finite products;
\item $f$ satisfies the property that the restriction $f^{\ast}|_{\cal C}:{\cal C}\to {\cal F}$ of $f^{\ast}$ to $\cal C$ is full and faithful;
\item for any family of arrows $\cal T$ in $\cal C$ with common codomain, if the image of $\cal T$ under $f^{\ast}$ is epimorphic in $\cal F$ then $\cal T$ is epimorphic in $\cal E$, and 
\item every subobject in $\cal F$ of an object of the form $f^{\ast}(c)$ where $c$ is an object of $\cal C$, is, up to isomorphism, of the form $f^{\ast}(m)$ (where $m$ is a subobject in ${\cal E}$).
\end{enumerate}
Then $f$ is hyperconnected.   
\end{lemma}

\begin{proofs}
We have to prove that, under the specified hypotheses, $f^{\ast}:{\cal E}\to {\cal F}$ is full and faithful and its image is closed under subobjects in $\cal F$. 

To prove that $f^{\ast}$ is faithful, we have to verify that for any arrows $h, k:u\to v$ in $\cal E$, $f^{\ast}(h)=f^{\ast}(k)$ implies $h=k$. Since $\cal C$ is a separating set for $\cal E$, we can suppose without loss of generality $u$ to lie in $\cal C$. Now, $h=k$ if and only if the equalizer $z:w\mono u$ of $h$ and $k$ is an isomorphism. The family of arrows from objects of $\cal C$ to $w$ is epimorphic, whence, as $f^{\ast}(z)$ is an isomorphism (since $f^{\ast}(h)=f^{\ast}(k)$), the family formed by the composition of these arrows with $z$ is epimorphic on $u$; indeed, the members of this latter family all lie in $\cal C$ and the image of this family under $f^{\ast}$ is epimorphic. It follows that $z$ is an epimorphism, equivalently an isomorphism, that is $h=k$, as required.     

To prove the fullness of $f^{\ast}$, we have to verify that for any objects $a$ and $b$ of $\cal E$ and any arrow $s:f^{\ast}(a)\to f^{\ast}(b)$ in $\cal F$, there exists a (unique) arrow $t:a\to b$ in $\cal E$ such that $f^{\ast}(t)=s$. As $\cal C$ is a separating set for $\cal E$ there exist epimorphic families $\{f_{i}:c_{i} \to b \textrm{ | } i\in I\}$ and $\{g_{j}:d_{j} \to a \textrm{ | } j\in J\}$ in $\cal E$ consisting of arrows whose domains lie in $\cal C$.

For any $i\in I$ and $j\in J$, consider the pullback square
\[  
\xymatrix {
r_{i, j} \ar[d]^{q_{j}} \ar[r]^{p_{i}} & f^{\ast}(c_{i}) \ar[d]^{f^{\ast}(f_{i})}\\
f^{\ast}(d_{j}) \ar[r]_{s\circ f^{\ast}(g_{j})} & f^{\ast}(b).} 
\]
By the universal property of pullbacks, the arrow $<p_{i}, q_{j}>:r_{i, j} \to f^{\ast}(c_{i})\times f^{\ast}(d_{j})\cong f^{\ast}(c_{i}\times d_{j})$ is a monomorphism. Since $c_{i}\times d_{j}$ lies in $\cal C$ by our hypothesis, $r_{i, j}$ lies in the image of $f^{\ast}$ and hence can be covered by arrows $h^{i, j}_{k}:f^{\ast}(e^{i, j}_{k})\to r_{i, j}$ (for $k\in K_{i, j}$) where each object $e^{i, j}_{k}$ lies in $\cal C$. Consider the arrows $p_{i}\circ h_{k}^{i, j}:f^{\ast}(e_{k}^{i, j}) \to f^{\ast}(c_{i})$ and $q_{j}\circ h_{k}^{i, j}:f^{\ast}(e_{k}^{i, j}) \to f^{\ast}(d_{j})$. Since the family $\{f^{\ast}(g_{j}\circ m_{k}^{i,j}) \textrm{ | } i\in I, j\in J,  k\in K_{i,j}\}$ is epimorphic in $\cal F$, our hypotheses ensure that the family $\{ g_{j}\circ m_{k}^{i,j}:e_{k}^{i, j} \to A \textrm{ | } i\in I, j\in J, k\in K_{i,j}\}$ is epimorphic in $\cal E$. To define an arrow $t:a\to b$ in $\cal E$ it is therefore equivalent to specify, for each $i\in I$, $j\in J$ and $k\in K_{i, j}$, an arrow $v_{k}^{i, j}:e_{k}^{i, j}\to b$ in such a way that the compatibility conditions of Corollary \ref{cordesc} are satisfied. Take $v_{k}^{i, j}$ equal to $f_{i}\circ l_{k}^{i, j}$. The compatibility conditions hold since the images of them under the functor $f^{\ast}$ do and, as we have proved above, $f^{\ast}$ is faithful. Thus we have a unique arrow $t:a\to b$ in $\cal E$ such that $t\circ g_{j}\circ m_{k}^{i, j}=f_{i}\circ l_{k}^{i, j}$. The fact that $f^{\ast}(t)=s$ follows at once.    

To conclude the proof of the lemma, it remains to show that the image of $f^{\ast}$ is closed under subobjects in $\cal F$. Let $a$ be an object of $\cal E$ and $k:z\mono f^{\ast}(a)$ a subobject in $\cal F$. Since $\cal C$ is a separating set for $\cal E$, the canonical arrow from the coproduct of all objects in $\cal C$ to $a$ is an epimorphism $y$; by pulling back $f^{\ast}(y)$ along the monomorphism $k$ we obtain a monomorphism $u$ from $z'$ to $dom(f^{\ast}(y))$ and an epimorphism $e:z'\to z$. Now, by our hypotheses, the pullbacks of $u$ along the coproduct arrows belong to the image of $f^{\ast}$; from the fact that $f^{\ast}$ preserves coproducts it thus follows that $u$ itself belongs to the image of $f^{\ast}$, it being of the form $f^{\ast}(m)$, where $m$ is a coproduct in $\cal E$ of subobjects in $\cal C$. Now, since $e$ is an epimorphism, $z$ is isomorphic to the coequalizer of its kernel pair $r:w\mono f^{\ast}(b)\times f^{\ast}(b)\cong f^{\ast}(b\times b)$, where $b=dom(m)$. Using the coproduct representation of $m$ and the fact that $\cal C$ is closed under finite products in $\cal E$, one can prove by considering the pullbacks of $r$ along the images under $f^{\ast}$ of the coproduct arrows to $b\times b$, that $r$ lies, up to isomorphism, in the image of $f^{\ast}$; therefore, as $f^{\ast}$ preserves coequalizers, $k$ is isomorphic to a subobject in the image of $f^{\ast}$, as required. 
\end{proofs}

\begin{theorem}\label{factor}
Let $\mathbb T$ be a geometric theory over a signature $\Sigma$ and ${\mathbb T}'$ a geometric expansion of $\mathbb T$ over a signature $\Sigma'$. If ${\mathbb T}'$ is a hyperconnected (resp. a localic) expansion of $\mathbb T$ then $p_{\mathbb T}^{{\mathbb T}'}:\Set[{\mathbb T}']\to \Set[{\mathbb T}]$ is a hyperconnected (resp. a localic) geometric morphism. 

In particular, the hyperconnected-localic factorization of the geometric morphism $p_{\mathbb T}^{{\mathbb T}'}:\Set[{\mathbb T}']\to \Set[{\mathbb T}]$ is given $p_{\mathbb T}^{{\mathbb T}''}\circ p_{{\mathbb T}''}^{{\mathbb T}'}$, where ${\mathbb T}''$ is the intermediate expansion of ${\mathbb T}$ obtained by adding to the signature $\Sigma$ of $\mathbb T$ no new sorts and all the relation and function symbols over $\Sigma'$ which only involve the sorts over $\Sigma$, and all the sequents over this extended signature which are provable in ${\mathbb T}'$. 
\end{theorem}

\begin{proofs}
Suppose that ${\mathbb T}'$ is localic over ${\mathbb T}$. We have to prove that every object of $\Set[{\mathbb T}']$ is a quotient of a subobject of an object of the form $f^{\ast}(H)$ where $H$ is an object of $\Set[{\mathbb T}]$. Since the signature of ${\mathbb T}'$ does not contain any sorts already present in the signature of $\mathbb T$ and the objects of the category ${\cal C}_{{\mathbb T}'}$ form a separating set for the topos $\Set[{\mathbb T}']$, we can conclude that the subobjects of the objects of the form ${p_{\mathbb T}^{{\mathbb T}'}}^{\ast}(\{\vec{x}. \top\})$ form a separating set for $\Set[{\mathbb T}']$. Therefore any object of $\Set[{\mathbb T}']$ is a quotient of a coproduct of subobjects of objects of the form ${p_{\mathbb T}^{{\mathbb T}'}}^{\ast}(\{\vec{x}. \top\})\cong \{\vec{x}. \top\}$; but the fact that ${p_{\mathbb T}^{{\mathbb T}'}}^{\ast}$ preserves coproducts immediately implies that any such coproduct is a subobject of an object in the image of ${p_{\mathbb T}^{{\mathbb T}'}}^{\ast}$.

Suppose instead that ${\mathbb T}'$ is hyperconnected over ${\mathbb T}$. We have to prove that ${p_{\mathbb T}^{{\mathbb T}'}}^{\ast}$ is hyperconnected. It suffices to apply Lemma \ref{hypercon} by taking $\cal E$ equal to $\Set[{\mathbb T}]$ and $\cal C$ equal to the syntactic category ${\cal C}_{\mathbb T}$; the fact that the hypotheses of the lemma are satisfied follows immediately from the fact that ${\mathbb T}'$ is hyperconnected over $\mathbb T$. 
\end{proofs}

\begin{theorem}\label{critclassifyingtopos}
Let $\mathbb T$ be a geometric theory, $\cal E$ a Grothendieck topos and $M$ a model of $\mathbb T$ in $\cal E$. Then $\cal E$ is the classifying topos for $\mathbb T$ and $M$ is a universal model for $\mathbb T$ if and only if the following conditions are satisfied:
\begin{enumerate}[(i)]
\item The family $\cal F$ of objects which can be built from the interpretations in $M$ of the sorts, function symbols and relation symbols over the signature of $\mathbb T$ by using geometric logic constructions (i.e. the objects given by the interpretations in $M$ of geometric formulae over the signature of $\mathbb T$) is separating for $\cal E$;

\item The model $M$ is conservative for $\mathbb T$; that is, for any geometric sequent $\sigma$ over the signature of $\mathbb T$, $\sigma$ is valid in $M$ if and only if it is provable in $\mathbb T$; 

\item Any arrow $k$ in $\cal E$ between objects $A$ and $B$ in the family $\cal F$ of point $(i)$ is definable; that is, if $A$ (resp. $B$) is equal to the interpretation of a geometric formula $\phi(\vec{x})$ (resp. $\psi(\vec{y})$) over the signature of $\mathbb T$, there exists a $\mathbb T$-provably functional formula $\theta$ from $\phi(\vec{x})$ to $\psi(\vec{x})$ such that the interpretation of $\theta$ in $M$ is equal to $k$.
\end{enumerate}  
\end{theorem}

\begin{proofs}
By the universal property of the geometric syntactic category ${\cal C}_{\mathbb T}$ of $\mathbb T$, the $\mathbb T$-model $M$ corresponds to a geometric functor $F_{M}:{\cal C}_{\mathbb T} \to {\cal E}$ assigning to each object (or arrow) of ${\cal C}_{\mathbb T}$ its interpretation in $M$. Condition $(ii)$ in the statement of the theorem is equivalent to the assertion that the functor $F_{M}$ is faithful, while condition $(iii)$ is equivalent to saying that $F_{M}$ is full. Therefore under conditions $(ii)$ and $(iii)$, ${\cal C}_{\mathbb T}$ embeds as a full subcategory of $\cal E$. Now, condition $(i)$ ensures that ${\cal C}_{\mathbb T}$ is dense in $\cal E$, whence the Comparison Lemma yields an equivalence ${\cal E}\simeq \Sh({\cal C}_{\mathbb T}, J)$, where $J$ is the Grothendieck topology on ${\cal C}_{\mathbb T}$ induced by the canonical topology on $\cal E$, that is the geometric topology on ${\cal C}_{\mathbb T}$. By the syntactic construction of classifying toposes and universal models, we can thus conclude that $\cal E$ is `the' classifying topos for $\mathbb T$ and $M$ is a universal model for $\mathbb T$. 
\end{proofs}

As immediate corollaries of Theorem \ref{critclassifyingtopos}, one recovers the following known results:

\begin{enumerate}[(i)]
\item Let $\cal C$ be a separating set of objects for a Grothendieck topos $\cal E$, and $\Sigma_{{\cal C}}$ the signature consisting of one sort $\name{c}$ for each object $c$ of $\cal C$ and one function symbol $\name{f}$ for each arrow $f$ in $\cal E$ between objects in $\cal C$. Then there exists a geometric theory $\mathbb T$ over the signature $\Sigma_{{\cal C}}$ classified by $\cal E$, whose universal model is given by the `tautological' $\Sigma_{{\cal C}}$-structure in $\cal E$ (cf. p. 837 \cite{El});  

\item Let $B$ be a pre-bound for $\cal E$ over $\Set$ (that is, an object such that the subobjects of its finite powers form a generating set for $\cal E$); then there exists a one-sorted geometric theory $\mathbb T$ classified by $\cal E$ and a universal model for $\mathbb T$ whose underlying object is $B$ (cf. Theorem D3.2.5 \cite{El}). 
 
\end{enumerate}
 
The first result can be obtained from Theorem \ref{critclassifyingtopos} by observing that the theory $\mathbb T$ of the tautological $\Sigma_{{\cal C}}$-structure obviously satisfies all the hypotheses of the theorem. That the first two conditions are satisfied is obvious, while the fact that the third holds can be proved as follows. Since every object of the syntactic category of $\mathbb T$ is a subobject of a finite product of objects of the form $\{x^{\name{c}}. \top\}$ (for $c\in {\cal C}$), it is enough to prove that every arrow $k:R\to c$ in $\cal E$ having as codomain an object $c$ in $\cal C$ is definable. Suppose that $R=[[\vec{x}. \phi]]_{M}$, where the sorts of the variables in $\vec{x}=(x_{1}, \ldots, x_{n})$ are respectively $\name{c_{1}'}, \ldots, \name{c_{n}'}$, and denote by $r:R\to c_{1}'\times \cdots \times c_{n}'$ the corresponding subobject. Since $\cal C$ is a separating set for $\cal E$, the family of arrows $\{f_{i}:c_{i}\to R \textrm{ | } i\in I\}$ from objects of $\cal C$ to $R$ is epimorphic; hence the geometric formula     
\[
\mathbin{\mathop{\textrm{ $\bigvee$}}\limits_{i\in I}} (\exists z_{i}^{\name{c_{i}}})(x_{1}=\name{g_{1}}(z_{i})\wedge \cdots \wedge x_{n}=\name{g_{n}}(z_{i})\wedge x^{\name{c}}=\name{k\circ f_{i}}(z_{i})),
\]
where $r\circ f_{i}=<g_{1}, \ldots, g_{n}>$, is $\mathbb T$-provably functional from $\{\vec{x}. \phi\}$ to $\{x^{\name{c}}. \top\}$ and its interpretation coincides with $k$. By Diaconescu's theorem, the theory $\mathbb T$ can be explicitly characterized as the theory of flat $J_{\cal E}^{\cal C}$-continuous functors on $\cal C$, where $J_{\cal E}^{\cal C}$ is the Grothendieck topology on $\cal C$ induced by the canonical topology on $\cal E$.

The second result can be deduced from Theorem \ref{critclassifyingtopos} by taking $\mathbb T$ to be the theory of the tautological structure over the one-sorted signature $\Sigma_{B}$ consisting of an $n$-ary relation 
symbol $\name{R}$ for each subobject $R\mono B^{n}$ in $\cal E$. The fact that $\mathbb T$ satisfies the first two conditions of the theorem is obvious, while the validity of the third condition follows from the fact that the graphs of morphisms $B^{n}\to B$ in $\cal E$ are interpretations of $(n+1)$-relation symbols over $\Sigma_{B}$.  

The following theorem provides a converse to Theorem \ref{factor}.

\begin{theorem}\label{expth}
Let $p:{\cal E}\to \Set[{\mathbb T}]$ be a geometric morphism to the classifying topos of a geometric theory $\mathbb T$. Then $p$ is, up to isomorphism, of the form $p_{\mathbb T}^{{\mathbb T}'}$ for some geometric expansion ${\mathbb T}'$ of ${\mathbb T}$. If $p$ is hyperconnected (resp. localic) then we can take ${\mathbb T}'$ to be a hyperconnected (resp. localic) expansion of ${\mathbb T}$.
\end{theorem}

\begin{proofs}
Choose a triplet ${\cal T}=({\cal C}_{ob}, {\cal C}_{arr}, {\cal C}_{rel})$ consisting of a set ${\cal C}_{ob}$ of objects of $\cal E$, of a set ${\cal C}_{arr}$ of arrows in $\cal E$ from finite products of objects of ${\cal C}_{ob}$ to objects of ${\cal C}_{ob}$ and of a set ${\cal C}_{rel}$ of subobjects of finite products of objects of ${\cal C}_{ob}$ with the property that the family of objects of $\cal E$ which can be built out of objects in ${\cal C}_{ob}$, arrows in ${\cal C}_{arr}$ and subobjects in ${\cal C}_{rel}$ by using geometric logic constructions is separating for $\cal E$. By definition of Grothendieck topos, such a triplet always exists. Let us define an expansion ${\mathbb T}_{{\cal T}}$ of $\mathbb T$ as follows: the signature $\Sigma_{{\mathbb T}_{{\cal T}}}$ of ${\mathbb T}_{{\cal T}}$ is obtained by adding to the signature of $\mathbb T$ one sort $\name{c}$ for each object $c$ in ${\cal C}_{ob}$ which is not of the form $f^{\ast}(H)$ for some object $H$ of ${\cal C}_{\mathbb T}\hookrightarrow \Set[{\mathbb T}]$, one function symbol $\name{f}$ for each arrow $f$ in ${\cal C}_{arr}$ whose domain or codomain is not of the form $f^{\ast}(H)$ (with the obvious sorts), one relation symbol for each subobject in ${\cal C}_{rel}$ (with the obvious sorts, the ones corresponding to an object of the form $f^{\ast}(\{\vec{x}.\phi\})$ being the sorts of the variables $\vec{x}$) and an additional relation symbol $\name{R}$ for any subobject $R\mono c_{1}\times \cdots \times c_{n}$ in $\cal E$ (where $c_{1}, \ldots, c_{n}$ are objects in ${\cal C}_{ob}$) which cannot be obtained from the data in $\cal T$ by means of geometric logic constructions (whose sorts are the obvious ones).

Consider the tautological ${\mathbb T}_{{\cal T}}$-structure $M$ in $\cal E$ obtained by interpreting each sort $\name{c}$ by the corresponding object $c$ (and similarly for the function and relation symbols added to the signature of $\mathbb T$), and each sort $A$ over the signature of $\mathbb T$ by the object $f^{\ast}(\{x^{A}. \top\})$ (and similarly for the function and relation symbols over the signature of $\mathbb T$). Define ${\mathbb T}'$ to be the theory of $M$ over the signature ${\mathbb T}_{{\cal T}}$. The theory ${\mathbb T}'$ satisfies the conditions of Theorem \ref{critclassifyingtopos}; the validity of the first two conditions is obvious, while the validity of the third follows from the fact that any subobject of a product of objects in ${\cal C}_{ob}$ is definable in ${\mathbb T}'$ and the model $M$ is conservative for ${\mathbb T}'$, whence the formula defining the graph of the given arrow is ${\mathbb T}'$-provably functional from the formula defining the domain to the formula defining the codomain. Therefore ${\mathbb T}'$ is classified by the topos $\cal E$ with universal model $M$. Notice that ${\mathbb T}'$ is an expansion of $\mathbb T$. This proves the first part of the theorem.

Suppose now that $f$ is localic. We can define a triplet ${\cal T}=({\cal C}_{ob}, {\cal C}_{arr}, {\cal C}_{rel})$ satisfying the conditions specified above as follows: we set ${\cal C}_{ob}$ equal to the set of objects of the form $f^{\ast}(H)$ where $H$ is an object of ${\cal C}_{\mathbb T}\hookrightarrow \Set[{\mathbb T}]$, ${\cal C}_{arr}$ equal to the empty set and ${\cal C}_{rel}$ equal to the set of all subobjects of (finite products of) objects in ${\cal C}_{ob}$. Alternatively, one can take ${\cal C}_{ob}$ to be the set of objects of the form $f^{\ast}(\{x^{A}. \top\})$ (where $A$ is any sort over the signature of $\mathbb T$), ${\cal C}_{arr}$ to be the set of arrows of the form $f^{\ast}([\theta])$, where $[\theta]$ is an arrow in ${\cal C}_{\mathbb T}$, and ${\cal C}_{rel}$ to be the set of all subobjects of (finite products of) objects in ${\cal C}_{ob}$. In either case, the signature of the resulting expansion will contain no new sorts with respect to the signature of $\mathbb T$ and hence it will be localic.

Suppose instead that $f$ is hyperconnected. Given any set of objects $\cal K$ of $\cal E$ which, together with the objects of the form $f^{\ast}(H)$ (where $H$ is an object of ${\cal C}_{\mathbb T}\hookrightarrow \Set[{\mathbb T}]$), form a separating set of $\cal E$, we can define a triplet ${\cal T}=({\cal C}_{ob}, {\cal C}_{arr}, {\cal C}_{rel})$ satisfying the required conditions by setting ${\cal C}_{ob}={\cal K}$, ${\cal C}_{arr}=\emptyset$ and ${\cal C}_{rel}=\emptyset$. Since $f$ is hyperconnected, $f^{\ast}$ is full and the image of $f^{\ast}$ is closed under subobjects; hence the signature of ${\mathbb T}_{{\cal T}}$ does not contain any relation or function symbol only involving the sorts of $\mathbb T$. On the other hand, any geometric sequent over the signature of $\mathbb T$ is provable in ${\mathbb T}'$ if and only if it is provable in $f^{\ast}(M_{\mathbb T})$, where $M_{\mathbb T}$ is the universal model of $\mathbb T$ lying in $\Set[{\mathbb T}]$ (since $f^{\ast}$ is full and faithful), i.e. if and only if it is provable in ${\mathbb T}'$. Hence ${\mathbb T}'$ is a hyperconnected expansion of $\mathbb T$, as required.
 
This completes our proof. 
\end{proofs}

Theorem \ref{expth} yields, in view of the equivalence between conditions $(iii)(b)$ and $(iii)(c)$ of Theorem \ref{main} and Remark \ref{conditionstar}(d), the following reformulation of a particular case of Theorem \ref{unif}.

\begin{theorem}
Let $\mathbb T$ be a theory of presheaf type such that the finitely presentable models of $\mathbb T$ coincide with the finitely presentable models of ${\mathbb T}_{m}$. Then there is a faithful expansion (in the sense of section \ref{genth}) of the injectivization of $\mathbb T$ which is classified by the topos $[{\textrm{f.p.} {\mathbb T}\textrm{-mod}(\Set)}_{m}, \Set]$. If moreover all the monic homomorphisms in ${{\mathbb T}\textrm{-mod}(\Set)}$ are all sortwise monic then this expansion can be taken to be fully faithful.
\end{theorem}\qed
 
A simple example of a theory satisfying both of the hypotheses of the theorem is the theory $\mathbb A$ of commutative rings with unit. Indeed, the finitely presented commutative rings with unit coincide precisely with the finitely generated ones, that is with the finitely presentable models of ${\mathbb A}_{m}$; also, the monic ring homomorphisms are precisely the injective ones. 

Another consequence of Theorem \ref{critclassifyingtopos} is the following criterion for a geometric theory to be of presheaf type.

\begin{theorem}\label{defpresheaf}
Let $\mathbb T$ be a geometric theory over a signature $\Sigma$. Then $\mathbb T$ is of presheaf type if and only if the following conditions are satisfied:
\begin{enumerate}[(i)]
\item Every finitely presentable model is presented by a geometric formula over $\Sigma$;

\item Every property of finite tuples of elements of a (finitely presentable) $\mathbb T$-model which is preserved by $\mathbb T$-model homomorphisms is definable by a geometric formula over $\Sigma$; 

\item The finitely presentable $\mathbb T$-models are jointly conservative for $\mathbb T$.

\end{enumerate}
\end{theorem} 

\begin{proofs}
The fact that any theory of presheaf type satisfies the given conditions was established in \cite{OCU}. It thus remains to prove the `if' part of the theorem. 

Consider the $\Sigma$-structure $U$ in the topos $[\textrm{f.p.} {\mathbb T}\textrm{-mod}(\Set), \Set]$ given by the forgetful functors at each sort. Clearly, $U$ is a $\mathbb T$-model.

To deduce our thesis, we shall verify that $\mathbb T$ satisfies the conditions of Theorem \ref{critclassifyingtopos} with respect to the model $U$. 

Since every finitely presentable $\mathbb T$-model is presented by a geometric formula over $\Sigma$, the first condition of the theorem is satisfied; indeed, any representable functor $Hom_{\textrm{f.p.}{\mathbb T}\textrm{-mod}(\Set)}(c, -)$ is isomorphic to the interpretation of a formula $\phi(\vec{x})$ in the model $U$ (take as $\phi(\vec{x})$ any formula presenting $c$). The second condition of the theorem follows immediately from the fact that the finitely presentable $\mathbb T$-models are jointly conservative for $\mathbb T$. It remains to show that the third condition of the theorem is satisfied. To this end, we observe that for any geometric formulae $\phi(\vec{x})$ and $\psi(\vec{y})$ over the signature of $\mathbb T$, the graph of any arrow $[[\vec{x}.\phi]]_{U}\to [[\vec{y}.\psi]]_{U}$ in the topos $[\textrm{f.p.} {\mathbb T}\textrm{-mod}(\Set), \Set]$ is a subobject of the product $[[\vec{x}.\phi]]_{U}\times [[\vec{y}.\psi]]_{U}$ in $[\textrm{f.p.} {\mathbb T}\textrm{-mod}(\Set), \Set]$ and hence it defines a property of tuples of elements of finitely presentable models of $\mathbb T$ which is preserved by $\mathbb T$-model homomorphisms; therefore, by our assumptions, there exists a formula $\theta(\vec{x}, \vec{y})$ over $\Sigma$ such that its interpretation in $U$ coincides with such subobject. Since $U$ is conservative and such subobject is the graph of an arrow in the topos $[\textrm{f.p.} {\mathbb T}\textrm{-mod}(\Set), \Set]$, the formula $\theta(\vec{x}, \vec{y})$ is $\mathbb T$-provably functional from $\phi(\vec{x})$ to $\psi(\vec{y})$, as required.  
\end{proofs}

\subsection{Expanding a geometric theory to a theory of presheaf type}\label{presheafcompletion}

In this section we shall discuss the problem of expanding a geometric theory $\mathbb T$ to a theory of presheaf type classified by the topos $[\textrm{f.p.} {\mathbb T}\textrm{-mod}(\Set), \Set]$. We shall say that a geometric theory is a \emph{presheaf completion} of a geometric theory $\mathbb T$ if it is an expansion of $\mathbb T$ such that the geometric morphism $p_{\mathbb T}^{{\mathbb T}'}:\Set[{\mathbb T}']\to \Set[{\mathbb T}]$ is isomorphic to the canonical geometric morphism $[\textrm{f.p.} {\mathbb T}\textrm{-mod}(\Set), \Set]\to \Set[{\mathbb T}]$.

The results of section \ref{genth} show that in order to obtain a presheaf completion of a given geometric theory $\mathbb T$ we can add a new sort $\name{c}$ for each finitely presentable $\mathbb T$-model $c$ which is not presented by a geometric formula over the signature of $\mathbb T$, and a relation symbol for each subobject of a finite product of objects of $[\textrm{f.p.} {\mathbb T}\textrm{-mod}(\Set), \Set]$ which are either of the form $Hom_{\textrm{f.p.} {\mathbb T}\textrm{-mod}(\Set)}(c, -)$ (where $c$ is not presented by any geometric formula over the signature of $\mathbb T$), or of the form $U_{A}$ (evaluation functor $\textrm{f.p.} {\mathbb T}\textrm{-mod}(\Set)\to \Set$ at the sort $A$) where $A$ is a sort over the signature of $\mathbb T$ (the sort corresponding to such a representable $Hom(c, -)$ being $\name{c}$ and to a functor $U_{A}$ being $A$). Indeed, by Theorem \ref{critclassifyingtopos}, the theory of the tautological structure over this extended signature will be an expansion of the theory $\mathbb T$ classified by the topos $[\textrm{f.p.} {\mathbb T}\textrm{-mod}(\Set), \Set]$.  Notice that if $\mathbb T$ satisfies the property that every finitely presentable model of $\mathbb T$ is presented by a geometric formula over its signature then this expansion of $\mathbb T$ is localic over $\mathbb T$.

Of course, as it is clear from the results of section \ref{genth}, there are in general many syntactic ways of `completing' a given geometric theory to a theory of presheaf type; the procedure described above represents just a particular choice which is by no means canonical. In fact, what is most interesting in practice is to obtain explicit axiomatizations of presheaf-type completions of a given theory $\mathbb T$ directly from the axioms of $\mathbb T$ (cf. section \ref{examples} below for some examples). Nonetheless, the results established above provide a useful guide in seeking such axiomatizations, as they indicate that in order to complete a geometric theory $\mathbb T$ to a theory classified by the topos $[\textrm{f.p.} {\mathbb T}\textrm{-mod}(\Set), \Set]$ one should expand the language of $\mathbb T$ in order to make each finitely presentable $\mathbb T$-model presented by a formula in the extended signature and every property of finite tuples of elements of finitely presentable models of $\mathbb T$ definable by a geometric formula in the extended signature; we shall see concrete applications of these general remarks in section \ref{examples}.    In fact, Theorem \ref{defpresheaf} ensures that making every finitely presentable model $\mathbb T$ finitely presented by a formula written in a possibly expanded signature is a necessary condition for the resulting theory to be classified by the topos $[\textrm{f.p.} {\mathbb T}\textrm{-mod}(\Set), \Set]$. Moreover, if one is able to prove that any finitely presentable $\mathbb T$-model is strongly finitely presented as a model of the theory $\mathbb T$, when the latter is considered over an extended signature $\Sigma'$, condition $(iii)$ of Theorem \ref{main} is automatically satisfied, while conditions $(i)$ and $(ii)$ can be made to hold at the cost of adding further axioms to $\mathbb T$ over the signature $\Sigma'$ (cf. Theorem \ref{crit}). On the other hand, by Theorem \ref{defpresheaf}, in order to complete $\mathbb T$ to a theory of presheaf type classified by the topos $[\textrm{f.p.} {\mathbb T}\textrm{-mod}(\Set), \Set]$, one should also expand the signature of $\mathbb T$ by adding relation symbols for making any property of finitely presentable $\mathbb T$-models which is preserved by homomorphisms of models definable over the extended signature (if it is not already definable over the original signature).

The following lemma shows that, under appropriate conditions, models which are finitely presented over a given signature remain finitely presented over a larger signature obtained from the former by adding relation symbols characterized
by disjunctive sequents of a certain form.
 
\begin{lemma}\label{lemmafp}
Let $\mathbb T$ be a geometric theory over a signature $\Sigma$, $\Sigma'$ a signature obtained from $\Sigma$ by only adding relation symbols $R$ and ${\mathbb S}$ a geometric theory over $\Sigma'$ obtained from $\mathbb T$ by adding pairs of axioms of the form $(\top \vdash_{\vec{z}} R(\vec{z}) \vee \mathbin{\mathop{\textrm{ $\bigvee$}}\limits_{i\in I}}\phi_{i}^{R}(\vec{z}))$ and $(R \wedge \mathbin{\mathop{\textrm{ $\bigvee$}}\limits_{i\in I}}\phi^{R}_{i}(\vec{z}) \vdash_{\vec{z}} \bot)$, where for each $i\in I$ $\phi_{i}^{R}$ is a geometric formula over $\Sigma$ such that there exists a geometric formula $\psi^{R}_{i}(\vec{z})$ over $\Sigma$ with the property that the sequents ($\phi^{R}_{i} \wedge \psi^{R}_{i} \vdash_{\vec{z}} \bot$) and $(\top \vdash_{\vec{z}} \phi^{R}_{i} \vee \psi^{R}_{i})$ are provable in $\mathbb S$.

Let ${\mathbb R}$ be the cartesianization of $\mathbb S$ (in the sense of Remark \ref{cartesianiz}(d)) and $\phi(\vec{x})=\phi(x_{1}, \ldots, x_{n})$ a ${\mathbb R}$-cartesian formula over $\Sigma'$ with the property that there exists a $\mathbb R$-model $M_{\phi}$ with $n$ generators $\vec{\xi}=(\xi_{1}, \ldots, \xi_{n})\in [[\vec{x}. \phi]]_{M_{\phi}}$ such that for any $\mathbb R$-model $N$, the elements of the interpretation $[[\vec{x}. \phi]]_{N}$ are in natural bijective correspondence with the $\Sigma$-structure homomorphisms $f:M_{\phi}\to N$ such that $f(\vec{\xi})\in [[\vec{x}. \phi]]_{N}$ through the assignment $f \to f(\vec{\xi})$. Then $M_{\phi}$ is finitely presented by the formula $\{\vec{x}.\phi\}$ as a $\mathbb R$-model.  
\end{lemma}

\begin{proofs}
We have to prove that, for any $\mathbb R$-model $N$ and any tuple $\vec{a}\in [[\vec{x}. \phi]]_{N}$, the unique $\Sigma$-structure homomorphism $f:M_{\phi}\to N$ such that $f(\vec{\xi})=\vec{a}$ preserves the satisfaction of all the relation symbols $R$ in $\Sigma'$, i.e. that for any such symbol $R$ of arity $m$ and any $m$-tuple $(y_{1}, \ldots, y_{m})$ of elements of $M_{\phi}$ satisfying $R_{M_{\phi}}$, the $m$-tuple $(f(y_{1}), \ldots, f(y_{m}))$ satisfies $R_{N}$.

As $M_{\phi}$ is by our hypothesis generated by the elements $\xi_{1}, \ldots, \xi_{n}$, each $y_{i}$ is equal to the interpretation in $M_{\phi}$ of a term $t_{i}(\xi_{1}, \ldots, \xi_{n})$ evaluated in the tuple $\vec{\xi}=(\xi_{1}, \ldots, \xi_{n})$. For each $i\in I$, consider the formula ${\psi'_{i}}^{R}(x_{1}, \ldots, x_{n})$ obtained by replacing each of the variables $z_{1}, \ldots, z_{m}$ in the formula $\psi_{i}^{R}$ with the corresponding term $t_{i}$. Let us show that the sequent $\phi \vdash_{\vec{x}} {\psi'}^{R}_{i}$ is valid in every $\mathbb R$-model. This amounts to showing that for any $\mathbb R$-model $P$ and any tuple $\vec{b}\in [[\vec{x}. \phi]]_{P}$, $\vec{b}$ satisfies the formula $\psi_{i}'^{R}$. To prove this, we observe that, since $M_{\phi}$ is by our hypotheses a $\mathbb S$-model, the $m$-tuple $(y_{1}, \ldots, y_{m})$ satisfies the formula $\psi_{i}^{R}$ (for each $i\in I$). Now, $(y_{1}, \ldots, y_{m})\in [[\vec{y}.\psi_{i}^{R}]]_{M_{\phi}}$ means that $(\xi_{1}, \ldots, \xi_{n})\in [[\vec{x}. {\psi_{i}'}^{R}]]_{M_{\phi}}$, which implies, since $f_{\vec{b}}:M_{\phi}\to P$ is a $\Sigma$-structure homomorphism, that $\vec{b}\in [[\psi_{i}'^{R}]]_{N}$, as required. Since ${\mathbb R}$ has enough models as it is cartesian, we can conclude that the sequent $\phi \vdash_{\vec{x}} {\psi'}^{R}_{i}$ is provable in $\mathbb R$ and hence in $\mathbb S$ (for each $i\in I$). It follows that the cartesian sequent $\phi \vdash_{\vec{x}} R(t_{1}\slash z_{1}, \ldots, t_{m}\slash z_{m})$ is provable in $\mathbb R$. By evaluating this sequent in the model $N$ at the tuple $f(\vec{\xi})\in [[\vec{x}. \phi]]_{N}$, and using the fact that $\Sigma$-structure homomorphisms commute with the interpretation of $\Sigma$-terms, we obtain that $(f(y_{1}), \ldots, f(y_{m}))=(f(t_{1}^{M_{\phi}}(\vec{\xi})), \ldots, f(t_{m}^{M_{\phi}}(\vec{\xi})))=(t_{1}^{N}(f(\vec{\xi})), \ldots, t_{m}^{N}(f(\vec{\xi})))$ satisfies $R_{N}$, as required.
\end{proofs}

Notice that if a model of a geometric theory $\mathbb T$ is finitely presentable as a model of its cartesianization then it is strongly finitely presented. In fact, any structure of the form $Hom_{\cal E}(E, M)$, where $M$ is a model of $\mathbb T$ in a topos $\cal E$, is a model of the cartesianization of $\mathbb T$, as it is obtained by applying a global section functor, which is cartesian, to a model of $\mathbb T$. Conversely, if a theory $\mathbb T$ is of presheaf type then any finitely presentable model of $\mathbb T$ is strongly finitely presented (cf. Corollary \ref{reprtopos}) and hence finitely presented relatively to whatever sub-theory $\mathbb S$ of $\mathbb T$ (that is, theory $\mathbb S$ of which $\mathbb T$ is a quotient) whose set-based models admit representations as global sections $Hom_{\cal E}(E, M)$ of models $M$ of $\mathbb T$ in Grothendieck toposes (pairs of theories satisfying these conditions are investigated for instance in \cite{cole} and \cite{coste}). This remark can be often profitably applied to the cartesianization of $\mathbb T$, even though it is not known in general if it is always the case that every model of it admits a sheaf representation of the above kind. 

Summarizing, we have the following theorem.

\begin{theorem}
Let $\mathbb T$ be a theory of presheaf type and $\mathbb S$ be a sub-theory of $\mathbb T$ such that every set-based model of $\mathbb S$ admits a representation as a structure of global sections $Hom_{\cal E}(1, M)$ of a model $M$ of $\mathbb T$ in a Grothendieck topos $\cal E$. Then every finitely presentable model of $\mathbb T$ is finitely presented as a model of $\mathbb S$.
\end{theorem}\qed

\subsection{Faithful expansions}\label{exp}

Let $\mathbb T$ be a geometric theory over a signature $\Sigma$ and ${\mathbb T}'$ a geometric expansion of $\mathbb T$ over a signature $\Sigma'$. 

Suppose that ${\cal C}$ is a subcategory of the category of finitely presentable ${\mathbb T}'$-models such that the canonical functor 
\[
(p_{\mathbb T}^{{\mathbb T}'})_{\Set}:{{\mathbb T}'}\textrm{-mod}(\Set) \to {\mathbb T}\textrm{-mod}(\Set)
\]
induced by the geometric morphism $p_{\mathbb T}^{{\mathbb T}'}:\Set[{\mathbb T}']\to \Set[{\mathbb T}]$ restricts to a functor 
\[
j:{\cal C} \to {\textrm{f.p.} {\mathbb T}\textrm{-mod}(\Set)}.
\]
 
Then we have a commutative diagram
\[  
\xymatrix {
\Sh({\cal C}_{{\mathbb T}'}, J_{{\mathbb T}'})  \ar[rrr]^{p_{\mathbb T}^{{\mathbb T}'}} & & & \Sh({\cal C}_{{\mathbb T}}, J_{{\mathbb T}}) \\
[{\cal C}, \Set] \ar[u]^{s_{\cal C}^{{\mathbb T}'}  } \ar[rrr]_{[j, \Set]} & & & [{\textrm{f.p.} {\mathbb T}\textrm{-mod}(\Set)}, \Set] \ar[u]^{t^{{\mathbb T}}},} 
\]            
where the geometric morphisms 
\[
s_{\cal C}^{{\mathbb T}'}:[{\cal C}, \Set] \to \Sh({\cal C}_{{\mathbb T}'}, J_{{\mathbb T}'})
\]
and 
\[
t^{{\mathbb T}}:[{\textrm{f.p.} {\mathbb T}\textrm{-mod}(\Set)}, \Set] \to \Sh({\cal C}_{{\mathbb T}}, J_{{\mathbb T}})
\]
are the canonical ones and $[j, \Set]$ is the geometric morphism canonically induced by the functor $j$. We shall refer to this diagram as to $(\ast)$.

\begin{theorem}\label{inheritancecontext}
Let $\mathbb T$ be a theory of presheaf type over a signature $\Sigma$ and ${\mathbb T}'$ a geometric expansion of $\mathbb T$ over a signature $\Sigma'$. Suppose that ${\mathbb T}'$ is classified by the topos $[{\cal C}, \Set]$, where $\cal C$ is a full subcategory of ${\textrm{f.p.} {\mathbb T}'\textrm{-mod}(\Set)}$ such that the functor $(p_{\mathbb T}^{{\mathbb T}'})_{\Set}:{{\mathbb T}'}\textrm{-mod}(\Set) \rightarrow {\mathbb T}\textrm{-mod}(\Set)$ restricts to a faithful functor $j:{\cal C} \rightarrow {\textrm{f.p.} {\mathbb T}\textrm{-mod}(\Set)}$. Then for any model $c$ of ${\mathbb T}'$ in $\cal C$ whose associated $\mathbb T$-model $j(c)$ is finitely presented by a geometric formula $\phi(\vec{x})$ over $\Sigma$, there exists a geometric formula $\psi(\vec{x})$ over $\Sigma'$ in the context $\vec{x}$ which presents the ${\mathbb T}'$-model $c$ and such that the sequent $\psi \vdash_{\vec{x}} \phi$ is provable in ${\mathbb T}'$. 
\end{theorem}

Before giving the proof of the theorem, we need to recall the following straightforward lemma, of which we give a proof for the reader's convenience.

\begin{lemma}\label{ad}
Let $R:{\cal A}\to {\cal B}$ and $L:{\cal B}\to {\cal A}$ be a pair of adjoint functors, where $R$ is the right adjoint and $L$ is the left adjoint. Let $\eta:1_{\cal B} \to R\circ L$ be the unit of the adjunction and $b$ an object of $\cal B$. Then $\eta(b)$ is monic if and only if for any arrows $f,g$ in $\cal B$ with codomain $b$, $Lf=Lg$ implies $f=g$. In particular, $\eta$ is pointwise monic if and only if the functor $L$ is faithful. 
\end{lemma}

\begin{proofs}
Let us denote by $\tau_{a,b}$ the bijection between the sets $Hom_{{\cal A}}(Lb, a)$ and $Hom_{{\cal B}}(b, Ra)$ given by the adjunction. By the naturality in $b$ of $\tau_{a, b}$, for any arrow $h:b'\to b$ in $\cal B$, the arrow $\eta_{b}\circ h$ corresponds under $\tau_{Lb, b'}$ to the arrow $Lh$. Therefore, as $\tau_{Lb, b'}$ is a bijection, $\eta_{b}\circ f=\eta_{b}\circ g$ if and only if $Lg=Lf$. From this the thesis follows at once. 
\end{proofs}

We can now prove the theorem.

\begin{proofs}
The geometric morphism 
\[
[j, \Set]:[{\cal C}, \Set] \to [{\textrm{f.p.} {\mathbb T}\textrm{-mod}(\Set)}, \Set]
\]
is essential, that is its inverse image $[j, \Set]^{\ast}$ admits a left adjoint $[j, \Set]_{!}$, namely the left Kan extension along the functor $j$, and the following diagram commutes:
\[  
\xymatrix {
[{\cal C}, \Set]  \ar[rrr]^{[j, \Set]_{!}} & & & [{\textrm{f.p.} {\mathbb T}\textrm{-mod}(\Set)}, \Set] \\
{\cal C}^{\textrm{op}} \ar[u]^{y_{\cal C}}   \ar[rrr]_{j^{\textrm{op}}} & & & {\textrm{f.p.} {\mathbb T}\textrm{-mod}(\Set)}^{\textrm{op}}, \ar[u]^{y}} 
\] 
where $y_{\cal C}$ and $y$ are the Yoneda embeddings.

The functor 
\[
[j, \Set]_{!}:[{\cal C}, \Set] \to [{\textrm{f.p.} {\mathbb T}\textrm{-mod}(\Set)}, \Set]
\]
satisfies the property that for any object $c$ of $\cal C$ and any arrows $\alpha, \beta:P \to y_{\cal C}(c)$, where $P$ is an object of $[{\cal C}, \Set]$, $[j, \Set]_{!}(\alpha)=[j, \Set]_{!}(\beta)$ implies $\alpha=\beta$. Indeed, this is clearly true for $P$ equal to a representable by the commutativity of the above diagram, the full and faithfulness of the Yoneda embeddings $y_{\cal C}$ and $y'$ and the fact that the functor $j$ is faithful by our hypotheses, and one can always reduce to this case by considering a covering of $P$ in $[{\cal C}, \Set]$ by representables.   

Now, by our hypotheses the geometric morphisms $s_{\cal C}^{{\mathbb T}'}$ and $t^{{\mathbb T}}$ defined above are equivalences. Lemma \ref{ad} thus implies that the geometric morphism $p_{\mathbb T}^{{\mathbb T}'}$ is essential and the unit of the adjunction between $(p_{\mathbb T}^{{\mathbb T}'})^{\ast}$ (right adjoint) and $(p_{\mathbb T}^{{\mathbb T}'})_{!}$ (left adjoint) is monic when evaluated at any object of the form $y_{{\cal C}_{{\mathbb T}'}}(\{\vec{y}. \chi\})$, where $\chi(\vec{y})$ is a formula presenting a ${\mathbb T}'$-model. Let $c$ be a ${\mathbb T}'$-model in $\cal C$. Since ${\mathbb T}'$ is classified by the topos $[{\cal C}, \Set]$, $c$ is finitely presented by a formula $\{\vec{x}. \chi\}$ over the signature $\Sigma'$ of ${\mathbb T}'$. The commutativity of the above square and of diagram $(\ast)$ thus implies that that $(p_{\mathbb T}^{{\mathbb T}'})_{!}(y_{{\cal C}_{{\mathbb T}'}}(\{\vec{y}. \chi\}))\cong y_{{\cal C}_{{\mathbb T}}}(\{\vec{x}. \phi\})$, where $\phi(\vec{x})$ is a formula over the signature $\Sigma$ which presents the model $j(c)$ and $y_{{\cal C}_{{\mathbb T}}}:{\cal C}_{\mathbb T}\to \Sh({\cal C}_{\mathbb T}, J_{\mathbb T})$, $y_{{\cal C}_{{\mathbb T}'}}:{\cal C}_{{\mathbb T}'}\to \Sh({\cal C}_{{\mathbb T}'}, J_{{\mathbb T}'})$ are the Yoneda embeddings. By definition of the geometric morphism $p_{\mathbb T}^{{\mathbb T}'}$, we have that $(p_{\mathbb T}^{{\mathbb T}'})^{\ast}(y_{{\cal C}_{{\mathbb T}}}(\{\vec{x}. \phi\}))=y_{{\cal C}_{{\mathbb T}'}}(\{\vec{x}. \phi\})$ (where the latter $\phi$ is considered as a formula over the signature $\Sigma'$). The unit of the adjunction between $(p_{\mathbb T}^{{\mathbb T}'})^{\ast}$ and $(p_{\mathbb T}^{{\mathbb T}'})_{!}$ thus yields a monic arrow $y_{{\cal C}_{{\mathbb T}'}}(\{\vec{y}. \chi\}) \mono {p_{\mathbb T}^{{\mathbb T}'}}^{\ast}({p_{\mathbb T}^{{\mathbb T}'}}_{!}(y_{{\cal C}_{{\mathbb T}'}}(\{\vec{y}. \chi\})))\cong y_{{\cal C}_{{\mathbb T}'}}(\{\vec{x}. \phi\})$ in the topos $\Sh({\cal C}_{{\mathbb T}'}, J_{{\mathbb T}'})$, in other words a monic arrow $\{\vec{y}. \chi\} \mono \{\vec{x}. \phi\}$ in the geometric syntactic category ${\cal C}_{{\mathbb T}'}$. Therefore $\{\vec{y}. \chi\}$ is isomorphic, as an object of ${\cal C}_{{\mathbb T}'}$, to an object $\{\vec{x}. \phi'\}$ such that the sequent $(\phi' \vdash_{\vec{x}} \phi)$ is provable in ${\mathbb T}'$; hence $\phi'(\vec{x})$ presents $c$ as a ${\mathbb T}'$-model, as required.  
\end{proofs}

\section{Examples}\label{examples}

In this section we shall discuss in detail various non-trivial examples of theories of presheaf type in light of the theory developed in the paper.

\subsection{Theories whose finitely presentable models are finite}

As an application of Corollary \ref{fincor}, one can recover at once the well-known results that the following theories are of presheaf type:

\begin{enumerate}
\item The theory of decidable objects (cf. \cite{johwra} and p. 907 of \cite{El});

\item The theory of decidable Boolean algebras (cf. Example D3.4.12 \cite{El});

\item The theory of linear orders;

\item The theory of total orders with endpoints (cf. section VIII.8 of \cite{MM})
\end{enumerate}

Notice that the former two theories are injectivizations of two cartesian theories, namely the empty theory over a signature consisting exactly of one sort and the algebraic theory of Boolean algebras. The latter two theories clearly satisfy the hypothesis of Corollary \ref{corinjectivizations}; hence their injectivizations are of presheaf type as well.

\subsubsection{The theory of abstract circles}

The theory $\mathbb C$ of abstract circles has been introduced by I. Moerdijk and shown in \cite{Moerdijk} to have the property that the points of Connes' topos of cyclic sets (cf. \cite{connes2}) can be identified with the set-based models of $\mathbb C$. In the same paper it is also stated that Connes' topos actually classifies $\mathbb C$, but the argument given therein seems incomplete. 

We shall prove in this section, by using Corollary \ref{corquotients}, that $\mathbb C$ is of presheaf type classified by Connes' topos. Specifically, we shall show that $\mathbb C$ satisfies the hypotheses of the corollary with respect to its Horn part (i.e., the Horn theory consisting of the collection of all Horn sequents which are provable in $\mathbb C$). We will also prove that the injectivization of $\mathbb C$ is of presheaf type as well.

The signature $\Sigma$ of the theory $\mathbb C$ consists of two sorts $P$ and $S$ (variables of type $P$ will be denoted by letters $x, y, \ldots$, while variables of type $S$ will be denoted by letters $a, b, \ldots$), two unary function symbols $0:P\to S$ and $1:P\to S$, two unary function symbols $\delta_{0}:S\to P$ and $\delta_{1}:S\to P$, one unary function symbol $*:S\to S$ and a ternary predicate $R$ of type $S$. The axioms of $\mathbb C$ can be formulated as follows:

\begin{enumerate}
\item \emph{Non-triviality axioms}:
\[
(\top \vdash_{[]} (\exists x)(x=x));
\]
\[
(\top \vdash_{x, y} (\exists a)(\delta_{0}(a)=x \wedge \delta_{1}(a)=y));
\]
\[
(0(x)=1(x) \vdash_{x} \bot)
\]

\item \emph{`Equational' axioms}:
\[
(\top \vdash_{a} {a^{\ast}}^{\ast}=a);
\]
\[
(\top \vdash_{a} \delta_{0}(a^{\ast})=\delta_{1}(a));
\]
\[
(\top \vdash_{x} \delta_{0}(0(x))=x \wedge \delta_{1}(0(x))=x);
\]
\[
(\top \vdash_{x} 0(x)^{\ast}=1(x));
\]
\[
(\delta_{0}(a)=x \wedge \delta_{1}(a)=x \vdash_{x, a} a=0(x) \vee a=1(x));
\]

\item \emph{Axioms for concatenation}:
\[
(R(a, b, c) \wedge R(a, b, c') \vdash_{a, b, c, c'} c=c');
\] 
\[
(R(a, b, c) \vdash_{a, b, c} \delta_{0}(c)=\delta_{0}(a) \wedge \delta_{1}(c)=\delta_{1}(b)); 
\]
\[
(R(a, b, c) \vdash_{a, b, c} R(c^{\ast}, a, b^{\ast}));
\] 
\[
R(a, b, d) \wedge  R(c,d, e) \vdash_{a, b,c,d,e} (\exists e')(R(b,c,e') \wedge R(a, e', e));  
\] 
\[
(R(a, b, 0(x)) \vdash_{a, b, x} a=0(x));
\] 
\[
(\delta_{0}(a)=x \vdash_{a,x} R(0(x), a, a);
\] 
\[
(\delta_{1}(a)=\delta_{0}(b) \vdash_{a, b}  (\exists c)R(a, b,c) \vee (\exists d)R(b^{\ast}, a^{\ast}, d)).
\]  
\end{enumerate}

A model of $\mathbb C$ in $\Set$ is said to be an \emph{abstract circle}. Any set of points $P$ of the circle $S^{1}$ defines an abstract circle $S_{P}$ whose segments $a\in S$ such that $\delta_{0}(a)=x$ and $\delta_{1}(a)=y$ are the oriented arcs on $S^{1}$ from the point $x$ to the point $y$. For any natural number $n\gt 0$, there is exactly one abstract circle, up to isomorphism, whose set of points has $n$ elements; we shall denote it by the symbol $C_{n}$. 

To prove that $\mathbb C$ is of presheaf type, we first notice that the sequent 
\[
(\delta_{0}(a)=\delta_{0}(b) \wedge \delta_{1}(a)=\delta_{1}(b) \vdash_{a, b} a=b)
\]
is provable in $\mathbb C$. We shall refer to it as to the `uniqueness axiom'. This sequent can be easily deduced as a consequence of the seventh axiom of group $3$ and the fifth axiom of group $2$.  

It follows that, modulo the non-triviality axioms and the uniqueness axiom, we can replace any expression of the form $(\exists c)\phi(c)$ arising in an axiom of $\mathbb C$ with the requirement that the unique segment $d$ such that $\delta_{0}(d)=\delta_{0}(c)$ and $\delta_{1}(d)=\delta_{1}(c)$ satisfies $\phi$. In particular, the fourth axiom of group $3$ is provably equivalent, modulo the non-triviality axioms and the uniqueness axiom, to the following sequent:
\[
(\delta_{0}(u)=\delta_{0}(b) \wedge \delta_{1}(u)=\delta_{1}(c) \wedge R(a, b, d) \wedge R(d,c,e) \vdash_{a,b,c,d,e,u} R(b, c, u)\wedge R(a, u,e)). 
\]

Similarly, the seventh axiom of group $3$ is provably equivalent, modulo the non-triviality axioms and the uniqueness axiom, to the following sequent:
\[
(\delta_{0}(c)=\delta_{0}(a) \wedge \delta_{1}(c)=\delta_{1}(b) \wedge \delta_{1}(a)=\delta_{0}(b) \vdash_{a,b,c} R(a, b, c) \vee R(b^{\ast}, a^{\ast}, c^{\ast})).
\]

From these remarks we see that $\mathbb C$ admits a presentation in which all the axioms do not contain quantifications except for the first and second of group $1$. 

Let us show that $\mathbb C$ satisfies the hypotheses of Corollary \ref{corquotients} with respect to its Horn part and the category of finite $\mathbb C$-models.
 
Given a homomorphism $f:c\to Hom_{\cal E}(E, M)$, where $c$ is a finite model of the Horn part of $\mathbb C$, $M$ is a model of $\mathbb C$ in a Grothendieck topos $\cal E$ and $E$ is an object of $\cal E$, we can `localize' $f$ (in the sense of section \ref{secquotients}) a finite number of times (once for the first axiom and once for each point of $c$ for the second axiom) to obtain $\Sigma$-substructure homomorphisms $f_{i}:c_{i}\hookrightarrow Hom_{\cal E}(E_{i}, M)$ such that the structures $c_{i}$ are finite models of the Horn part of $\mathbb C$ satisfying the non-triviality axioms (notice that any substructure of a structure of the form $Hom_{\cal E}(E_{i}, M)$, where $M$ is a model of $\mathbb C$ in $\cal E$, satisfies all the Horn sequents provable in $\mathbb C$). We can clearly further localize each of these homomorphisms so to obtain the satisfaction of all the other axioms of $\mathbb C$; since this can be done without modifying their domains, the final result will be a family of homomorphisms whose domains are structures which satisfy all the axioms of $\mathbb C$. 
  
Next, we notice that for any set-based model $M$ of the Horn part of $\mathbb C$, any point $x\in MP$ and any segment $a\in MS$ there exists a finite substructure $N$ of $M$ such that $NP$ contains $x$ and $NS$ contains $a$ (take $N$ equal to the substructure of $M$ given by the sets
\[
NP=\{x, \delta_{0}(a), \delta_{1}(a)\}
\]
and 
\[
NS=\{a, 0(x), 1(x), 0(\delta_{0}(a)), 0(\delta_{1}(a)),1(\delta_{0}(a)),1(\delta_{1}(a)))\}).
\]
Moreover, any two finite substructures $N_{1}$ and $N_{2}$ of $M$ are contained in a common substructure $N$ of $M$ (take $N$ equal to the substructure of $M$ given by $NP=N_{1}P \cup N_{2}P$ and $N_{1}S=N_{1}S \cup N_{2}S$). 

This discussion, combined with the argument above (specialized to the case ${\cal E}=\Set$), shows that every model of $\mathbb C$ in $\Set$ is a directed union of finite models of $\mathbb C$. It follows in particular that every finitely presentable $\mathbb C$-model is finite (it being a retract of a finite model). On the other hand, every finite model of $\mathbb C$ is finitely presentable as a model of the Horn part of $\mathbb C$ (cf. Theorem 6.4 \cite{OCG}). We can thus conclude that the hypotheses of Corollary \ref{corquotients} are satisfied, whence $\mathbb C$ is of presheaf type classified by the topos of covariant set-valued functors on the category of finite models of $\mathbb C$.   

Due to the presence of conjunctions in the premises of some axioms of $\mathbb C$, we cannot directly apply Corollary \ref{corinjectivizations} to conclude that the the injectivization ${\mathbb C}_{m}$ of $\mathbb C$ is of presheaf type. We shall instead apply Corollary \ref{injcorollary}. Since every $\mathbb C$-model in $\Set$ is a directed union of finite $\mathbb C$-models, the finitely presentable ${\mathbb C}_{m}$-models are exactly the finite $\mathbb C$-models. Moreover, by Proposition \ref{freemono}, the monic $\mathbb C$-model homomorphisms in $\Set$ are precisely the homomorphisms which are sortwise injective; indeed, the formulae $\{x^{P}. \top\}$ and $\{x^{S}. \top\}$ strongly present respectively the $\mathbb C$-models $C_{1}$ and $C_{2}$. To show that the hypotheses of Corollary \ref{injcorollary} are satisfied, it remains to verify that for any Grothendieck topos $\cal E$, object $E$ of $\cal E$ and $\Sigma$-structure homomorphism $x:c\to Hom_{\cal E}(E, M)$, where $c$ is a finite ${\mathbb C}$-model and $M$ is a sortwise decidable $\mathbb C$-model, there exists an epimorphic family $\{e_{i}:E_{i}\to E \textrm{ | } i\in I\}$ in $\cal E$ and for each $i\in I$ a $\mathbb C$-model homomorphism $f_{i}:c\to c_{i}$ of finite $\mathbb C$-models and a sortwise disjunctive $\Sigma$-structure homomorphism $x_{i}:c_{i} \to Hom_{{\cal E}}(E_{i}, M)$ such that $x_{i} \circ f_{i}=Hom_{\cal E}(e_{i}, M)\circ x$ for all $i\in I$.

In order to apply Proposition \ref{findecid}, we make $\mathbb C$ into a one-sorted theory by identifying points $x$ with the segments $0_{x}$ and rewriting the axioms appropriately. The proposition yields an epimorphic family $\{e_{i}:E_{i}\to E \textrm{ | } i\in I\}$ in $\cal E$ and for each $i\in I$ a sortwise surjective homomorphism $q_{i}:c\to c_{i}$, where $c_{i}$ is a finite $\Sigma$-structure, and a disjunctive $\Sigma$-structure homomorphism (in the sense of Lemma \ref{explicit2}) $J_{i}:c_{i}\mono Hom_{\cal E}(E_{i}, M)$ such that $J_{i}\circ q_{i}=Hom_{\cal E}(e_{i}, M)\circ f$ for all $i\in I$. Now, since $c$ is a $\mathbb C$-model and $q_{i}$ is sortwise surjective, the structure $c_{i}$ satisfies the non-triviality axioms. We can clearly suppose $E_{i}\ncong 0$ without loss of generality, and hence the arrows $J_{i}$ to be injective. If we consider the image factorizations of the $\Sigma$-structure homomorphisms $J_{i}$ (in the sense of Lemma \ref{lemmafact}), we thus obtain substructures $c_{i}'\mono Hom_{{\cal E}}(E_{i}, M)$ whose underlying sets are the same as those of $c_{i}$ (since the $J_{i}$ are injective) and $\Sigma$-structure homomorphisms $q_{i}':c\to c_{i}'$. Since the $c_{i}'$ have the same underlying set as $c_{i}$, they all satisfy the non-triviality axiom, and, at the cost of refining the epimorphic family $\{e_{i}:E_{i}\to E \textrm{ | } i\in I\}$, we can suppose them to satisfy all the other axioms of $\mathbb C$ (cf. the argument given above for showing that $\mathbb C$ satisfies the hypotheses of Corollary \ref{corquotients} with respect to its Horn part and the category of finite $\mathbb C$-models). So the hypotheses of Corollary \ref{injcorollary} are satisfied, and we can conclude that ${\mathbb C}_{m}$ is of presheaf type classified by the topos of covariant set-valued functors on the category of finite models of $\mathbb C$ and injective homomorphisms between them.

\subsubsection{The geometric theory of finite sets}

In this section we shall revisit, from the point of view of the theory developed in the paper, a well-known example of a non-trivial theory of presheaf type, namely the geometric theory $\mathbb T$ of finite sets. Recall from \cite{El} (Example D1.1.7(k)) that the signature $\Sigma$ of $\mathbb T$ consists of one sort $A$ and a $n$-ary relation symbol $R_{n}$ for each $n\gt 0$. The axioms of $\mathbb T$ are the following: for each $n$, one has the axiom 
\[
\sigma_{n}:=(R_{n}(x_{1}, \ldots, x_{n}) \vdash_{\vec{x}, y} \mathbin{\mathop{\textrm{ $\bigvee$}}\limits_{1\leq i\leq n}}y=x_{i}),  
\]
expressing the requirement that if an $n$-tuple of individuals satisfies the relation $R_{n}$ then it exhausts the members of (the set interpreting) the sort $A$. (The case $n = 0$ of this 
axiom is $(R_{0} \vdash_{[]} \bot)$, which says that if $R_{0}$ holds then the interpretation of the sort $A$ must be empty.) We also have, for each $n\gt 0$, the axiom 
\[
(\top \vdash_{[]}  \mathbin{\mathop{\textrm{ $\bigvee$}}\limits_{1\leq i\leq n}}(\exists x_{1})\cdots (\exists x_{n})R_{n}(x_{1}, \ldots, x_{n})).
\]
Finally, to ensure that the interpretations of the $R_{n}$ are uniquely determined by that of the sort $A$ (i.e. that $R_{n}$ holds for all $n$-tuples which exhaust the elements of the interpretation of sort $A$, and not just for some of them), one adds the axioms 
\[
(R_{n}(x_{1},\ldots, x_{n})) \vdash_{x_{1}, \ldots, x_{n}} R_{m}(x_{f(1)}, \ldots, x_{f(m)}))
\]
whenever $f: \{1,2,... ,m\} \to \{1,2,... ,n\}$ is a surjection, and 
\[
(R_{n}(x_{1}, \ldots, x_{n}) \wedge x_{i}=x_{j} \vdash_{x_{1}, \ldots, x_{n}} R_{n-1}(x_{1}, \ldots, x_{i-1}, x_{i+1}, \ldots, x_{n}))
\]
whenever $1 \leq i \lt j \leq n$. 

We can deduce that $\mathbb T$ is of presheaf type as an application of Corollary \ref{corquotients}. 

The models of $\mathbb T$ in $\Set$ can be identified with the finite sets, while the $\mathbb T$-model homomorphism are the precisely the surjective functions between them.

We can regard $\mathbb T$ as a quotient of its Horn part (i.e. of the theory consisting of all the Horn sequents over the signature of $\mathbb T$ which are provable in $\mathbb T$). The criterion for finite presentability given by Lemma 6.2 \cite{OCG} ensures that every finite set, regarded as a model of $\mathbb T$, is finitely presentable as a model of the Horn part of $\mathbb T$. Indeed, the finiteness of the structure immediately implies that the second condition of the lemma is satisfied (cf. the proof of Theorem 6.4 \cite{OCG}), while the satisfaction of the first condition follows from the fact that every function from a finite model of cardinality $n$ of the Horn part of $\mathbb T$ to a set-based model $M$ of the Horn part of $\mathbb T$ which preserves the predicate $R_{n}$ preserves the predicate $R_{m}$ for any $m\lt n$ (by the `introduction' and `elimination' rules expressed by the last two groups of axioms for $\mathbb T$). 

To apply Corollary \ref{corquotients}, it thus remains to prove that the second condition in the statement of the corollary is satisfied. Given a homomorphism $f:a\to Hom_{\cal E}(E, M)$, where $a$ is a finite model of the Horn part of $\mathbb T$ and $M$ is a $\mathbb T$-model in a Grothendieck topos $\cal E$, if the cardinality of $aA$ is $n$ then for any $m\gt n$, the sequent $\sigma_{m}$ is satisfied in $a$ provided that $\sigma_{n}$ is. Indeed, if $m\geq n$ then for any $m$-tuple $(y_{1}, \ldots, y_{m})$ of elements in $aA$ there exists a sub-$n$-tuple of elements of $aA$ obtained by removing $m-n$ repetitions, which satisfies the relation $R_{n}$ if the $m$-tuple $(y_{1}, \ldots, y_{m})$ satisfies the relation $R_{m}$, since $a$, as a model of the Horn part of $\mathbb T$, satisfies the Horn sequent expressing the `elimination' rule (i.e., the last group of axioms for $\mathbb T$). Similarly, by invoking the `introduction' rules, one can prove that if the cardinality of $aA$ is $n$ then for any $k\lt n$, the sequent $\sigma_{k}$ is satisfied in $a$ provided that $\sigma_{n}$ is. 

Let us show that we can inductively localize the morphism $f$ to eventually arrive at $\Sigma$-structure homomorphisms $f_{i}:a_{i}\to Hom_{\cal E}(E_{i}, M)$ defined on $\Sigma$-structures $a_{i}$ which satisfy all the axioms of $\mathbb T$. Notice that if the $\Sigma$-structure homomorphisms $h_{i}:a\to a_{i}$ in such localization are quotient maps (in the sense that a tuple in $a_{i}$ satisfies a relation if and only if it is the image under $h_{i}$ of a tuple satisfying that relation in $a$) then the fact that the $a_{i}$ satisfy the last two groups of axioms of $\mathbb T$ will follow automatically from the fact that $a$ does. It will thus be enough to show, by the above considerations, that each $a_{i}$ satisfies the sequent $\sigma_{n_{i}}$, where $n_{i}$ is the cardinality of $a_{i}A$, and the second axiom of $\mathbb T$.

Starting from a $\Sigma$-structure homomorphism $f:a\to Hom_{\cal E}(E, M)$, where $a$ is a finite model of the Horn part of $\mathbb T$, from the fact that $M$ is a model of $\mathbb T$ it follows that there exists an epimorphic family $\{e_{i}:E_{i}\to E \textrm{ | } i\in I\}$ and for each $i\in I$ a natural number $n_{i}$ and a $n_{i}$-tuple $(\xi_{1},\ldots, \xi_{n_{i}})$ of generalized elements $E_{i}\to MA$ such that $<\xi_{1},\ldots, \xi_{n_{i}}>$ factors through the interpretation of $R_{n_{i}}$ in $M$. By taking $a_{i}$ to be the $\Sigma$-substructure of $Hom_{\cal E}(E_{i}, M)$ on the finite subset consisting of the elements $\xi_{1},\ldots, \xi_{n_{i}}$ plus all the elements in the image of the homomorphism $Hom_{\cal E}(e_{i}, M)\circ f$, we clearly obtain a structure satisfying the second axiom of $\mathbb T$. This structure will also satisfy the last two groups of axioms for $\mathbb T$ (the validity of Horn sequents is inherited by substructures). Now, in order to obtain from this family of $\Sigma$-structure homomorphisms a localization such that the domains of its homomorphisms satisfy all the axioms of $\mathbb T$, for each $i\in I$, it suffices to localize each $f_{i}$ a finite number of times (one for each $(n_{i}+1)$-tuple $\vec{u}$ of elements of $a_{i}$, where $n_{i}$ is the cardinality of $a_{i}$), endowing the sets $d_{i}$ arising in the surjective-injective factorizations of the homomorphisms with the quotient structure induced by the domains $a_{i}$ via the relevant quotient map $q_{i}$ (in the sense that the interpretation of each relation symbol $R$ over the signature of $\mathbb T$ in such a set $d_{i}$ is defined to be equal to the image of the interpretation of $R$ in $a_{i}$ under the quotient map $q_{i}$). By the above remarks, these structures will satisfy all the axioms of $\mathbb T$.

This completes the proof of the fact that the hypotheses of Corollary \ref{corquotients} are satisfied by the theory $\mathbb T$ with respect to its Horn part and the category $\cal C$ of finite $\mathbb T$-models; therefore the theory $\mathbb T$ is of presheaf type classified by the topos $[{\cal C}, \Set]$.

\subsection{The theory of Diers' fields}\label{diers}

As observed in \cite{diers}, the coherent theory $\mathbb T$ of fields is not of presheaf type. In fact, one can easily identify two properties which are preserved by homomorphisms of fields but which are not definable by geometric formulae in the signature of the coherent theory of fields: the property of a field to have characteristic $0$, and the property of a tuple $(x_{1}, \ldots, x_{n}, x_{n+1})$ of elements of a field that the element $x_{n+1}$ is transcendental over the subfield generated by the elements $x_{1}, \ldots, x_{n}$. This can be easily seen by arguing as follows. Assuming the axiom of choice, the theory $\mathbb T$ has enough $\Set$-based models (it being coherent). Hence if the property of having characteristic $0$ were definable by a geometric sentence $\phi$ over the signature of  fields then the sequent $(\top \vdash_{[]} \phi \vee \mathbin{\mathop{\textrm{ $\bigvee$}}\limits_{p\in {\mathbb P}}}\phi_{p})$, where $\mathbb P$ is the set of prime numbers and $\phi_{p}$ (for each $p\in {\mathbb P}$) is the sentence $p.1=0$ expressing the property of having characteristic $p$, would be provable in $\mathbb T$. But, $\mathbb T$ being coherent, the infinitary disjunction on the right-hand side of the sequent would then be provably equivalent to a finite subdisjunction, which is absurd as it would imply that the set of all possible characteristics of a field is finite. A similar argument works for the other property, which, like the former, is the complement of a property definable by a strictly infinitary geometric formula.

In order to make such properties definable and possibly obtain a presheaf completion of the theory $\mathbb T$, it is thus necessary to enlarge the signature of $\mathbb T$ with new relation symbols. In fact, Johnstone introduces in \cite{diers} a $0$-ary predicate $R_{0}$, expressing the property of a field to have characteristic $0$, and for each natural number $n\geq 0$ a $n+1$-predicate $R_{n+1}(x_{1}, \ldots, x_{n}, x_{n+1})$ expressing the property of $x_{n+1}$ of being transcendental over the subfield generated by the elements $x_{1}, \ldots, x_{n}$. Let $\Sigma'$ be the resulting signature. Formally, one has to impose the following axioms over $\Sigma'$ to ensure that these predicates have indeed the required meaning (below we use the abbreviation $Inv(z)$ for the formula $(\exists x)(x\cdot z=1 \wedge z\cdot x=1)$):
\[
(\top \vdash R_{0} \vee \mathbin{\mathop{\textrm{ $\bigvee$}}\limits_{p\in {\mathbb P}}}p.1=0 ); 
\] 
\[
(R_{0} \wedge p.1=0 \vdash  \bot)
\]
(for each $p\in {\mathbb P}$);
\[
(\top \vdash_{x_{1}, \ldots, x_{n+1}} R_{n}(x_{1}, \ldots, x_{n}, x_{n+1}) \vee (\mathbin{\mathop{\textrm{ $\bigvee$}}\limits_{m\in {\mathbb N}, \vec{c_{1}^{n}}, \ldots, \vec{c_{m}^{n}}}}  (\mathbin{\mathop{\textrm{ $\sum$}}\limits^{m}\limits_{j=0}} P_{\vec{c^{n}_{j}}}x_{n+1}^{j}=0)                       \wedge \mathbin{\mathop{\textrm{ $\bigvee$}}\limits_{i\in \{0, 1, \ldots, m\} }} Inv(P_{\vec{c^{n}_{i}}} )))
\]
(for each natural number $n\geq 0$), where the former disjunction is taken over all the natural numbers $m\geq 0$ and all the tuples $\vec{c_{i}^{n}}$ (for $i\in \{0, \ldots, m\}$) of integer coefficients (i.e., coefficients of the form $1+\cdots +1$ an integer number of times) of polynomials in $n$ variables $x_{1},\ldots, x_{n}$ of degree $\leq m$ and the expression $P_{\vec{c^{n}_{i}}}$ (for each $i\in \{0, \ldots, m\}$) denotes the polynomial term in the variables $x_{1},\ldots, x_{n}$ corresponding to the tuple $\vec{c^{n}_{i}}$, and 
\[
(R_{n}(x_{1}, \ldots, x_{n}, x_{n+1}) \wedge (\mathbin{\mathop{\textrm{ $\bigvee$}}\limits_{m\in {\mathbb N}, \vec{c_{1}^{n}}, \ldots, \vec{c_{m}^{n}}}} (\mathbin{\mathop{\textrm{ $\sum$}}\limits^{m}\limits_{j=0}} P_{\vec{c^{n}_{j}}}x_{n+1}^{j}=0) \wedge \mathbin{\mathop{\textrm{ $\bigvee$}}\limits_{i\in \{0, 1, \ldots, m\} }} Inv(P_{\vec{c^{n}_{i}}} ))) \vdash_{x_{1}, \ldots, x_{n+1}} \bot ).
\] 

Let $\mathbb D$ be the theory, called in \cite{diers} of \emph{Diers fields}, obtained from $\mathbb T$ by adding these new predicates and the above-mentioned axioms. 

Note that the geometric formula
\[
\mathbin{\mathop{\textrm{ $\bigvee$}}\limits_{m\in {\mathbb N}, \vec{c_{1}^{n}}, \ldots, \vec{c_{m}^{n}}}} (\mathbin{\mathop{\textrm{ $\sum$}}\limits^{m}\limits_{j=0}} P_{\vec{c^{n}_{j}}}x_{n+1}^{j}=0) \wedge \mathbin{\mathop{\textrm{ $\bigvee$}}\limits_{i\in \{0, 1, \ldots, m\} }} Inv(P_{\vec{c^{n}_{i}}} ))
\]
is $\mathbb D$-provably equivalent to a disjunction of geometric formulae, namely, the formulae $(\mathbin{\mathop{\textrm{ $\sum$}}\limits^{m}\limits_{j=0}} P_{\vec{c^{n}_{j}}}x_{n+1}^{j}=0) \wedge Inv(P_{\vec{c^{n}_{i}}})$ (for each $m$, $\vec{c_{1}^{n}}, \ldots, \vec{c_{m}^{n}}$ and $i\in \{0, \ldots, m\}$), each of which has a $\mathbb D$-provable complement, namely the formula $Inv(\mathbin{\mathop{\textrm{ $\sum$}}\limits^{m}\limits_{j=0}} P_{\vec{c^{n}_{j}}}x_{n+1}^{j}=0) \vee P_{\vec{c^{n}_{i}}}=0$. 

Notice also that, for any tuples $\vec{c_{i}^{n}}$ (for $i\in \{0, \ldots, m\}$) of integer coefficients of polynomial expressions $P_{\vec{c^{n}_{i}}}$ in the variables $x_{1}, \ldots, x_{n}$, the cartesian sequent
\[
(\ast)\textrm{ } (R_{n}(x_{1}, \ldots, x_{n}, x_{n+1}) \wedge  Inv(P_{\vec{c^{n}_{i}}}) \vdash_{x_{1}, \ldots, x_{n+1}} Inv(\mathbin{\mathop{\textrm{ $\sum$}}\limits^{m}\limits_{j=0}} P_{\vec{c^{n}_{j}}}x_{n+1}^{j}=0))
\]
is provable in $\mathbb D$.

Following \cite{diers}, we observe that the theory $\mathbb D$ satisfies the property that every finitely presentable model of it, i.e. any finitely generated field, is finitely presented as a model of its cartesianization. This will follow from Lemma \ref{lemmafp} once we have proved that every finitely generated field $F$ is presented by a finite set of generators, in the weak sense of the lemma, by a formula over the signature of $\mathbb D$. To this end, we regard $\mathbb T$ as axiomatized in the signature of von Neumann regular rings, which contains a unary function for the operation of pseudoinverse; indeed, over this signature, every finitely generated field (in the sense of field theory) becomes finitely generated (in the sense of model theory).

Notice that, by the above remarks, $\mathbb D$ satisfies the first set of hypotheses of Lemma \ref{lemmafp} with respect to the theory $\mathbb T$. 

We shall prove that every finitely generated field $F$ is (weakly) presented (as a model of the cartesianization of $\mathbb D$) by a geometric formula over the signature of $\mathbb D$ by induction on the number $n$ of generators of $F$. If $n=0$ then $F$ is equal to its prime field; so, either $F$ has characteristic $p$, in which case it is equal to ${\mathbb Z}_{p}$, or $F$ has characteristic $0$, in which case it is equal to $\mathbb Q$. Now, the field ${\mathbb Z}_{p}$ is clearly (weakly) presented by the formula $p.1=0$, while the field $\mathbb Q$ is (weakly) presented by the formula $R_{0}$ since the sequent $R_{0} \vdash Inv(n.1)$ is provable in the cartesianization of $\mathbb D$ for each non-zero natural number $n$. Now, consider a field $F$ generated by $n+1$ elements $x_{1}, \ldots, x_{n}, x_{n+1}$. If we denote by $F_{0}$ its prime field, we have that $F=F_{0}(x_{1}, \ldots, x_{n})(x_{n+1})$.    

Suppose that $F_{0}(x_{1}, \ldots, x_{n})$ is (weakly) presented by a formula in $n$ variables $\phi(\vec{x})$ with the elements $x_{1}, \ldots, x_{n}$ as generators.

There are two cases: either the element $x_{n+1}$ is transcendental over the field $F_{0}(x_{1}, \ldots, x_{n})$ or not. 

In the first case, $F$ is isomorphic to the field of rational functions in one variable with coefficients in $F_{0}(x_{1}, \ldots, x_{n})$, and it is (weakly) presented by the formula $\phi \wedge R_{n}(x_{1}, \ldots, x_{n}, x_{n+1})$ with generators $x_{1}, \ldots, x_{n}, x_{n+1}$; in other words, for any model $(A, \{(R_{n})_{A} \textrm{ } n\in \mathbb N\})$ of the cartesianization of $\mathbb D$, the function which assigns a ring homomorphism $f:F\to A$ with the property that $f(x_{1}, \ldots, x_{n}, x_{n+1})\in (R_{n+1})_{A}\cap [[\vec{x}. \phi]]_{A}$ to the element $f(x_{1}, \ldots, x_{n}, x_{n+1})$ is injective and surjective on $(R_{n+1})_{A}\cap [[\vec{x}. \phi]]_{A}$. This can be proved as follows. Since $F$ is generated by the elements $x_{1}, \ldots, x_{n+1}$, the injectivity is clear, so it remains to prove the surjectivity, i.e. that for any $(n+1)$-tuple $(a_{1}, \ldots, a_{n}, a_{n+1})\in (R_{n+1})_{A}\cap [[\vec{x}. \phi]]_{A}$ there exists a ring homomorphism $F\to A$ which sends $x_{i}$ to $a_{i}$ for each $i\in \{1, \ldots, n+1\}$. By the induction hypothesis, since $(a_{1}, \ldots, a_{n})\in [[\vec{x}. \phi]]_{A}$, there exists a unique ring homomorphism $g:F_{0}(x_{1}, \ldots, x_{n}) \to A$ such that $g(x_{i})=a_{i}$ for each $i\in \{1, \ldots, n\}$. By definition of $F$, there exists a ring homomorphism $f:F\to A$ which extends $g$ and sends $x_{n+1}$ to $a_{n+1}$ if and only if for every polynomial $P$ with a non-zero coefficient $P_{\vec{c^{n}_{i}}}(x_{1}, \ldots, x_{n})$ in $F_{0}(x_{1}, \ldots, x_{n})$, $\mathbin{\mathop{\textrm{ $\sum$}}\limits^{m}\limits_{j=0}} g(P_{\vec{c^{n}_{j}}})a^{j}$ is invertible in $A$. But since $P_{\vec{c^{n}_{i}}}(x_{1}, \ldots, x_{n})$ is non-zero (equivalently, invertible) in the field $F_{0}(x_{1}, \ldots, x_{n})$, its image $g(P_{\vec{c^{n}_{i}}}(x_{1}, \ldots, x_{n}))=P_{\vec{c^{n}_{i}}}(g(x_{1}), \ldots, g(x_{n}))$ under the homomorphism $g$ is invertible in $A$; therefore, since the sequent $(\ast)$ holds in $A$, the condition $(f(x_{1}), \ldots, f(x_{n}), a)\in (R_{n})_{A}$ entails the fact that $\mathbin{\mathop{\textrm{ $\sum$}}\limits^{m}\limits_{j=0}} g(P_{\vec{c^{n}_{j}}})a^{j}$ is invertible in $A$, as required.   

In the second case, consider the minimal polynomial $P$ for $x_{n+1}$ over $F_{0}(x_{1}, \ldots, x_{n})$; then $F$ is isomorphic to the quotient of $F_{0}(x_{1}, \ldots, x_{n})$ by the ideal generated by the polynomial $P$. It is immediate to see that $F$ is (weakly) presented by the formula in $n+1$ variables $\phi \wedge P(x_{n+1})=0$ with generators $x_{1}, \ldots, x_{n}, x_{n+1}$.

These arguments show that the theories $\mathbb T$ and $\mathbb D$ satisfy the hypotheses of Lemma \ref{lemmafp}. It follows that all the finitely generated fields are finitely presented models of the cartesianization of $\mathbb D$, as required. Condition $(iii)$ of Theorem \ref{main} is therefore satisfied (cf. Proposition \ref{criteriongeometricity}(i)). Alternatively, we could have deduced the fact that the theory $\mathbb D$ satisfies condition $(iii)$ of Theorem \ref{main} from Corollary \ref{corcondition3} and Theorem \ref{condition3}. 

In \cite{diers}, Johnstone shows that $\mathbb D$ is a theory of presheaf type by assuming a form of the axiom of choice to ensure that $\mathbb D$ has enough set-based models. Theorem \ref{main} allows to prove that $\mathbb D$ is of presheaf type directly, without assuming any non-constructive principles. Having already proved that condition $(iii)$ of the theorem is satisfied, it remains to see that conditions $(i)$ and $(ii)$ hold as well.   

The fact that condition $(ii)(a)$ of Theorem \ref{thmcond2} holds follows immediately from the above-mentioned discussion in view of Remark \ref{mon}(b). Condition $(ii)(c)$ of Theorem \ref{thcond1} is automatically satisfied while condition $(ii)(a)$ of Theorem \ref{thmcond2} follows from condition $(ii)(a)$ of Theorem \ref{thmcond2} (cf. Remarks \ref{remcond1}(b)-(c)), By Remark \ref{mon}(a), to prove that condition $(ii)(b)$ of Theorem \ref{thmcond2} holds, it suffices to verify that condition $(ii)(b)$ of Theorem \ref{thcond1} does. It thus remains to verify that condition $(ii)(b)$ of Theorem \ref{thcond1} holds, i.e. that for any finitely generated fields $c$ and $d$, $\mathbb D$-model $M$ in a Grothendieck topos $\cal E$ and $\Sigma'$-structure homomorphisms $x:c\to Hom_{\cal E}(E, M)$ and $y:d\to Hom_{\cal E}(E, M)$, there exists an epimorphic family $\{e_{i}:E_{i} \to E \textrm{ | } i\in I\}$ in $\cal E$ and for each $i\in I$ a finitely generated field $b_{i}$, field homomorphisms $u_{i}:c\to b_{i}$, $v_{i}:d\to b_{i}$ and a $\Sigma'$-structure homomorphism $z_{i}:b_{i}\to Hom_{\cal E}(E_{i}, M)$ such that $Hom_{\cal E}(e_{i}, M)\circ x=z_{i}\circ u_{i}$ and $Hom_{\cal E}(e_{i}, M)\circ y=z_{i}\circ v_{i}$. 

We can prove this by induction on the sum $n$ of the minimal number of generators of $c$ and of $d$. 

Before proceeding with the proof, it is convenient to remark the following fact: for any non-zero model $(A, \{(R_{n})_{A} \textrm{ } n\in \mathbb N\})$ of the cartesianization of $\mathbb D$ (for instance, a $\Sigma'$-structure of the form $Hom_{\cal E}(E, M)$, where $M$ is a model of $\mathbb D$ in the topos $\cal E$ and $E$ is a non-zero object of $\cal E$) and any field $e$, considered as a model of $\mathbb D$, all the $\Sigma'$-structure homomorphisms $e\to A$ reflect the satisfaction of the relations $R_{n}$, i.e. $f(x_{1}, \ldots, x_{n})\in (R_{n})_{A}$ implies $(x_{1}, \ldots, x_{n})\in (R_{n})_{e}$. This easily follows from the disjunctive axioms of $\mathbb D$ defining $R_{n}$ and the cartesian axiom $(\ast)$. Note also that such homomorphisms are always injective (since their domain is a field and their codomain is a non-zero ring).  

In proving our claim, we can suppose without loss of generality all the objects $E$ arising in $\Sigma'$-structure homomorphisms to structures of the form $Hom_{\cal E}(E, M)$ to be non-zero (since removing zero arrows from an epimorphic family leaves the family epimorphic). 

If $n=0$ (that is, if both $c$ and $d$ are equal to their prime fields) then $c$ and $d$ have the same characteristic; indeed, equalities of the form $n.1=0$ are preserved and reflected by the homomorphisms $x$ and $y$ (cf. the above remarks). So $c$ and $d$ are isomorphic, whence the claim is trivially satisfied. Let us now assume that the condition is true for all $k\leq n$ and prove it for $n+1$. We can represent $c$ as $c'(x)$ in such a way that a set of generators for $c$ can be obtained by adding $x$ to a set of generators for $c'$; then by the induction hypothesis there exists an epimorphic family $\{e_{i}:E_{i}\to E \textrm{ | } i\in I\}$ and for each $i\in I$ a finitely generated field $u_{i}$, field homomorphisms $f_{i}:c'\to u_{i}$ and $g_{i}:d\to u_{i}$ and a $\Sigma'$-structure homomorphism $r_{i}:u_{i}\to Hom_{\cal E}(E_{i}, M)$ such that $r_{i}\circ f_{i}=Hom_{\cal E}(E, M)\circ f|_{c'}$ and $r_{i}\circ g_{i}=Hom_{\cal E}(E, M)\circ g$. Now, consider for each $i\in I$ the element $f(x)\circ e_{i}\in Hom_{\cal E}(E_{i}, M)$. Suppose that $c'$ has $m$ generators $\chi_{1}, \ldots, \chi_{m}$; then for each $i\in I$, by the disjunctive axiom of $\mathbb D$ involving $R_{m}$, there exists an epimorphic family $\{f_{i, j}:F_{i, j} \to  E_{i} \textrm{ | } j\in J_{i}\}$ such that for any $j\in J_{i}$, either $R_{m}((r_{i}\circ f_{i})(\chi_{1})\circ e_{i}\circ f_{i, j}, \ldots, (r_{i}\circ f_{i})(\chi_{m})\circ e_{i}\circ f_{i, j}, f(x)\circ e_{i}\circ f_{i, j})$ or $f(x)\circ e_{i}\circ f_{i, j}$ is the root of a non-zero polynomial $P$ with coefficients belonging to the von Neumann regular sub-ring of $Hom_{\cal E}(F_{i, j}, M)$ generated by the elements $(r_{i}\circ f_{i})(\chi_{1})\circ e_{i}\circ f_{i, j}, \ldots, (r_{i}\circ f_{i})(\chi_{m})\circ e_{i}\circ f_{i, j}$. In the first case, $f(x)\circ e_{i}\circ f_{i, j}$ is transcendental over $c'$ via the embedding $Hom_{\cal E}(f_{i, j}, M)\circ r_{i}\circ f_{i}$ and over $u_{i}$ via the embedding $Hom_{\cal E}(f_{i, j}, M)\circ r_{i}$ (apply the above remarks to these two embeddings); it follows that the homomorphism $f_{i}$ extends to a homomorphism from $c=c'(x)$ to $u_{i}$.

In the second case, the homomorphism $Hom_{\cal E}(e_{i}\circ f_{i, j}, M) \circ f$ being injective, there exists a non-zero polynomial $P$ with coefficients in $c'$ such that $f(x)\circ e_{i}\circ f_{i, j}$ is a root of the image of $P$ under $Hom_{\cal E}(f_{i, j}, M)\circ r_{i}\circ f_{i}=Hom_{\cal E}(e_{i}\circ f_{i, j}, M) \circ f|_{c'}$. It follows that $x$ is a root of $P$ in $c$. We can thus clearly suppose $P$ to be irreducible without loss of generality, and represent $c=c'(x)$ as $c=c'[z]\slash P(z)$, via an isomorphism sending $z$ to $x$; denoting $P'$ the image of $P$ under the homomorphism $f_{i}$, we thus obtain that the arrow $f_{i}$ yields an arrow $c=c'[z]\slash P(z) \to u_{i}[w]\slash P'(w)$ and the homomorphism $Hom_{\cal E}(f_{i, j}, M) \circ r_{i}$ factors through the quotient map $u_{i}\to u_{i}[w]\slash P'(w)$ yielding a ring homomorphism which is in fact a $\Sigma'$-structure homomorphism (by Lemma \ref{lemmafp}, cf. the argument given above for proving that every finitely generated field is finitely presented as a model of the cartesianization of $\mathbb D$). 

From these remarks it is now straightforward to obtain a set of data satisfying the requirements of our condition.       

We could have alternatively proved that $\mathbb D$ satisfies condition $(ii)(b)$ of Theorem \ref{thcond1} either by using Remark \ref{remcond1}(d) or by applying Theorem \ref{flatcond}.

\subsection{The theory of algebraic extensions of a given field}

Let $F$ be a field. We define the theory ${\mathbb T}_{F}$ of \emph{algebraic extensions} of $F$ as the expansion of the coherent theory of fields obtained by adding one constant symbol $\overline{a}$ for each element $a\in F$ and the following axioms:
\[
(\top \vdash \overline{1_{F}}=1);
\]
\[
(\top \vdash \overline{0_{F}}=0);
\]
\[
(\top \vdash \overline{a}+\overline{b}=\overline{a+_{F} b}),
\]
for any elements $a, b\in F$ (the symbol $+_{F}$ denotes the addition operation in $F$);
\[
(\top \vdash \overline{a}\cdot \overline{b}=\overline{a\cdot_{F} b}),
\]
for any elements $a, b\in F$ (the symbol $\cdot_{F}$ denotes the multiplication operation in $F$), plus the algebraicity axiom
\[
(\top \vdash_{x}  \mathbin{\mathop{\textrm{ $\bigvee$}}\limits_{n\in {\mathbb N}, a_{0}, \ldots, a_{n-1}, a_{n} \in F}} \overline{a_{n}}x^{n}+\overline{a_{n-1}}x^{n-1}+ \cdots + \overline{a_{0}}=0).
\]

We shall prove that ${\mathbb T}_{F}$ is of presheaf type. Clearly, the finitely presentable models of ${\mathbb T}_{F}$ are exactly the finitely generated algebraic extensions of $F$, that is the finite extensions of $F$. One can prove, by adapting the argument used in the proof of the result that every finitely generated field is finitely presented as a model of the cartesianization of the theory of Diers' fields established in section \ref{diers}, that every finite extension of $F$ is finitely presented as a model of the cartesianization of ${\mathbb T}_{F}$. Specifically, every finite extension  $F(x_{1}, \ldots, x_{n})$ of $F$ is presented by the conjunction of the formulae of the form $P_{i}(x_{1}, \ldots, x_{i+1})(x_{i})=0$ (for $i=0, \ldots, n-1$), where $P$ is the minimal polynomial for the element $x_{i+1}$ over the field $F(x_{1}, \ldots, x_{i})$.
 
Condition $(iii)$ of Theorem \ref{main} is thus satisfied by the theory ${\mathbb T}_{F}$ with respect to the category of finite extensions of $F$ (cf. Proposition \ref{criteriongeometricity}(i)). 

In verifying that conditions $(i)$ and $(ii)$ of Theorem \ref{main} are satisfied, one is reduced as in section \ref{diers} to check that condition $(ii)(b)$ of Theorem \ref{thcond1} holds; this can be done again by adapting the  argument given in section \ref{diers} to this case. We can thus conclude that ${\mathbb T}_{F}$ is of presheaf type classified by the topos of covariant set-valued functors on the category of finite extensions of $F$.

Next, let us consider the theory ${\mathbb S}_{F}$ of ${\mathbb T}_{F}$ of \emph{separable extensions} of $F$, that is the quotient of ${\mathbb T}_{F}$ obtained by adding the following sequent:
\[
(\top \vdash_{x}  \mathbin{\mathop{\textrm{ $\bigvee$}}\limits_{n\in {\mathbb N}, (a_{0}, \ldots, a_{n-1}, a_{n}) \in {\cal S}_{F}^{n}}} \overline{a_{n}}x^{n}+\overline{a_{n-1}}x^{n-1}+ \cdots + \overline{a_{0}}=0),
\]
where ${\cal S}_{F}^{n}$ is the set of $n$-tuples of elements $a_{1},\ldots, a_{n}$ of $F$ such that the polynomial $a_{n}Z^{n}+a_{n-1}Z^{n-1}+\cdots, +a_{0} \in F[Z]$ is irreducible and separable.

Clearly, the finitely presentable ${\mathbb S}_{F}$-models are precisely the finite separable extensions of $F$. In particular, every finitely presentable ${\mathbb S}_{F}$-model is finitely presentable as a ${\mathbb T}_{F}$-model. In fact, by Artin's primitive element theorem, every finite separable extension of $F$ is presented by a formula of the form $\{x. \overline{a_{n}}x^{n}+\overline{a_{n-1}}x^{n-1}+ \cdots + \overline{a_{0}}=0\}$.

The hypotheses of Corollary \ref{corquotients} are trivially satisfied; indeed, for any ${\mathbb T}_{F}$-model homomorphism $M\to N$ (in an arbitrary Grothendieck topos), if $N$ is separable (i.e., a model of ${\mathbb S}_{F}$) then $M$ is \emph{a fortiori} separable as well. We can thus conclude that also the theory ${\mathbb S}_{F}$ is of presheaf type, classified by the category of covariant set-valued functors on the category of finite separable extensions of $F$. 

\begin{remark}
The theory of fields of a fixed finite characteristic $p$ which are algebraic over their prime field, which was proved in \cite{OCJ} to be of presheaf type classified by the category of covariant set-valued functors on the category of finite fields of characteristic $p$, is (trivially) Morita-equivalent to the theory ${\mathbb T}_{{\mathbb Z}_{p}}$ introduced above.
\end{remark}

\subsection{Groups with decidable equality}

In this section we shall study the injectivization $\mathbb G$ of the (algebraic) theory of groups. 

Clearly, the finitely presentable $\mathbb G$-models are precisely the finitely generated groups.

Even if $\mathbb G$ satisfies condition $(iii)$ of Theorem \ref{main} with respect to the category of finitely generated groups and injective homomorphisms between them (by Theorem \ref{condition3} and Corollary \ref{corcondition3}), $\mathbb G$ is not of presheaf type. To see this, consider the property of an element $x$ of a group $G$ to be non-nilpotent. This property is clearly preserved by injective homomorphisms of groups, so if $\mathbb G$ were of presheaf it would be definable by a geometric formula $\phi(x)$ over the signature of $\mathbb G$. Then the sequent  
\[
(\top \vdash_{x}  \mathbin{\mathop{\textrm{ $\bigvee$}}\limits_{n\in {\mathbb N}}} \phi \vee (x^{n}=1)),
\]
would be provable in $\mathbb G$ (since $\mathbb G$ has enough set-based models and this sequent is valid in every set-based $\mathbb G$-model by definition of $\phi$) and hence, as $\mathbb G$ is coherent, the disjunction on the right hand side would be $\mathbb G$-provably equivalent to a finite sub-disjunction; but this is absurd since it implies that there exists a natural number $n$ such that every nilpotent element $x$ of a group satisfies $x^{n}=1$.

To obtain a presheaf completion of the theory $\mathbb G$, we add to the signature of $\mathbb G$ a relation symbol $R_{N}^{n}$ for each natural number $n$ and any normal subgroup of the free group $F_{n}$ on $n$ generators, and the following axioms (where the symbol $\neq$ denotes the predicate of $\mathbb G$ which is $\mathbb G$-provably complemented to the equality relation):
\[
(\top \vdash_{\vec{x}} R_{N}^{n}(\vec{x}) \vee (\mathbin{\mathop{\textrm{ $\bigvee$}}\limits_{w(\vec{x}),w'(\vec{x})\in F_{n} \textrm{ | } ww'^{-1}\in N}} w\neq w') \vee (\mathbin{\mathop{\textrm{ $\bigvee$}}\limits_{w(\vec{x}),w'(\vec{x})\in F_{n} \textrm{ | } ww'^{-1}\notin N}}w=w'))
\]
and
\[
(R_{N}^{n}(\vec{x}) \wedge (\mathbin{\mathop{\textrm{ $\bigvee$}}\limits_{w(\vec{x}),w'(\vec{x})\in F_{n} \textrm{ | } ww'^{-1}\in N}} w\neq w' \vee \mathbin{\mathop{\textrm{ $\bigvee$}}\limits_{w(\vec{x}),w'(\vec{x})\in F_{n} \textrm{ | } ww'^{-1}\notin N}}w=w') \vdash_{\vec{x}} \bot)
\]
(for any $n\in {\mathbb N}$ and any normal subgroup $N$ of $F_{n}$), and
\[
(\top \vdash_{\vec{x}} \mathbin{\mathop{\textrm{ $\bigvee$}}\limits_{N\in {\cal N}_{n}}} R_{N}^{n}(\vec{x})), 
\]
where ${\cal N}_{n}$ is the set of normal subgroups of $F_{n}$ (for any $n\in {\mathbb N}$). 

Let ${\mathbb G}_{p}$ be the resulting theory; we shall prove that it is of presheaf type. 

Notice that the theories ${\mathbb G}$ and ${\mathbb G}_{p}$ satisfy the first set of hypotheses of Lemma \ref{lemmafp}. 

Let us first show that every finitely generated group is finitely presented as a model of the cartesianization of ${\mathbb G}_{p}$. By Lemma \ref{lemmafp}, to prove that $F_{n}\slash N$ is presented by the formula $R_{N}^{n}(\vec{x})$ as a model of the cartesianization of ${\mathbb G}_{p}$, it suffices to verify that for any set-based model $(G', \{(R_{N}^{n})_{G} \textrm{ | } n\in {\mathbb N}, N\in {\cal N}_{n} \})$ of the cartesianization of ${\mathbb G}_{p}$, denoting by $\vec{\xi}=(\xi_{1}, \ldots, \xi_{n})$ the generators of the group $F_{n}$, the group homomorphisms $f:F_{n}\slash N \to G'$ such that $f(\vec{\xi})\in [[\vec{x}.R_{N}^{n}]]_{G'}$ correspond exactly to the $n$-tuples $\vec{y}$ of elements of $G'$ which belong to the interpretation in $G'$ of the formula $R_{N}^{n}$ (via the assignment $f\to f(\vec{\xi})$). Given a $n$-tuple $\vec{y}$ of elements of $G'$ which belong to the interpretation of $R_{N}^{n}$ in $G'$, we can define a function $f:F_{n}\slash N \to G'$ by setting $f([w])=w_{G'}(\vec{y})$. This is a well-defined injective group homomorphism since the following sequents are provable in the cartesianization of ${\mathbb G}_{p}$ and hence are valid in $G'$:
\[
(R_{N}^{n} \vdash_{\vec{x}} w=w' )
\]
for any $w, w'\in F_{n}$ such that $ww'^{-1}\in N$, and
\[
(R_{N}^{n} \vdash_{\vec{x}} w\neq w' )
\]
for any $w, w'\in F_{n}$ such that $ww'^{-1}\notin N$. Since every finitely generated group is, up to isomorphism, of the form $F_{n}\slash N$ for some natural number $n$ and a normal subgroup $N$ of $F_{n}$, we can conclude that the theory ${\mathbb G}_{p}$ satisfies condition $(iii)$ of Theorem \ref{main} (cf. Proposition \ref{criteriongeometricity}(i)). An alternative way to prove this would have been to invoke Corollary \ref{corcondition3} and Theorem \ref{condition3}.  

As in the case of Diers fields treated in section \ref{diers}, in order to verify that the theory ${\mathbb G}_{p}$ satisfies conditions $(i)$ and $(ii)$ of Theorem \ref{main}, one is reduced to show that condition $(ii)(b)$ of Theorem \ref{thcond1} is satisfied; but this follows from Remark \ref{remcond1}(d). The theory ${\mathbb G}_{p}$ is thus of presheaf type classified by the topos $[\textrm{f.g.} \textbf{Grp}, \Set]$, where $\textrm{f.g.} \textbf{Grp}$ is the category of finitely generated groups and injective homomorphisms between them.

The category $\textrm{f.g.} \textbf{Grp}$ is cocartesian; indeed, it has an initial object (namely, the trivial group) and pushouts (given by the free product with amalgamation construction, cf. \cite{schreier}). It follows the topos $[\textrm{f.g.} \textbf{Grp}, \Set]$ is coherent. The theory ${\mathbb G}_{p}$, in spite of being infinitary, is thus classified by a coherent topos. This fact has various implications for ${\mathbb G}_{p}$. For instance, the coherence of the classifying topos for ${\mathbb G}_{p}$ implies that every formula over the signature of ${\mathbb G}_{p}$ presenting a finitely generated group (regarded as a ${\mathbb G}_{p}$-model) is not only ${\mathbb G}_{p}$-irreducible, but also a coherent object of the classifying topos; in particular, for any ${\mathbb G}_{p}$-compact formula $\{\vec{x}. \phi\}$ over the signature of ${\mathbb G}_{p}$, the formula $\{x,y.\phi(x) \wedge \phi(y)\}$ is also ${\mathbb G}_{p}$-compact (recall that a formula $\{\vec{z}. \chi\}$ over the signature $\Sigma$ of a geometric theory $\mathbb T$ is $\mathbb T$-compact if whenever $\{\vec{z}. \chi\}$ $\mathbb T$-provably entails a disjunction of geometric formulas over $\Sigma$, $\{\vec{z}. \chi\}$ $\mathbb T$-provably entails a finite sub-disjunction of it). Semantically, a formula $\{\vec{x}. \phi\}$ is ${\mathbb G}_{p}$-compact if whenever $\{S_{i} \textrm{ | } i\in I\}$ is a family of assignments $G\to S^{G}_{i}\subseteq G^{n}$ sending each finitely generated group $G$ to a subset $S^{G}_{i}\subseteq [[\vec{x}. \phi]]_{G}$ in such a way that every injective homomorphism $f:G\to G'$ of groups sends tuples in $S^{G}_{i}$ to tuples in $S^{G'}_{i}$, if $[[\vec{x}. \phi]]_{G}=\mathbin{\mathop{\textrm{ $\bigcup$}}\limits_{i\in I}} S^{G}_{i}$ for all $G$ then there exists a finite subset $J\subseteq I$ such that $[[\vec{x}. \phi]]_{G}=\mathbin{\mathop{\textrm{ $\bigcup$}}\limits_{i\in J}} S^{G}_{i}$. The fact that the classifying topos of ${\mathbb G}_{p}$ is coherent also implies, by Deligne's theorem (assuming the axiom of choice), the existence of set-based models for any non-contradictory quotient of ${\mathbb G}_{p}$ whose associated Grothendieck topology on ${\textrm{f.g.} \textbf{Grp}}^{\textrm{op}}$ is of finite type.

\subsection{Locally finite groups}
 
Let $\mathbb A$ be the algebraic theory of groups. Since every finite group is finitely presented as a $\mathbb A$-model (cf. Theorem 6.4 \cite{OCG}) and $\mathbb A$ is of presheaf type, Theorem \ref{crit} ensures that there exists a quotient $\mathbb U$ of $\mathbb A$ classified by the topos $[{\cal C}, \Set]$, where $\cal C$ is the category of finite groups and homomorphisms between them, which can be characterized as the set of all geometric sequents over the signature of $\mathbb A$ that are valid in every finite group. On the other hand, since $\mathbb U$ is classified by the topos $[{\cal C}, \Set]$, the set-based models of $\mathbb U$ are exactly the groups which can be expressed as filtered colimits of finite groups, such groups are exactly the groups which validate all the geometric sequents over the signature of $\mathbb A$ which hold in every finite group. As a by-product, we obtain the following characterization of locally finite groups (equivalently, of the groups which can be expressed as filtered colimits of finite groups).

\begin{proposition}
The locally finite groups are exactly the groups which satisfy all the geometric sequents over the signature of the theory of groups which hold for all finite groups.
\end{proposition}

\begin{proofs}
In view of the above remarks, it remains to verify that a group is locally finite (in the sense that all its finitely generated subgroups are finite) if and only if it is a filtered colimit of finite groups. This can be proved as follows. If a group is locally finite then it is the filtered union of all its finitely generated (and hence finite) subgroups. Conversely, suppose that  $G$ is a filtered colimit of finite groups. Then $G$ is the directed union of the images in $G$ of these finite subgroups, which are again finite. It follows that every finitely generated subgroup $H$ of $G$ is contained in one of them (notice that, since the union is filtered, there exists one of them which contains all the generators of $H$) and hence it is \emph{a fortiori} finite, as required. 
\end{proofs}

The injectivization of $\mathbb A$ is also of presheaf type (by Corollary \ref{corinjectivizations}) and can be characterized as the quotient of $\mathbb G$ consisting of all the geometric sequents over the signature of the injectivization $\mathbb G$ of the theory of groups which hold in all finite groups.

\subsection{Vector spaces}

Let ${\mathbb V}_{K}$ be the expansion of the algebraic theory of vector spaces over a field $K$ obtained by adding a $n$-ary predicate $R_{n}$ for each natural number $n$ expressing the property of a $n$-tuple of elements to be linearly independent, i.e. the following sequents:
\[
(\top \vdash_{\vec{x}} R_{n}(\vec{x}) \vee \mathbin{\mathop{\textrm{ $\bigvee$}}\limits_{\substack{(k_{1},\ldots, k_{n})\in K^{n}\\ k_{i}\neq 0 \textrm{ for some $i$}}}} k_{1}x_{1}+ \cdots + k_{n}x_{n}=0), 
\]  
and
\[
(R_{n}(\vec{x}) \wedge (\mathbin{\mathop{\textrm{ $\bigvee$}}\limits_{\substack{(k_{1},\ldots, k_{n})\in K^{n} \\ k_{i}\neq 0 \textrm{ for some $i$}}}} k_{1}x_{1}+ \cdots + k_{n}x_{n}=0) \vdash_{\vec{x}} \bot). 
\]

The category of models of the theory ${\mathbb V}_{K}$ in $\Set$ has as objects the vector spaces over $K$ and as arrows the injective homomorphisms between them. The finitely presentable ${\mathbb V}_{K}$-models are precisely the finite-dimensional vector spaces over $K$. 

By using techniques analogous to those employed in section \ref{diers}, one can prove that ${\mathbb V}_{K}$ is of presheaf type.

Also, by using Corollary \ref{corquotients}, one can easily prove that for any fixed natural number $n$, both the expansion of the theory of vector spaces over a field $K$ and of the theory ${\mathbb V}_{K}$ obtained by adding the sequent
\[
\top \vdash_{(x_{1}, \ldots, x_{n+1})} \mathbin{\mathop{\textrm{ $\bigvee$}}\limits_{\substack{(k_{1},\ldots, k_{n+1})\in K^{n+1} \\ k_{i}\neq 0 \textrm{ for some $i$}}}} k_{1}x_{1}+ \cdots + k_{n+1}x_{n+1}=0)
\]
are of presheaf type.
 
The models in $\Set$ of the latter theory are precisely the vector spaces over $K$ of dimension $\leq n$.

\subsection{The theory of abelian $l$-groups with strong unit}

Recall that an abelian $l$-group with strong unit is a lattice-ordered group $(G, 0, \leq)$ with a distinguished element $u$, called the unit of the group, such that for any element $x\in G$ such that $x \geq 0$ there exists a natural number $n$ such that $x\leq nu$, where $nu=u+\cdots+u$ $n$ times. We refer the reader to chapter $2$ of \cite{Steinberg} as an introduction to the theory of lattice-ordered groups.

We can axiomatize the theory ${\mathbb L}_{u}$ of abelian $l$-groups with strong unit over a signature $\Sigma$ consisting of four binary function symbols $+$, $-$, $inf$, $sup$, two constants $0$ and $u$ and a binary relation symbol $\leq$, by using Horn sequents to formalize the notion of abelian $l$-group and the following geometric sequent to express the property of strong unit:
\[
(x\geq 0\vdash_{x} \bigvee_{n\in\mathbb{N}}x\leq nu).
\]    

The following lemma will be useful in showing that the theory ${\mathbb L}_{u}$ is of presheaf type.

\begin{lemma}\label{lemmaunit}
Let $G$ be an abelian group with a distinguished element $u$ and generators $x_{1}, \ldots, x_{n}$. If for every $i\in \{1, \ldots, n\}$ there exists a natural number $k_{i}$ such that $|x_{i}|\leq k_{i}u$ then $u$ is a strong unit for $G$.
\end{lemma}
  
\begin{proofs}
Recall that the absolute value $|x|$ of an element $x$ of an abelian $l$-group with unit $(G, 0, +, -, \leq, inf, sup)$ is the element $sup(x, -x)$. For any $x\in G$, $|x|\geq 0$ $|x|=|-x|$, and for any $x, y\in G$, the triangular inequality $|x+y|\leq |x|+|y|$ holds.

Since $G$ is generated by elements $x_{1}, \ldots, x_{n}$, every element $x$ of $G$ can be expressed as the interpretation $t(x_{1}, \ldots, x_{n})$ of a term $t$ over the signature $\Sigma$. We shall prove that there exists a natural number $n$ such that $|x|\leq nu$ by induction on the structure of $t$. This will clearly imply our thesis, since if $x\geq 0$ then $|x|=x$. If $t$ is a variable then the claim is clearly true by our hypothesis. If $x=x'+x''$ with $|x'|\leq n'u$ and $|x''|\leq n''u$ then by the triangular inequality we have $|x|\leq n'+n''$, and similarly for the subtraction. The $inf$ and $sup$ cases are similarly straightforward.  
\end{proofs}  

Let us now verify that the theory ${\mathbb L}_{u}$ satisfies the hypotheses of Corollary \ref{corquotients} with respect to its Horn part $\mathbb H$.

We have to prove that for any finitely presentable $\mathbb H$-model $c$, any model $G$ of ${\mathbb L}_{u}$ in a Grothendieck topos $\cal E$, any object $E$ of $\cal E$ and any $\Sigma$-structure homomorphism $f:c\to Hom_{\cal E}(E, G)$, there exists an epimorphic family $\{e_{i}:E_{i} \to E \textrm{ | } i\in I\}$ in $\cal E$ and for each $i\in I$ a finitely presentable model $c_{i}$ of ${\mathbb L}_{u}$ and $\Sigma$-structure homomorphisms $f_{i}:c\to c_{i}$ and $u_{i}:c_{i}\to Hom_{\cal E}(E_{i}, G)$ such that $Hom_{\cal E}(e_{i}, G)\circ f=u_{i}\circ f_{i}$ for all $i\in I$. Let us suppose that $c$ is presented as a $\mathbb H$-model by a cartesian formula $\phi(\vec{y})$ with generators $x_{1}, \ldots, x_{n}$.  Since $G$ is a $l$-group with strong unit, there exists an epimorphic family $\{e_{i}:E_{i} \to E \textrm{ | } i\in I\}$ in $\cal E$ and for each $k\in \{1, \ldots, n\}$ and $i\in I$ a natural number $m_{k, i}$ such that $f(|x_{i}|)\circ e_{i} \leq m_{k, i}u_{Hom_{\cal E}(E_{i}, G)}$ (where $u_{Hom_{\cal E}(E_{i}, G)}$ denotes the unit of the $\ell$-group $Hom_{\cal E}(E_{i}, G)$). For each $i\in I$, let $c_{i}$ the $\mathbb H$-model presented by the cartesian formula $\phi(x_{1}, \ldots, x_{n}) \wedge |x_{1}|\leq m_{1, i} \wedge \ldots \wedge |x_{n}|\leq m_{n, i}$. For each $i\in I$, we have a natural quotient homomorphism $f_{i}:c\to c_{i}$ through which $Hom_{\cal E}(e_{i}, G)\circ f$ factors; the resulting factorization $u_{i}$ satisfies the required property $Hom_{\cal E}(e_{i}, G)\circ f=u_{i}\circ f_{i}$. Since $\mathbb H$ is a Horn theory, each $c_{i}$ is generated by the $n$-tuple $f_{i}(x_{1}), \ldots, f_{i}(x_{n})$ which presents it as a $\mathbb H$-model (cf. Remark \ref{remblasce}(a)). Therefore the $c_{i}$ are $l$-groups with strong unit by Lemma \ref{lemmaunit}.

This argument also shows that the finitely presentable ${\mathbb L}_{u}$-models are exactly the finitely presented $\mathbb H$-models whose unit is strong (cf. Theorem \ref{rigidity}).

Therefore all the hypotheses of Corollary \ref{corquotients} are satisfied and we can conclude that the theory ${\mathbb L}_{u}$ is of presheaf type. In fact, ${\mathbb L}_{u}$ is Morita-equivalent to the (algebraic) theory of MV-algebras (cf. \cite{OCAR}).


{\small

\begin{thebibliography}{10}

\bibitem{grothendieck}  M. Artin, A. Grothendieck and J. L. Verdier, \emph{Th\'eorie des topos et cohomologie \'etale des sch\'emas}, S\'eminaire de G\'eom\'etrie Alg\'ebrique du Bois-Marie, ann\'ee 1963-64; second edition published as Lecture Notes in Math., vols. 269, 270 and 305(Springer-Verlag, 1972).

\bibitem{beke} T. Beke, Theories of presheaf type, \emph{J. Symbolic Logic} 69 (3), 923-934 (2004).

\bibitem{blasce} A. Blass and A. \v{S}\v{c}edrov, Classifying topoi and finite forcing, \emph{Journal of Pure and Applied Algebra} 28 (2), 111-140 (1983).

\bibitem{borceux} F. Borceux, \emph{Handbook of categorical algebra} vol. 1, (Cambridge University Press, 1994).

\bibitem{OC} O. Caramello, Yoneda representations of flat functors and classifying toposes, \emph{Theory and Applications of Categories} 26 (21), 538-553 (2012).

\bibitem{OCG} O. Caramello, Topological Galois Theory, \emph{arXiv:math.CT/1301.0300}.

\bibitem{OC2} O. Caramello, Fra\"iss\'e's construction from a topos-theoretic perspective, \emph{arXiv:math.CT/0805.2778}, to appear in \emph{Logica Universalis}. 

\bibitem{OCS} O. Caramello, Syntactic characterizations of properties of classifying toposes, \emph{Theory and Applications of Categories} 26 (6), 176-193 (2012).

\bibitem{OC6} O. Caramello, Lattices of theories, \emph{arXiv:math.CT/0905.0299}. 

\bibitem{OC10} O. Caramello, The unification of Mathematics via Topos Theory,\\ \emph{arXiv:math.CT/1006.3930v1}.

\bibitem{OCJ} O. Caramello and P. T. Johnstone, De Morgan's law and the theory of fields, \emph{Advances in Mathematics} 222 (6), 2145-2152 (2009).

\bibitem{OCAR} O. Caramello and A. C. Russo, The Morita-equivalence between MV-algebras and abelian $\ell$-groups with strong unit, \emph{arXiv:math.CT/arXiv:1312.1272}.

\bibitem{OCU} O. Caramello, Universal models and definability, \emph{Math. Proc. Cam. Phil. Soc.} 152 (2), 279-302 (2012).

\bibitem{cole} J. Cole, The bicategory of toposes, and spectra, unpublished manuscript.

\bibitem{connes2} A. Connes, Cohomologie cyclique et foncteur $Ext^{n}$, \emph{C.R. Acad. Sci. Paris, Sér. I Math}
296, 953-958 (1983).

\bibitem{coste} M. Coste, Localisation, spectra and sheaf representation, \emph{Applications of sheaves}, Lecture Notes in Math. vol. 753 (Springer-Verlag, 1979), 212-238. 

\bibitem{Coumans} D. Coumans, Generalising canonical extension to the categorical setting, \emph{Annals of Pure and Applied Logic}, 163 (12), 1940-1961 (2012). 

\bibitem{Giraud} J. Giraud, Analysis situs, \emph{S\'eminaire Bourbaki}, expos\'e 256 (1963).

\bibitem{Hakim} M. Hakim, \emph{Topos annel\'es et sch\'emas relatifs} (Springer-Verlag, 1972).

\bibitem{Hodges} W. Hodges, \emph{Model theory}, Encyclopedia of Mathematics and its Applications vol. 42 (Cambridge University Press, 1993).

\bibitem{El} P. T. Johnstone, \emph{Sketches of an Elephant: a topos theory compendium, Vols. 1 and 2}, vols. 43 and 44 of \emph{Oxford Logic Guides} (Oxford University Press, 2002).

\bibitem{diers} P. T. Johnstone, A syntactic approach to Diers' localizable categories, in \emph{Applications of sheaves}, Lecture Notes in Math. vol. 753 (Springer-Verlag, 
1979), 466-478. 

\bibitem{johwra} P. T. Johnstone and G. Wraith, Algebraic theories in toposes, in \emph{Indexed categories and
their applications}, Lecture Notes in Mathematics (Springer-Verlag) vol. 661 141-242 (1978)

\bibitem{JW} A. Joyal and G. Wraith, Eilenberg-MacLane toposes and cohomology, \emph{Mathematical applications of category theory}, Contemporary Mathematics, vol. 30 (American Mathematical Society, 1984), 117-131.

\bibitem{Lawvere} W. Lawvere, Functorial Semantics of Algebraic Theories, 1963, republished in \emph{Reprints in Theory and Applications of Categories} 5, 1-121 (2004). 

\bibitem{MM} S. Mac Lane and I. Moerdijk, \emph{Sheaves in geometry and logic: a first introduction to topos theory} (Springer-Verlag, 1992).

\bibitem{MR} M. Makkai and G. Reyes, \emph{First-order categorical logic},  Lecture Notes in Mathematics (Springer-Verlag) vol. 611 (1977).

\bibitem{Moerdijk} I. Moerdijk, Cyclic sets as a classifying topos (Utrecht, 1996), available at \texttt{http://ncatlab.org/nlab/files/MoerdijkCyclic.pdf}.

\bibitem{PS} R. Par\'e and D. Schumacher, Abstract families and the Adjoint Functor Theorems, in \emph{Indexed categories and their applications}, Lecture Notes in Math. vol. 661 (Springer-Verlag, 1978), 1-125.

\bibitem{schreier} O. Schreier, Die Untergruppen der freien Gruppen, \emph{Abh. Math. Sem. Univ. Hamburg} vol 5, 161-183 (1927).

\bibitem{Steinberg}
S. A. Steinberg, \emph{Lattice-ordered Rings and Modules}, Springer (2010).

\bibitem{vickers} S. Vickers, Strongly algebraic = SFP (topically), \emph{Mathematical Structures in Computer Science} 11 (6), 717-742 (2001).

\end{thebibliography}
}
\end{document}